\documentclass[11pt,reqno,a4paper]{amsart}
\usepackage{amsmath,amssymb}
\usepackage{mathrsfs}
\usepackage{mathtools}
\usepackage[hidelinks]{hyperref}
\usepackage{comment}
\usepackage{enumitem}
\usepackage{cases}
\usepackage{array}

\numberwithin{equation}{section}

\newtheorem{theorem}{Theorem}
\newtheorem{proposition}{Proposition}[section]
\newtheorem{lemma}[proposition]{Lemma}
\newtheorem{corollary}[proposition]{Corollary}

\theoremstyle{remark}
\newtheorem{remark}[proposition]{Remark}
\newtheorem*{remark*}{Remark}

\renewcommand{\epsilon}{\varepsilon}

\begin{document}
\raggedbottom

\title[The vortex filament conjecture for Gross--Pitaevskii]
{On the vortex filament conjecture for the Gross-Pitaevskii equation}

\author[M.~del Pino]{Manuel del Pino}
\address{\noindent M.~del Pino: Department of Mathematical Sciences, University of Bath, Bath BA2 7AY, United Kingdom}
\email{mdp59@bath.ac.uk}

\author[R.~Juneman]{Rowan Juneman}
\address{\noindent R.~Juneman: Department of Mathematics, Chinese University of Hong Kong, Shatin, NT, Hong Kong.}
\email{rowanjuneman@cuhk.edu.hk}

\author[M.~Musso]{Monica Musso}
\address{\noindent M.~Musso: Department of Mathematical Sciences, University of Bath, Bath BA2 7AY, United Kingdom}
\email{mm2683@bath.ac.uk}

\author[J.~Wei]{Juncheng Wei}
\address{\noindent J.~Wei: Department of Mathematics, Chinese University of Hong Kong, Shatin, NT, Hong Kong.}
\email{wei@math.cuhk.edu.hk}

\begin{abstract}
We establish one form of the vortex filament conjecture for the three-dimensional Gross–Pitaevskii equation. Given any smooth closed embedded binormal flow of curves on a compact time interval, we construct, in the small-core limit \(\epsilon\to0\), a family of exact solutions whose degree-one vortex filaments converge uniformly to the prescribed flow. Near the filament the solutions have the standard planar vortex profile, while away from it their phase gradients converge to the associated Biot–Savart field. We also derive a refined modulation law for the vortex filament.
\end{abstract}

\maketitle

\section{Introduction}

For \(0<\epsilon<1\), we consider the three-dimensional
Gross--Pitaevskii equation
\begin{equation}\label{gpepsintro}
    |\log\epsilon|\,i\partial_tu
    +\Delta_xu+\frac1{\epsilon^2}(1-|u|^2)u=0,
    \qquad (x,t)\in\mathbb R^3\times[0,T].
\end{equation}
Equation \eqref{gpepsintro} arises in quantum physics as a mean-field model for superfluids
and Bose--Einstein condensates.  In this description, the solution represents
a complex order parameter whose codimension-two defects are vortex filaments:
curves along which \(u\) vanishes and around which its phase has non-zero
degree.  The vortex filament conjecture predicts that, while a thin filament remains smooth and separated from other filaments, its leading-order
motion is governed by the binormal curvature flow
\begin{equation}\label{bcfintro}
    \partial_t\gamma
    =\partial_s\gamma\times\partial_{ss}\gamma,
    \qquad |\partial_s\gamma|=1.
\end{equation}

Despite extensive formal and analytical work, the relation between the Gross--Pitaevskii equation and binormal flow has remained largely open; see e.g. \cite{jerrard2014icm} for a survey. The available rigorous results include conditional derivations based on an assumed vortex-core description \cite{lin2000}, concentration results and constructions in highly symmetric settings \cite{bethuelorlandismets2004,chiron2005,jerrard2002}, and adiabatic results under small-curvature assumptions \cite{zhang2022}. Related weak-solution theory for binormal flow itself was developed in \cite{jerrardsmets2015}, while recent asymptotic and numerical studies of the correspondence with Gross--Pitaevskii dynamics can be found in \cite{arrayasetal2025}.

In this paper we give a constructive result for general filament geometries.
Given any smooth binormal flow of embedded closed curves on a compact time
interval, we construct a family of solutions to \eqref{gpepsintro}
carrying a single degree-one vortex filament which converges uniformly to the
prescribed flow as \(\epsilon\to0\). The construction determines the
asymptotic profile of the solution both near and far from the filament, and
yields a refined modulation equation for the centre of the vortex core.

\subsection{Setup and statement of the main result}

The Gross--Pitaevskii equation \eqref{gpepsintro} is a Hamiltonian evolution,
with Hamiltonian given by the Ginzburg--Landau energy
\begin{equation*}
    E_\epsilon(u)
    :=\int_{\mathbb R^3}
        \left\{\frac12|\nabla u|^2
        +\frac1{4\epsilon^2}(1-|u|^2)^2\right\}.
\end{equation*}
In particular, \(E_\epsilon(u(t))\) is independent of \(t\).  The corresponding
natural energy space is
\[
    \mathcal E(\mathbb R^3)
    :=\bigl\{u\in H^1_{\mathrm{loc}}(\mathbb R^3;\mathbb C):
    \nabla u\in L^2(\mathbb R^3),\
    1-|u|^2\in L^2(\mathbb R^3)\bigr\}.
\]
Thus finite-energy configurations approach the vacuum manifold \(S^1\) at
spatial infinity in an averaged sense, while still allowing non-trivial phase
circulation around curves on which \(u\) vanishes.  For each fixed
\(\epsilon>0\), the Cauchy problem is globally well posed in
\(\mathcal E(\mathbb R^3)\); see, for example, \cite{gerard2006}.  The equation
is also invariant under constant phase rotations and Euclidean transformations.

The asymptotic description of a degree-one filament involves two complementary
scales.  At the scale of the ambient space, the filament carries energy
\(\pi|\log\epsilon|\) per unit length to leading order, so that the normalised
energy density converges formally to \(\pi\) times arclength on the limiting
curve.  The factor \(|\log\epsilon|\) in front of \(\partial_tu\) in
\eqref{gpepsintro} selects the slow time scale on which this curve evolves by
binormal flow.

At the scale of the vortex core, the solution is modelled on the standard
planar Ginzburg--Landau vortex.  Writing
\(y=(r\cos\theta,r\sin\theta)\), this profile takes the equivariant form
\begin{equation*}
    W(y)=w(r)e^{i\theta}
\end{equation*}
and solves
\begin{equation*}
    \Delta_yW+(1-|W|^2)W=0
    \qquad\text{in }\mathbb R^2.
\end{equation*}
The solution \(W\) has degree one at infinity and satisfies
\(|W(y)|\to1\) as \(|y|\to\infty\).  Its radial modulus is the unique positive
solution of
\begin{gather*}
    w''+\frac1r w'-\frac1{r^2}w+(1-w^2)w=0,\\[4pt]
    w(0)=0,\qquad
    w(r)\longrightarrow1\quad\text{as }r\longrightarrow\infty.
\end{gather*}
Moreover, \(w\) is strictly increasing and
\begin{align*}
    w(r)&=w'(0)r+O(r^3)
        &&\text{as }r\rightarrow0^+,\\
    w(r)&=1-\frac1{2r^2}+O(r^{-4})
        &&\text{as }r\rightarrow\infty;
\end{align*}
see, for example, \cite{chenelliottqi1994,herveherve1994}.

We next fix the geometric conventions used throughout the paper.  Let
\(\Gamma\subset\mathbb R^3\) be a smooth embedded closed curve of length
\(L\).  We represent \(\Gamma\) by an \(L\)-periodic arclength
parametrisation \(\gamma=\gamma(s)\), \(s\in\mathbb R\); thus
\begin{equation*}
    \gamma(s+L)=\gamma(s),
    \qquad |\partial_s\gamma(s)|=1.
\end{equation*}
Write \(\vec{\mathsf t}=\partial_s\gamma\) for the unit tangent.  Starting
from an oriented orthonormal basis of the normal plane at one point, normal
parallel transport produces a frame along \([0,L]\).  In general, however,
the normal vectors return rotated after one circuit; see
\cite{bishop1975}.  If \(\Theta\) denotes this normal holonomy angle,
rotating the parallel frame through the angle \(-\Theta s/L\) produces a
smooth \(L\)-periodic oriented orthonormal frame
\((\vec{\mathsf t},\vec e_1,\vec e_2)\).  Equivalently, this construction
reflects the fact that the oriented normal bundle of a closed curve in
\(\mathbb R^3\) is trivial.

For this choice of frame, the structure equations take the form
\begin{equation}\label{periodicframeintro}
\begin{split}
    \partial_s\vec{\mathsf t}
       &=\kappa_1\vec e_1+\kappa_2\vec e_2,\\
    \partial_s\vec e_1&=-\kappa_1\vec{\mathsf t}
       +\mathfrak a\vec e_2,\\
    \partial_s\vec e_2&=-\kappa_2\vec{\mathsf t}
       -\mathfrak a\vec e_1,
\end{split}
\end{equation}
where \(\kappa_1\) and \(\kappa_2\) are the curvature components of
\(\gamma\) in the chosen normal frame, while
\(\mathfrak a=-\Theta/L\) is the normal connection coefficient.  Since
\(\gamma\) is parametrised by arclength and
\((\vec e_1,\vec e_2)\) is orthonormal, the usual curvature
\(\kappa=|\partial_{ss}\gamma|\) satisfies
\[
    \kappa^2=\kappa_1^2+\kappa_2^2.
\]
We identify a normal vector \(q_1\vec e_1+q_2\vec e_2\) with the complex
number \(q=q_1+iq_2\).  In this notation, normal covariant differentiation
is represented by
\begin{equation*}
    D_sq=(\partial_s+i\mathfrak a)q.
\end{equation*}
Changing the lift of \(\Theta\) changes \(\mathfrak a\) by an integral
multiple of \(2\pi/L\) and produces a gauge-equivalent periodic frame.  For
a smooth family of closed curves on a time interval, the holonomy angle and
the corresponding periodic frame may be chosen smoothly in time.  We make
one such choice throughout.

Since \(\Gamma\) is embedded, there is a number \(\rho_\Gamma>0\) such that
\begin{equation*}
    x=\gamma(s)+z_1\vec e_1(s)+z_2\vec e_2(s),
    \qquad |z|<\rho_\Gamma,
\end{equation*}
gives periodic tubular coordinates in a neighbourhood of \(\Gamma\).  Both
the frame and the normal variable are now genuinely periodic in \(s\).

We apply this notation to the prescribed binormal flow.  Thus
\(\gamma^0=\gamma^0(s,t)\) is smooth, \(L\)-periodic in \(s\), parametrised
by arclength and satisfies \eqref{bcfintro}.  We write
\[
    \Gamma_t^0:=\{\gamma^0(s,t):s\in\mathbb R\}.
\]
The binormal flow preserves arclength, so \(L\) is independent of \(t\).
We assume that \(\Gamma_t^0\) is embedded for every \(t\in[0,T]\).  By
smoothness and compactness, the corresponding tubular radii may then be
chosen uniformly in time; we fix one such radius \(\rho_0>0\).

Finally, for compact subsets \(A,B\subset\mathbb R^3\), let
\(d_{\rm H}(A,B)\) denote their Hausdorff distance.  Wherever a
complex-valued function \(u\) does not vanish, we write its phase gradient in
the intrinsic form
\begin{equation*}
    \nabla\arg u=\frac{\operatorname{Im}(\overline u\nabla u)}{|u|^2}.
\end{equation*}
We can now state our main result.

\begin{theorem}\label{mainthm}
Let \(\gamma^0\) be a smooth solution of the binormal flow
\eqref{bcfintro} on \([0,T]\), and suppose that \(\Gamma_t^0\) is an
embedded closed curve for every \(t\in[0,T]\).  Then there exist
\(\epsilon_0>0\) and, for every \(0<\epsilon<\epsilon_0\), a smooth
finite-energy solution \(u_\epsilon\) of \eqref{gpepsintro} with the
following properties.

For every \(t\in[0,T]\), the zero set of \(u_\epsilon(\cdot,t)\) in the
\(\rho_0/2\)-tubular neighbourhood of \(\Gamma_t^0\) is a smooth embedded
closed curve \(\Gamma_{\epsilon,t}\), of degree one on each oriented normal
disc, and
\begin{equation*}
    \sup_{0\leq t\leq T}
    d_{\rm H}(\Gamma_{\epsilon,t},\Gamma_t^0)
    \leq \frac{C}{|\log\epsilon|}.
\end{equation*}
Moreover, \(\Gamma_{\epsilon,t}\) has a parametrisation
\(\gamma_\epsilon(s,t)\), with a smooth periodic oriented normal frame
\((\vec e_{1,\epsilon},\vec e_{2,\epsilon})\), and there is a real-valued
phase \(\beta_\epsilon(s,t)\) such that, for every fixed \(R>0\),
\begin{equation*}
\begin{split}
 &u_\epsilon\big(\gamma_\epsilon(s,t)
   +\epsilon y_1\vec e_{1,\epsilon}(s,t)
   +\epsilon y_2\vec e_{2,\epsilon}(s,t),t\big)\\
 &\hspace{8em}=e^{i\beta_\epsilon(s,t)}
   \big(W(y)+r_\epsilon(y,s,t)\big),
\end{split}
\end{equation*}
where
\begin{equation*}
    \sup_{t\in[0,T]}\sup_s
    \|r_\epsilon(\cdot,s,t)\|_{C^1(B_R)}\longrightarrow0
    \qquad\text{as }\epsilon\longrightarrow0.
\end{equation*}
Finally, away from the space--time filament traced by \(\gamma^0\),
\(|u_\epsilon|\to1\) locally uniformly and
\begin{equation*}
    \frac{\operatorname{Im}(\overline{u_\epsilon}\nabla u_\epsilon)}
         {|u_\epsilon|^2}
    \longrightarrow
    \frac12\int_0^L
    \frac{\partial_s\gamma^0(s,t)\times
    (x-\gamma^0(s,t))}{|x-\gamma^0(s,t)|^3}\,ds
\end{equation*}
locally uniformly.
\end{theorem}

The local well-posedness theory for binormal flow gives a unique smooth
evolution from every smooth embedded closed initial curve, and embeddedness
persists for a positive time.  Thus Theorem~\ref{mainthm} applies, starting
from any such curve, as long as embeddedness persists.  The proof moreover
contains considerably more information
about the solution profile and filament motion than is recorded above.  In
particular,
after choosing compatible parametrisations, the filament curve admits the
sharper description
\[
 \gamma_\epsilon
 =\gamma_\epsilon^{(1)}
  +O\big(\epsilon|\log\epsilon|^C\big)
\]
uniformly on \([0,T]\), with the error controlled in every fixed Sobolev norm.
Here \(\gamma_\epsilon^{(1)}\) is a distinguished centre curve used in the
construction, which evolves according to the following nonlocal geometric
equation.  If \(\ell_\epsilon(t)\) denotes its length, \(\sigma\) is arclength,
\(\vec{\mathsf t}_\epsilon=\partial_\sigma\gamma_\epsilon^{(1)}\), and
\(\vec K_\epsilon=\partial_\sigma^2\gamma_\epsilon^{(1)}\), then this centre
curve is chosen to satisfy the refined evolution equation
\eqref{geometriccurveequation}, namely
\begin{equation*}
\begin{split}
 \partial_t\gamma_\epsilon^{(1)}
 ={}&a_\epsilon\vec{\mathsf t}_\epsilon
   +\vec{\mathsf t}_\epsilon\times\vec K_\epsilon\\
 &+\frac1{|\log\epsilon|}
 \left\{
 \left(\log\ell_\epsilon+\mathfrak c_\epsilon\right)
       \vec{\mathsf t}_\epsilon\times\vec K_\epsilon
 +\Pi_{\vec{\mathsf t}_\epsilon}^{\perp}
       \mathcal N(\gamma_\epsilon^{(1)})
 \right\}.
\end{split}
\end{equation*}
Here \(a_\epsilon\) fixes a constant-speed parametrisation and does not alter
the image of the curve, while \(\mathfrak c_\epsilon\) is the vortex-core
constant described in \eqref{coreconstantasymptotics}.  The nonlocal term is
the finite-part Biot--Savart field
\begin{equation}\label{nonlocalfieldintro}
\begin{split}
 \mathcal N(\gamma_\epsilon^{(1)})(\sigma)
 :=\int_0^{\ell_\epsilon}
 \bigg\{&
 \frac{\vec{\mathsf t}_\epsilon(\sigma')\times
  \big(\gamma_\epsilon^{(1)}(\sigma)
       -\gamma_\epsilon^{(1)}(\sigma')\big)}
 {\big|\gamma_\epsilon^{(1)}(\sigma)
       -\gamma_\epsilon^{(1)}(\sigma')\big|^3}\\
 &-\frac{\vec{\mathsf t}_\epsilon(\sigma)\times
              \vec K_\epsilon(\sigma)}
 {2d_{\ell_\epsilon}(\sigma,\sigma')}
 \bigg\}\,d\sigma',
\end{split}
\end{equation}
where
\[
 d_{\ell}(\sigma,\sigma')
 :=\min_{q\in\mathbb Z}|\sigma-\sigma'+q\ell|,
 \qquad
 \Pi_{\vec{\mathsf t}}^\perp\vec v
 :=\vec v-(\vec{\mathsf t}\cdot\vec v)\vec{\mathsf t},
 \quad \vec v\in\mathbb R^3.
\]
Thus binormal flow is the leading law, whereas the renormalised curvature
coefficient and the projected finite-part interaction give the first
correction, of relative size \( |\log\epsilon|^{-1}\).  The proof also
constructs an expansion of the solution, and the corresponding further
adjustments of the centre curve, to arbitrary algebraic order in
\(\epsilon\).  The main steps of the argument are described in
Section~\ref{proofoverviewsection}.

\subsection{Related literature}
\leavevmode\par
\medskip
\noindent\textit{Vortex filaments in classical fluids.}
Vortex filaments for the Gross--Pitaevskii equation have a close classical
analogue in the incompressible Euler equations.  The latter describe an ideal
fluid and admit flows whose vorticity is concentrated in thin tubes.  In the
small-core limit, both models give rise to the same Biot--Savart field and
hence to the same formal binormal flow dynamics \cite{jerrard2014icm}.  The
classical theory begins with Helmholtz's study of vortex tubes and Da Rios'
derivation of binormal curvature flow for their centre lines
\cite{darios1906}.  Subsequent rigorous work has produced special solutions
concentrated near circles and helices
\cite{averkioumusso2026,daviladelpinomussowei2022,fraenkel1970,fraenkelberger1974}, together
with a conditional justification of binormal flow for Euler vorticity
satisfying suitable concentration and energy assumptions
\cite{jerrardseis2017}.  We note that the general constructive problem remains
open for 3D Euler.

For the Navier--Stokes equations, viscosity spreads the vorticity into a
Lamb--Oseen-type core.  The vanishing-viscosity limit is understood for
circular rings \cite{gallaysverak2024}, and well-posedness is global near a
straight Oseen column \cite{bedrossiangermainharropgriffiths2023}.  In a small
Reynolds-number regime, some short-time results are available in
\cite{fontelosispizuavega2026,fontelosvega2023}.

\medskip

\noindent\textit{Multiple filaments and vortex reconnection.}
The connection between classical and quantum vortex dynamics is not limited to
isolated filaments.  For configurations of several nearly parallel filaments,
the formal Klein--Majda--Damodaran system
\cite{kleinmajdadamodaran1995} arising from three-dimensional Euler flow
has been rigorously recovered by Jerrard and Smets as an effective
Gross--Pitaevskii dynamics \cite{jerrardsmets2021}. Travelling-wave solutions
associated with this system were subsequently constructed near interacting
helical filaments \cite{daviladelpinomedinarodiac2022}, while interactions
between vortex rings have been studied in
\cite{jerrardsmets2018leapfrogging}. More general interactions may change the
topology of the vortex set \cite{encisoperaltasalas2021}.

\medskip

\noindent\textit{Parabolic Ginzburg--Landau dynamics.}
A complementary geometric limit arises when the Hamiltonian
Gross--Pitaevskii evolution is replaced by the parabolic Ginzburg--Landau flow
\begin{equation*}
    \partial_tu
    =\Delta u
    +\frac1{\epsilon^2}(1-|u|^2)u,
    \qquad
    u:\mathbb{R}^n\times[0,\infty)\longrightarrow\mathbb{C}.
\end{equation*}
In this setting, codimension-two defect sets \(\{u=0\}\) converge to
mean-curvature flow, as first proved for well-prepared initial data by Lin
\cite{lin1998} and Jerrard--Soner \cite{jerrardsoner1999}. Convergence beyond
the lifespan of the smooth mean-curvature flow was later justified in Brakke's
weak formulation, first by Ambrosio and Soner \cite{ambrosiosoner1997} under a
density hypothesis, and subsequently by Bethuel, Orlandi and Smets
\cite{bethuelorlandismets2006} under a natural energy bound.

\medskip

\noindent\textit{Hyperbolic and gauged problems.}
For a class of hyperbolic problems including the Ginzburg--Landau wave equation
\begin{equation}\label{waveGLintro}
    \partial_t^2u
    =\Delta u
    +\frac1{\epsilon^2}(1-|u|^2)u,
    \quad
    u:\mathbb{R}^n\times[0,\infty)\longrightarrow\mathbb{C},\quad 3\leq n\leq4,
\end{equation}
Jerrard \cite{jerrardwave2011} proved the existence of suitably prepared
solutions whose normalised energy concentrates on prescribed smooth timelike
minimal surfaces of codimension two.
More recently, Geevechi and Jerrard \cite{geevechijerrard2025} considered the
critical hyperbolic Abelian Higgs model (a gauged version of
\eqref{waveGLintro}) and constructed slowly moving, nearly parallel vortex
filaments whose leading-order dynamics are governed by a wave map into the
moduli space of planar Abelian Higgs vortices; analogous parabolic results yield
harmonic map heat flow into the same moduli space.

\medskip

\noindent\textit{Stability.}
The results of this paper apply to suitably constructed solutions.  More
general questions concerning the long-time dynamics of \eqref{gpepsintro} and
filament motion for arbitrary initial data remain open.  The difficulty of
controlling perturbations over long times is already present for the planar
degree-one vortex: its orbital stability is known
\cite{gravejatpacheriesmets2022}, but its asymptotic stability remains
unresolved; see
\cite{collotgermainpacherie2025,luhrmannschlagshahshahani2025} for progress in
this direction.

Our restriction to degree-one filaments also reflects an important stability
distinction.  Whereas the planar degree-one vortex is variationally stable,
the equivariant vortices $W_d(r,\theta)=w_d(r)e^{id\theta}$ with $|d|\geq2$ are
saddle points of the renormalised Ginzburg--Landau energy
\cite{ovchinnikovsigal1997}.  For the degree-two and degree-three vortices, the
respective Morse indices have recently been computed in
\cite{delpinoliumussoweiyang2026}.  The resulting negative modes are expected
to split a higher-degree filament into unit-degree filaments, although a
rigorous dynamical description of this splitting mechanism is not yet
available.

\section{Overview of the proof}
\label{proofoverviewsection}

The proof of Theorem \ref{mainthm} is based on the construction of an accurate
approximate solution followed by linearised estimates for the remainder,
broadly following the strategy developed by the first three authors to
describe the dynamics of interacting vortices in the plane
\cite{delpinojunemanmusso2025}.  Passing from point vortices to a filament,
however, introduces substantial geometric complications: the
finite-dimensional modulation equations are replaced by an evolution PDE, the
effective dynamics contain genuinely nonlocal operators, and the translational
modes form an infinite-dimensional family.  The principal steps of the
argument are described below.

\subsection{First ansatz}

A classical construction from electromagnetism associates every oriented
closed curve with a solid-angle function, defined modulo $4\pi$ on the
complement of the curve \cite{binysh2018,maxwell1873}.  Taking half of the
solid angle determines an $S^1$-valued phase of degree one around the curve,
thus encoding the circulation of a unit vortex filament.

To make this construction explicit, let
\[
    \Gamma_t=\{\gamma(s,t):0\leq s<L\}
\]
be a smooth embedded oriented closed curve, where
$s\mapsto\gamma(s,t)$ is an $L$-periodic parametrisation compatible with its
orientation.  Choose a smooth compact oriented Seifert surface
$\Sigma_t\subset\mathbb R^3$ with $\partial\Sigma_t=\Gamma_t$.  We orient
$\Sigma_t$ so that its induced boundary orientation is opposite to that of
$\gamma$, and denote its unit normal by $\vec n_{\Sigma_t}$.

For $x\notin\Sigma_t$, the signed solid angle subtended by $\Sigma_t$ at $x$
is defined as the flux through $\Sigma_t$ of the inverse-square radial field
centred at $x$.  The associated half-angle phase is then given by
\begin{equation}\label{solidanglephaseoverview}
    \varphi(x,t)
    :=\frac12\int_{\Sigma_t}
      \frac{(x-y)\cdot\vec n_{\Sigma_t}(y)}
           {|x-y|^3}\,dS_y.
\end{equation}
The traces of the solid angle on the two sides of $\Sigma_t$ differ by
$4\pi$.  Hence $\varphi$ jumps by $2\pi$, and $e^{i\varphi}$ extends across
the interior of $\Sigma_t$ to a single-valued smooth $S^1$-valued function on
$\mathbb R^3\setminus\Gamma_t$.  This function has degree one around
$\Gamma_t$, and is independent of the particular Seifert surface chosen
(replacing $\Sigma_t$ by another such surface changes $\varphi$ only by
integer multiples of $2\pi$).

The significance of \eqref{solidanglephaseoverview} is twofold.  On the one
hand, a calculation using Stokes' theorem gives
\[
    \nabla_x\varphi(x,t)
    =\frac12\int_0^L
      \frac{\partial_s\gamma(s',t)\times
            (x-\gamma(s',t))}
           {|x-\gamma(s',t)|^3}\,ds',
\]
which is precisely the Biot--Savart velocity field of a classical vortex
filament of circulation $2\pi$.  On the other hand, the degree-one winding of
$e^{i\varphi}$ makes it the natural global phase for our first ansatz.  Let
$(\vec e_1,\vec e_2)$ be a periodic oriented orthonormal frame for the normal
bundle of $\Gamma_t$, and write
\[
    x=\gamma(s,t)+z_1\vec e_1(s,t)+z_2\vec e_2(s,t)
\]
in tubular coordinates.
For $\delta>0$ sufficiently small that these
coordinates are well-defined for $|z|<2\delta$, we choose a smooth cut-off
$\eta(z)$ satisfying
\[
    \eta(z)=1\quad\text{for }|z|\leq\delta,
    \qquad
    \eta(z)=0\quad\text{for }|z|\geq2\delta.
\]
Writing
\[
    W(y)=w(|y|)e^{i\theta}
\]
for the standard planar degree-one vortex, we then define our first
approximate solution to \eqref{gpepsintro} by
\begin{equation}\label{firstansatzoverview}
    u_\epsilon^{(0)}(x,t)
    =
    \left[
      \eta(z)w\!\left(\frac{|z|}{\epsilon}\right)+1-\eta(z)
    \right]e^{i\varphi(x,t)},
\end{equation}
inside the tubular neighbourhood, and set
$u_\epsilon^{(0)}=e^{i\varphi}$ outside it. This ansatz inserts a planar vortex of core radius $\epsilon$ into each normal
plane to $\Gamma_t$, while the solid-angle phase provides its global
continuation.  Near the curve, $\varphi$ has the same leading angular
behaviour as the polar angle $\theta$, and the relative phase
$\varphi-\theta$ extends continuously to $\Gamma_t$.  Consequently,
$u_\epsilon^{(0)}$ extends continuously across $\Gamma_t$ and vanishes
precisely there.

After multiplying by $\epsilon^2$, substitution of
\eqref{firstansatzoverview} into \eqref{gpepsintro} shows that the leading core
terms cancel since $W$ solves the planar vortex equation.  The remaining
inner error begins at order $\epsilon$, while the cut-off and outer errors are
$O(\epsilon^2|\log\epsilon|)$.  A major difficulty is to eliminate the
components of the inner error along the planar translation modes, a requirement
that leads to the nonlocal evolution equation for $\gamma$ discussed next.

\subsection{The corrected curve equation}

To improve the ansatz \eqref{firstansatzoverview} near the filament, we
construct an inner correction by inverting the planar linearised
Ginzburg--Landau operator in each normal plane.  A suitable planar inverse
exists only when the source satisfies the appropriate solvability conditions
against the translation kernel, spanned by $\partial_{y_1}W$ and
$\partial_{y_2}W$; otherwise the inverse develops growing components that
cannot be matched to the outer approximation.  Requiring these two projections
to vanish yields an evolution equation for $\gamma$.

More precisely, writing $\sigma$ for arclength, let
\[
    \vec{\mathsf t}=\partial_\sigma\gamma,
    \qquad
    \vec K=\partial_\sigma^2\gamma
\]
denote the unit tangent and curvature vector, and set
$\gamma_t^\perp=\Pi_{\vec{\mathsf t}}^\perp\partial_t\gamma$.  The two
solvability conditions take the form
\begin{equation}\label{correctedcurveoverview}
 -\gamma_t^\perp
 +\left(1+O\!\left(|\log\epsilon|^{-1}\right)\right)
      \vec{\mathsf t}\times\vec K
 +\frac1{|\log\epsilon|}
      \Pi_{\vec{\mathsf t}}^\perp\mathcal N(\gamma)
 =0,
\end{equation}
where $\mathcal N(\gamma)$ is the finite-part Biot--Savart field defined in
\eqref{nonlocalfieldintro}.  Thus binormal flow gives the leading evolution.
The correction is not, however, of lower order: $\mathcal N(\gamma)$ contains
a nonlocal term of logarithmic order two.

To solve problem \eqref{correctedcurveoverview}, Section
\ref{normalgraphsection} writes $\gamma$ as a normal graph over the prescribed
binormal-flow solution $\gamma^0$ and chooses the tangential velocity so that
the parametrisation remains of constant speed.  The constant-speed constraint
allows the dependent graph and length variables to be recovered smoothly and
tamely from a set of independent graph variables.  In this gauge, the normal
evolution takes the form of a nonlocal quasilinear Schr\"odinger equation.  The
principal operator has logarithmic order two and is self-adjoint before
composition with the normal complex structure $J_\gamma$; it therefore
contributes only a skew-adjoint term to the evolution.  This prevents
derivative loss at the highest-order level, while the remaining
curvature-dependent logarithmic term, variable coefficients, and commutators
satisfy tame linearised estimates with a fixed finite loss of derivatives.

The derivative loss prevents a direct contraction argument, but is
accommodated by a Nash--Moser iteration.  Theorem
\ref{uniformlocalcurvetheorem} consequently provides, for all sufficiently
small $\epsilon$, a unique corrected curve on the full interval $[0,T]$, with
the same initial curve as $\gamma^0$.  After choosing compatible
parametrisations, for each fixed $m\geq6$ this curve satisfies
\[
 \sup_{0\leq t\leq T}
 \|\gamma_\epsilon(\cdot,t)-\gamma^0(\cdot,t)\|_{H^m}
 \leq \frac{C_{m,T}}{|\log\epsilon|}.
\]
The corresponding forced construction depends smoothly and tamely on its
data; this is used to remove the higher-order translation obstructions
arising in the subsequent approximation scheme.

\subsection{Inner--outer construction of an arbitrary-order approximation}

Once the first translation projections vanish, the remaining inner source is
inverted using the two-dimensional linearised Ginzburg--Landau operator.  The inverse is
constructed Fourier mode by Fourier mode as in
\cite{delpinojunemanmussolinear2025}, with the integration constants chosen
to give regularity at the origin and the asymptotic behaviour required for
matching.  This gives a local improvement near the filament.  To obtain a
global improvement, one must also cancel the outer error already present in
the initial ansatz and the additional terms generated by cutting off the
slowly decaying phase components of the inner correction.  Section
\ref{firstglobalimprovementsection} treats these contributions together: the
inner correction is localised near the filament, the leading combined phase
error is cancelled by solving a three-dimensional wave equation, and the
corresponding amplitude correction is determined algebraically.

Section \ref{arbitraryapproxsection} iterates this construction.  At each
stage, the leading translation contribution is cancelled at main order by
applying the tame inverse of the curve equation linearised at a suitable
reference curve.  The remaining inner error is successively reduced by a
finite batch of planar inversions, which are allowed to have controlled
polynomial growth in the normal variable.  The resulting global error then
supplies the source for an outer correction, whose phase is obtained from a
three-dimensional wave equation and whose amplitude is determined
algebraically.  Earlier outer waves are kept fixed as physical functions and
the reference curve linearisation is reused, so its logarithmic order-two
principal part remains the same skew-adjoint operator treated by the
Nash--Moser theory.  Section \ref{threeDwaveestimatesection} provides the
uniform wave estimates needed throughout the iteration.  The conclusion is
Theorem \ref{arbitraryapproximationtheorem}, which yields, to any prescribed
algebraic order $M$, a smooth family of approximate solutions parametrised by
its initial curve, for which the residual and its first two initial-curve
derivatives are $O(\epsilon^M)$ at any prescribed finite Sobolev order.

\subsection{Linear theory around the accurate approximation}

Let $U_\epsilon$ be an accurate approximate solution produced by the
inner--outer construction.  To correct its residual, Section
\ref{curvedlinearsection} develops inverse estimates for the full
time-dependent Gross--Pitaevskii operator linearised at $U_\epsilon$.  Passing
to the global core variable, we set
\[
    X=\frac{x}{\epsilon},
    \qquad
    V_\epsilon(X,t)=U_\epsilon(\epsilon X,t),
\]
and consider
\[
\begin{cases}
\begin{aligned}
    \epsilon^2|\log\epsilon|\,i\partial_t\Phi
    +\Delta_X\Phi+(1-|V_\epsilon|^2)\Phi
    -2\operatorname{Re}(\overline{V_\epsilon}\Phi)V_\epsilon
    &=F,\\
    \Phi(\cdot,0)&=0.
\end{aligned}
\end{cases}
\]
The objective is an estimate for $\Phi$ in terms of the general forcing $F$,
uniformly on the full time interval up to a fixed polynomial loss in
$\epsilon$.

The first step is a unitary transformation to material coordinates moving
with the filament, together with a relative phase that fixes the leading
vortex profile.  Writing $q=\mathscr S_t\Phi$, where $\mathscr S_t$ denotes
this transformation, the transformed equation involves a time-dependent
real-linear elliptic operator $H_t$, symmetric with respect to the real
$L^2$ inner product.  We denote the corresponding quadratic form by $B_t[q]$.

The form $B_t$ is not nonnegative.  Near the filament, angular decomposition
in each normal plane controls the zero-mode phase gradient and amplitude,
and the nonzero modes modulo the two planar translations.  Let
$\tau_1,\tau_2$ be suitably adjusted cut-offs of the two normal derivatives
of the complete approximation, expressed in the moving coordinates and
cut off at a fixed physical distance from the filament.  At each position
along the filament, labelled by $\sigma_0$, we choose coefficients
$b_t(q)=(b_1,b_2)$ so that, in the normal-plane representation,
$q-b_1\tau_1-b_2\tau_2$ is orthogonal to $i\tau_1$ and $i\tau_2$ with respect
to the real $L^2$ inner product.  We then define
\[
 Q_t[q]
 =B_t[q]+C_{\rm tr}c_\epsilon
       \|b_t(q)\|_{L^2_{\sigma_0}}^2
       +\gamma c_\epsilon^2\|q\|_{L^2}^2,
 \qquad c_\epsilon=\epsilon^2|\log\epsilon|.
\]

For sufficiently large fixed $C_{\rm tr},\gamma>0$ and sufficiently small
$\epsilon$, this is a positive quadratic functional satisfying
\[
 c_\epsilon^2\|q\|_{H^1}^2\leq C Q_t[q],
 \qquad Q_t[q]\leq C\|q\|_{H^1}^2.
\]
The translation term controls the two nearly neutral directions, while
the $L^2$ term absorbs the remaining negative zeroth-order contributions.
These coercivity bounds follow from angular decomposition, a radial
estimate for the scalar phase, and cancellation of the boundary fluxes
between the inner and exterior regions.

Controlling the time derivative of $Q_t[q(t)]$ is delicate, since
contributions from the moving coordinates, the material Jacobian, and
the time-dependent phase and amplitude would be too large if estimated
separately.  Writing the perturbation relative to the complete approximation
reveals cancellations: the density and modulus variations cancel the
corresponding scalar transport terms, while the phase rotation cancels
the corresponding phase-gradient contribution from the current.  Further
cancellations among the modulus terms follow from the continuity identity
for the approximation, including its residual.  Combining these identities
with the evolution of the translation coefficients gives
\[
 \frac{d}{dt}Q_t[q(t)]
 \leq C_0|\log\epsilon|\,Q_t[q(t)]
 +\frac{C}{c_\epsilon}
       Q_t[q(t)]^{1/2}\|F(t)\|_{H^1_X}.
\]
Together with the coercivity bounds and the zero initial condition,
Gronwall's inequality then yields
\[
 \|\Phi\|_{L_t^\infty H^1_X}
 \leq C\epsilon^{-4-C_0T}
       \|F\|_{L_t^\infty H^1_X}.
\]

Finally, Section~\ref{fullnonlinearsection} applies this inverse with
\[
    F=-R_\epsilon-\mathcal N_\epsilon(\Phi),
\]
where $R_\epsilon$ is the residual of $V_\epsilon$ and
$\mathcal N_\epsilon(\Phi)$ is quadratic. Since the approximation
can be constructed with arbitrarily small algebraic residual, the polynomial
loss in the inverse is absorbed and the resulting fixed-point map is
contractive. The conclusions of Theorem \ref{mainthm} readily follow.

\section{Ansatz and error of approximation}\label{ansatzsection}

We work with the Gross--Pitaevskii equation
\begin{equation}\label{scaledGP}
 |\log\epsilon|\,i\partial_tu+\Delta_xu
 +\frac{1}{\epsilon^2}(1-|u|^2)u=0,
 \qquad (x,t)\in\mathbb R^3\times[0,T].
\end{equation}
Equivalently, after multiplying by $\epsilon^2$,
\begin{equation}\label{scaledGPoperatorform}
 S_\epsilon(u):=
 \epsilon^2|\log\epsilon|\,i\partial_tu
 +\epsilon^2\Delta_xu+(1-|u|^2)u=0.
\end{equation}
The sign and time scaling in \eqref{scaledGP} are chosen so that a
positively oriented degree-one
filament evolves at leading order by
\begin{equation}\label{bcfintrosection}
 \partial_t\gamma
 =\partial_s\gamma\times\partial_{ss}\gamma.
\end{equation}

\subsection{Periodic normal coordinates}

Let $\gamma=\gamma(s,t)$ be a smooth embedded closed curve, parametrised by
arclength $s$ at each fixed time.  As in the introduction, we choose a smooth
periodic oriented orthonormal frame
$(\vec{\mathsf t},\vec e_1,\vec e_2)$ and use tubular coordinates
\begin{equation*}
 x=\gamma(s,t)+z_1\vec e_1(s,t)+z_2\vec e_2(s,t),
\end{equation*}
where $\vec{\mathsf t}=\partial_s\gamma$.  The frame equations are
\begin{align*}
 \partial_s\vec{\mathsf t}
 &=\kappa_1\vec e_1+\kappa_2\vec e_2,\\
 \partial_s\vec e_1&=-\kappa_1\vec{\mathsf t}
    +\mathfrak a\vec e_2,\\
 \partial_s\vec e_2&=-\kappa_2\vec{\mathsf t}
    -\mathfrak a\vec e_1.
\end{align*}
Here $\mathfrak a=-\Theta/L$ is constant in $s$, where $\Theta$ is a
continuous lift of the normal holonomy angle.  The integer in this lift
is kept fixed throughout the time interval and the local curve chart.  The
construction preceding \eqref{periodicframeintro} then gives a periodic frame with
constant connection, smooth in $t$ and in the local curve parameters.

We identify a normal vector $q_1\vec e_1+q_2\vec e_2$ with
$q=q_1+iq_2$.  In this trivialisation the normal covariant derivative is
\begin{equation*}
 D_sq=(\partial_s+i\mathfrak a)q.
\end{equation*}
If $\kappa=|\partial_{ss}\gamma|$ and $\vec n,\vec b,\tau$ are the Frenet
normal, binormal and torsion (defined where $\kappa\neq0$), then
\begin{equation*}
 (\kappa_1,\kappa_2)=\kappa(\cos\alpha,\sin\alpha),
 \qquad \partial_s\alpha+\mathfrak a=\tau,
\end{equation*}
and
\begin{align*}
 \vec n&=\cos\alpha\,\vec e_1+\sin\alpha\,\vec e_2,\\
 \vec b&=-\sin\alpha\,\vec e_1+\cos\alpha\,\vec e_2.
\end{align*}
In real components we write
\begin{equation*}
 \vec\kappa=(\kappa_1,\kappa_2),
 \qquad
 J\vec\kappa=(-\kappa_2,\kappa_1),
\end{equation*}
where $J(a_1,a_2)=(-a_2,a_1)$.  The normal vector represented by
$J\vec\kappa$ is
\[
 -\kappa_2\vec e_1+\kappa_1\vec e_2
 =\vec{\mathsf t}\times\partial_{ss}\gamma,
\]
which equals $\kappa\vec b$ wherever $\kappa\neq0$.  Thus the component
vector of $D_s\vec\kappa$ is
$\partial_s\vec\kappa+\mathfrak aJ\vec\kappa$.

Put $\partial_\theta=Jz\cdot\nabla_z$.  The coordinate vector in the
$s$-direction and the corresponding horizontal derivative are
\begin{equation*}
 \partial_sx=h\vec{\mathsf t}+\mathfrak aJz,
 \qquad
 \mathscr D_s:=\partial_s-\mathfrak a\partial_\theta,
\end{equation*}
where
\begin{equation*}
 h=1-\kappa_1z_1-\kappa_2z_2.
\end{equation*}
For an evolving frame we also write
\begin{equation*}
 \mathfrak c=\partial_t\vec e_1\cdot\vec e_2,
 \qquad
 D_t=\partial_t+i\mathfrak c,
 \qquad
 \mathscr D_t=\partial_t-\mathfrak c\partial_\theta.
\end{equation*}
The remaining terms in the coordinate time derivative are generated by the
velocity of the centre curve and by its tangential reparametrisation.  The
operator $\mathscr D_t$ includes the normal rotation of the periodic frame.
For a scalar function $q=q(s,z)$,
\begin{align}
 \nabla_xq
 &=h^{-1}(\mathscr D_sq)\vec{\mathsf t}
   +(\partial_{z_1}q)\vec e_1+(\partial_{z_2}q)\vec e_2,
 \label{gradientperiodicframe}\\
 \Delta_xq
 &=\Delta_zq-\frac{\vec\kappa}{h}\cdot\nabla_zq
   +\frac1h\mathscr D_s
      \left(h^{-1}\mathscr D_sq\right)\label{laplaceperiodicframe}\\
 &=\Delta_zq-\frac{\vec\kappa}{h}\cdot\nabla_zq
   +\frac1{h^2}\mathscr D_s^2q
   +\frac{(D_s\vec\kappa)\cdot z}{h^3}\mathscr D_sq.
 \nonumber
\end{align}
These identities follow directly from the coordinate metric in
$(s,z_1,z_2)$, for which $\det(g_{ij})=h^2$.  Radial
functions of $z$ are annihilated by $\partial_\theta$, whereas angular
Fourier modes see the shifted derivative
$\partial_s-im\mathfrak a$.

\subsection{The solid-angle phase and the first ansatz}

Let $\Sigma_t$ be an oriented Seifert surface with
$\partial\Sigma_t=\Gamma_t$, where
\begin{equation*}
 \Gamma_t=\{\gamma(s,t):s\in[0,L]\}.
\end{equation*}
We orient $\Sigma_t$ so that its induced boundary orientation is opposite to
the parametrisation $s\mapsto\gamma(s,t)$.  This convention makes the
singular part of the phase below equal to $+\theta$ for the positively
oriented periodic frame.
For $x\notin\Sigma_t$, define
\begin{equation}\label{globalphase}
 \varphi(x,t)
 =\frac12\int_{\Sigma_t}
   \frac{(x-y)\cdot\vec n_{\Sigma_t}(y)}
        {|x-y|^3}\,dS_y.
\end{equation}
\begin{lemma}[Solid-angle phase]\label{solidanglephaselemma}
The phase $\varphi$ jumps by $2\pi$ across $\Sigma_t$, while
$e^{i\varphi}$ is single-valued and smooth on
$\mathbb R^3\setminus\Gamma_t$.  Its gradient is the Biot--Savart field
\begin{equation}\label{solidanglegradient}
 \nabla_x\varphi(x,t)
 =\frac12\int_0^L
 \frac{\partial_s\gamma(s',t)\times
       (x-\gamma(s',t))}
      {|x-\gamma(s',t)|^3}\,ds'.
\end{equation}
In particular, $\Delta_x\varphi=0$ away from $\Gamma_t$.
\end{lemma}

\begin{proof}
Put \(K_x(y)=(x-y)/|x-y|^3\).  Away from \(x\), one has
\(\operatorname{div}_yK_x=0\).  The surface Stokes identity
\[
 \int_{\partial\Sigma_t}K_x\times d\ell
 =-\int_{\Sigma_t}
   \big(\nabla_yK_x+(\operatorname{div}_yK_x)I\big)
   \vec n_{\Sigma_t}\,d\Sigma_t
\]
and \(\nabla_xK_x=-\nabla_yK_x\) therefore give, with the boundary
orientation fixed above,
\[
 \nabla_x\left(\int_{\Sigma_t}K_x\cdot\vec n_{\Sigma_t}\,d\Sigma_t\right)
 =\int_0^L
 \frac{\partial_s\gamma(s',t)\times(x-\gamma(s',t))}
        {|x-\gamma(s',t)|^3}\,ds'.
\]
Dividing by two proves \eqref{solidanglegradient}.

It remains to check the jump at an interior point of $\Sigma_t$, in a
ball disjoint from $\Gamma_t$.  Flattening the surface, the local
plane integral has one-sided limits $2\pi$ and $-2\pi$; the factor
$1/2$ in \eqref{globalphase} therefore gives a jump of $2\pi$ for
$\varphi$.  The Biot--Savart gradient already proved above is smooth
through this ball.  Subtracting the constant jump from one branch
therefore joins the branches to a smooth local phase, so
$e^{i\varphi}$ extends smoothly across the interior of $\Sigma_t$.
This proves smoothness and single-valuedness on
$\mathbb R^3\setminus\Gamma_t$.  Finally, for
fixed \(s'\), the integrand in \eqref{solidanglegradient} is the curl of
\[
 \frac{\partial_s\gamma(s',t)}{|x-\gamma(s',t)|}.
\]
It is consequently divergence-free away from \(\Gamma_t\).  Taking the
divergence in \eqref{solidanglegradient} therefore gives
\(\Delta_x\varphi=0\) there.  These classical solid-angle identities are
also discussed in
\cite{binysh2018,maxwell1873}.
\end{proof}

Write
\begin{equation*}
 W(y)=w(r)e^{i\theta},
 \qquad y=(r\cos\theta,r\sin\theta),
\end{equation*}
for the standard degree-one planar vortex.  Fix
$\delta>0$ smaller than one quarter of the tubular radius and take a smooth
cutoff $\eta=\eta(z)$ satisfying
\begin{equation*}
 \eta=1\quad\text{for }|z|\leq\delta,
 \qquad
 \eta=0\quad\text{for }|z|\geq2\delta.
\end{equation*}
As a first approximate solution to \eqref{scaledGPoperatorform}, we take, in the tubular
neighbourhood of $\Gamma_t$,
\begin{equation}\label{firstansatz}
 u_\epsilon^{(0)}(x,t)
 =\left[
   \eta(z)w\left(\frac{|z|}{\epsilon}\right)
   +(1-\eta(z))
  \right]e^{i\varphi(x,t)}.
\end{equation}
Away from the tubular neighbourhood we set
$u_\epsilon^{(0)}=e^{i\varphi}$.  Since $\eta=0$ for $|z|\geq2\delta$,
the two definitions agree wherever both apply and hence define a global
approximation.  Thus \eqref{firstansatz} describes the first approximation
near the filament, with the vortex profile cut off at a fixed physical
distance; in the inner variable $y=z/\epsilon$, the transition occurs at
$|y|\simeq\delta/\epsilon$.

\subsection{Exact expression for the first error}
\label{firsterrorsubsect}

Our next goal is to compute the error associated to \eqref{firstansatz}.  For
real-valued $f$ and $\phi$, direct differentiation gives
\begin{equation}\label{amplitudephaseidentity}
\begin{split}
 S_\epsilon(fe^{i\phi})
 =ife^{i\phi}\Bigg[&
 \epsilon^2|\log\epsilon|\frac{f_t}{f}
 +2\epsilon^2\frac{\nabla f}{f}\cdot\nabla\phi
 +\epsilon^2\Delta\phi\\
 &+i\left(
 \epsilon^2|\log\epsilon|\phi_t
 +\epsilon^2|\nabla\phi|^2
 -\epsilon^2\frac{\Delta f}{f}
 -(1-f^2)\right)\Bigg].
\end{split}
\end{equation}
Although \eqref{amplitudephaseidentity} is written where $f>0$, its
unfactored form extends across the zero set.

In the region $|z|\leq\delta$, put
\begin{equation*}
 y=\frac z\epsilon,\qquad r=|y|.
\end{equation*}
Using \eqref{laplaceperiodicframe}, \eqref{amplitudephaseidentity}, and
$\Delta_x\varphi=0$, we obtain
\begin{equation}\label{exactinnerfirsterror}
\begin{split}
 S_\epsilon(u_\epsilon^{(0)})
 =iu_\epsilon^{(0)}\Bigg[&
 \epsilon\frac{\nabla_yw}{w}\cdot
 \left(|\log\epsilon|\,\partial_tz
       +2\nabla_z\varphi\right)\\
 &+i\Bigg(
 \epsilon^2|\log\epsilon|\varphi_t
 +\epsilon^2|\nabla_z\varphi|^2
 +\frac{\epsilon^2}{h^2}(\mathscr D_s\varphi)^2
 +\epsilon\frac{\vec\kappa}{h}\cdot
       \frac{\nabla_yw}{w}\\
 &\hspace{7em}
 -\frac{\Delta_yw}{w}-(1-w^2)
 \Bigg)\Bigg].
\end{split}
\end{equation}
Here $\partial_tz$ is taken at fixed $x$, and
\begin{equation}\label{ztleading}
 \partial_tz
 =-(\gamma_t\cdot\vec e_1,\gamma_t\cdot\vec e_2)
 +O(|z|).
\end{equation}

The remaining regions produce lower-order errors.  Indeed, the vortex
asymptotics \(w(r)=1-\frac1{2r^2}+O(r^{-4})\) imply, in the transition annulus
\(\delta\leq |z|\leq2\delta\), that the cut-off amplitude
\[
 f_\epsilon(z):=\eta(z)w(|z|/\epsilon)+1-\eta(z)
\]
satisfies \(f_\epsilon-1=O(\epsilon^2)\), together with the same bound for
any fixed number of spatial and time derivatives.  Hence
\begin{equation*}
 S_\epsilon(u_\epsilon^{(0)})=O(\epsilon^2|\log\epsilon|)
 \qquad\text{for }\delta\leq|z|\leq2\delta.
\end{equation*}
For \(|z|\geq2\delta\), where \(u_\epsilon^{(0)}=e^{i\varphi}\), one has the
exact formula
\begin{equation}\label{outerfirsterror}
 S_\epsilon(u_\epsilon^{(0)})
 =-\epsilon^2e^{i\varphi}
 \left(|\log\epsilon|\,\varphi_t+|\nabla_x\varphi|^2\right).
\end{equation}
In particular, \eqref{outerfirsterror} is
\(O(\epsilon^2|\log\epsilon|)\) on compact sets separated from the filament.

\subsection{Local expansions of the phase}

We next record some useful expansions for terms appearing in
\eqref{exactinnerfirsterror}.  We recall that
\[
 z=|z|(\cos\theta,\sin\theta),\qquad
 \kappa\vec b=\vec{\mathsf t}\times\partial_{ss}\gamma,
 \qquad
 (\kappa_1,\kappa_2)=\kappa(\cos\alpha,\sin\alpha).
\]

\begin{lemma}[Local Biot--Savart expansion]
\label{localbiotsavartexpansionlemma}
Uniformly for \(|z|\) small,
\begin{equation}\label{nablaphase}
\begin{split}
 \nabla_x\varphi
 ={}&\frac{-z_2\vec e_1+z_1\vec e_2}{|z|^2}
 +\frac{\kappa}{2}
   \left(\log\frac{L}{|z|}-1\right)\vec b\\
 &+\frac{\kappa}{2}\cos(\theta-\alpha)
   \frac{-z_2\vec e_1+z_1\vec e_2}{|z|}
 +\frac12\mathcal N(\gamma)[s]
 +O(|z|\log|z|).
\end{split}
\end{equation}
Here \(|s-s'|_L:=\min_{q\in\mathbb Z}|s-s'+qL|\), and the finite-part
operator is
\begin{equation}\label{ngammadef}
\begin{split}
 \mathcal N(\gamma)[s]
 :={}&\int_0^L\bigg\{
 \frac{\partial_s\gamma(s')\times
       (\gamma(s)-\gamma(s'))}
      {|\gamma(s)-\gamma(s')|^3}
 -\frac{\kappa(s)\vec b(s)}{2|s-s'|_L}
 \bigg\}\,ds'.
\end{split}
\end{equation}
\end{lemma}

\begin{proof}
Write \(s'=s+h\), using the symmetric representative
\(h\in(-L/2,L/2]\), and split the integral in
\eqref{solidanglegradient} into \(|h|<h_0\) and \(|h|\geq h_0\).  In the local
part insert
\[
 \gamma(s+h)=\gamma(s)+h\vec{\mathsf t}
 +\frac{h^2}{2}\kappa\vec n+O(h^3),
 \qquad
 \partial_s\gamma(s+h)
 =\vec{\mathsf t}+h\kappa\vec n+O(h^2).
\]
Integrating the homogeneous terms in \((h,z)\) gives respectively the
polar field, the logarithmic binormal term, and the angular curvature
term in \eqref{nablaphase}.  After subtraction of
\(\kappa\vec b/(2|h|)\), the value at \(z=0\) of the remaining local
integral combines with the far integral to give precisely
\(\frac12\mathcal N(\gamma)[s]\).  Taylor's formula with integral
remainder gives \(O(|z|\log|z|)\), uniformly with derivatives on a
bounded smooth embedded family.  This is the standard local-induction
expansion; see \cite{callegariting1978,margeritbrancher2001} for the complete coefficient
calculation.
\end{proof}

The time derivative of the solid angle must be taken with the opposite
cross-product order from the one used for the spatial Biot--Savart field.
Transporting $\Sigma_t$ and applying Stokes' theorem gives
\begin{equation}\label{phasetexact}
 \partial_t\varphi(x,t)
 =\frac12\int_0^L
 \frac{
 \big(\partial_t\gamma(s',t)\times
      (x-\gamma(s',t))\big)
 \cdot\partial_s\gamma(s',t)}
 {|x-\gamma(s',t)|^3}\,ds'.
\end{equation}
\begin{lemma}[Local expansion of the time derivative of the phase]
\label{phasetlocallemma}
Uniformly for $|z|$ small,
\begin{equation}\label{phasetlocal}
 \partial_t\varphi(x,t)
 =-\nabla_z\theta\cdot\gamma_t^\perp(s,t)
 +O\big(1+|\log|z||\big).
\end{equation}
Here
\[
 \nabla_z\theta
 =\frac{-z_2\vec e_1+z_1\vec e_2}{|z|^2},
 \qquad
 \gamma_t^\perp
 =\Pi_{\vec{\mathsf t}}^\perp\gamma_t.
\]
The remainder is uniform on bounded smooth embedded families.
\end{lemma}

\begin{proof}
Fix $(s,t)$, write $s'=s+q$ with $q\in(-L/2,L/2]$, and set
\[
 \rho=|z|,\qquad
 \vec{\mathsf t}_0=\partial_s\gamma(s,t),\qquad
 V_0=\partial_t\gamma(s,t),\qquad
 K_0=\partial_{ss}\gamma(s,t).
\]
Split \eqref{phasetexact} into $|q|<q_0$ and $|q|\geq q_0$, where
$q_0$ is smaller than a fixed embeddedness scale.  The latter part is
$O(1)$.  In the local part, with
\[
 R(q):=x-\gamma(s+q),\qquad
 T(q):=\partial_s\gamma(s+q),\qquad
 V(q):=\partial_t\gamma(s+q),
\]
Taylor expansion gives
\begin{equation}\label{localcurvetaylorphaset}
\begin{split}
 R(q)&=z-q\vec{\mathsf t}_0-\frac12q^2K_0+O(|q|^3),\\
 T(q)&=\vec{\mathsf t}_0+qK_0+O(q^2),\\
 V(q)&=V_0+O(|q|).
\end{split}
\end{equation}
Consequently,
\begin{equation}\label{localdistancephaset}
 |R(q)|^2\asymp\rho^2+q^2,\qquad
 \left||R(q)|^2-(\rho^2+q^2)\right|
 \leq C(\rho q^2+q^4),
\end{equation}
and, writing
\[
 N(q,z)=\big(V(q)\times R(q)\big)\cdot T(q),\qquad
 N_0(z)=(V_0\times z)\cdot\vec{\mathsf t}_0,
\]
one has
\[
 N(q,z)=N_0(z)+O(\rho|q|+q^2).
\]
It follows, with $d=(\rho^2+q^2)^{1/2}$, that
\[
 \left|
 \frac{N(q,z)}{|R(q)|^3}
 -\frac{N_0(z)}{d^3}
 \right|
 \leq C\left[
 \frac{\rho|q|+q^2}{d^3}
 +\frac{\rho(\rho q^2+q^4)}{d^5}
 \right].
\]
The right-hand side has integral $O(1+|\log\rho|)$; the logarithm comes
only from $q^2d^{-3}$.  Hence
\[
 \partial_t\varphi(x,t)
 =
 \frac{N_0(z)}2\int_{-\infty}^{\infty}
 \frac{dq}{(\rho^2+q^2)^{3/2}}
 +O(1+|\log\rho|)
 =\frac{N_0(z)}{\rho^2}+O(1+|\log\rho|).
\]
Finally,
\[
 \frac{N_0(z)}{\rho^2}
 =-\frac{-z_2\vec e_1+z_1\vec e_2}{|z|^2}
   \cdot V_0^\perp
 =-\nabla_z\theta\cdot\gamma_t^\perp(s,t).
\]
This proves \eqref{phasetlocal}.
\end{proof}

\subsection{The complete order-\texorpdfstring{\(\epsilon\)}{epsilon} inner error}

At the inner scale $z=\epsilon y$, choose the branch of $\varphi$ compatible
with the polar angle $\theta$ and define
\begin{equation}\label{nablaphistar}
\begin{split}
 \varphi^*(s,z,t)
 &:={}
   \varphi\big(\gamma(s,t)+z_\alpha\vec e_\alpha(s,t),t\big)-\theta(z),\\
 \Phi_\epsilon(s,y,t)&:=\varphi^*(s,\epsilon y,t).
\end{split}
\end{equation}
Here and in what follows, repeated normal indices $\alpha$ are summed over
$\alpha=1,2$.  We also write
\[
 e_r=(\cos\theta,\sin\theta),\qquad
 e_\theta=(-\sin\theta,\cos\theta),\qquad
 \nabla_y\theta=\frac1r e_\theta.
\]
At $z=\epsilon y$ the normal-coordinate gradient satisfies
\[
 \nabla_z\varphi
 =\epsilon^{-1}\nabla_y\theta+\nabla_z\varphi^*.
\]
The normal part of the expansion in
Lemma~\ref{localbiotsavartexpansionlemma} therefore gives
\begin{equation}\label{phistargradientexpansion}
\begin{split}
 2\nabla_z\varphi^*
 ={}&\kappa\left(\log\frac{L}{\epsilon r}-1\right)\vec b
 +\kappa\cos(\theta-\alpha)e_\theta
 +\Pi_{\vec{\mathsf t}}^\perp\mathcal N(\gamma)[s]\\
 &+O(\epsilon r|\log(\epsilon r)|).
\end{split}
\end{equation}

At the same scale, Lemma~\ref{phasetlocallemma} gives
\begin{equation}\label{phasetinnerscaled}
 \epsilon^2\varphi_t
 =-\epsilon\nabla_y\theta\cdot\gamma_t^\perp
 +O\big(\epsilon^2(1+|\log(\epsilon r)|)\big).
\end{equation}

The vortex equation
\begin{equation*}
 w''+\frac1r w'-\frac1{r^2}w+(1-w^2)w=0
\end{equation*}
cancels all order-one terms in \eqref{exactinnerfirsterror}.  Define
\begin{equation}\label{Vepsilonfirsterror}
 V_\epsilon(s,y,t)
 :=-|\log\epsilon|\,\gamma_t^\perp(s,t)
   +2\nabla_z\varphi^*(s,\epsilon y,t).
\end{equation}
Combining \eqref{ztleading}, \eqref{phasetinnerscaled} and
\eqref{exactinnerfirsterror}, we find
\begin{equation}\label{epsilonerror}
\begin{split}
 S_\epsilon(u_\epsilon^{(0)})
 =\epsilon e^{i\Phi_\epsilon}H_\epsilon^{(1)}
 +\text{terms of order at least }\epsilon^2,
\end{split}
\end{equation}
where
\begin{equation}\label{firstinnersource}
\begin{split}
 H_\epsilon^{(1)}
 =iW\Bigg[
 &\frac{w'}w e_r\cdot V_\epsilon\\
 &+i\left\{
 \frac1r e_\theta\cdot V_\epsilon
 +\frac{w'}w\vec\kappa\cdot e_r
 \right\}\Bigg].
\end{split}
\end{equation}
To improve the approximation near the filament, we intend to invert the
planar linearised Ginzburg--Landau operator with leading source
$H_\epsilon^{(1)}$.  The two
translation modes
\[
 Z_j=\partial_{y_j}W,\qquad j=1,2,
\]
are obstructions to this inversion.  We shall therefore require
$H_\epsilon^{(1)}$ to be orthogonal to $Z_1$ and $Z_2$ on a suitable disc.
The precise truncated solvability conditions are stated and computed in
\S\ref{translation-solvability-subsection} below.

Before imposing these conditions, we make the remainder in
\eqref{epsilonerror} precise in the form needed for the solvability argument.
Write
\[
 \Lambda_\epsilon=|\log\epsilon|,
\]
and let $D_s=\partial_s+\mathfrak aJ$ and
$D_t=\partial_t+\mathfrak cJ$ act on pairs of translation coefficients.
We shall also use the following equivalent form of
\eqref{firstinnersource}:
\begin{equation}\label{firstinnersourcegeometric}
 H_\epsilon^{(1)}
 =iV_\epsilon\cdot\nabla_yW
 -e^{i\theta}w'(r)\vec\kappa\cdot e_r.
\end{equation}

\begin{lemma}[Uniform remainder in the translation projection]
\label{first-error-projection-remainder-lemma}
Fix an integer $q\geq0$.  Let $\gamma$ range over a family of smooth
embedded curves which is bounded in $C^{q+5}$ in $(s,t)$, whose lengths lie
in a compact subinterval of $(0,\infty)$, and which has a common chord--arc
constant and tubular radius.  Choose the periodic normal frames with
connection coefficients uniformly bounded in $C^{q+4}$, and fix $\delta$
below the common tubular radius.  For $R_\epsilon=\delta/\epsilon$, define
\begin{equation*}
 \mathcal P_{\epsilon,j}^{\rm full}(s,t)
 :=\operatorname{Re}\int_{B_{R_\epsilon}}
 e^{-i\Phi_\epsilon(s,y,t)}
 S_\epsilon(u_\epsilon^{(0)})
 \big(\gamma(s,t)+\epsilon y_\alpha\vec e_\alpha(s,t),t\big)
 \overline{Z_j(y)}\,dy.
\end{equation*}
Then
\begin{equation*}
 \mathcal P_{\epsilon,j}^{\rm full}
 =\epsilon\operatorname{Re}\int_{B_{R_\epsilon}}
 H_\epsilon^{(1)}\overline{Z_j}\,dy
 +\mathfrak R_{\epsilon,j},
 \qquad j=1,2,
\end{equation*}
where, uniformly for $0<\epsilon<\epsilon_q$,
\begin{equation}\label{uniform-first-error-projection-bound}
 \sup_{a+b\leq q}\sup_{s,t}
 \left|D_s^aD_t^b
   (\mathfrak R_{\epsilon,1},\mathfrak R_{\epsilon,2})\right|
 \leq C_q\epsilon^2\Lambda_\epsilon.
\end{equation}
In particular, for $m+b\leq q$,
\begin{equation}\label{uniform-first-error-projection-Sobolev-bound}
 \sup_{0\leq t\leq T}
 \left\|D_t^b
   (\mathfrak R_{\epsilon,1},\mathfrak R_{\epsilon,2})(t)
 \right\|_{H_s^m}
 \leq C_q\epsilon^2\Lambda_\epsilon.
\end{equation}
Equivalently, after division by the leading factor $\epsilon$, the error in
the two compatibility conditions is $O_{C^q}(\epsilon\Lambda_\epsilon)$.
The constants $C_q$ and $\epsilon_q$ depend only on $q$, $\delta$, and the
uniform geometric and frame bounds just specified.
\end{lemma}

\begin{proof}
All calculations take place in $|z|<\delta$, where the cut-off in
\eqref{firstansatz} is identically one.  Put $z=\epsilon y$, $\rho=|z|$,
and introduce
\begin{equation*}
\begin{split}
 \mathfrak z(s,z,t)
 &:=\partial_tz\big|_x+\gamma_t^\perp(s,t),\\
 p_t(s,z,t)
 &:=\varphi_t\big(\gamma(s,t)+z_\alpha\vec e_\alpha(s,t),t\big)
    +\nabla_z\theta\cdot\gamma_t^\perp(s,t).
\end{split}
\end{equation*}
Thus \eqref{ztleading} says, with the same estimate after $q$ derivatives,
\begin{equation}\label{frame-motion-remainder-bound}
 |\mathfrak z(s,z,t)|\leq C_q\rho.
\end{equation}
Expanding the two squares in \eqref{exactinnerfirsterror}, using the vortex
equation, and subtracting \eqref{firstinnersourcegeometric} gives the exact
identity
\begin{equation}\label{exact-first-error-remainder-identity}
\begin{split}
 \mathcal E_\epsilon
 :={}&e^{-i\Phi_\epsilon}S_\epsilon(u_\epsilon^{(0)})
       \big(\gamma+\epsilon y_\alpha\vec e_\alpha,t\big)
       -\epsilon H_\epsilon^{(1)}\\
 ={}&i\epsilon\Lambda_\epsilon e^{i\theta}w'(r)
       e_r\cdot\mathfrak z(s,\epsilon y,t)\\
 &-\epsilon^2W(y)\left[
       \Lambda_\epsilon p_t(s,\epsilon y,t)
       +|\nabla_z\varphi^*(s,\epsilon y,t)|^2
       +h^{-2}|\mathscr D_s\varphi(s,\epsilon y,t)|^2
       \right]\\
 &-\epsilon e^{i\theta}w'(r)
       \left(\frac1h-1\right)\vec\kappa\cdot e_r.
\end{split}
\end{equation}

We next record the angular cancellations which prevent the large disc
$B_{R_\epsilon}$ from producing a loss.  A refinement of Lemma
\ref{phasetlocallemma} gives
\begin{equation}\label{refined-phase-time-expansion}
 p_t(s,\rho e_r,t)
 =a_t(s,t)\log\frac1\rho+b_t(s,\rho,\theta,t),
 \qquad
 a_t=-\frac12\gamma_t^\perp\cdot J_\gamma\vec K.
\end{equation}
Here $b_t$ is uniformly bounded for $0<\rho<\delta$, together with every
required covariant and angular derivative.  In particular, the only
unbounded coefficient is independent of $\theta$.  We include the verification.
Use the notation in \eqref{localcurvetaylorphaset} and set
$d=(\rho^2+q^2)^{1/2}$.  Taylor expansion through degree two gives
\begin{equation}\label{refined-time-numerator}
 N(q,z)=N_0(z)+qN_1(z)
 -\frac{q^2}{2}\gamma_t^\perp\cdot J_\gamma\vec K
 {}+O_q(\rho q^2+|q|^3),
 \qquad |N_1(z)|\leq C_q\rho.
\end{equation}
Indeed, the part of the numerator which is independent of $z$ and
quadratic in $q$ is
\[
 \big(\gamma_t\times(-\tfrac12q^2\vec K)\big)
       \cdot\vec{\mathsf t}
 +\big(\gamma_t\times(-q\vec{\mathsf t})\big)
       \cdot(q\vec K)
 =-\frac{q^2}{2}\gamma_t^\perp\cdot J_\gamma\vec K.
\]
The distance expansion in \eqref{localdistancephaset}, one order further,
shows that replacing $|R(q)|^{-3}$ by $d^{-3}$ in the first three terms of
\eqref{refined-time-numerator} produces a uniformly integrable remainder.
The term $qN_1(z)d^{-3}$ is odd and integrates to zero on
$(-q_0,q_0)$.  The remainder in \eqref{refined-time-numerator} satisfies
\[
 \int_{-q_0}^{q_0}
 \frac{\rho q^2+|q|^3}{(\rho^2+q^2)^{3/2}}\,dq\leq C_q,
\]
while
\[
 \frac12\int_{-q_0}^{q_0}
 \frac{q^2}{(\rho^2+q^2)^{3/2}}\,dq
 =\log\frac1\rho+O(1).
\]
After the polar term $N_0(z)/\rho^2$ is cancelled by the second term in
the definition of $p_t$, the part with $|q|\geq q_0$ is smooth in $z$.
This proves \eqref{refined-phase-time-expansion}.  Differentiating the
subtracted kernels proves the same bounds after the stated derivatives;
the estimates above are unchanged because angular differentiation acts
only on factors of $z=\rho e_r$.
Consequently the logarithmic term has zero first angular moment and
\begin{equation}\label{phase-time-first-angular-moment}
 \left|\int_0^{2\pi}p_t(s,\rho e_r,t)e_r\,d\theta\right|
 \leq C_q.
\end{equation}

Similarly, \eqref{nablaphase}--\eqref{phistargradientexpansion} give
\begin{equation*}
 \nabla_z\varphi^*(s,\rho e_r,t)
 =A(s,t)\log\frac1\rho+B(s,\theta,t)
   +O_q\big(\rho(1+|\log\rho|)\big),
\end{equation*}
where $A$ is independent of $\theta$ and $B$ is uniformly bounded with all
the required derivatives.  The term
$|A|^2\log^2(1/\rho)$ in the square is therefore angularly constant and has
zero projection against $e_r$.  The cross term has only one logarithm.
Moreover, by \eqref{gradientperiodicframe},
$\mathscr D_s\varphi=h\vec{\mathsf t}\cdot\nabla_x\varphi$; all the
singular terms displayed in \eqref{nablaphase} are normal, so this quantity
has an angularly constant leading term and an
$O_q(\rho(1+|\log\rho|))$ remainder.  Since
$h^{-2}=1+O_q(\rho)$, we obtain
\begin{equation}\label{quadratic-phase-first-angular-moment}
\begin{split}
 \left|\int_0^{2\pi}\left(
 |\nabla_z\varphi^*|^2
 +h^{-2}|\mathscr D_s\varphi|^2\right)e_r\,d\theta\right|
 \leq C_q(1+|\log\rho|).
\end{split}
\end{equation}
Estimates \eqref{phase-time-first-angular-moment} and
\eqref{quadratic-phase-first-angular-moment} remain valid after applying
$D_s^aD_t^b$, $a+b\leq q$.

We now project \eqref{exact-first-error-remainder-identity}.  Direct
differentiation of $W=w(r)e^{i\theta}$ gives the two identities
\begin{equation}\label{translation-projection-pointwise-identities}
\begin{split}
 \operatorname{Re}\big(W\overline{Z_j}\big)
 &=ww'e_{r,j},\\
 \operatorname{Re}\big(ie^{i\theta}w'c\,\overline{Z_j}\big)
 &=\frac{ww'}r c\,e_{\theta,j}
 \qquad(c\in\mathbb R).
\end{split}
\end{equation}
For the last term in \eqref{exact-first-error-remainder-identity}, the real
part is $w'^2c e_{r,j}$.  The standard vortex asymptotics imply
\begin{equation}\label{radial-integrals-first-error-projection}
\begin{split}
 &\int_0^{R_\epsilon}|ww'|r
       (1+|\log(\epsilon r)|)\,dr\leq C\Lambda_\epsilon,\\
 &\int_0^{R_\epsilon}|ww'|r\,dr
   +\int_0^{R_\epsilon}(w')^2r^2\,dr\leq C.
\end{split}
\end{equation}
Indeed, $w(r)=O(r)$ and $w'(r)=O(1)$ at the origin, while
$w'(r)=O(r^{-3})$ at infinity; splitting at $r=1$ proves both bounds.

The first line of \eqref{exact-first-error-remainder-identity}, together
with \eqref{frame-motion-remainder-bound} and the second identity in
\eqref{translation-projection-pointwise-identities}, contributes at most
\[
 C\epsilon\Lambda_\epsilon
 \int_0^{R_\epsilon}|ww'|\epsilon r\,dr
 \leq C\epsilon^2\Lambda_\epsilon.
\]
The $p_t$ term is bounded, using
\eqref{phase-time-first-angular-moment}, by
\[
 C\epsilon^2\Lambda_\epsilon
 \int_0^{R_\epsilon}|ww'|r\,dr
 \leq C\epsilon^2\Lambda_\epsilon.
\]
The two quadratic phase terms contribute at most
$C\epsilon^2\Lambda_\epsilon$ by
\eqref{quadratic-phase-first-angular-moment} and the first line of
\eqref{radial-integrals-first-error-projection}.  Finally,
$|h^{-1}-1|\leq C\epsilon r$ and the last radial integral in
\eqref{radial-integrals-first-error-projection} bound the curvature
correction by $C\epsilon^2$.

These estimates prove \eqref{uniform-first-error-projection-bound} when
$q=0$.  Applying $D_s^aD_t^b$ does not change any radial weight, and all
coefficient and angular-moment estimates above were stated in differentiated
form.  The same calculation therefore proves
\eqref{uniform-first-error-projection-bound} for $a+b\leq q$.
The Sobolev estimate
\eqref{uniform-first-error-projection-Sobolev-bound} follows by integrating
the differentiated pointwise bounds in $s$.
\end{proof}

\subsection{Translation solvability conditions}
\label{translation-solvability-subsection}

Set
\begin{equation*}
 R_\epsilon=\frac{\delta}{\epsilon}.
\end{equation*}
The two solvability conditions for inverting the planar linearised operator
are
\begin{equation}\label{translationorthogonalitysectiontwo}
 \operatorname{Re}\int_{B_{R_\epsilon}}
 H_\epsilon^{(1)}
 \overline{\partial_{y_j}W}\,dy=0,
 \qquad j=1,2.
\end{equation}
No condition against $iW$ is required for the weighted inverse used later.

We now compute \eqref{translationorthogonalitysectiontwo}.  For any normal
vector field $a=a(r,\theta)$,
\begin{equation*}
\begin{split}
 &\begin{pmatrix}
 \langle ia\cdot\nabla W,\partial_{y_1}W\rangle_{B_R}\\
 \langle ia\cdot\nabla W,\partial_{y_2}W\rangle_{B_R}
 \end{pmatrix}
 =J\int_0^R\int_0^{2\pi}ww'a(r,\theta)\,d\theta\,dr,
\end{split}
\end{equation*}
where $\langle f,g\rangle=\operatorname{Re}\int f\overline g$.  Also,
\begin{equation}\label{curvatureprojectionidentity}
\begin{pmatrix}
 \langle-e^{i\theta}w'\vec\kappa\cdot e_r,
         \partial_{y_1}W\rangle_{B_R}\\
 \langle-e^{i\theta}w'\vec\kappa\cdot e_r,
         \partial_{y_2}W\rangle_{B_R}
\end{pmatrix}
=-\pi\vec\kappa\int_0^R(w')^2r\,dr.
\end{equation}
Finally,
\begin{equation}\label{angularaverageidentity}
 \int_0^{2\pi}\cos(\theta-\alpha)e_\theta\,d\theta
 =\pi\frac{J\vec\kappa}{\kappa}.
\end{equation}

Define the vortex-core constant
\begin{equation*}
 \mathfrak c_\epsilon
 :=-\frac{2}{w^2(R_\epsilon)}
 \int_0^{R_\epsilon}ww'\log r\,dr
 +\frac{1}{w^2(R_\epsilon)}
 \int_0^{R_\epsilon}(w')^2r\,dr
 -\frac12.
\end{equation*}
Equations \eqref{phistargradientexpansion}, \eqref{Vepsilonfirsterror},
\eqref{firstinnersourcegeometric} and
\eqref{translationorthogonalitysectiontwo}--\eqref{angularaverageidentity}
show that the order-$\epsilon$ translation conditions are equivalent to
\begin{equation}\label{firstgammaeq}
 -|\log\epsilon|\,\gamma_t^\perp
 +c_\epsilon(L)\,\partial_s\gamma\times\partial_{ss}\gamma
 +\Pi_{\vec{\mathsf t}}^\perp\mathcal N(\gamma)=0,
\end{equation}
where
\begin{equation*}
 c_\epsilon(L)
 :=\log\left(\frac{L}{\epsilon}\right)+\mathfrak c_\epsilon.
\end{equation*}
In \(\mathfrak c_\epsilon\), the logarithmic moment comes from
\eqref{phistargradientexpansion}, the term $-\frac12$ is the sum of the
$-1$ in the local-induction expansion and the $+\frac12$ contributed by
the angular term there, and the integral involving $(w')^2r$ comes from
\eqref{curvatureprojectionidentity}.  If the full, rather than truncated,
phase expansion is used in the projection, the left-hand side of
\eqref{firstgammaeq} has an additional $O(\epsilon|\log\epsilon|)$ normal
remainder.  This affects only the next improvement of the curve.
In addition, Lemma \ref{first-error-projection-remainder-lemma} shows that
replacing the coefficient $\epsilon H_\epsilon^{(1)}$ by the complete
residual in the compatibility conditions changes the left-hand side, after
division by $\epsilon$, by another
$O_{C^q}(\epsilon|\log\epsilon|)$ normal term for every fixed $q$.

For fixed $\delta$,
\begin{equation}\label{coreconstantasymptotics}
 \mathfrak c_\epsilon
 =\mathfrak c_0+O(\epsilon^2|\log\epsilon|),
 \qquad
 \mathfrak c_0
 :=-2\int_0^\infty ww'\log r\,dr
   +\int_0^\infty(w')^2r\,dr-\frac12.
\end{equation}
Consequently,
\begin{equation*}
 \frac{c_\epsilon(L)}{|\log\epsilon|}
 =1+\frac{\log L+\mathfrak c_0}
          {|\log\epsilon|}
 +o\left(\frac1{|\log\epsilon|}\right).
\end{equation*}
Solving \eqref{firstgammaeq} for the normal velocity therefore gives
\begin{equation}\label{binormalleadingfromorthogonality}
 \gamma_t^\perp
 =\partial_s\gamma\times\partial_{ss}\gamma
 +O\left(\frac1{|\log\epsilon|}\right).
\end{equation}
Formally letting $\epsilon\to0$, one recovers the binormal flow
\eqref{bcfintrosection}, with the orientation fixed by \eqref{scaledGP}.
The remainder in \eqref{binormalleadingfromorthogonality}, however, is not
lower order in derivatives: it contains the nonlocal correction
\[
 \frac{1}{|\log\epsilon|}
 \Pi_{\vec{\mathsf t}}^\perp\mathcal N(\gamma),
\]
whose linearisation is of logarithmic order two and therefore marginally
higher order than the binormal term.  Issues related to this and the
solvability of \eqref{firstgammaeq} are discussed in the next section.

\begin{remark}[Time-dependent length]
For notational simplicity, the calculations in this section were written for
arclength-parametrised curves of constant length \(L\).  For curves of
time-dependent length \(\ell(t)\), one instead interprets all spatial
derivatives and periodic distances in intrinsic arclength, replaces \(L\) by
\(\ell(t)\), and in particular uses \(c_\epsilon(\ell(t))\) and
\(\mathfrak a=-\Theta/\ell(t)\).  To work on a fixed parameter circle, one
writes \(|\partial_s\gamma|=\ell(t)/L\) and introduces the tangential velocity
that preserves this constant-speed parametrisation.  Since the solid-angle and
Biot--Savart formulae are reparametrisation invariant and the tangential
velocity does not contribute to \eqref{phasetexact}, the local expansions and
normal translation solvability conditions are unchanged.
\end{remark}

\section{The curve equation in a normal graph gauge}\label{normalgraphsection}

The goal of this section is to solve the corrected curve equation
\eqref{firstgammaeq}.  Writing the unknown curve as a graph over the
prescribed binormal flow $\gamma^0$ reduces the problem to a nonlocal
quasilinear Schr\"odinger equation whose principal part is skew-adjoint and
admits uniform tame estimates with a fixed loss of derivatives.  A
Nash--Moser iteration then yields a smooth solution which remains
$O(|\log\epsilon|^{-1})$-close to $\gamma^0$ in Sobolev norms.

\subsection{The corrected curve equation and statement of the result}

We first record intrinsically the curve equation to be solved.  Unlike
binormal flow, \eqref{firstgammaeq} does not in general preserve the length of
the curve, so an arclength parametrisation cannot be maintained on the fixed
parameter circle $\mathbb R/L\mathbb Z$.  We therefore let
$\gamma=\gamma(s,t)$ have constant speed
$|\partial_s\gamma|=\lambda(t)$, where
$\ell(t)=L\lambda(t)$ is its time-dependent length, and choose the tangential
velocity so that this condition is propagated.  Writing $\sigma$ for
arclength, set $\vec{\mathsf t}=\partial_\sigma\gamma$ and
$\vec K=\partial_\sigma\vec{\mathsf t}$.  For $\vec v\in\mathbb R^3$, write
\[
 \Pi_{\vec{\mathsf t}}^\perp\vec v
 :=\vec v-(\vec{\mathsf t}\cdot\vec v)\vec{\mathsf t}.
\]
Set
\begin{equation*}
    c_\epsilon(\ell)
    :=\log\bigg(\frac{\ell}{\epsilon}\bigg)+\mathfrak c_\epsilon
\end{equation*}
and define the normal velocity
\begin{equation}\label{normalvelocitydefinition}
    \vec V_\epsilon[\gamma]
    :=\frac{1}{|\log\epsilon|}
    \left\{
    c_\epsilon(\ell)\,\vec{\mathsf t}\times\vec K
    +\Pi_{\vec{\mathsf t}}^\perp\mathcal N(\gamma)
    \right\}.
\end{equation}
Here $\mathcal N(\gamma)$ is parametrised intrinsically by
$\sigma\in\mathbb R/\ell(t)\mathbb Z$.  The two translation solvability
conditions in \eqref{firstgammaeq} prescribe the normal velocity, so
the full evolution has the form
\begin{equation}\label{geometriccurveequation}
    \partial_t\gamma=a\vec{\mathsf t}+\vec V_\epsilon[\gamma]
\end{equation}
for a tangential velocity $a=a(s,t)$.

To determine $a$, differentiate the constant-speed condition.  Since
$\partial_s\gamma=\lambda\vec{\mathsf t}$ and
$\vec V_\epsilon\cdot\vec{\mathsf t}=0$, equation
\eqref{geometriccurveequation} gives
\[
\begin{split}
 \partial_t\lambda
 &=\vec{\mathsf t}\cdot\partial_s
   (a\vec{\mathsf t}+\vec V_\epsilon)\\
 &=\partial_sa-\lambda\vec K\cdot\vec V_\epsilon.
\end{split}
\]
The left-hand side is independent of $s$.  Periodicity of $a$ therefore
gives
\begin{equation}\label{lambdatimeequation}
    \partial_t\lambda
    =-\frac{\lambda}{L}\int_0^L
      \vec K\cdot\vec V_\epsilon\,ds,
\end{equation}
and the exact constant-speed gauge is
\begin{equation}\label{atangentialequation}
    \partial_sa
    =\lambda\left(
    \vec K\cdot\vec V_\epsilon
    -\frac1L\int_0^L\vec K\cdot\vec V_\epsilon\,ds
    \right).
\end{equation}
We fix the remaining freedom to add a function of time by imposing
\begin{equation}\label{ameanzero}
    \int_0^L a(s,t)\,ds=0.
\end{equation}
Equivalently, the length satisfies
\begin{equation*}
    \frac{d\ell}{dt}
    =-\int_0^{\ell(t)}\vec K\cdot\vec V_\epsilon\,d\sigma.
\end{equation*}
Thus the mean-zero part of $\vec K\cdot\vec V_\epsilon$ is absorbed by the
tangential reparametrisation, while its mean gives the length variation
anticipated above.  In particular, the binormal term makes no contribution
because $\vec K\cdot(\vec{\mathsf t}\times\vec K)=0$, whereas the nonlocal
term need not integrate to zero.

We will solve \eqref{geometriccurveequation} by writing $\gamma$ as a graph
over the prescribed binormal flow.  Thus, let $\gamma^0$ satisfy
\begin{equation}\label{backgroundbf}
    \partial_t\gamma^0=\partial_s\gamma^0\times\partial_{ss}\gamma^0,
    \qquad |\partial_s\gamma^0|=1,
\end{equation}
where $s\in\mathbb R/L\mathbb Z$.  Let
$(\vec{\mathsf t}^{\,0},\vec e_1^{\,0},\vec e_2^{\,0})$ be a smooth periodic
oriented orthonormal frame along $\gamma^0$.  It may be chosen so that the
spatial normal connection $\mathfrak a^0$ is constant in both $s$ and $t$.
Indeed, the holonomy construction makes it constant in $s$, while the frame
compatibility equations for binormal flow give
$\partial_t\mathfrak a^0=0$.  Thus
\begin{equation}\label{backgroundperiodicframe}
\begin{split}
    \partial_s\vec{\mathsf t}^{\,0}
        &=k_1^0\vec e_1^{\,0}+k_2^0\vec e_2^{\,0},\\
    \partial_s\vec e_1^{\,0}&=-k_1^0\vec{\mathsf t}^{\,0}
       +\mathfrak a^0\vec e_2^{\,0},\\
    \partial_s\vec e_2^{\,0}&=-k_2^0\vec{\mathsf t}^{\,0}
       -\mathfrak a^0\vec e_1^{\,0}.
\end{split}
\end{equation}
We use the abbreviation
\begin{equation*}
 D_s^0=\partial_s+i\mathfrak a^0,
 \qquad k^0=k_1^0+ik_2^0.
\end{equation*}
We seek a solution of the orthogonality equations in the form
\begin{equation}\label{normalgraphansatz}
    \gamma=\gamma^0+f_3\vec{\mathsf t}^{\,0}
    +f_1\vec e_1^{\,0}+f_2\vec e_2^{\,0}.
\end{equation}
Here all the functions depend on $(s,t)$ and are $L$-periodic in $s$.

The result proved in this section is the following.

\begin{theorem}[Existence on the prescribed interval]
\label{uniformlocalcurvetheorem}
Let
\[
 \gamma^0:\mathbb R/L\mathbb Z\times[0,T]\longrightarrow\mathbb R^3
\]
be a smooth embedded solution of \eqref{backgroundbf}, parametrised by
arclength, and assume that it has a uniform tubular radius on $[0,T]$.
There is $\epsilon_0=\epsilon_0(\gamma^0,T)\in(0,e^{-2})$ such that for
every $0<\epsilon\leq\epsilon_0$ the equation
\eqref{geometriccurveequation}, with
\eqref{normalvelocitydefinition}, \eqref{atangentialequation} and
\eqref{ameanzero}, has a unique smooth solution on the full interval
$[0,T]$ satisfying
\begin{equation*}
    \gamma_\epsilon(s,0)=\gamma^0(s,0).
\end{equation*}
The curves $\gamma_\epsilon(\cdot,t)$ are closed, embedded and have constant
speed $\lambda_\epsilon(t)$.  For every integer $m\geq6$,
\begin{equation}\label{uniformgraphestimate}
\begin{split}
    \sup_{0\leq t\leq T}
    \left(
       \|f_\epsilon(t)\|_{H^m}
       +\|f_{3,\epsilon}(t)\|_{H^m}
       +|\lambda_\epsilon(t)-1|
    \right)
    \leq\frac{C_{m,T}}{|\log\epsilon|}.
\end{split}
\end{equation}
Here $f_\epsilon=f_{1,\epsilon}+if_{2,\epsilon}$ and
$f_{3,\epsilon}$ are the coefficients in \eqref{normalgraphansatz}.
They solve \eqref{exactcomplexfequation}, \eqref{exactf123system} and the
exact constant-speed constraint \eqref{f3spatialexact}, with
\begin{equation}\label{zerographinitialdata}
    f_\epsilon(0)=0,\qquad f_{3,\epsilon}(0)=0,
    \qquad\lambda_\epsilon(0)=1.
\end{equation}
\end{theorem}

\subsection{Reduction to the fixed normal graph variables}

We now reduce \eqref{geometriccurveequation} to an evolution equation for the
graph variables introduced above and identify the constraints used in later
linearised estimates.

\subsubsection{The constant-speed constraint}

Differentiating \eqref{normalgraphansatz} with respect to $s$ and using
\eqref{backgroundperiodicframe}, we find
\begin{equation*}
    \partial_s\gamma
    =p_0\vec{\mathsf t}^{\,0}+p_1\vec e_1^{\,0}+p_2\vec e_2^{\,0},
\end{equation*}
where
\begin{equation}\label{pidefinition}
\begin{split}
    p_0&=1+\partial_sf_3-k_1^0f_1-k_2^0f_2,\\
    p_1&=\partial_sf_1-\mathfrak a^0f_2+k_1^0f_3,\\
    p_2&=\partial_sf_2+\mathfrak a^0f_1+k_2^0f_3.
\end{split}
\end{equation}
Equivalently, $p_1+ip_2=D_s^0f+k^0f_3$, where $f=f_1+if_2$.
The constant-speed condition introduced above is
\begin{equation*}
    |\partial_s\gamma(s,t)|=\lambda(t),
\end{equation*}
and hence
\begin{equation*}
    \operatorname{Length}(\gamma(\cdot,t))
    =\int_0^L|\partial_s\gamma|\,ds=L\lambda(t)=\ell(t).
\end{equation*}
Consequently,
\begin{equation}\label{constantconstraintpi}
    p_0^2+p_1^2+p_2^2=\lambda^2.
\end{equation}
For a sufficiently small perturbation we have $p_0>0$, and hence this
constraint is equivalently written as the equation
\begin{equation}\label{f3spatialexact}
\begin{split}
\partial_sf_3={}&k_1^0f_1+k_2^0f_2-1\\
&+\sqrt{\lambda^2-
    (\partial_sf_1-\mathfrak a^0f_2+k_1^0f_3)^2-
    (\partial_sf_2+\mathfrak a^0f_1+k_2^0f_3)^2}.
\end{split}
\end{equation}
In particular, at the linearised level about $f_1=f_2=f_3=0$ and
$\lambda=1$, equation \eqref{f3spatialexact} reduces to
\begin{equation}\label{f3linearconstraintnew}
    \partial_sf_3=k_1^0f_1+k_2^0f_2+(\lambda-1).
\end{equation}
The periodicity condition for $f_3$ determines the constant mode of
$\lambda-1$ at this order.  In the exact equation the same role is played by
the integral over one period of the right-hand side of
\eqref{f3spatialexact}.

The unit tangent and the intrinsic arclength derivative of the perturbed
curve are
\begin{equation}\label{unknownunitangent}
    \vec{\mathsf t}
    =\frac{1}{\lambda}
    \big(p_0\vec{\mathsf t}^{\,0}+p_1\vec e_1^{\,0}+p_2\vec e_2^{\,0}\big),
    \qquad
    \partial_\sigma=\frac{1}{\lambda}\partial_s,
\end{equation}
where $\sigma$ denotes arclength on the unknown curve.  We write
\begin{equation*}
    \vec K:=\partial_\sigma\vec{\mathsf t}
    =\kappa_1\vec e_1+\kappa_2\vec e_2
\end{equation*}
for its curvature vector; here $(\vec e_1,\vec e_2)$ is any smooth periodic
oriented normal frame of the unknown curve.

\subsubsection{Recovery of the velocity in graph coordinates}

The formulae above recover $\ell$, $\vec{\mathsf t}$ and $\vec K$ from the
graph variables.  Substitution in \eqref{normalvelocitydefinition} gives the
normal velocity, while \eqref{atangentialequation}--\eqref{ameanzero}
recover the unique tangential velocity in the chosen gauge.  Thus every
coefficient in \eqref{geometriccurveequation} is now an explicit geometric
functional of $(f_1,f_2,f_3,\lambda-1)$.

\subsubsection{Graph evolution}

For the background binormal-flow solution, put
\begin{equation*}
 q_1^0=\partial_sk_1^0-\mathfrak a^0k_2^0,
 \qquad
 q_2^0=\partial_sk_2^0+\mathfrak a^0k_1^0,
\end{equation*}
so that $D_s^0k^0=q_1^0+iq_2^0$.  The time derivatives of the periodic
frame have the form
\begin{equation}\label{backgroundframetime}
\begin{split}
    \partial_t\vec{\mathsf t}^{\,0}
       &=-q_2^0\vec e_1^{\,0}+q_1^0\vec e_2^{\,0},\\
    \partial_t\vec e_1^{\,0}
       &=q_2^0\vec{\mathsf t}^{\,0}
         +\omega^0\vec e_2^{\,0},\\
    \partial_t\vec e_2^{\,0}
       &=-q_1^0\vec{\mathsf t}^{\,0}
         -\omega^0\vec e_1^{\,0},
\end{split}
\end{equation}
where
\begin{equation}\label{omegabackground}
    \omega^0(s,t)
    =-\frac12\big((k_1^0)^2+(k_2^0)^2\big)+\omega_0(t).
\end{equation}
The compatibility relation in the normal connection is
$\partial_t\mathfrak a^0-\partial_s\omega^0
=\frac12\partial_s|k^0|^2$.  Since $\mathfrak a^0$ is constant in $s$,
periodicity gives $\partial_t\mathfrak a^0=0$ and then
\eqref{omegabackground}.
Differentiating \eqref{normalgraphansatz} in time and using
\eqref{backgroundbf} and \eqref{backgroundframetime}, we obtain
\begin{equation}\label{graphtimederivativecomponents}
\begin{split}
\partial_t\gamma\cdot\vec{\mathsf t}^{\,0}
   &=\partial_tf_3+f_1q_2^0-f_2q_1^0,\\
\partial_t\gamma\cdot\vec e_1^{\,0}
   &=-k_2^0+\partial_tf_1-f_3q_2^0-\omega^0f_2,\\
\partial_t\gamma\cdot\vec e_2^{\,0}
   &=\phantom{-}k_1^0+\partial_tf_2+f_3q_1^0+\omega^0f_1.
\end{split}
\end{equation}
For brevity, let
\begin{equation*}
    \vec G_\epsilon[f_1,f_2,f_3]
    :=a\vec{\mathsf t}+\vec V_\epsilon[\gamma],
\end{equation*}
where $\gamma$, $\vec{\mathsf t}$, $\vec K$, $\lambda$ and $a$ are determined
by \eqref{normalgraphansatz}, \eqref{pidefinition},
\eqref{constantconstraintpi}, \eqref{unknownunitangent},
\eqref{normalvelocitydefinition}, \eqref{lambdatimeequation},
\eqref{atangentialequation} and \eqref{ameanzero}.  Equating
\eqref{graphtimederivativecomponents} with the components of
\eqref{geometriccurveequation}, we arrive at the exact system
\begin{equation}\label{exactf123system}
\boxed{
\begin{aligned}
\partial_tf_1
   &=\vec G_\epsilon[f_1,f_2,f_3]\cdot\vec e_1^{\,0}
     +k_2^0+f_3q_2^0+\omega^0f_2,\\[3pt]
\partial_tf_2
   &=\vec G_\epsilon[f_1,f_2,f_3]\cdot\vec e_2^{\,0}
     -k_1^0-f_3q_1^0-\omega^0f_1,\\[3pt]
\partial_tf_3
   &=\vec G_\epsilon[f_1,f_2,f_3]\cdot\vec{\mathsf t}^{\,0}
     -f_1q_2^0+f_2q_1^0.
\end{aligned}}
\end{equation}
This system is supplemented by the exact constant-speed constraint
\eqref{f3spatialexact} and the evolution equation
\eqref{lambdatimeequation} for $\lambda$.

Writing
\begin{equation*}
    f=f_1+if_2,
\end{equation*}
the first two equations of \eqref{exactf123system} may equivalently be
expressed as
\begin{equation}\label{exactcomplexfequation}
\begin{split}
\partial_tf={}&\vec G_\epsilon[f_1,f_2,f_3]\cdot
       (\vec e_1^{\,0}+i\vec e_2^{\,0})
       -ik^0-if_3D_s^0k^0-i\omega^0f.
\end{split}
\end{equation}

Set
\begin{equation*}
 U=(f_1,f_2,f_3,\lambda-1).
\end{equation*}
The constraint \eqref{f3spatialexact}, the equations
\eqref{lambdatimeequation} and \eqref{atangentialequation}, and the system
\eqref{exactf123system} constitute a closed evolution system for $U$.
Consequently, the existence problem in Theorem
\ref{uniformlocalcurvetheorem} reduces to solving this system with
$U(0)=0$.  Conversely, every sufficiently small periodic solution $U$
recovers through \eqref{normalgraphansatz} a solution of
\eqref{geometriccurveequation} in the constant-speed gauge.  If the
initial data satisfy the constant-speed constraint, then the choice
\eqref{atangentialequation} propagates it.  Periodicity of the graph functions
also guarantees that the perturbed curve remains closed.  For sufficiently
small graph functions, embeddedness follows from the existence of a uniform
tubular neighbourhood of the reference curve.

\subsection{The linearised graph equation and its nonlocal operator}
\label{periodic-nonlocal-subsection}

\subsubsection{The constrained linearisation at the background curve}

To solve the graph system \eqref{exactf123system} supplemented by
\eqref{lambdatimeequation}, \eqref{f3spatialexact}, we need uniform tame
estimates for the corresponding linearised initial-value problems.  We first examine the
linearisation at the background curve $U=0$, where the principal structure
is most transparent.  Writing
\[
 H=(h_1,h_2,h_3,\mu),\qquad h=h_1+ih_2,
\]
and recalling that the background curve has length $L$, the linearised
constraint \eqref{f3linearconstraintnew} becomes
\[
 \partial_sh_3=k_1^0h_1+k_2^0h_2+\mu,
 \qquad
 \mu=-\frac1L\int_0^L(k_1^0h_1+k_2^0h_2)\,ds.
\]
Thus $\mu$ and the mean-zero part of $h_3$ are determined by the normal
variation $h$; the mean of $h_3$ is the remaining tangential coordinate.

At highest order, the normal curvature variation is $(D_s^0)^2h$.  In
addition, the principal nonlocal contribution arising from the finite-part
field \eqref{ngammadef} is (see \S\ref{nonlocalfieldstructure} below)
\[
 \frac{i}{2|\log\epsilon|}
 \mathcal T_L^{\mathfrak a^0}(D_s^0)^2h,
 \qquad
 \mathcal T_L^{\mathfrak a^0}[g](s)
 :=\int_{-L/2}^{L/2}
 \frac{e^{i\mathfrak a^0r}g(s+r)-g(s)}{|r|}\,dr.
\]
Consequently the complex normal equation in the linearised $U$-system
formally has the structure
\begin{equation}\label{background-linearised-principal-form}
\begin{split}
 i\partial_th+(D_s^0)^2h
 +\frac{1}{|\log\epsilon|}\left\{
 (\log L+\mathfrak c_\epsilon)(D_s^0)^2h
 +\frac12\mathcal T_L^{\mathfrak a^0}(D_s^0)^2h
 \right\}
 &=\mathscr E_\epsilon^0[H].
\end{split}
\end{equation}
Here $\mathscr E_\epsilon^0[H]$ contains the local terms of order at most
one, the remaining $|\log\epsilon|^{-1}$ terms of order at most two, and
the constrained couplings to $h_3$ and $\mu$.  

The Schr\"odinger structure in
\eqref{background-linearised-principal-form} is the graph-coordinate
counterpart of the classical Hasimoto representation.  Wherever the Frenet
frame is defined, Hasimoto's transformation
\[
 q(\sigma,t)=\kappa(\sigma,t)
 \exp\left(i\int_0^\sigma\tau(\sigma',t)\,d\sigma'\right)
\]
converts binormal flow into the focusing cubic nonlinear Schr\"odinger
equation, up to a time-dependent gauge \cite{hasimoto1972}.  Applied
formally to \eqref{geometriccurveequation}, the same transformation gives,
modulo tangential transport and a time-dependent gauge, an equation of the
schematic form
\[
\begin{aligned}
 i\partial_tq+\partial_{\sigma\sigma}q+\frac12|q|^2q
 +\frac1{|\log\epsilon|}\bigg\{&
 \big(\log\ell+\mathfrak c_\epsilon\big)
 \left(\partial_{\sigma\sigma}q+\frac12|q|^2q\right)\\
 &+\frac12\partial_{\sigma\sigma}\mathcal T_\ell[q]
 +\mathcal Q_\ell(q)\bigg\}=0.
\end{aligned}
\]
Here $\partial_{\sigma\sigma}\mathcal T_\ell$ is the Hasimoto-gauge
counterpart of the logarithmic operator
$\mathcal T_L^{\mathfrak a^0}(D_s^0)^2$ in
\eqref{background-linearised-principal-form}, and $\mathcal Q_\ell$ collects
quasilinear geometric terms.  Although the
correction carries the small factor $|\log\epsilon|^{-1}$, the multiplier
of $\partial_{\sigma\sigma}\mathcal T_\ell$ grows like
$|\xi|^2\log(2+|\xi|)$, so this is a logarithmic order-two perturbation of
cubic Schr\"odinger flow.

\subsubsection{The intrinsic operator and its flat model}
\label{nonlocalfieldstructure}

We now rigorously analyse the structure of the nonlocal field
$\mathcal N(\gamma)$ in preparation for the tame estimates used later.  We
isolate the logarithmic order-two term appearing formally in
\eqref{background-linearised-principal-form} and show, first in the flat
model, that the remaining field has classical order two.

Let $\ell$ be the length of a smooth closed embedded curve $\gamma$, write
$\sigma\in\mathbb R/\ell\mathbb Z$ for arclength, and use the symmetric
representative $h\in(-\ell/2,\ell/2]$.  Choose the periodic normal frame
with constant connection $\mathfrak a_\gamma=-\Theta_\gamma/\ell$ and set
\begin{equation*}
 D_\sigma=\partial_\sigma+i\mathfrak a_\gamma.
\end{equation*}
For brevity, throughout the remainder of this section we write
$\mathfrak a=\mathfrak a_\gamma$ in operator superscripts.

Let
$P_{\sigma+h\to\sigma}^\perp$ denote normal parallel transport along the
curve from $\sigma+h$ back to $\sigma$, and for a normal vector field $X$
along $\gamma$ define
\begin{equation*}
 \mathcal T_\ell^\perp X(\sigma)
 :=\int_{-\ell/2}^{\ell/2}
 \frac{P_{\sigma+h\to\sigma}^\perp X(\sigma+h)-X(\sigma)}{|h|}\,dh.
\end{equation*}
Define the normal remainder by
\begin{equation}\label{Ndecompositionsectiontwo}
 \widetilde{\mathcal N}(\gamma)
 :=\Pi_{\vec{\mathsf t}}^\perp\mathcal N(\gamma)
   -\frac12J_\gamma\mathcal T_\ell^\perp\vec K,
 \qquad J_\gamma X=\vec{\mathsf t}\times X.
\end{equation}

If $g=g_1+ig_2$ denotes the periodic-frame components of a normal field
$X$, then normal parallel transport is represented by multiplication by
$e^{i\mathfrak a h}$, and hence
\begin{equation}\label{periodicframeTsectiontwo}
 \mathcal T_{\ell,\mathfrak a}[g](\sigma)
 :=\int_{-\ell/2}^{\ell/2}
 \frac{e^{i\mathfrak a h}g(\sigma+h)-g(\sigma)}{|h|}\,dh
\end{equation}
represents $\mathcal T_\ell^\perp X$.  Thus the singular term in
\eqref{Ndecompositionsectiontwo} has complex component
$\frac i2\mathcal T_{\ell,\mathfrak a}[\kappa_1+i\kappa_2]$.
We henceforth denote the component operator in
\eqref{periodicframeTsectiontwo} by $\mathcal T_\ell^{\mathfrak a}$.

The operator $\mathcal T_\ell^{\mathfrak a}$ is self-adjoint and commutes
with $D_\sigma$.  On the Fourier mode
$e^{i\xi_k\sigma}$, $\xi_k=2\pi k/\ell$, it has multiplier
\begin{equation}\label{Tfouriermultiplier}
    \widehat{\mathcal T_\ell^{\mathfrak a}}(\xi_k)
    =-2\operatorname{Cin}\left(
       \frac{\ell|\xi_k+\mathfrak a_\gamma|}{2}\right),
    \qquad
    \operatorname{Cin}(r):=\int_0^r\frac{1-\cos q}{q}\,dq.
\end{equation}
In particular its multiplier grows logarithmically, so it is unbounded on
each fixed Sobolev space.  Since the principal normal curvature variation is
$(D_s^0)^2h$ and $J_{\gamma^0}$ is multiplication by $i$, the extracted term
in \eqref{Ndecompositionsectiontwo} gives precisely the principal nonlocal
term in \eqref{background-linearised-principal-form}.

It remains to understand the remainder term
$\widetilde{\mathcal N}(\gamma)$.  We first freeze the geometry at a straight
filament.  Let
$J:\mathbb R^2\to\mathbb R^2$ be the rotation
$J(z_1,z_2)=(-z_2,z_1)$.  For a smooth periodic
$z:\mathbb R/\ell\mathbb Z\to\mathbb R^2$ define
\begin{equation}\label{flatBoperator}
\begin{split}
    \mathcal B_\ell^{\mathfrak a}z(\sigma)
    :=\operatorname{p.v.}\int_{-\ell/2}^{\ell/2}
    \frac{1}{|h|^3}\Big\{
       &z(\sigma)-e^{i\mathfrak a_\gamma h}
       \bigl(z(\sigma+h)-hD_\sigma z(\sigma+h)\bigr)\\
       &\quad-\tfrac12h^2D_\sigma^2z(\sigma)
    \Big\}\,dh
\end{split}
\end{equation}
and
\begin{equation}\label{flatRoperator}
    \mathcal R_\ell^{\mathfrak a}
    :=\mathcal B_\ell^{\mathfrak a}
      -\tfrac12\mathcal T_\ell^{\mathfrak a}D_\sigma^2.
\end{equation}
The operator $J\mathcal B_\ell^{\mathfrak a}$ is the
linearisation, in normal graph
coordinates, of the finite-part Biot--Savart operator
$\mathcal N$ at a straight filament.  Correspondingly,
$J\mathcal R_\ell^{\mathfrak a}$ is the flat linearisation of
$\widetilde{\mathcal N}$.

\begin{lemma}[Flat remainder multiplier]\label{flatmultiplierlemma}
The operator $\mathcal R_\ell^{\mathfrak a}$ is a self-adjoint
Fourier multiplier.
Its symbol has the form
\begin{equation*}
    r_\ell^{\mathfrak a}(\xi)
    =(\xi+\mathfrak a_\gamma)^2\rho\left(
       \frac{\ell|\xi+\mathfrak a_\gamma|}{2}\right),
\end{equation*}
where
\begin{equation}\label{rhoexplicitformula}
 \rho(Q)=\frac{\sin Q}{Q}-\frac{1-\cos Q}{Q^2}-\frac12,
 \qquad \rho(0):=0.
\end{equation}
In particular $\rho$ is smooth and bounded on $[0,\infty)$ and
$\rho_\infty=-\frac12$.  On the Fourier lattice
$\xi_k=2\pi k/\ell$,
\begin{equation}\label{rhosymbolbounds}
    r_\ell^{\mathfrak a}(\xi_k)
    =-\frac12(\xi_k+\mathfrak a_\gamma)^2+q_{\ell,k},
    \qquad |q_{\ell,k}|\leq C_\ell(1+|\xi_k|).
\end{equation}
Consequently $\mathcal R_\ell^{\mathfrak a}$ has order two and
$J\mathcal R_\ell^{\mathfrak a}$ is skew-adjoint on
$L^2(\mathbb R/\ell\mathbb Z;\mathbb R^2)$.
\end{lemma}

\begin{proof}
For $z(\sigma)=e^{i\xi\sigma}z_0$, set
$\zeta=\xi+\mathfrak a_\gamma$.  Symmetry of the interval in
\eqref{flatBoperator} shows that the imaginary odd part vanishes.  The real
multiplier is
\begin{equation*}
    b_{\ell,\mathfrak a_\gamma}(\xi)
    =2\zeta^2\int_0^{\ell|\zeta|/2}
      \frac{1-\cos q-q\sin q+\tfrac12q^2}{q^3}\,dq.
\end{equation*}
On the other hand, by \eqref{Tfouriermultiplier}, the multiplier of
$\tfrac12\mathcal T_\ell^{\mathfrak a}D_\sigma^2$ is
$\zeta^2\operatorname{Cin}(\ell|\zeta|/2)$.  Therefore
\begin{equation*}
\begin{split}
    \rho(Q)
    =2\int_0^Q
      \frac{1-\cos q-q\sin q+\tfrac12q^2}{q^3}\,dq
      -\operatorname{Cin}(Q).
\end{split}
\end{equation*}
The numerator in the first integral is $q^4/8+O(q^6)$ at the
origin.  Moreover,
\begin{equation*}
    \rho'(Q)
    =\frac{\cos Q}{Q}-\frac{2\sin Q}{Q^2}
      +\frac{2(1-\cos Q)}{Q^3}.
\end{equation*}
The last expression is the derivative of
\[
 \frac{\sin Q}{Q}-\frac{1-\cos Q}{Q^2}.
\]
Since $\rho(0)=0$, integration gives \eqref{rhoexplicitformula}.
Consequently $\rho_\infty=-\frac12$ and
$\rho(Q)+\frac12=O((1+Q)^{-1})$.  This proves
\eqref{rhosymbolbounds}.  The oscillatory order-one remainder comes from the
sharp endpoint $|h|=\ell/2$ in the symmetric periodic representation; its
pointwise lattice bound gives $H^{m+1}\to H^m$, and it remains a real scalar
multiplier.  Since
$J^*=-J$ and $J$ commutes with these scalar operators,
$J\mathcal R_\ell^{\mathfrak a}$ is skew-adjoint.
\end{proof}

Thus, in the flat model, the complete order-two linearisation of the normal
nonlocal field is
\begin{equation*}
 J\mathcal B_\ell^{\mathfrak a}
 =\frac12J\mathcal T_\ell^{\mathfrak a}D_\sigma^2
  +J\mathcal R_\ell^{\mathfrak a}.
\end{equation*}
The first term has logarithmic order two and gives the term displayed
explicitly in \eqref{background-linearised-principal-form}; the second has
classical order two and is included there in $\mathscr E_\epsilon^0[H]$.
Both are skew-adjoint.  We next show that this flat principal structure
persists for a curved filament and control the remaining variable-coefficient
terms.

\subsection{Exact linearisation and the curved-filament symbol}

We now pass from the flat calculation above to the complete linearisation on
a curved filament, without freezing the length, the normal connection, or the
normal frame.  Put
$\mathbb T=\mathbb R/\mathbb Z$ and represent a constant-speed curve by
\begin{equation}\label{fixed-circle-curve}
 \Gamma(\theta)=\gamma(\ell\theta),\qquad
 |\Gamma_\theta|=\ell,\qquad \theta\in\mathbb T.
\end{equation}
Let $d_{\mathbb T}$ denote periodic distance on $\mathbb T$.  If
$\Gamma_\tau$ is a family in the constant-speed manifold, put
$X=\partial_\tau\Gamma_\tau|_{\tau=0}$ and
$\dot\ell=\partial_\tau\ell_{\Gamma_\tau}|_{\tau=0}$.  Differentiating
$|\Gamma_\theta|^2=\ell^2$ gives the exact tangent-space identities
\begin{equation}\label{constant-speed-tangent-space}
 \vec{\mathsf t}\cdot X_\theta=\dot\ell
 \quad\hbox{is independent of $\theta$},
 \qquad
 \dot{\vec{\mathsf t}}
 =\frac{X_\theta-\dot\ell\vec{\mathsf t}}{\ell},
 \qquad
 \dot{\vec K}=\frac{X_{\theta\theta}}{\ell^2}
       -2\frac{\dot\ell}{\ell}\vec K.
\end{equation}
Conversely, if $Y\perp\vec{\mathsf t}$ is prescribed, its unique
mean-zero tangential lift to this tangent space is
\begin{equation}\label{constant-speed-normal-lift}
 X_Y=Y+\phi_Y\vec{\mathsf t},\qquad
 \partial_\sigma\phi_Y
 =\langle\vec K,Y\rangle
  -\frac1\ell\int_0^\ell\langle\vec K,Y\rangle\,d\sigma,\qquad
 \int_0^\ell\phi_Y\,d\sigma=0.
\end{equation}
In particular, $Y\mapsto\phi_Y$ has order $-1$.

To compare operators associated with nearby curves on the fixed circle, we
require a tame choice of the length, holonomy and periodic normal frame.

\begin{lemma}[A tame periodic normal trivialisation]
\label{tameperiodicnormaltrivialisation}
Let $m\geq4$ and let $\Gamma_*$ be a smooth embedded constant-speed curve.
There is a neighbourhood $\mathscr U$ of $\Gamma_*$, fixed at a
low Sobolev level, in the
constant-speed manifold and a choice
\begin{equation}\label{localperiodicframechoice}
 \Gamma\longmapsto
 \big(\ell_\Gamma,\vartheta_\Gamma,
       \vec e_{1,\Gamma},\vec e_{2,\Gamma}\big)
\end{equation}
with the following properties.  The frame is periodic and oriented,
$\vartheta_\Gamma=\ell_\Gamma\mathfrak a_\Gamma$ is a fixed real lift of
the normal holonomy, and the map \eqref{localperiodicframechoice} is smooth
tame.  More precisely, for every fixed $q\geq0$,
\begin{equation}\label{framevariationtamebound}
\begin{split}
 |D^q\ell_\Gamma[X_1,\ldots,X_q]|
 &\leq C_q\prod_{\mu=1}^q\|X_\mu\|_{H^{3}},\\
 |D^q\vartheta_\Gamma[X_1,\ldots,X_q]|
 &\leq C_q\prod_{\mu=1}^q\|X_\mu\|_{H^{3}},\\
 \|D^q\vec e_{\alpha,\Gamma}[X_1,\ldots,X_q]
       \|_{H^{r+1}}
 &\leq C_{r,q}\bigg[
 (1+\|\Gamma\|_{H^{r+2}})
       \prod_{\mu=1}^q\|X_\mu\|_{H^3}\\
 &\hspace{3em}+
 \sum_{\nu=1}^q\|X_\nu\|_{H^{r+2}}
       \prod_{\mu\ne\nu}\|X_\mu\|_{H^3}\bigg]
\end{split}
\end{equation}
for $2\leq r\leq m$.  The same estimates hold for differences of the
frame, length and holonomy, with one factor replaced by
$\|\Gamma-\bar\Gamma\|$ at the corresponding Sobolev level.
\end{lemma}

\begin{proof}
Choose a unit vector $v_*$ which is not parallel to
$\vec{\mathsf t}_*(0)$.  After decreasing the neighbourhood, projection
of $v_*$ on $N_{\Gamma(0)}\Gamma$ and normalisation give a smooth initial
normal vector $E_1(0)$; put $E_2(0)=\vec{\mathsf t}(0)\times E_1(0)$.
Solve the normal parallel-transport equations
\begin{equation}\label{paralleltransportfixedcircle}
 \partial_\theta E_\alpha+
  \big(E_\alpha\cdot\partial_\theta\vec{\mathsf t}\big)
       \vec{\mathsf t}=0.
\end{equation}
The endpoint map on the oriented normal plane is a rotation.  In the
chosen neighbourhood it has a unique smooth angle lift
$\Theta_\Gamma$ near a fixed lift at $\Gamma_*$.  Its integer is constant
on this chart; the choice introduces no additional variation term.  Set
\begin{equation*}
 \begin{pmatrix}\vec e_1\\ \vec e_2\end{pmatrix}
 =R(-\theta\Theta_\Gamma)
   \begin{pmatrix}E_1\\ E_2\end{pmatrix},
 \qquad
 \vartheta_\Gamma=-\Theta_\Gamma.
\end{equation*}
This is periodic and has constant connection
$\mathfrak a_\Gamma=\vartheta_\Gamma/\ell_\Gamma$.

It remains to verify the estimates.  Since the curve has constant speed,
\[
 \ell_\Gamma=|\partial_\theta\Gamma|,
 \qquad
 \vec{\mathsf t}=\ell_\Gamma^{-1}\partial_\theta\Gamma.
\]
Differentiating these identities and
\eqref{paralleltransportfixedcircle} gives linear ordinary differential
equations for every variation of $E_\alpha$.  Their coefficients contain
at most two derivatives of $\Gamma$.  One-dimensional Sobolev
multiplication and Gronwall's inequality give
\eqref{framevariationtamebound}, with at most one factor at the high
Sobolev level.  The endpoint rotation and its local angle are smooth
finite-dimensional functions of $E_\alpha(1)$, so the same estimates hold
for $\Theta_\Gamma$.  Applying the fundamental theorem of calculus to the
path between two curves in the local constant-speed chart proves the
difference assertions.
\end{proof}

The fixed-circle unitary trivialisation used below is
\begin{equation}\label{fixed-circle-unitary-trivialisation}
 \mathscr U_\Gamma Y(\theta)
 :=\ell^{1/2}\big(Y(\ell\theta)\cdot\vec e_1(\ell\theta)
       +iY(\ell\theta)\cdot\vec e_2(\ell\theta)\big).
\end{equation}
It carries $\nabla_\sigma^\perp$ to
\begin{equation*}
 \mathscr D_\Gamma=\ell^{-1}(\partial_\theta+i\vartheta_\Gamma),
 \qquad \vartheta_\Gamma=\ell\mathfrak a_\Gamma,
\end{equation*}
and carries $J_\Gamma$ to multiplication by $i$.  Thus operators associated
with different lengths act on the same Hilbert space.  In these coordinates
\eqref{ngammadef} becomes
\begin{equation}\label{N-normalised-circle}
 \mathcal N(\Gamma)(\theta)
 =\int_{\mathbb T}\left\{
 \frac{\Gamma_{\theta'}(\theta')\times
       (\Gamma(\theta)-\Gamma(\theta'))}
      {|\Gamma(\theta)-\Gamma(\theta')|^3}
 -\frac{J_\Gamma\vec K(\theta)}{2d_{\mathbb T}(\theta,\theta')}
 \right\}\,d\theta'.
\end{equation}
The integral is understood as the limit of symmetric truncations.  Notice
that the factors of $\ell$ cancel in the counterterm.  This is the reason
for using the normalised circle when differentiating the operator.

\begin{lemma}[Full first variation of the finite part]
\label{full-first-variation-N}
Let $\Gamma$ be a smooth embedded constant-speed curve satisfying the
chord--arc condition
\[
 |\Gamma(\theta)-\Gamma(\theta')|
 \geq c_{\rm ca}\ell\,d_{\mathbb T}(\theta,\theta'),
 \qquad \theta,\theta'\in\mathbb T,
\]
for some $c_{\rm ca}>0$, and let $X$ satisfy
\eqref{constant-speed-tangent-space}. Set
\[
 R=\Gamma(\theta)-\Gamma(\theta'),\qquad
 \delta X=X(\theta)-X(\theta').
\]
Then the derivative of the finite-part field is the symmetric finite part
of
\begin{equation}\label{exact-DN-formula}
\begin{split}
 D\mathcal N(\Gamma)[X](\theta)
 ={}&\int_{\mathbb T}\bigg\{
 \frac{X_{\theta'}(\theta')\times R}{|R|^3}
 +\frac{\Gamma_{\theta'}(\theta')\times\delta X}{|R|^3}\\
 &\hspace{21mm}
 -3\frac{(\Gamma_{\theta'}(\theta')\times R)
          (R\cdot\delta X)}{|R|^5}\\
 &\hspace{21mm}
 -\frac{\dot{\vec{\mathsf t}}\times\vec K
          +\vec{\mathsf t}\times\dot{\vec K}}
        {2d_{\mathbb T}(\theta,\theta')}
 \bigg\}\,d\theta'.
\end{split}
\end{equation}
Consequently, if
$F(\Gamma)=\Pi_{\vec{\mathsf t}}^\perp\mathcal N(\Gamma)$, then its full
linearisation in the constant-speed gauge is
\begin{equation}\label{exact-DF-formula}
 DF(\Gamma)[X]
 =\Pi_{\vec{\mathsf t}}^\perp D\mathcal N(\Gamma)[X]
 -(\vec{\mathsf t}\cdot\mathcal N(\Gamma))\dot{\vec{\mathsf t}}
 -(\dot{\vec{\mathsf t}}\cdot\mathcal N(\Gamma))\vec{\mathsf t}.
\end{equation}
These formulae include the variations of the length, tangent, curvature,
normal projection and finite-part counterterm.  They therefore apply to
the lift \eqref{constant-speed-normal-lift}.
\end{lemma}

\begin{proof}
For $\delta>0$, restrict \eqref{N-normalised-circle} to
$d_{\mathbb T}(\theta,\theta')>\delta$.  The resulting integral has a
smooth kernel, so differentiation under the integral is legitimate.  The
three terms in the first two lines of \eqref{exact-DN-formula} follow from
\[
 \frac{d}{d\tau}\frac{A_\tau\times R_\tau}{|R_\tau|^3}
 =\frac{\dot A\times R+A\times\dot R}{|R|^3}
  -3\frac{(A\times R)(R\cdot\dot R)}{|R|^5},
\]
with $A=\Gamma_{\theta'}(\theta')$ and $\dot R=\delta X$.
The last term follows from
$\frac d{d\tau}(J_\Gamma\vec K)
=\dot{\vec{\mathsf t}}\times\vec K
+\vec{\mathsf t}\times\dot{\vec K}$ and
\eqref{constant-speed-tangent-space}.  Taylor expansion at
$\theta'=\theta$ through fourth order gives the leading term
$(\dot{\vec{\mathsf t}}\times\vec K+
\vec{\mathsf t}\times\dot{\vec K})/(2|h|)$.
Subtracting the differentiated counterterm leaves an integrable
remainder, with an integrable majorant for its tangential derivatives
uniform on the bounded chord--arc family.  Dominated
convergence proves \eqref{exact-DN-formula} first for smooth curves and in
$C^{m-3}$ for a $C^m$ family; the Sobolev version used below follows by
density.  Finally
\[
 \big(D\Pi_{\vec{\mathsf t}}^\perp[X]\big)\vec v
 =-(\vec{\mathsf t}\cdot\vec v)\dot{\vec{\mathsf t}}
  -(\dot{\vec{\mathsf t}}\cdot\vec v)\vec{\mathsf t},
 \qquad \vec v\in\mathbb R^3,
\]
and \eqref{exact-DF-formula} follows.
\end{proof}

Reparametrisation covariance gives the exact identity
\begin{equation}\label{tangential-linearisation-identity}
 DF(\gamma)[\psi\vec{\mathsf t}]
 =\psi\partial_\sigma F(\gamma).
\end{equation}
Together, \eqref{constant-speed-normal-lift},
\eqref{exact-DN-formula}, \eqref{exact-DF-formula} and
\eqref{tangential-linearisation-identity} determine the full Fr\'echet
derivative in any parametrisation.

We next put the nonlocal part of this derivative into the form needed for
the linear estimates.  Proposition \ref{curved-principal-symbol} separates
it into a scalar logarithmic order-two operator, a transport term, a
curvature-dependent logarithmic term and a tame order-zero remainder.
Before proving that proposition, we record in Lemma
\ref{finite-part-calculus-lemma} the kernel estimates used to control the
remainder.

For this purpose we use the fixed circle throughout and the real
\(L^2\) inner product for two-component normal fields.  An admissible
neighbourhood has a fixed positive chord--arc constant, lengths in a
compact subinterval of \((0,\infty)\), and a bounded \(H^{16}\) norm.
Higher norms of the curve are displayed in the estimates, rather than
included in their constants.  All comparisons of normal fields use
\eqref{fixed-circle-unitary-trivialisation}.

\begin{lemma}[Bounded remainder kernels]
\label{finite-part-calculus-lemma}
Let $\mathbb T=\mathbb R/\mathbb Z$.  All Sobolev norms below are
taken on $\mathbb T$, componentwise for vector- and matrix-valued
functions.  Fix an even cutoff
$\chi\in C_c^\infty((-1/4,1/4))$ equal to one near zero, and a
translation $h_0\in\mathbb T$.

Let $c\in C^\infty(\mathbb T;\mathbb R^{2\times2})$ and let
$K:\mathbb T\times[-1/2,1/2]\to\mathbb R^{2\times2}$ be a measurable
kernel satisfying
\[
 |K|_{m,1}
 :=\int_{-1/2}^{1/2}
       \|K(\cdot,h)\|_{H^m}\,dh<\infty
 \qquad\text{for every }m\geq0.
\]
For $v\in C^\infty(\mathbb T;\mathbb R^2)$, consider the linear
operators whose actions are
\begin{equation}\label{exact-order-zero-kernel-class}
\begin{split}
 &c(\theta)\operatorname{p.v.}\int_{-1/2}^{1/2}
       \frac{\chi(h)}h v(\theta+h)\,dh,
 \qquad
 \int_{-1/2}^{1/2}K(\theta,h)v(\theta+h)\,dh,\\
 &c(\theta)v(\theta),\qquad
 c(\theta)v(\theta+h_0),\qquad
 L\partial_\theta^{-1}(I-P_0)v.
\end{split}
\end{equation}
Here addition in the argument of $v$ is taken modulo one,
$\operatorname{p.v.}$ denotes the symmetric principal value at $h=0$,
and
\[
 P_0v:=\int_0^1v(\theta)\,d\theta,
 \qquad
 L=\log(1-\partial_\theta^2)^{1/2}.
\]
The mean $P_0v$ is regarded as a constant function, $I$ is the
identity operator, and $\partial_\theta^{-1}(I-P_0)v$ is the unique
zero-mean periodic function whose derivative is $v-P_0v$.
The operator $L$ acts componentwise.

Each operator in \eqref{exact-order-zero-kernel-class} extends to a
bounded operator on $H^m(\mathbb T;\mathbb R^2)$ for every $m\geq0$.
The same holds for any fixed finite sum or composition of these
operators, with independently chosen coefficients and kernels.
In particular, for $m\geq2$,
\[
 \left\|
   \int_{-1/2}^{1/2}K(\theta,h)v(\theta+h)\,dh
 \right\|_{H^m}
 \leq C_m\left(
       |K|_{2,1}\|v\|_{H^m}
       +|K|_{m,1}\|v\|_{H^2}\right).
\]

Now let $\Gamma$ range over an admissible neighbourhood
$\mathscr U$ of smooth constant-speed curves, as specified above.
Let $R_\Gamma$ be a fixed finite sum of finite compositions of the
operators in \eqref{exact-order-zero-kernel-class}, whose coefficients
$c_\Gamma$ and kernels $K_\Gamma$ depend smoothly on $\Gamma$ in
the respective Sobolev and kernel seminorms.  The cutoff,
translations, and choice of sums and compositions are independent
of $\Gamma$.

Write $D^q$ for the $q$th derivative with respect to $\Gamma$ in the
chosen local curve chart, and let $X_1,\ldots,X_q$ be smooth variation
directions.  For every $m\geq8$ and $q\geq0$, assume that each of the
quantities
\[
 \|D^qc_\Gamma[X_1,\ldots,X_q]\|_{H^m},
 \qquad
 |D^qK_\Gamma[X_1,\ldots,X_q]|_{m,1}
\]
is bounded by
\[
 C_{m,q}\left\{
 (1+\|\Gamma\|_{H^{m+8}})
       \prod_{j=1}^q\|X_j\|_{H^{16}}
 +\sum_{i=1}^q\|X_i\|_{H^{m+8}}
       \prod_{j\ne i}\|X_j\|_{H^{16}}
 \right\}.
\]
These assumptions apply to every coefficient and kernel occurring
in $R_\Gamma$.  For $q=0$, the derivative is the coefficient or
kernel itself, an empty product equals one, and an empty sum equals
zero.

Then, for every $m\geq8$,
\begin{equation}\label{order-zero-remainder-bound}
 \|R_\Gamma v\|_{H^m}
 \leq C_m\left\{
       \|v\|_{H^m}
       +(1+\|\Gamma\|_{H^{m+8}})\|v\|_{H^8}
       \right\}.
\end{equation}
Moreover, for every $q\geq1$,
\begin{equation}\label{remainder-family-tame-bound}
\begin{split}
 \|D^qR_\Gamma[X_1,\ldots,X_q]v\|_{H^m}
 \leq C_{m,q}\bigg\{&
 \|v\|_{H^m}\prod_{j=1}^q\|X_j\|_{H^{16}}\\
 &+(1+\|\Gamma\|_{H^{m+8}})\|v\|_{H^8}
       \prod_{j=1}^q\|X_j\|_{H^{16}}\\
 &+\|v\|_{H^8}
       \sum_{i=1}^q\|X_i\|_{H^{m+8}}
       \prod_{j\ne i}\|X_j\|_{H^{16}}
 \bigg\}.
\end{split}
\end{equation}
The constants depend on the indicated indices, the cutoff, the
fixed admissible neighbourhood, the coefficient bounds, and the
chosen finite expression for $R_\Gamma$.  Difference estimates
follow by integrating the corresponding next derivative along a
path in the same curve chart.
\end{lemma}

\begin{proof}
The principal-value convolution in
\eqref{exact-order-zero-kernel-class} has Fourier multiplier
\[
 2i\int_0^{1/4}\chi(h)\frac{\sin(2\pi nh)}h\,dh,
 \qquad n\in\mathbb Z,
\]
which is uniformly bounded in $n$.  This follows from the bounded
primitive of $\sin x/x$ and integration by parts against the smooth
cutoff.  It is therefore bounded on every $H^m$.
The last operator in \eqref{exact-order-zero-kernel-class} has
multiplier
\[
 \frac{\log(1+4\pi^2n^2)}{4\pi i n}
 \quad\text{for }n\ne0,
 \qquad 0\quad\text{for }n=0,
\]
which is also bounded.

Translations are isometries on every $H^m$.  For $m\geq2$, the
Sobolev product estimate
\begin{equation}\label{onedproductestimate}
 \|fg\|_{H^m}
 \leq C_m\left(
       \|f\|_{H^m}\|g\|_{H^2}
       +\|f\|_{H^2}\|g\|_{H^m}\right)
\end{equation}
controls multiplication by $c$.  Applying the same estimate at
each $h$ and then using Minkowski's inequality gives
\[
\begin{split}
 \left\|
   \int_{-1/2}^{1/2}K(\theta,h)v(\theta+h)\,dh
 \right\|_{H^m}
 &\leq \int_{-1/2}^{1/2}
       \|K(\cdot,h)v(\cdot+h)\|_{H^m}\,dh\\
 &\leq C_m\left(
       |K|_{2,1}\|v\|_{H^m}
       +|K|_{m,1}\|v\|_{H^2}\right).
\end{split}
\]
At $m=0$, boundedness follows from
\[
 \int_{-1/2}^{1/2}\|K(\cdot,h)\|_{L^\infty}\,dh
 \leq C|K|_{1,1}
\]
and the analogous bound for $c$.  Interpolation gives boundedness
for $0<m<2$.  Finite sums and compositions are consequently bounded
on every $H^m$.

For the curve-dependent assertions, first consider a single operator
in \eqref{exact-order-zero-kernel-class}.  Since the cutoff and
translations are fixed, differentiation with respect to $\Gamma$
acts only on its coefficient or kernel.  The assumed coefficient
bounds and the preceding product estimates therefore give
\eqref{order-zero-remainder-bound} and
\eqref{remainder-family-tame-bound} for each such operator.

At the low level $m=8$, the uniform $H^{16}$ bound on $\Gamma$ gives
\[
 \|D^qR_\Gamma[X_1,\ldots,X_q]v\|_{H^8}
 \leq C_q\|v\|_{H^8}
       \prod_{j=1}^q\|X_j\|_{H^{16}}.
\]
For a finite composition, the Leibniz rule distributes the curve
derivatives among its factors.  Repeated application of the
high-level estimates and this low-level bound retains at most one
high Sobolev norm in each term.  This proves
\eqref{order-zero-remainder-bound} and
\eqref{remainder-family-tame-bound} for the stated finite sums and
compositions.  Integrating the next derivative along a path in the
curve chart gives the difference estimates.
\end{proof}

\begin{proposition}[The curvature-dependent logarithmic term]
\label{curved-principal-symbol}
Let \(Y\) be normal, let \(X_Y\) be its constant-speed lift from
\eqref{constant-speed-normal-lift}, and put \(w=\mathscr U_\Gamma Y\).
Write \(k=(\kappa_1,\kappa_2)^T\) for the curvature components in the
chosen normal frame, evaluated at \(\sigma=\ell\theta\), and set
\begin{equation}\label{curvature-logarithmic-matrix}
 M_\Gamma=J\left(\frac12kk^T-\frac34|k|^2I\right),
 \qquad J=\begin{pmatrix}0&-1\\1&0\end{pmatrix}.
\end{equation}
Then
\begin{equation}\label{full-linearised-symbol-decomposition}
 \mathscr U_\Gamma\Pi_\Gamma^\perp DF(\Gamma)[X_Y]
 =i\mathscr B_{\ell,\vartheta}w
       -n_\Gamma\mathscr D_\Gamma w+M_\Gamma Lw+R_\Gamma^0w,
\end{equation}
where \(n_\Gamma=(\vec{\mathsf t}\cdot\mathcal N)\circ(\ell\,\cdot)\),
\(\mathscr D_\Gamma=\ell^{-1}(\partial_\theta+i\vartheta)\), and
\(\mathscr B_{\ell,\vartheta}=\mathscr U_\Gamma
\mathcal B_\ell^{\mathfrak a}\mathscr U_\Gamma^{-1}\).
The remainder satisfies \eqref{order-zero-remainder-bound} and
\eqref{remainder-family-tame-bound}.  In particular,
\begin{equation}\label{logarithmic-symmetric-matrix}
 \frac{M_\Gamma+M_\Gamma^T}{2}
 =\frac14(Jkk^T-kk^TJ),
 \qquad
 \left\|\frac{M_\Gamma+M_\Gamma^T}{2}\right\|_{\rm op}
 =\frac14|k|^2.
\end{equation}
The remainder after the scalar principal operator and transport have
been removed therefore has logarithmic order zero.
\end{proposition}

\begin{proof}
We first work in arclength and identify the normal plane at \(\sigma\)
with \(\mathbb R^2\).  Put \(t_0=t(\sigma)\), let \(E_0\) be its
oriented orthonormal normal basis, and let \(E_h\) be its parallel
transport to \(\sigma+h\).  In this proof
\[
\begin{gathered}
 \pi_0=E_0^T,\qquad T_h=t(\sigma+h),\qquad
 Q_h=\frac{\Gamma(\sigma+h)-\Gamma(\sigma)}h,\\
 g_j=|Q_h|^{-j},\qquad Z_h=\pi_0(T_h\times Q_h).
\end{gathered}
\]
All these coefficients extend smoothly across \(h=0\).  The chord--arc
condition bounds \(|Q_h|\) away from zero for \(|h|\leq\ell/2\).
Write \(k_h=E_h^TK(\sigma+h)\), \(k=k_0\), and
\[
 u=E_0^TY(\sigma),\quad v(h)=E_h^TY(\sigma+h),\quad
 \phi_h=\phi_Y(\sigma+h),\quad
 \eta=\ell^{-1}\int_0^\ell K\cdot Y\,d\sigma.
\]
Primes on \(v(h)\), \(\phi_h\) and \(k(h)=k_h\) denote
\(h\) differentiation; \(k'\) without an argument is its value
at zero.  We write \(u''=E_0^T(\nabla_\sigma^\perp)^2Y(\sigma)\).
Two-component outputs are identified with normal vectors by \(E_0\).
The lift and speed variation give
\[
 X_Y'(\sigma+h)=E_h(v'+\phi_hk_h)-\eta T_h,\qquad
 \pi_0\delta K=u''+kk^Tu+\phi_0k'.
\]

\smallskip
\emph{Exact coefficients and their Taylor polynomials.}
For a column \(a\in\mathbb R^2\), define the three matrix-valued
coefficients and three vector-valued coefficients
\begin{equation}\label{exact-six-normal-kernel-coefficients}
\begin{aligned}
 \mathsf A a&=-h\,\pi_0((E_ha)\times Q_h)g_3,\\
 \mathsf B a&=-\pi_0(T_h\times E_ha)g_3
                         +3Z_h(Q_h\cdot E_ha)g_5,\\
 \mathsf C a&=\pi_0(T_h\times E_0a)g_3
                         -3Z_h(Q_h\cdot E_0a)g_5,\\
 \mathsf d&=-h\,\pi_0((E_hk_h)\times Q_h)g_3
                         +3Z_h(Q_h\cdot T_h)g_5,\\
 \mathsf e&=\pi_0(T_h\times t_0)g_3
                         -3Z_h(Q_h\cdot t_0)g_5,\qquad
 \mathsf f=hZ_hg_3.
\end{aligned}
\end{equation}
Thus the normal integrand of \(D\mathcal N[X_Y]\), including its
counterterm, is exactly
\begin{equation}\label{exact-normal-integrand-coefficients}
 \mathcal I(h)=
 \frac{\mathsf A v'+\mathsf Bv+\mathsf Cu+
                 \mathsf d\phi_h+\mathsf e\phi_0+\mathsf f\eta}{|h|^3}
       -\frac{J(u''+kk^Tu+\phi_0k')}{2|h|}.
\end{equation}
For example, the terms with \(g_5\) are the differentiated radial
denominator in \eqref{exact-DN-formula}; the term \(\mathsf f\eta\)
comes from the speed variation in \(X_Y'\).

The identities \(T'=E_hk_h\), \(E_h'=-T_hk_h^T\) and
\(Q_h=\int_0^1T_{\rho h}\,d\rho\) give
\[
\begin{split}
 T_h&=t_0+hE_0k+\tfrac12h^2(E_0k'-|k|^2t_0)+O(h^3),\\
 E_h&=E_0-ht_0k^T
       -\tfrac12h^2(t_0(k')^T+E_0kk^T)+O(h^3),\\
 Q_h&=t_0+\tfrac12hE_0k
       +\tfrac16h^2(E_0k'-|k|^2t_0)+O(h^3),\\
 g_3&=1+\tfrac18|k|^2h^2+O(h^3),\qquad
 g_5=1+\tfrac5{24}|k|^2h^2+O(h^3).
\end{split}
\]
Substitution into \eqref{exact-six-normal-kernel-coefficients} gives
the following exact factorisations:
\begin{equation}\label{six-kernel-Taylor-factorisations}
\begin{aligned}
 \mathsf A&=hJ+h^3A_3,\\
 \mathsf B&=-J+h^2B_2+h^3B_3,&
 B_2&=\tfrac14Jkk^T+\tfrac38|k|^2J,\\
 \mathsf C&= J+h^2C_2+h^3C_3,&
 C_2&=\tfrac34Jkk^T-\tfrac38|k|^2J,\\
 \mathsf d&=-\tfrac12hJk+h^3d_3,\\
 \mathsf e&=\tfrac12hJk+\tfrac12h^2Jk'+h^3e_3,\qquad&
 \mathsf f&=-\tfrac12h^2Jk+h^3f_3.
\end{aligned}
\end{equation}
Here each subscript-three coefficient is the integral Taylor
remainder of its corresponding coefficient:
\[
 a_3(\sigma,h)=\frac12\int_0^1(1-\rho)^2
                 \partial_h^3\mathsf a(\sigma,\rho h)\,d\rho.
\]
This defines the coefficients also at \(h=0\).  In particular, the
quadratic coefficient of \(\mathsf d\) is zero, whereas that of
\(\mathsf e\) cancels the \(\phi_0k'\) counterterm.

\smallskip
\emph{The logarithmic term and an explicit remainder.}
Let \(\chi_\ell(h)=\chi(h/\ell)\); below write it as \(\chi\).
The flat integrand and its full integral are
\[
 \mathcal I_{\rm fl}(h)
   =J\frac{u-v+hv'-\tfrac12h^2u''}{|h|^3},\qquad
 \operatorname{f.p.}\int_{-\ell/2}^{\ell/2}\mathcal I_{\rm fl}(h)\,dh
   =J\mathcal B_\ell^{\mathfrak a}Y.
\]
Subtracting this integrand in
\eqref{exact-normal-integrand-coefficients} leaves the singular terms
\begin{equation}\label{even-logarithmic-kernel-calculation}
 -\frac12Jk\frac{h(\phi_h-\phi_0)}{|h|^3}
 -\frac{\eta Jk}{2|h|}
 +\frac{Jkk^T(u+v)}{4|h|}
 -\frac{3|k|^2J(u-v)}{8|h|}.
\end{equation}
The remaining terms are precisely
\[
 \operatorname{sgn}(h)
       (A_3v'+B_3v+C_3u+d_3\phi_h+e_3\phi_0+f_3\eta).
\]
Integrate the first term in
\eqref{even-logarithmic-kernel-calculation} by parts on
\(\{\delta<|h|<\ell/2\}\), using
\(h/|h|^3=-\partial_h(|h|^{-1})\) and
\(\phi_h'=k_h^Tv-\eta\).  Its inner boundary term tends to
\(-Jk(k^Tu-\eta)\).  Its \(\eta\) integral cancels the second
term of \eqref{even-logarithmic-kernel-calculation}.  Put
\[
 C_\Gamma=J\left(\tfrac14kk^T-\tfrac38|k|^2I\right).
\]
The remaining singular integral is
\begin{equation}\label{extracted-even-logarithmic-kernel}
 C_\Gamma\int\chi(h)\frac{u-v(h)}{|h|}\,dh.
\end{equation}
One more integration by parts, on the term \(A_3v'\), gives the
following identity for the entire near remainder:
\begin{equation}\label{explicit-near-curve-remainder}
\begin{split}
 \mathcal R_{\rm near}Y={}&-2A_3(\sigma,0)u-Jk(k^Tu-\eta)\\
 &-\int\operatorname{sgn}(h)\partial_h(\chi A_3)v\,dh\\
 &+\int\chi\operatorname{sgn}(h)
              (B_3v+C_3u+d_3\phi_h+e_3\phi_0+f_3\eta)\,dh\\
 &-\frac12Jk\int
       \frac{\chi'(\phi_h-\phi_0)+\chi(k_h-k)^Tv}{|h|}\,dh .
\end{split}
\end{equation}
Every integral in this formula is absolutely convergent.
In the last line, \(\chi'\) is supported away from zero and
\[
 \frac{k_h-k}{|h|}
       =\operatorname{sgn}(h)\int_0^1k'(\rho h)\,d\rho.
\]
The two displayed multiplication terms in
\eqref{explicit-near-curve-remainder} retain both inner boundary
contributions.

For the far remainder use the exact expression
\[
 \mathcal R_{\rm far}Y
   =\int(1-\chi)(\mathcal I-\mathcal I_{\rm fl})\,dh.
\]
The \(u''\) terms cancel.  The only derivative of the distinguished
field has coefficient
\(a_{\rm far}=(1-\chi)(\mathsf A-hJ)/|h|^3\).
Replace its integral against \(v'\) by
\[
 a_{\rm far}(\ell/2)v(\ell/2)
 -a_{\rm far}(-\ell/2)v(-\ell/2)-\int a_{\rm far}'v\,dh.
\]
This exhibits the far part as smooth integral operators, endpoint
translations, and multiplication or integral operators applied to
\(\phi_Y\) and \(\eta\).

\smallskip
\emph{Uniform bounds and passage to the fixed frame.}
All estimates are taken after \(h=\ell\xi\), so that the integration
interval, cutoff and endpoint translations are fixed when the curve
varies.  In the integral Taylor formulas, \(\partial_h^j\) becomes
\(\ell^{-j}\partial_\xi^j\), and the curve is evaluated at the fixed
normalised point \(\theta+\rho\xi\).  Parameter differentiation
therefore leaves every translation argument fixed.
The parallel components satisfy
\[
 v(h)=e^{\mathfrak a hJ}E_{\rm per}(\sigma+h)^TY(\sigma+h).
\]
Here \(E_{\rm per}=(\vec e_1,\vec e_2)\) is the chosen periodic
normal frame.  Thus these components are a fixed translation of the periodic components followed
by a smooth rotation.  The coefficients in
\eqref{explicit-near-curve-remainder} involve the Taylor remainders
in \eqref{six-kernel-Taylor-factorisations} and
\(\partial_hA_3\).  The latter uses four \(h\) derivatives of
\(\mathsf A\); since \(Q_h\) is an averaged tangent, this uses at most
five derivatives of the curve.  The same bound covers \(d_3\), whose
formula includes \(T_h'\).  All other displayed coefficients have
no larger count.  Smooth reciprocal powers of \(|Q_h|\), and their
parameter derivatives, are controlled by the chord--arc lower bound.

We now verify the coefficient hypotheses of Lemma
\ref{finite-part-calculus-lemma}.  For every integral kernel
\(K_\Gamma\) and multiplication or endpoint coefficient \(c_\Gamma\)
obtained above, we have
\begin{equation}\label{explicit-curve-kernel-seminorm-bound}
\begin{split}
 &|D^pK_\Gamma[X_1,\ldots,X_p]|_{m,1}
       +\|D^pc_\Gamma[X_1,\ldots,X_p]\|_{H^m}\\
 &\quad\leq C_{m,p}\bigg\{
 (1+\|\Gamma\|_{H^{m+8}})\prod_i\|X_i\|_{H^{16}}
 +\sum_i\|X_i\|_{H^{m+8}}\prod_{j\ne i}\|X_j\|_{H^{16}}
                                                        \bigg\}.
\end{split}
\end{equation}
Indeed, the integral Taylor formula bounds each coefficient in
\(H^m\), uniformly in \(\xi\), by the corresponding differentiated
curve and frame norms.  The bounded factors
\(\operatorname{sgn}(\xi)\), the cutoff derivatives, and the divided
difference in the last line of
\eqref{explicit-near-curve-remainder} give integrable majorants.
Two further derivatives in estimating these coefficients and one
Sobolev embedding fit within the reserve eight.  Products and
reciprocal powers obey the one-high-factor product bound.  This proves
\eqref{explicit-curve-kernel-seminorm-bound} for any fixed number of
parameter derivatives; differences follow by integration along the
chart path.  The lift
\(\phi_Y=\partial_\sigma^{-1}(I-P_0)(K\cdot Y)\) and the scalar
\(\eta\) are bounded maps of the same type.

The tangential coefficient can also be estimated without using the
first variation.  In ambient coordinates set
\[
 G(\sigma,h)=-h(T_h\times Q_h)|Q_h|^{-3}
       =\tfrac12h^2(t_0\times K)+h^3G_3(\sigma,h),
\]
where \(G_3\) is its integral Taylor remainder.  The defining
counterterm then gives the exact formula
\begin{equation}\label{exact-undifferentiated-finite-part-remainder}
 \mathcal N(\Gamma)(\sigma)
       =\int_{-\ell/2}^{\ell/2}\operatorname{sgn}(h)G_3(\sigma,h)\,dh.
\end{equation}
The same coefficient estimates apply to this absolutely convergent
integral and to \(n_\Gamma=t\cdot\mathcal N\).

In the periodic frame the scalar operator in
\eqref{extracted-even-logarithmic-kernel} has multiplier
\[
 \tau_\vartheta(n)=
   2\int_0^{1/2}\chi(\xi)
       \frac{1-\cos((2\pi n+\vartheta)\xi)}{\xi}\,d\xi.
\]
Splitting at \(\langle n\rangle^{-1}\) and integrating the remaining
oscillatory integral by parts gives
\(\tau_\vartheta(n)-2\log\langle2\pi n\rangle=O(1)\).
Every positive-order \(\vartheta\) derivative is uniformly bounded,
since its integral has a factor \(\xi^{p-1}\), \(p\geq1\).
Let \(\mathsf T_{\chi,\vartheta}\) denote this multiplier.
The remainder in the proposition is therefore exactly
\[
\begin{split}
 R_\Gamma^0w={}&
 \mathscr U_\Gamma(\mathcal R_{\rm near}+\mathcal R_{\rm far})Y
       +C_\Gamma(\mathsf T_{\chi,\vartheta}-2L)w\\
 &\hspace{3em}-n_\Gamma\mathscr U_\Gamma(\phi_YK).
\end{split}
\]
The last term is the bounded part of
\(\Pi^\perp D\Pi^\perp[X_Y]\mathcal N
 =-n_\Gamma(\nabla_\sigma^\perp Y+\phi_YK)\).
Lemma \ref{finite-part-calculus-lemma} and
\eqref{explicit-curve-kernel-seminorm-bound} now prove the stated
remainder and family bounds.  Since \(M_\Gamma=2C_\Gamma\), this also
proves \eqref{full-linearised-symbol-decomposition}.  Finally,
\(J^T=-J\) gives \eqref{logarithmic-symmetric-matrix}; in a frame with
\(k=(|k|,0)^T\), the symmetric matrix has off-diagonal entries
\(|k|^2/4\).
\end{proof}

\subsection{Tame estimates for the finite-part Biot--Savart operator}

We now turn the structural decomposition of Proposition
\ref{curved-principal-symbol} into the tame Sobolev and energy estimates
needed for the Nash--Moser argument.  For \(\Gamma\in\mathscr U\), write
\(\ell=\ell_\Gamma\) and \(\vartheta=\vartheta_\Gamma\) for the length and
real holonomy parameter chosen in Lemma
\ref{tameperiodicnormaltrivialisation}.  The scalar operator
\(\mathscr B_{\ell,\vartheta}\) appearing in
\eqref{full-linearised-symbol-decomposition} has Fourier multiplier \(b_n\),
while the flat remainder \(\mathcal R_\ell^{\mathfrak a}\) defined in
\eqref{flatRoperator}, expressed on the fixed circle, has multiplier \(r_n\).
For \(p_n=2\pi n+\vartheta\), these multipliers are
\begin{equation}\label{fixed-circle-principal-multipliers}
\begin{split}
 b_n&=\ell^{-2}p_n^2\{\operatorname{Cin}(|p_n|/2)+\rho(|p_n|/2)\},\\
 r_n&=\ell^{-2}p_n^2\rho(|p_n|/2).
\end{split}
\end{equation}
Thus \(b_n\) has logarithmic order two, whereas \(r_n\) has classical
order two.  The explicit formula \eqref{rhoexplicitformula} gives, for
\(j,k\geq0\),
\begin{equation*}
\begin{split}
 |\partial_\ell^j b_n|&\leq C_j\langle n\rangle^2
                                      (1+\log\langle n\rangle),\\
 |\partial_\ell^j\partial_\vartheta^k b_n|
 &\leq C_{j,k}\langle n\rangle(1+\log\langle n\rangle),
                                                    \quad k\geq1,\\
 |\partial_\ell^j r_n|&\leq C_j\langle n\rangle^2,
 \qquad
 |\partial_\ell^j\partial_\vartheta^k r_n|
       \leq C_{j,k}\langle n\rangle,\quad k\geq1.
\end{split}
\end{equation*}
All constants are uniform for \(\Gamma\) in the admissible neighbourhood
\(\mathscr U\) fixed above.  The endpoint contributions at
\(h=\pm\ell/2\) are retained; on the fixed circle they are translations by
\(\pm1/2\) composed with operators of order at most one.

\begin{proposition}[Tame field bounds and logarithmic quadratic estimate]
\label{tameskewremainder}
In a fixed constant-speed chart the full field
\(\mathscr F(\Gamma)=\mathscr U_\Gamma\Pi_\Gamma^\perp\mathcal N(\Gamma)\)
is smooth tame of a fixed finite degree.  One may use the reserve eight:
\begin{equation}\label{higher-variation-tame-full-N}
\begin{split}
 &\|D^q\mathscr F(\Gamma)[X_1,\ldots,X_q]\|_{H^m}\\
 &\quad\leq C_{m,q}\bigg\{
 (1+\|\Gamma\|_{H^{m+8}})\prod_j\|X_j\|_{H^{16}}
 +\sum_i\|X_i\|_{H^{m+8}}
                       \prod_{j\ne i}\|X_j\|_{H^{16}}\bigg\}.
\end{split}
\end{equation}
The same finite-reserve assertion holds for the tangential component
\(n_\Gamma\) and for \(\widetilde{\mathcal N}\).  Differences are
obtained by integrating the next variation along a chart segment.

Put \(A=(1-\partial_\theta^2)^{1/2}\) and \(L=\log A\), and define
\[
 \mathfrak q_m(w)
 :=\operatorname{Re}\left\langle A^mw,
 A^m\bigl(i\mathscr B_{\ell,\vartheta}w
       -n_\Gamma\mathscr D_\Gamma w+M_\Gamma Lw+R_\Gamma^0w\bigr)
                         \right\rangle_{L^2}.
\]
Thus \(\mathfrak q_m(w)\) is the contribution of the linearised operator in
\eqref{full-linearised-symbol-decomposition} to the derivative of
\(\frac12\|w\|_{H^m}^2\).  Then
\begin{equation}\label{skewlinearisedestimate}
\begin{split}
 |\mathfrak q_m(w)|\leq{}&
 \frac14\|k\|_\infty^2\|L^{1/2}A^mw\|_2^2+C_m\|w\|_{H^m}^2\\
 &+C_m(1+\|\Gamma\|_{H^{m+10}})
                      \|w\|_{H^8}\|w\|_{H^m}.
\end{split}
\end{equation}
In particular, the coefficient of the unbounded logarithmic norm is
independent of \(m\).  The associated polarised form satisfies the analogous
estimate, with the product of the two logarithmic norms in place of their
square.
\end{proposition}

\begin{proof}
Formula \eqref{exact-undifferentiated-finite-part-remainder}
expresses the full field as an absolutely convergent integral of the
factored coefficient \(G_3\).  Its integral Taylor formula uses at
most four derivatives of the curve.  Differentiating this formula,
after rescaling to the fixed circle, produces products of divided
differences and smooth reciprocal chord quotients.  The product
estimate of Lemma \ref{finite-part-calculus-lemma} places one curve
or direction factor at the high level.  Normal projection and the
output frame have the bounds of
Lemma \ref{tameperiodicnormaltrivialisation}.  Composition with the
constant-speed chart requires solving the constraint \eqref{f3spatialexact}
for the dependent graph and length variables.  The required tame bounds
follow from the differentiated constraint \eqref{chart-recovery-ODE},
analysed at the beginning of the next subsection.  These comparisons and
Sobolev products fit within the reserve eight.
This proves \eqref{higher-variation-tame-full-N}, including \(q=0\),
and the analogous estimate for \(n_\Gamma\), directly from the defining
integral.  Parameter differentiation is justified by the uniform
integrable coefficient bounds, so it commutes with this integral.

The term subtracted in defining \(\widetilde{\mathcal N}\) is a
logarithmic multiplier applied to curvature, with smoothly varying
transport.  Its multiplier and parameter derivatives have the bounds
in \eqref{fixed-circle-principal-multipliers}; the logarithm costs
one derivative.  Thus it has the same finite-reserve bound.
Differences follow by integrating the next variation along a chart
segment.

It remains to prove \eqref{skewlinearisedestimate}, which controls the
symmetric part of the linearised operator in the \(H^m\) inner product.
The term \(i\mathscr B_{\ell,\vartheta}w\) makes no contribution to
\(\mathfrak q_m(w)\), since \(\mathscr B_{\ell,\vartheta}\) is a
self-adjoint scalar Fourier multiplier and commutes with \(A\).

For the transport term, set \(a=-n_\Gamma\).  Since \(a\) is real,
\[
 2\operatorname{Re}\langle A^mw,A^m(a\mathscr D_\Gamma w)\rangle
 =-\int(\partial_\sigma a)|A^mw|^2
   +2\operatorname{Re}\langle A^mw,[A^m,a]\mathscr D_\Gamma w\rangle.
\]
The standard one-dimensional Sobolev commutator estimate bounds the second
term without placing more than \(m\) derivatives on \(w\).

For the curvature-dependent logarithmic term, put \(v=A^mw\).  Then
\[
 A^mM_\Gamma Lw=M_\Gamma Lv+[A^m,M_\Gamma]Lw.
\]
The Fourier difference of the weights, or the fractional product
commutator estimate, gives
\begin{equation}\label{logarithmic-commutator-bound}
 \|[A^m,M_\Gamma]Lw\|_2
 \leq C_m\bigl(\|M_\Gamma\|_{H^3}\|w\|_{H^m}
                       +\|M_\Gamma\|_{H^m}\|w\|_{H^3}\bigr).
\end{equation}
Indeed the difference of the \(m\)th-order weights gains one
derivative, and \(\log\langle n\rangle\leq C\langle n\rangle\).
Also \([L^{1/2},M_\Gamma]L^{1/2}\) is bounded on \(L^2\), with a bound
depending on \(\|M_\Gamma\|_{H^3}\); this follows from the difference of
the two logarithmic Fourier weights and the Fourier Schur test.  Writing
\(M_{\Gamma,{\rm sym}}=(M_\Gamma+M_\Gamma^T)/2\), we consequently have
\[
 \operatorname{Re}\langle v,M_\Gamma Lv\rangle
 =\langle L^{1/2}v,M_{\Gamma,{\rm sym}}L^{1/2}v\rangle
                          +O(\|M_\Gamma\|_{H^3}\|v\|_2^2).
\]
By \eqref{logarithmic-symmetric-matrix}, the first term on the right is
bounded in absolute value by
\[
 \frac14\|k\|_\infty^2\|L^{1/2}A^mw\|_2^2.
\]
The tame bound for \(R_\Gamma^0\), the preceding transport estimate, and
the frame and coefficient estimates control the remaining terms.  The
\(H^{m+10}\) norm of \(\Gamma\) displayed in
\eqref{skewlinearisedestimate} provides the derivatives required for
\(n_\Gamma\), \(M_\Gamma\), and the frame.  Combining these estimates proves
\eqref{skewlinearisedestimate} for every real \(m\geq8\).  Applying the same
estimates to the associated bilinear form proves the polarised assertion.
\end{proof}

\subsection{The reduced equation and its tame linearisation}

We now return to the normal-graph evolution system
\eqref{exactf123system}, together with the constant-speed constraint
\eqref{f3spatialexact} and the length and tangential-gauge equations
\eqref{lambdatimeequation}--\eqref{ameanzero}.  Using the estimates of
the preceding subsection, we express this system in independent
variables and establish uniform tame bounds for its linearisation.
Proposition \ref{NashMoser-linear-inverse-proposition} in the next
subsection will then construct the tame inverse required by the
Nash--Moser argument.

We first choose coordinates for the constrained system.  Write
\(f_3=b+g\), where
\[
 b=\overline{f_3}:=\frac1L\int_0^L f_3(s,t)\,ds,
 \qquad \int_0^L g(s,t)\,ds=0,
\]
and put
\[
 z=(f_1,f_2,b).
\]
The constant-speed constraint \eqref{f3spatialexact}, together with
periodicity, determines \((g,\lambda-1)\) locally from \(z\).
Indeed, at the zero graph its derivative in
\((g,\mu)=(g,\lambda-1)\) is
\[
 (g,\mu)\longmapsto\partial_sg-\mu:
 H^m_0\times\mathbb R\longrightarrow H^{m-1},
\]
which is an isomorphism; here \(H^m_0\) denotes the zero-mean
\(L\)-periodic Sobolev space.  Denote the resulting recovery map,
the corresponding curve, and its differential by
\begin{equation}\label{NashMoser-recovery-map}
 \mathfrak R_t(z)=(f_1,f_2,f_3,\lambda-1),
 \qquad \mathfrak C_t(z)=\Gamma_z,
 \qquad \mathscr A_t(z)=D_z\mathfrak C_t(z).
\end{equation}
These maps are smooth tame on a fixed small neighbourhood of zero,
uniformly for \(t\in[0,T]\).

We next establish the constraint estimates that justify this tame
dependence and control the linearisation at arbitrary background
paths, whether or not they solve the evolution equation.
Put \(P=D_s^0f+k^0f_3\) and \(p_0=(\lambda^2-|P|^2)^{1/2}>0\).
For a direction \(y\), write \(\eta=D_zf_3[y]\) and
\(\dot\lambda=D_z\lambda[y]\).  Differentiation gives
\begin{equation}\label{chart-recovery-ODE}
 \eta'+\frac{P\cdot k^0}{p_0}\eta
      -\frac\lambda{p_0}\dot\lambda
 =k^0\cdot y_\perp-\frac P{p_0}\cdot D_s^0y_\perp,
 \qquad \overline\eta=y_b.
\end{equation}
The operator in \((\eta-\overline\eta,\dot\lambda)\) is a small
perturbation of \((\partial_s,-1)\).  It has a uniformly bounded
periodic inverse, obtained by an integrating factor and the mean
condition.  Integrate the term containing \(y_\perp'\) once by
parts in this formula.  It becomes multiplication by \(-P/p_0\)
on \(y_\perp\), followed by an integral with one derivative on its
coefficient.  This proves that the recovery differential is order
zero in the distinguished direction.  The nonlinear recovery follows
also by the implicit-function theorem from \(H^m\) to \(H^{m-1}\),
since the equation contains only first derivatives and the inverse
gains one derivative.

Further differentiation of \eqref{chart-recovery-ODE} gives the tame
bounds for the recovery and its inverse.  In particular, for \(m\geq8\),
\begin{equation}\label{chart-second-differential-bound}
\begin{split}
 \|D_z^2\mathfrak C_t(z)[r,y]\|_{H^m}
 \leq C_m\bigl\{&\|r\|_{H^{10}}\|y\|_{H^m}
                     +\|r\|_{H^{m+2}}\|y\|_{H^8}\\
 &+(1+\|z\|_{H^{m+2}})
                         \|r\|_{H^{10}}\|y\|_{H^8}\bigr\}.
\end{split}
\end{equation}
To check the count, a second differentiation adds products containing
at most one spatial derivative of each direction.  In the inverse
formula, move the derivative on the distinguished direction to its
coefficient.  This uses at most two derivatives of the other direction
and at most two of the background.  Apply
\eqref{onedproductestimate}; all remaining undifferentiated factors
are bounded on the low neighbourhood.  The same argument after any
fixed number of parameter differentiations puts one factor at the high
level.  No assumption that \(z\) solves an evolution equation is used.

Let \(\mathfrak X_T\) be the smooth functions
\([0,T]\to C^\infty(\mathbb R/L\mathbb Z;\mathbb R^2)\times\mathbb R\),
graded by
\[
 \|z\|_{\mathfrak X_T^m}
 =\sum_{2j+q\leq m}\sup_{0\leq t\leq T}
                         \|\partial_t^j\partial_s^qz(t)\|_{L^2}.
\]
The spatial norm includes the absolute value of the scalar component.
Put \(\mathfrak H=C^\infty(\mathbb R/L\mathbb Z;\mathbb R^2)
\times\mathbb R\), \(\mathfrak Y_T=\mathfrak X_T\times\mathfrak H\),
and
\[
 \|(h,y^{\rm in})\|_{\mathfrak Y_T^m}
 =\|h\|_{\mathfrak X_T^m}+\|y^{\rm in}\|_{\mathfrak H^m}.
\]
For the tame calculus use the even subgrading
\(\|z\|_{(a)}=\|z\|_{\mathfrak X_T^{2a}}\), and similarly for
\(\mathfrak Y_T\).  It defines the same smooth topology.  Odd
seminorms are bounded by the next even seminorm; interpolation will
always use even endpoint grades.  Extension across the time endpoints
and Fourier cutoffs with temporal weight two give, for even
\(m\geq r\geq0\), smoothing operators satisfying
\begin{equation}\label{NashMoser-smoothing-estimates}
 \|S_\theta v\|_m\leq C_{m,r}\theta^{m-r}\|v\|_r,
 \qquad
 \|(I-S_\theta)v\|_r\leq C_{m,r}\theta^{r-m}\|v\|_m.
\end{equation}

Write \(\Lambda_\epsilon=|\log\epsilon|\).  The normal velocity is
\(J_\Gamma K_\Gamma+\Lambda_\epsilon^{-1}\mathcal R_\epsilon[\Gamma]\),
where
\begin{equation}\label{Repsilongeometric}
 \mathcal R_\epsilon[\Gamma]
   =(\log\ell_\Gamma+\mathfrak c_\epsilon)J_\Gamma K_\Gamma
                          +\Pi_\Gamma^\perp\mathcal N(\Gamma).
\end{equation}
Include the mean-zero tangential velocity from
\eqref{atangentialequation}--\eqref{ameanzero}, and call the full
velocity \(\mathscr V_\epsilon\).  Define the reduced vector field by
\begin{equation}\label{reduced-vector-field-identity}
 \mathscr A_t(z)\mathscr G_\epsilon(t,z)
 =\mathscr V_\epsilon(\Gamma_z)-\partial_t^{\rm ch}\mathfrak C_t(z),
 \qquad
 \Phi_\epsilon(z)
 =\bigl(\partial_tz-\mathscr G_\epsilon(t,z),z(0)\bigr).
\end{equation}
The superscript \({\rm ch}\) holds \(z\) fixed.  The estimates above
show that \(\Phi_\epsilon\) is smooth tame of fixed finite degree,
uniformly for \(0<\epsilon\leq e^{-2}\).

For later coefficient estimates we record an explicit way to count
derivatives of this vector field.  Normal projection of
\eqref{reduced-vector-field-identity} determines its two normal
components.  The tangential mean is determined by the mean condition
on the gauge.  Thus the inverse derivative in the tangential recovery
does not add to the differential order of the normal velocity.  The
bound \eqref{higher-variation-tame-full-N}, the curvature formula, and
\eqref{chart-recovery-ODE} give a reserve eight for
\(\mathscr G_\epsilon\).  A reserve twelve suffices for all spatial,
time-weighted and polarised estimates for
\(\mathscr L_z=D_z\mathscr G_\epsilon(t,z)\).  The time derivatives
fall on the coefficients and on \(z\), with weight two, and the
number of parameter derivatives does not increase this fixed reserve.
In particular, for even \(m\geq40\), fixed \(j\geq0\), and base directions
\(v_1,\ldots,v_j\),
\begin{equation}\label{curve-coefficient-family-tame}
\begin{split}
 &\|D_z^j\mathscr L_z[v_1,\ldots,v_j]Y\|_{\mathfrak X_T^m}\\
 &\quad\leq C_{m,j}\bigg\{
 \|Y\|_{\mathfrak X_T^{m+12}}\prod_i\|v_i\|_{\mathfrak X_T^{40}}\\
 &\qquad +(1+\|z\|_{\mathfrak X_T^{m+12}})
       \|Y\|_{\mathfrak X_T^{40}}\prod_i\|v_i\|_{\mathfrak X_T^{40}}\\
 &\qquad +\|Y\|_{\mathfrak X_T^{40}}
       \sum_i\|v_i\|_{\mathfrak X_T^{m+12}}
                           \prod_{l\ne i}\|v_l\|_{\mathfrak X_T^{40}}
 \bigg\}.
\end{split}
\end{equation}
This follows by differentiating the recovered vector field, applying
\eqref{higher-variation-tame-full-N} and
\eqref{chart-recovery-ODE}, and using the anisotropic product estimate.
In each term the total temporal weight is the number of time
 derivatives multiplied by two; the spatial reserve is unchanged.

\begin{lemma}[The linear equation at an arbitrary background]
\label{complete-lower-block-decomposition}
Let \(z(t)\) be any smooth path in a sufficiently small fixed low
neighbourhood, and let
\(Y_t-\mathscr L_zY=h\).  Set
\begin{equation*}
\begin{gathered}
 X=\mathscr A_t(z)Y=Y^\perp+\phi t_z,
 \qquad w=\mathscr U_{\Gamma_z}Y^\perp,\\
 \beta=\frac1{\ell_z}\int_0^{\ell_z}\phi\,d\sigma,
 \qquad W=\mathscr S_zY=(w,\beta).
\end{gathered}
\end{equation*}
The coordinate maps \(\mathscr S_z\) and \(\mathscr S_z^{-1}\) are
order zero and satisfy
\begin{equation}\label{Eqz-tame-symmetriser-bound}
 \|\mathscr S_zY\|_{H^m\times\mathbb R}
 \leq C_m\{\|Y\|_{H^m}
           +(1+\|z\|_{H^{m+2}})\|Y\|_{H^8}\},
 \qquad m\geq8,
\end{equation}
with the analogous estimate for the inverse and their polarised
variations.  Only the fixed low-level norm equivalence is uniform
without a high background term.

In these coordinates the exact linear equation is
\begin{equation}\label{exact-fixed-space-linear-system}
\begin{split}
 W_t={}&\begin{pmatrix}i\mathscr H_{\epsilon,z}&0\\0&0\end{pmatrix}W
                       +\mathscr Q_{\epsilon,z}W+\mathscr S_zh\\
 &+\begin{pmatrix}(a_z-\Lambda_\epsilon^{-1}n_z)\mathscr D_z
              +\Lambda_\epsilon^{-1}M_zL&0\\0&0\end{pmatrix}W.
\end{split}
\end{equation}
where \(a_z\) is the tangential velocity, \(n_z=t_z\cdot\mathcal N(\Gamma_z)\)
in fixed-circle coordinates, and
\begin{equation*}
 \mathscr H_{\epsilon,z}
 =\mathscr D_z^2+\Lambda_\epsilon^{-1}
       \{\mathscr B_{\ell_z,\vartheta_z}
          +(\log\ell_z+\mathfrak c_\epsilon)\mathscr D_z^2\},
 \qquad
 \mathscr D_z=\ell_z^{-1}(\partial_\theta+i\vartheta_z).
\end{equation*}
It is self-adjoint, commutes with \(A\), and commutes with multiplication
by \(i\).  The matrix \(M_z\) is \eqref{curvature-logarithmic-matrix}.

Put
\begin{equation}\label{corrected-coefficient-reserve}
 \mathcal Z_m(t)=\|z(t)\|_{H^{m+12}}
                         +\|z_t(t)\|_{H^{m+10}}.
\end{equation}
The remaining coefficients obey
\begin{equation}\label{complete-transport-coefficient-bound}
 \|a_z\|_{H^{m+1}}+\Lambda_\epsilon^{-1}\|n_z\|_{H^{m+1}}
          +\|M_z\|_{H^m}
 \leq C_m(1+\mathcal Z_m(t)),
\end{equation}
and the bounded block satisfies
\begin{equation}\label{complete-order-zero-operator-bound}
 \|\mathscr Q_{\epsilon,z}W\|_{H^m\times\mathbb R}
 \leq C_m\{\|W\|_{H^m\times\mathbb R}
                +(1+\mathcal Z_m(t))\|W\|_{H^8\times\mathbb R}\}.
\end{equation}
These bounds hold for real \(m\geq8\), and have polarised and
difference versions with one high factor.  Their constants depend only
on a fixed \(\mathfrak X_T^{40}\) bound for the background, in
addition to the prescribed flow and the indicated indices.
\end{lemma}

\begin{proof}
The constant-speed tangent identities give
\begin{equation}\label{linearised-chart-constraint-exact}
 \phi=\beta+\partial_\sigma^{-1}(I-P_0)(k_z\cdot Y^\perp),
 \qquad
 \dot\ell[Y]=-\int_0^{\ell_z}k_z\cdot Y^\perp\,d\sigma.
\end{equation}
Together with \eqref{chart-recovery-ODE} and the frame bounds these
give \eqref{Eqz-tame-symmetriser-bound}, including its inverse.
They also show that every block containing \(\beta\) is a bounded
finite-dimensional or finite-rank map.

There is an additional term when the background is not a solution.
Differentiating \eqref{reduced-vector-field-identity} and then using
the linear equation gives the identity
\begin{equation}\label{off-solution-chart-identity}
 X_t=D_\Gamma\mathscr V_\epsilon(\Gamma_z)X
       +D_z^2\mathfrak C_t(z)[r_z,Y]+\mathscr A_t(z)h,
 \qquad r_z=z_t-\mathscr G_\epsilon(t,z).
\end{equation}
Indeed the two explicit time derivatives of the chart cancel and the
two second chart differentials combine to the displayed defect.
Equation \eqref{chart-second-differential-bound} proves that this defect
term is order zero in \(Y\).  Since
\(\mathscr G_\epsilon\) has reserve eight, its coefficients use at
most \(z\) at level \(m+10\), or \(z_t\) at level \(m+2\), before
the frame comparison.  Both are bounded by
\eqref{corrected-coefficient-reserve}.

Now project \eqref{off-solution-chart-identity} normally.  The normal
variation of the local binormal term is
\[
 J_z\{(\nabla_\sigma^\perp)^2Y^\perp
                +(k_z\cdot Y^\perp)k_z+\phi\nabla_\sigma^\perp k_z\}.
\]
The normal variation of \(a_zt_z\) is
\(a_z(\nabla_\sigma^\perp Y^\perp+\phi k_z)\); the variation of
\(a_z\) itself is tangent.  Proposition \ref{curved-principal-symbol}
gives the full nonlocal contribution.  Variations of
\(\log\ell_z+\mathfrak c_\epsilon\) multiply fixed curvature and
are finite-rank terms.  The time derivatives of normal projection and
of \(\mathscr U_{\Gamma_z}\) are multiplication operators on the
normal component, or multiplication followed by the inverse derivative
in \eqref{linearised-chart-constraint-exact}.  These observations give
the first row of \eqref{exact-fixed-space-linear-system}.

For the mean tangential row, use
\(\beta=\int_{\mathbb T}t_z\cdot X\,d\theta\).
The tangential gauge has zero mean in these coordinates, so its
variation has zero mean too.  Every derivative on \(Y^\perp\) in
the remaining tangential row can be integrated by parts against its
smooth coefficient.  The resulting row is bounded at the fixed low
level.  This accounts for the second row, including the defect term.

For clarity, \(\mathscr Q_{\epsilon,z}\) is defined by subtracting
the two displayed diagonal blocks in
\eqref{exact-fixed-space-linear-system} from
\((\partial_t\mathscr S_z)\mathscr S_z^{-1}
 +\mathscr S_z\mathscr L_z\mathscr S_z^{-1}\).
The preceding identities express it as the bounded remainder
\(\Lambda_\epsilon^{-1}R_{\Gamma_z}^0\), the local curvature
multiplications, the lift and mean terms, the time-frame terms, and
the defect in \eqref{off-solution-chart-identity}.

Here is a derivative count for these coefficients.  The kernel
coefficients satisfy the bounds with a reserve of eight derivatives
established in \eqref{explicit-curve-kernel-seminorm-bound}; Lemma
\ref{finite-part-calculus-lemma} transfers these bounds to the
corresponding operators.
The chart and frame comparisons use at most two further derivatives.
The defect uses the reserve eight of the vector field and two
derivatives from \eqref{chart-second-differential-bound}.  One
derivative of the transport is needed in energy.  The time-frame
coefficient uses at most two spatial derivatives of \(\Gamma_t\),
and the time derivative has weight two.  The reserve twelve in
\eqref{corrected-coefficient-reserve} bounds each of these counts.
Products are estimated with one high factor; compositions of the
order-zero coordinate maps are estimated first at level eight and
then at level \(m\).  The low coefficients needed this way are
controlled below level forty.  This proves
\eqref{complete-transport-coefficient-bound} and
\eqref{complete-order-zero-operator-bound} without placing a high
background norm in their constants.

The same factored kernels and differentiated constraint equations
prove the parameter versions.  A time derivative differentiates these
coefficients, including \(r_z\), so its additional cost is exactly
its anisotropic weight.  All estimates are uniform for
\(0<\epsilon\leq e^{-2}\).
\end{proof}

\subsection{A tame inverse with finite derivative loss}

We now use the linearised equation and coefficient bounds of Lemma
\ref{complete-lower-block-decomposition} to construct the smooth tame
inverse required by the Nash--Moser argument.  The key point is to
obtain estimates uniform in \(\epsilon\), with a fixed finite loss of
derivatives independent of the Sobolev index.

Fix
\begin{equation*}
 m_*=40,\qquad \mu=20.
\end{equation*}
These reserves are not optimised.  The coefficient reserve, the chart
comparisons and time reconstruction are included in \(\mu\).

\begin{proposition}[Tame inverse for the linearised curve equation]
\label{NashMoser-linear-inverse-proposition}
There are \(r_*>0\) and \(\epsilon_0>0\), depending only on the
prescribed flow on \([0,T]\), such that, if
\begin{equation}\label{NashMoser-low-ball}
 \|z\|_{\mathfrak X_T^{m_*}}\leq r_*,
 \qquad 0<\epsilon\leq\epsilon_0,
\end{equation}
then every smooth pair \((h,y^{\rm in})\) has a unique smooth solution
of
\begin{equation}\label{NashMoser-linear-IVP}
 Y_t-\mathscr L_zY=h,\qquad Y(0)=y^{\rm in}.
\end{equation}
For \(m\geq m_*\),
\begin{equation}\label{NashMoser-tame-inverse-estimate}
\begin{split}
 \|Y\|_{\mathfrak X_T^m}
 \leq C_{m,T}\bigl\{&
 \|(h,y^{\rm in})\|_{\mathfrak Y_T^{m+\mu}}\\
 &+(1+\|z\|_{\mathfrak X_T^{m+\mu}})
       \|(h,y^{\rm in})\|_{\mathfrak Y_T^{m_*+\mu}}\bigr\}.
\end{split}
\end{equation}
The solution operator \(\Psi_\epsilon(z)(h,y^{\rm in})=Y\) is smooth
tame.  Each fixed number of its base variations has a fixed finite
reserve, independent of \(m\) and \(\epsilon\); the reserves for
different numbers of variations need not equal \(\mu\).
\end{proposition}

\begin{proof}
The proof combines the geometric energy method for skew mean-curvature
flow with the periodic quasilinear Schr\"odinger energy method; see
\cite{feolaiandoli2019,song2021}.  Uniformity in \(\epsilon\) relies on
the self-adjoint scalar structure of the principal operator
\(\mathscr H_{\epsilon,z}\) in
\eqref{exact-fixed-space-linear-system}, whose logarithmic component
comes from \(\mathcal T_\ell^{\mathfrak a}\).  We divide the proof into
spatial estimates, construction and uniqueness, time reconstruction
based on an iterated-operator estimate, and smooth tame dependence.

\emph{Spatial energy.}
Use \eqref{exact-fixed-space-linear-system}.  Choose a constant
\(\lambda_*\) larger than
\[
 K_*:=\frac14\sup_{z,t,\theta}|k_z(\theta,t)|^2
\]
on the low ball, and reduce \(\epsilon_0\) so that
\begin{equation}\label{finite-loss-epsilon-choice}
 \Lambda_\epsilon\geq\max\{2,\lambda_*T\},
 \qquad 0<\epsilon\leq\epsilon_0.
\end{equation}
For any real \(m\geq8\), define
\begin{equation*}
 q_m(t)=m+1-\lambda_*t/\Lambda_\epsilon,
 \qquad
 E_m(t)=\|A^{q_m(t)}w(t)\|_2^2+|\beta(t)|^2.
\end{equation*}
Then \(m\leq q_m(t)\leq m+1\).  Differentiating the exponent produces
\[
 -2\lambda_*\Lambda_\epsilon^{-1}
                 \|L^{1/2}A^{q_m(t)}w\|_2^2.
\]
The principal block contributes zero, even at frequencies where its
combined symbol vanishes.  Apply
\eqref{logarithmic-commutator-bound} and
\eqref{logarithmic-symmetric-matrix} to the logarithmic block, and
the transport and bounded-block estimates to the other terms.  This
gives
\begin{equation}\label{NashMoser-spatial-linear-energy}
\begin{split}
 \frac12E_m'
 &+\frac{\lambda_*-K_*}{\Lambda_\epsilon}
                    \|L^{1/2}A^{q_m(t)}w\|_2^2\\
 &\leq C_m E_m+C_m\bigl\{
       \|\mathscr S_zh\|_{H^{m+1}\times\mathbb R}
       +(1+\mathcal Z_{m+1}(t))\|W\|_{H^8\times\mathbb R}
                                      \bigr\}E_m^{1/2}.
\end{split}
\end{equation}
The constants in the bounded terms may depend on \(m\); the
coefficient \(K_*\) in the logarithmic term does not.  This distinction
is what makes the loss independent of the output index.

For \(m=8\), the coefficient \(\mathcal Z_9\) is controlled by
the fixed \(\mathfrak X_T^{40}\) ball.  Thus the low estimate closes
by Gronwall.  Using the low coordinate equivalence, it gives
\begin{equation}\label{NashMoser-low-spatial-estimate}
 \sup_t\|Y(t)\|_{H^8}
 \leq C_T\left\{\|y^{\rm in}\|_{H^9}
                         +\|h\|_{L^1_tH^9}\right\}.
\end{equation}
Put \(D_j=\|y^{\rm in}\|_{H^j}+\|h\|_{L^1_tH^j}\).
At a higher level the last term of
\eqref{NashMoser-spatial-linear-energy} is a known source, bounded by
the low estimate.  Apply Gronwall to \(E_m^{1/2}\), then use the
tame coordinate estimates for the initial data, forcing and output.
For \(m\geq8\) this yields
\begin{equation}\label{NashMoser-high-spatial-estimate}
 \sup_t\|Y(t)\|_{H^m}
 \leq C_{m,T}\{D_{m+1}
                 +(1+\|z\|_{\mathfrak X_T^{m+13}})D_9\}.
\end{equation}
At the low level we use \eqref{NashMoser-low-spatial-estimate};
the higher background norms enter through the source term.

\smallskip
\emph{Construction and uniqueness.}
For \(0<\delta\leq1\), add
\(\delta\operatorname{diag}(\partial_\theta^4,1)W\) to the left
side of \eqref{exact-fixed-space-linear-system}, and apply orthogonal
Fourier--Galerkin truncation to every term.  The projectors commute
with \(A\), \(L\) and the principal block.  The additional energy
contribution is
\[
 \delta\|\partial_\theta^2A^{q_m(t)}w\|_2^2+\delta|\beta|^2.
\]
It is nonnegative.  The spatial estimates are therefore uniform in
both the truncation and \(\delta\); the finite-dimensional solutions
exist on \([0,T]\).  For smooth data they are bounded at every
spatial level.  The regularised equation bounds one time derivative
in \(H^q\) from a bound in \(H^{q+4}\).  Rellich compactness and
Arzel\`a--Ascoli, followed by a diagonal subsequence, first remove the
Galerkin cutoff and then \(\delta\), with convergence in
\(C_tH^{q+2}\) for every fixed \(q\).

All nonlocal terms are continuous under this convergence: the
principal multiplier maps \(H^{q+3}\) to \(H^q\), the logarithmic
matrix maps \(H^{q+1}\) to \(H^q\), and the remainder is bounded
as proved above.  The viscosity term tends to zero using its uniform
\(H^{q+4}\) bound.  The limit solves
\eqref{NashMoser-linear-IVP} and retains the spatial estimates.
Equation \eqref{NashMoser-spatial-linear-energy}, at the low level
and with zero data, proves uniqueness.  Higher time derivatives are
reconstructed from the limiting equation in the next step.

\smallskip
\emph{Time reconstruction.}
For \(r\geq1\) and \(\boldsymbol\nu=(\nu_1,\ldots,\nu_r)\), put
\[
 \mathscr W_{\boldsymbol\nu}
 =(\partial_t^{\nu_1}\mathscr L_z)\cdots
                       (\partial_t^{\nu_r}\mathscr L_z),
 \qquad d=q+2r+2|\boldsymbol\nu|.
\]
Here the product denotes composition, with the rightmost operator
acting first.  Each time derivative differentiates the coefficients
of \(\mathscr L_z\).
For even integers \(q\geq0\), this composition satisfies
\begin{equation}\label{iterated-linearised-composition-bound}
 \|\mathscr W_{\boldsymbol\nu}V\|_{H^q}
 \leq C_{q,r,\boldsymbol\nu,T}
       \left\{\|V\|_{H^{d+1}}
          +(1+\|z\|_{\mathfrak X_T^{d+18}})\|V\|_{H^8}\right\}.
\end{equation}
For an odd output index, apply this estimate at \(q+1\).
Both high indices contain \(2|\boldsymbol\nu|\).  Although the
operator has spatial order at most \(2r\) with \(r\) logarithms,
time derivatives of a coefficient must be included when distributing
the high norm between that coefficient and \(V\).

Here are the product estimates that justify this assertion.  Expand
this operator composition using \eqref{exact-fixed-space-linear-system} and the
coordinate maps.  Its factors are scalar multipliers, differential
operators, multiplication operators, and the bounded kernels in
\eqref{explicit-near-curve-remainder} and their far counterparts.
For an integral kernel, the spatial Fourier matrix is dominated by
the convolution sequence
\[
 c_k(t)=\int|\widehat K(t,k,\xi)|\,d\xi,
 \qquad
 \left|\int\widehat K(t,k,\xi)e^{2\pi i n\xi}\,d\xi\right|
       \leq c_k(t).
\]
This reduces integral compositions to convolution estimates.
Translations satisfy the same bounds.  Coefficient derivatives satisfy
\eqref{explicit-curve-kernel-seminorm-bound}.
The scalar multipliers obey
\eqref{fixed-circle-principal-multipliers}, including derivatives of
their length and holonomy parameters.

The total additional derivative weight in an operator composition is at most \(d\),
where a time derivative has weight two.  Every intermediate spatial
frequency is bounded by the sum of the input and coefficient
frequencies.  The entire logarithmic product is absorbed once, using
\[
 (1+\Lambda_\epsilon^{-1}\log\langle n\rangle)^r
       \leq C_r\langle n\rangle .
\]
This costs one spatial derivative on whichever factor receives the
largest frequency.  The coefficient reserve twelve from
Lemma \ref{complete-lower-block-decomposition}, with two further
derivatives for the coordinate comparisons, bounds the coefficient
seminorms before distributing these additional derivatives.  The
coordinate maps themselves have reserve two; a principal spatial
factor differentiates them at most twice, which is within this count.
One derivative for the logarithmic product and one for the spatial
Fourier \(\ell^1\) bound give reserve sixteen.  Rounding an
intermediate background grade up to an even integer uses at most one
more derivative.  Thus, writing \(\lceil a\rceil_2=2\lceil a/2\rceil\),
a term with additional weights \(d_0,d_1,\ldots,d_b\), of sum at
most \(d\), is bounded by
\[
 C\,\|V\|_{H^{d_0+1}}
     \prod_{i=1}^b(1+\|z\|_{\mathfrak X_T^{\lceil d_i+16\rceil_2}}).
\]
The even high background grade \(d+18\) covers these reserves.
The same description applies to the integral formulas for the chart
maps: their fixed-interval integrations are bounded, and
differentiating their coefficients has the counts just stated.

We spell out the interpolation needed to reduce this product to
one high factor.  The smoothing estimates
\eqref{NashMoser-smoothing-estimates}, with the cutoff optimised,
give interpolation between any two of the indicated grades.
The low background norm is bounded at level forty.  If \(d>22\),
the interpolation exponents of the background factors satisfy
\[
 \alpha_i=\frac{(\lceil d_i+16\rceil_2-40)_+}{d-22}
               \leq\frac{(d_i-23)_+}{d-22},
 \qquad
 \alpha_0=\frac{(d_0-7)_+}{d-7}
\]
for the field factor.  If some background factor is above its low
grade, the sum is at most
\[
 \frac{d-d_0-23}{d-22}
                  +\frac{(d_0-7)_+}{d-7}\leq1.
\]
If every background factor is low, only \(\alpha_0\leq1\)
is needed.  For \(d\leq22\) all background factors are low;
indices on \(V\) below eight are bounded by its \(H^8\) norm.
Young's inequality now proves
\eqref{iterated-linearised-composition-bound}.  After any fixed number of
parameter differentiations the same allocation gives one high
background or direction factor, with all other direction factors
at their fixed low grades.

Repeated substitution of \(Y_t=\mathscr L_zY+h\) gives
\[
 \partial_t^jY
 =\sum_{r+|\boldsymbol\nu|=j}c_{\boldsymbol\nu}
                         \mathscr W_{\boldsymbol\nu}Y
  +\sum_{r+|\boldsymbol\nu|+a=j-1}d_{\boldsymbol\nu,a}
                         \mathscr W_{\boldsymbol\nu}\partial_t^ah.
\]
The first sum has \(r\geq1\).  The second sum also includes the identity
operator, corresponding to \(r=0\).
Fix even \(m\geq40\) and take \(q=m-2j\).  The first sum uses
\(Y\) in \(H^{m+1}\) and the background at grade at most \(m+18\).
The lower spatial derivatives are controlled by this same Sobolev norm.
For a term with \(r\geq1\) in the forcing sum, set
\(d=q+2r+2|\boldsymbol\nu|\).  Then
\[
 d\geq2,\qquad d+2a=q+2j-2\leq m-2.
\]
The high field term in the composition estimate is bounded by
\(\|h\|_{\mathfrak X_T^{m-1}}\).  Its other term is bounded by
\[
 (1+\|z\|_{\mathfrak X_T^{d+18}})
                         \|h\|_{\mathfrak X_T^{2a+8}}.
\]
Interpolation between the even background grades
\(40,m+18\) and forcing grades \(16,m+4\) gives
\begin{equation}\label{time-forcing-product-interpolation}
\begin{split}
 &(1+\|z\|_{\mathfrak X_T^{d+18}})
                         \|h\|_{\mathfrak X_T^{2a+8}}\\
 &\qquad\leq C_m\left\{
 \|h\|_{\mathfrak X_T^{m+4}}
 +(1+\|z\|_{\mathfrak X_T^{m+18}})
                         \|h\|_{\mathfrak X_T^{16}}\right\}.
\end{split}
\end{equation}
Indeed, when both factors exceed their low grades, the exponents are
\[
 \frac{d-22}{m-22},\qquad \frac{2a-8}{m-12}.
\]
Their sum is at most \((m-32)/(m-22)<1\).  If one factor is
low, the other is within its displayed high grade:
\(d+18\leq m+16\) and \(2a+8\leq m+4\).
This proves \eqref{time-forcing-product-interpolation} in all cases.
The term with \(r=0\) gives \(\partial_t^{j-1}h\) directly and has no
background factor.

We have therefore proved
\begin{equation}\label{time-reconstruction-from-spatial-estimates}
\begin{split}
 \|Y\|_{\mathfrak X_T^m}\leq C_{m,T}\bigl\{&
       \sup_t\|Y(t)\|_{H^{m+1}}+\|h\|_{\mathfrak X_T^{m+4}}\\
 &+(1+\|z\|_{\mathfrak X_T^{m+18}})
       (\sup_t\|Y(t)\|_{H^8}+\|h\|_{\mathfrak X_T^{16}})\bigr\}.
\end{split}
\end{equation}
This estimate has been proved for even \(m\geq40\).
The same substitution at \(t=0\) controls the initial time jets.
Combining it with \eqref{NashMoser-low-spatial-estimate} and
\eqref{NashMoser-high-spatial-estimate} at level \(m+1\) uses
data at most at grade \(m+4\), background at grade \(m+18\),
and low data at grade at most sixteen.  For odd \(m\), use the
result at \(m+1\): the data grade is at most \(m+5\) and the
background grade at most \(m+19\).  Hence
\eqref{NashMoser-tame-inverse-estimate} holds at every integer
\(m\geq40\), with \(m_*=40\) and \(\mu=20\).

\smallskip
\emph{Smooth tame dependence.}
For two backgrounds, uniqueness gives the resolvent identity
\begin{equation*}
 (\Psi_\epsilon(z_1)-\Psi_\epsilon(z_0))d
 =\Psi_\epsilon(z_1)
       ((\mathscr L_{z_1}-\mathscr L_{z_0})\Psi_\epsilon(z_0)d,0).
\end{equation*}
The tame bounds already proved give continuity at every output level
from a fixed higher input level.  Apply this identity to a difference
quotient, use the differentiated coefficient estimates, and take the
limit in that higher topology.  This proves the first derivative.
For a nonempty set of directions \(I\), the subsequent derivatives solve
\begin{equation*}
 (\partial_t-\mathscr L_z)Y^{[I]}
 =\sum_{\varnothing\ne J\subseteq I}
       D_z^{|J|}\mathscr L_z[z_j:j\in J]Y^{[I\setminus J]},
 \qquad Y^{[I]}(0)=0.
\end{equation*}
Induction, using the resolvent identity for the difference quotients,
proves the higher derivative formulae.

Each inverse loses the fixed reserve twenty and each differentiated
coefficient map has reserve at most twelve by
\eqref{curve-coefficient-family-tame}.  We use the even subgrading
in these compositions.  Tame interpolation
reduces products of background norms to one high background norm.
An explicit generous choice for \(j\) inverse variations is
\begin{equation}\label{polarised-inverse-reserves}
 \mu_0=20,\qquad \mu_{j+1}=2\mu_j+40.
\end{equation}
Indeed the high source term uses at most the preceding reserve plus
twenty and twelve.  In the low source term multiplied by a high
background norm, interpolation adds at most one more copy of the
preceding reserve; the extra forty bounds the two operator reserves
and the fixed trace and product shifts.  This proves a tame estimate
with all data and direction norms at levels \(m+\mu_j\) and
\(m_*+\mu_j\), one high factor, and the background term at
\(m+\mu_j\).  An odd output grade is rounded up once; the gap
between the operator costs \(20+12\) and the added forty covers
this extra shift.  The constants may depend on \(j\), but the reserve
does not depend on \(m\) or \(\epsilon\).  This is smooth tameness
on the fixed low ball and completes the proof.
\end{proof}

\subsection{Proof of Theorem \ref{uniformlocalcurvetheorem}}

\begin{proof}[Proof of Theorem \ref{uniformlocalcurvetheorem}]
Put $\Lambda_\epsilon=|\log\epsilon|$.  We apply the Nash--Moser theorem to
the map \eqref{reduced-vector-field-identity} on the fixed time interval $[0,T]$.
Its derivative is
\begin{equation}\label{NashMoser-derivative-map}
 D\Phi_\epsilon(z)Y
 =\big(\partial_tY-D_z\mathscr G_\epsilon(t,z)Y,Y(0)\big).
\end{equation}
Proposition \ref{NashMoser-linear-inverse-proposition} says precisely that
\eqref{NashMoser-derivative-map} has the smooth tame inverse
$\Psi_\epsilon(z)$ required by the Nash--Moser inverse function theorem.
The neighbourhood \eqref{NashMoser-low-ball} and all tame constants may be
chosen independently of $0<\epsilon\leq\epsilon_0$, with
\eqref{finite-loss-epsilon-choice} imposed first.

At the origin of the chart the prescribed curve satisfies binormal flow, so
the local part of the graph equation vanishes.  The corrected part need not
vanish, and
\begin{equation*}
 \Phi_\epsilon(0)
 =\left(-\frac1{\Lambda_\epsilon}p_\epsilon^0,0\right),
\end{equation*}
where $p_\epsilon^0$ is obtained by evaluating
\eqref{Repsilongeometric}, the tangential gauge, and the length equation on
$\gamma^0$.  Smoothness of $\gamma^0$ and
\eqref{coreconstantasymptotics} give, for every $m$,
\begin{equation}\label{NashMoser-origin-residual-bound}
 \|\Phi_\epsilon(0)\|_{\mathfrak Y_T^m}
 \leq\frac{C_{m,T}}{\Lambda_\epsilon}.
\end{equation}

For clarity, the corresponding Newton--Moser scheme starts from $z_0=0$.
At each stage it smooths the current base point and residual, applies
$\Psi_\epsilon$ to solve the resulting linearised initial-value problem,
and smooths the correction according to a rapidly increasing sequence of
cutoff parameters.  The smoothing estimates
\eqref{NashMoser-smoothing-estimates} recover the fixed loss $\mu$ in
\eqref{NashMoser-tame-inverse-estimate}, while the accumulated substitution
errors and the quadratic Newton errors are summable.  This is the standard
construction in the Nash--Moser theorem
\cite{hamilton1982}; all of its hypotheses, including smooth tame dependence
of the inverse family, have been verified above.

Choose $\epsilon_0$ so small that the right-hand side of
\eqref{NashMoser-origin-residual-bound} at the finite low index required by
the iteration lies in the uniform Nash--Moser neighbourhood.  The iteration
then converges in $\mathfrak X_T$ to a smooth $z_\epsilon$ satisfying
\begin{equation}\label{NashMoser-solved-map}
 \Phi_\epsilon(z_\epsilon)=(0,0).
\end{equation}
The quantitative tame estimate in the iteration and
\eqref{NashMoser-origin-residual-bound} give, for every $m$,
\begin{equation}\label{NashMoser-solution-bound}
 \|z_\epsilon\|_{\mathfrak X_T^m}
 \leq\frac{C_{m,T}}{\Lambda_\epsilon}.
\end{equation}
Recover $(f_1,f_2,f_3,\lambda-1)$ using
\eqref{NashMoser-recovery-map}.  Its tame estimates turn
\eqref{NashMoser-solution-bound} into \eqref{uniformgraphestimate}, and the
initial-data component of \eqref{NashMoser-solved-map} gives
\eqref{zerographinitialdata}.

The recovered curve satisfies the exact constant-speed constraint at every
time.  Equivalently, \eqref{lambdatimeequation} and
\eqref{atangentialequation} give
\[
 \partial_t|\partial_s\gamma|
 =\partial_sa-\lambda\vec K\cdot\vec V_\epsilon
 =-\frac\lambda L\int_0^L\vec K\cdot\vec V_\epsilon\,ds
 =\partial_t\lambda.
\]
Periodicity preserves closure.  For small $\epsilon$ the low-norm case of
\eqref{NashMoser-solution-bound}, Sobolev embedding, and the uniform tubular
radius of $\gamma^0$ preserve embeddedness on $[0,T]$.

The Nash--Moser inverse is locally one-to-one.  This gives local
uniqueness for the nonlinear initial-value problem as well.  More
explicitly, restart the same inverse argument on a short interval at
any time at which two smooth solutions agree.  The equation determines
the same initial time jets for both solutions.  Their difference is
therefore small in the finite low space--time norm on a sufficiently
short interval, and local injectivity applies there.

To pass from this local statement to the uniqueness asserted in the
theorem, let $\bar\gamma_\epsilon$ be any smooth solution on $[0,T]$ with
the same initial curve and satisfying the same constant-speed gauge.  Its
normal-graph coordinate $\bar z_\epsilon$ belongs to the low ball on a
short initial interval, since $\bar z_\epsilon(0)=0$ and the equation
determines the same initial time jets as for $z_\epsilon$.  The preceding
local uniqueness therefore gives $\bar z_\epsilon=z_\epsilon$ there.  If
the interval of equality had a first finite endpoint, then
\eqref{NashMoser-solution-bound} would place both solutions strictly inside
the low ball at that endpoint.  Smoothness and the local uniqueness
argument would extend their equality beyond it, a contradiction.  Thus
$\bar\gamma_\epsilon=\gamma_\epsilon$ on all of $[0,T]$.  This completes
the proof.
\end{proof}

\subsection{Forced tame solution families}

For later use, we record a forced, parameter-dependent extension of
Theorem \ref{uniformlocalcurvetheorem}, together with the corresponding
tame difference and evaluation-map properties.

\begin{corollary}[Forced tame solution families]
\label{forcedNashMoserfamilycorollary}
There are fixed finite integers \(s_0\geq m_*\), \(\nu\geq\mu\)
and a number \(\eta_0>0\), independent of sufficiently small
\(\epsilon\), with the following property.  For smooth data
\(d=(r,z^{\rm in})\) satisfying
\(\|d\|_{\mathfrak Y_T^{s_0+\nu}}\leq\eta_0\), the problem
\[
 \partial_tz=\mathscr G_\epsilon(t,z)+r,
 \qquad z(0)=z^{\rm in}
\]
has a unique smooth solution in the local solution neighbourhood.
It depends smoothly and tamely on \(d\).  For two such data pairs
\(d,\bar d\), after decreasing \(\eta_0\) if necessary, their
solutions satisfy, for \(m\geq s_0\),
\begin{equation}\label{forcedNashMoser-difference-estimate}
\begin{split}
 \|z-\bar z\|_{\mathfrak X_T^m}
 \leq C_{m,T}\bigl\{&\|d-\bar d\|_{\mathfrak Y_T^{m+\nu}}\\
 &+(1+\|d\|_{\mathfrak Y_T^{m+\nu}}
          +\|\bar d\|_{\mathfrak Y_T^{m+\nu}})
             \|d-\bar d\|_{\mathfrak Y_T^{s_0+\nu}}\bigr\}.
\end{split}
\end{equation}
Every fixed number of parameter variations has a finite enlarged tame
reserve.  The nonlinear solution-map reserve \(\nu\) is allowed to
exceed the linear inverse reserve \(\mu\).
\end{corollary}

\begin{proof}
Apply the Nash--Moser inverse function theorem to \(\Phi_\epsilon\),
using Proposition \ref{NashMoser-linear-inverse-proposition}.  Its local
inverse is smooth tame, and its domain contains a convex ball specified
by finitely many data seminorms.  Enlarge the highest such index to
\(s_0+\nu\).  The origin residual tends to zero in every fixed
seminorm by \eqref{NashMoser-origin-residual-bound}, so a fixed smaller
data ball about zero is available for all sufficiently small
\(\epsilon\).

To prove the difference estimate, use the data segment
\(d_\tau=\bar d+\tau(d-\bar d)\), and let \(z_\tau\) be its
image under the local inverse.  Then
\[
 \partial_\tau z_\tau=\Psi_\epsilon(z_\tau)(d-\bar d).
\]
The tame bound for the local inverse, followed by
\eqref{NashMoser-tame-inverse-estimate}, gives
\eqref{forcedNashMoser-difference-estimate} after increasing the
fixed reserve \(\nu\).  The data norms on this segment are bounded
by the endpoint data norms.  Integrating in \(\tau\) completes the
proof.  Differentiation of this identity gives the parameter
assertion with finite reserves for each fixed number of variations.
\end{proof}

\subsection{Conservation of vector area}

For a smooth oriented closed curve we define its vector area by
\begin{equation}\label{vectorareadefinition}
    \mathcal A_{\rm vec}(\gamma)
    :=\frac12\int_0^{\ell}\gamma(\sigma)\times
       \partial_\sigma\gamma(\sigma)\,d\sigma
    =\frac12\oint_\gamma\gamma\times d\gamma.
\end{equation}
This quantity is invariant under reparametrisation and translation of the
curve.

\begin{lemma}[Vector-area conservation]\label{vectorareaconservationlemma}
Let $\gamma$ be a smooth embedded closed solution of the corrected curve
equation \eqref{geometriccurveequation}, with normal velocity given by
\eqref{normalvelocitydefinition}.  Then
\begin{equation}\label{vectorareaconservation}
    \frac{d}{dt}\mathcal A_{\rm vec}(\gamma(t))=0.
\end{equation}
In particular, this identity is unaffected by variation of the length or by
the tangential velocity used to preserve the constant-speed gauge.
\end{lemma}

\begin{proof}
We work at a fixed time and use intrinsic arclength $\sigma$.  Differentiating
\eqref{vectorareadefinition} in any fixed periodic parametrisation, and then
returning to arclength, gives
\begin{equation}\label{vectorareavariation}
    \frac{d}{dt}\mathcal A_{\rm vec}(\gamma(t))
    =\int_0^\ell\partial_t\gamma\times\vec{\mathsf t}\,d\sigma.
\end{equation}
The tangential velocity makes no contribution.  Moreover,
\[
    (\vec{\mathsf t}\times\vec K)\times\vec{\mathsf t}
    =\vec K,
    \qquad
    \int_0^\ell\vec K\,d\sigma
    =\int_0^\ell\partial_\sigma\vec{\mathsf t}\,d\sigma=0.
\]
Since $c_\epsilon(\ell)$ is independent of $\sigma$, it remains only to show
\begin{equation}\label{nonlocalvectorareacancellation}
    \int_0^\ell
    \mathcal N(\gamma)(\sigma)\times\vec{\mathsf t}(\sigma)\,d\sigma=0.
\end{equation}
The normal projection in \eqref{normalvelocitydefinition} may be omitted
here, because the tangential part of $\mathcal N(\gamma)$ vanishes after
crossing with $\vec{\mathsf t}$.

Let $d_\ell(\sigma,\sigma')$ denote periodic arclength distance, and set
\[
    X=\gamma(\sigma)-\gamma(\sigma'),\qquad
    \vec{\mathsf t}=\vec{\mathsf t}(\sigma),\qquad
    \vec{\mathsf t}'=\vec{\mathsf t}(\sigma').
\]
For $\delta>0$ define the symmetrically truncated integral
\begin{equation}\label{Jdeltavectorarea}
\begin{split}
    J_\delta
    :=\iint_{d_\ell(\sigma,\sigma')>\delta}
    \frac{(\vec{\mathsf t}'\times X)\times\vec{\mathsf t}}
         {|X|^3}\,d\sigma'\,d\sigma.
\end{split}
\end{equation}
The counterterm in \eqref{ngammadef} is
\[
    -\frac{\kappa\vec b}{2d_\ell(\sigma,\sigma')}
    =-\frac{\vec{\mathsf t}\times\vec K}
            {2d_\ell(\sigma,\sigma')}.
\]
After crossing with $\vec{\mathsf t}$, its truncated integral in $\sigma'$
is a scalar independent of $\sigma$ times $\vec K(\sigma)$.  Its subsequent
integral in $\sigma$ is therefore zero.  Hence the left-hand side of
\eqref{nonlocalvectorareacancellation} equals
$\lim_{\delta\downarrow0}J_\delta$.

Put
\[
    F(\sigma,\sigma')
    :=\frac{(\vec{\mathsf t}'\times X)\times\vec{\mathsf t}}
            {|X|^3}.
\]
The vector triple-product identity, followed by differentiation of
$|X|^{-1}$, yields
\begin{equation}\label{symmetrisedvectorareakernel}
\begin{split}
    F(\sigma,\sigma')+F(\sigma',\sigma)
    &=
    \partial_\sigma\left(\frac{\vec{\mathsf t}'}{|X|}\right)
    +\partial_{\sigma'}\left(\frac{\vec{\mathsf t}}{|X|}\right).
\end{split}
\end{equation}
Since the domain in \eqref{Jdeltavectorarea} is symmetric under interchange
of $\sigma$ and $\sigma'$, equation \eqref{symmetrisedvectorareakernel} and
integration over the boundary of the deleted diagonal strip give
\begin{equation}\label{Jdeltaboundaryvectorarea}
    J_\delta
    =\int_0^\ell
      \frac{\vec{\mathsf t}(\sigma+\delta)
            -\vec{\mathsf t}(\sigma)}
           {|\gamma(\sigma+\delta)-\gamma(\sigma)|}\,d\sigma.
\end{equation}
As $\delta\downarrow0$, the integrand in
\eqref{Jdeltaboundaryvectorarea} converges uniformly to $\vec K(\sigma)$.
Consequently,
\[
    \lim_{\delta\downarrow0}J_\delta
    =\int_0^\ell\vec K(\sigma)\,d\sigma=0,
\]
which proves \eqref{nonlocalvectorareacancellation}.  Substituting this into
\eqref{vectorareavariation} proves \eqref{vectorareaconservation}.
\end{proof}

\section{First global improvement of the approximation}
\label{firstglobalimprovementsection}

The aim of this section is to construct a first global improvement of the
initial ansatz \eqref{firstansatz}.  We take the centre curve to be the
solution of the corrected curve equation \eqref{geometriccurveequation},
which enforces the leading-order orthogonality conditions against the two
translation modes.  An inner correction, obtained by inverting the planar
linearised Ginzburg--Landau operator, then eliminates the part of the source
orthogonal to these modes.  We then solve a three-dimensional wave equation
to construct an outer correction that removes the leading error away from
the filament.  The construction follows the inner--outer strategy of
\cite{delpinojunemanmusso2025}, adapted here to a closed filament in three
space dimensions.

Throughout this section, let $\gamma$ be the corrected curve supplied by Theorem
\ref{uniformlocalcurvetheorem} and set
\begin{equation*}
    \Lambda_\epsilon:=|\log\epsilon|,
    \qquad R_\epsilon:=\frac{\delta}{\epsilon},
\end{equation*}
where $\delta>0$ is chosen so that $|z|<4\delta$ lies in a fixed embedded
tubular neighbourhood of $\gamma(t)$ for every $t\in[0,T]$.  All constants
are uniform for $0<\epsilon<\epsilon_0$ and $t\in[0,T]$.

We shall also need the construction uniformly for initial curves in a tame
constant-speed chart ball $\mathscr O_{\rm in}$ about the prescribed initial
curve, chosen sufficiently small that Theorem
\ref{uniformlocalcurvetheorem} supplies a common chord--arc constant and
tubular radius.  For $\Gamma_{\rm in}\in\mathscr O_{\rm in}$, write
$\gamma[\Gamma_{\rm in}]$ for the corresponding corrected curve and
$D_{\rm in}$ for differentiation with respect to $\Gamma_{\rm in}$.  The
periodic frame, tubular projection, solid-angle phase and cutoffs are chosen
uniformly on this chart ball.  Whenever tubular coordinates are used, $s$
denotes intrinsic arclength on $\mathbb R/\ell(t)\mathbb Z$.

\subsection{The global gluing ansatz}

To construct the global gluing ansatz, we first regularise the relative phase
near the filament.  To do this consistently around the closed curve, we
choose a periodic normal frame in which the difference between the
solid-angle phase and the polar angle has a single-valued real
representative on the punctured tube.  The following lemma provides the
required choice of frame.  In the lemma, $r$ denotes physical normal
distance.

\begin{lemma}[A periodic frame and a real relative phase]
\label{longitudinal-framing-lemma}
Let $\gamma(s,t)$ be the smooth family of embedded closed curves under
consideration, parametrised by arclength
$s\in\mathbb R/\ell(t)\mathbb Z$, and let $e^{i\varphi}$ be the
solid-angle phase used in the initial ansatz.  There exists a smooth
periodic oriented orthonormal normal frame $(\vec e_1,\vec e_2)$ whose
normal connection coefficient
\[
 \mathfrak a(t):=\partial_s\vec e_1(s,t)\cdot\vec e_2(s,t)
\]
is independent of $s$, and for which the relative phase has a smooth,
single-valued real representative on the punctured tube.

More precisely, writing
\[
 x(s,r,\theta,t)
 =\gamma(s,t)
 +r\bigl(\cos\theta\,\vec e_1(s,t)
        +\sin\theta\,\vec e_2(s,t)\bigr),
\]
there exists a smooth real-valued function
$\Phi_\gamma(s,r,\theta,t)$, defined for $0<r<4\delta$, such that
\[
 B_\gamma(s,r,\theta,t)
 :=\exp\!\bigl(i\varphi(x(s,r,\theta,t),t)-i\theta\bigr)
 =e^{i\Phi_\gamma(s,r,\theta,t)},
\]
with
\[
 \Phi_\gamma(s+\ell(t),r,\theta,t)
 =\Phi_\gamma(s,r,\theta+2\pi,t)
 =\Phi_\gamma(s,r,\theta,t).
\]
The frame and the real representative can be chosen smoothly in time
and in the curve parameters on a sufficiently small contractible chart.
For a fixed frame, choosing the value of $\Phi_\gamma$ at one base point
fixes its additive $2\pi\mathbb Z$ ambiguity.
\end{lemma}

\begin{proof}
Fix $t$ and write $L=\ell(t)$.  Begin with a smooth $L$-periodic
oriented normal frame whose connection coefficient
$\mathfrak a_0(t)$ is independent of $s$, obtained by the construction
preceding \eqref{periodicframeintro}.  Let $B_0$ be its relative phase,
defined by the same formula as $B_\gamma$.  For fixed $s$ and $r$, the
winding number of $\theta\mapsto B_0(s,r,\theta,t)$ is zero, since
$e^{i\varphi}$ and $e^{i\theta}$ both have degree one around a normal
circle.  Let $n\in\mathbb Z$ be the winding number of
$s\mapsto B_0(s,r,\theta,t)$ for fixed $r$ and $\theta$.  It is
independent of these choices because the corresponding loops are
homotopic within the punctured tube.

Rotate this frame by the angle
\[
 \alpha(s)=-\frac{2\pi n}{L}s.
\]
The rotated frame is periodic because
$\alpha(s+L)-\alpha(s)=-2\pi n$.  The old polar angle is
$\theta+\alpha(s)$, so the relative phase in the rotated frame is
\[
 B_\gamma(s,r,\theta,t)
 =e^{i\alpha(s)}B_0(s,r,\theta+\alpha(s),t).
\]
Its angular winding number is still zero.  Its longitudinal winding
number is $n-n=0$, since the angular winding number of $B_0$ is zero.
The normal connection coefficient becomes
\[
 \mathfrak a(t)=\mathfrak a_0(t)+\partial_s\alpha
 =\mathfrak a_0(t)-\frac{2\pi n}{L},
\]
which remains independent of $s$.

The two generators of the fundamental group of the punctured tube,
corresponding to one circuit in $s$ and one in $\theta$, therefore
have zero winding under $B_\gamma$.  The lifting criterion for
$\mathbb R\to S^1$ gives a smooth real representative $\Phi_\gamma$
with both asserted periodicities.  Any two such representatives differ
by a constant in $2\pi\mathbb Z$, so a value at one base point fixes
this ambiguity.  This frame is also in the zero self-linking, or
Seifert, class: for fixed $r$ and $\theta$, the nearby closed curve
$s\mapsto x(s,r,\theta,t)$ has linking number zero with $\Gamma_t$,
because its period of $d\varphi/(2\pi)$ is the longitudinal winding
number just computed; see \eqref{solidanglegradient}.

The preliminary periodic frame can be chosen smoothly in time and in
the curve parameters by the normal parallel-transport construction
used in Lemma \ref{tameperiodicnormaltrivialisation}.  The integer $n$
is locally constant under these variations, so the rotation above also
depends smoothly on the parameters.  After shrinking the curve chart
to be contractible, the chart and time interval introduce no further
periods.  Choosing the real representative at one base point and
lifting along paths therefore gives a smooth family of
$\Phi_\gamma$ on the whole punctured tube.
\end{proof}

We now use the frame and real representative supplied by Lemma
\ref{longitudinal-framing-lemma} to regularise the phase near the
filament.  In inner coordinates $y=z/\epsilon$, put $r=|y|$ and set
\begin{equation*}
 \Phi_\epsilon(s,y,t)
 :=\Phi_\gamma(s,\epsilon|y|,\theta(y),t),
 \qquad y\neq0.
\end{equation*}
The periodicity in $\theta$ makes this independent of the choice of
polar angle.  This fixes the real representative of the relative phase
in \eqref{nablaphistar}.  It is smooth for $y\neq0$.  Integrating
\eqref{phistargradientexpansion} in $y$ gives
\begin{equation}\label{phasejetfirstimprovement}
\begin{split}
 \Phi_\epsilon(s,y,t)
 ={}&\beta_\epsilon(s,t)
 +\frac{\epsilon}{2}
      \big(J\vec\kappa(s,t)\cdot y\big)
      \log\frac{\ell(t)}{\epsilon|y|}
 +\frac{\epsilon}{2}
      \Pi_{\vec{\mathsf t}}^\perp\mathcal N(\gamma)[s,t]\cdot y\\
 &+O\big(\epsilon^2|y|^2(1+|\log(\epsilon|y|)|)\big).
\end{split}
\end{equation}
Here
\[
 \beta_\epsilon(s,t):=\lim_{y\to0}\Phi_\epsilon(s,y,t)
\]
is independent of the direction of approach, periodic in $s$, and smooth
in time and in the curve parameters.  The term
$(J\vec\kappa\cdot y)\log|y|$ need not be smooth at $y=0$.
We therefore replace $\Phi_\epsilon$ by a smooth phase that agrees
with it for $r\geq2$.

Choose $\eta_0\in C_c^\infty([0,\infty))$ and a radial function
$\eta_{\rm in}\in C_c^\infty(\mathbb R^2)$ such that
\begin{equation}\label{cutoffsfirstimprovement}
\begin{aligned}
 &\eta_0(r)=1\quad\hbox{for }0\leq r\leq1,
 &&\eta_0(r)=0\quad\hbox{for }r\geq2,\\
 &\eta_{\rm in}(z)=1\quad\hbox{for }|z|\leq\tfrac12\delta,
 &&\eta_{\rm in}(z)=0\quad\hbox{for }|z|\geq\tfrac34\delta.
\end{aligned}
\end{equation}
Choose a smooth radial function $\zeta$ which is zero for $r\leq1$
and one for $r\geq2$.  Define the regularised phase by
\begin{equation}\label{regularphaseextension}
 \widehat\Phi_\epsilon
 :=\beta_\epsilon+\zeta(r)(\Phi_\epsilon-\beta_\epsilon).
\end{equation}
Thus $\widehat\Phi_\epsilon=\beta_\epsilon$ on $B_1$ and
$\widehat\Phi_\epsilon=\Phi_\epsilon$ for $r\geq2$.
It is smooth in Cartesian normal coordinates, periodic in $s$, and smooth
in time and in the curve parameters.
For later calculations put
\begin{equation}\label{core-phase-scaled-variables}
 P_\epsilon:=\epsilon^{-1}(\Phi_\epsilon-\beta_\epsilon),
 \qquad p_\epsilon:=\zeta P_\epsilon,
 \qquad q_\epsilon:=p_\epsilon-P_\epsilon.
\end{equation}
The last function is supported in $B_2$ and includes the complete
logarithmic jet.  The phase
$\widehat\Phi_\epsilon=\beta_\epsilon+\epsilon p_\epsilon$ is smooth,
and all fixed-order derivatives of $p_\epsilon$ on $B_3$ are bounded by
$C_m\Lambda_\epsilon^{C_m}$.  Radial cutoffs preserve angular averages
and commute with parameter differentiation at fixed normal coordinates.

Let $\phi_{1,\epsilon}=\phi_{1,\epsilon}(s,y,t)$ denote an inner correction
to be constructed in the next subsection.  In the tube $|z|<\delta$, define
\begin{equation}\label{innergluedprofile}
\begin{split}
 U_\epsilon^{\rm in}(x,t)
 :=e^{i\widehat\Phi_\epsilon}\Big[&
 \eta_0(r)\big(W(y)+\epsilon\eta_{\rm in}(z)
                         \phi_{1,\epsilon}(s,y,t)\big)\\
 &+(1-\eta_0(r))W(y)
   \exp\left(\epsilon\eta_{\rm in}(z)
                    \frac{\phi_{1,\epsilon}(s,y,t)}{W(y)}\right)
 \Big].
\end{split}
\end{equation}
For $r\geq2$ the quotient $\phi_{1,\epsilon}/W$ is regular.  Near
$r=0$ the additive expression is used instead.  Since
$\widehat\Phi_\epsilon=\Phi_\epsilon$ for $r\geq2$ and
$\eta_{\rm in}=0$ near $|z|=\delta$, formula
\eqref{innergluedprofile} agrees there with $e^{i\Phi_\epsilon}W$.
We therefore extend $U_\epsilon^{\rm in}$ to all of $\mathbb R^3$ by the
first ansatz $u_\epsilon^{(0)}$.  This extension is smooth and is exactly
$e^{i\varphi}$ for $|z|\geq2\delta$.

The first globally improved approximation will have the form
\begin{equation}\label{firstglobalgluingansatz}
 u_\epsilon^{(1)}(x,t)
 :=U_\epsilon^{\rm in}(x,t)
   \exp\big(i\psi_\epsilon^{\rm out}(x,t)\big),
 \qquad
 \psi_\epsilon^{\rm out}
 =\psi_{1,\epsilon}^{\rm out}+i\psi_{2,\epsilon}^{\rm out}.
\end{equation}
Thus $\psi_{1,\epsilon}^{\rm out}$ is a genuine phase correction, whereas
$e^{-\psi_{2,\epsilon}^{\rm out}}$ is the small accompanying amplitude
correction forced by the Gross--Pitaevskii system.

\subsection{The planar inner inverse}

The planar degree-one vortex $W$ satisfies
\begin{equation*}
 \Delta_yW+(1-|W|^2)W=0.
\end{equation*}
Linearising this equation at $W$ gives rise to the real-linear operator
\begin{equation*}
 \mathcal L_W[\phi]
 :=\Delta_y\phi+(1-|W|^2)\phi
       -2\operatorname{Re}(\overline W\phi)W.
\end{equation*}
The bounded kernel of $\mathcal L_W$ is spanned by the Jacobi fields
\begin{equation*}
    Z_1=\partial_{y_1}W,\qquad
    Z_2=\partial_{y_2}W,\qquad Z_0=iW,
\end{equation*}
see, for example, \cite[Theorem~3.2]{pacardriviere2000}.  Only solvability
conditions against $Z_1, Z_2$ will be required for the improvement procedure.

Let $D_t^x$ denote differentiation at fixed Cartesian $x$, pulled back
to the moving tube, and let $\nabla_{\gamma,\epsilon}$ and
$\Delta_{\gamma,\epsilon}$ be the operators in
\eqref{gradientperiodicframe}--\eqref{laplaceperiodicframe} after
$z=\epsilon y$.  For a real phase $\omega$ define the exact conjugated
residual
\[
 \mathscr C_{\epsilon,\gamma,\omega}[v]
 :=e^{-i\omega}S_\epsilon(e^{i\omega}v).
\]
Choose a radial $\chi\in C_c^\infty(B_3)$ which is one on $B_2$.
The source to be eliminated is
\begin{equation}\label{regularisedfirstinnersource}
\begin{split}
 H_\epsilon^{\rm c}
 &:=\epsilon^{-1}
       \mathscr C_{\epsilon,\gamma,\widehat\Phi_\epsilon}[W],\\
 \widehat H_\epsilon^{(1)}
 &:=\chi H_\epsilon^{\rm c}+(1-\chi)H_\epsilon^{(1)}.
\end{split}
\end{equation}
In particular, $\widehat H_\epsilon^{(1)}$ is a smooth Cartesian
source and agrees with \eqref{firstinnersource} for $r\geq3$.
For $\mathcal E_\epsilon
:=\mathscr C_{\epsilon,\gamma,\Phi_\epsilon}[W]
  -\epsilon H_\epsilon^{(1)}$, its defining identity is
\begin{equation}\label{exact-core-source-identity}
 \mathscr C_{\epsilon,\gamma,\widehat\Phi_\epsilon}[W]
 -\epsilon\widehat H_\epsilon^{(1)}
 =(1-\chi)\mathcal E_\epsilon.
\end{equation}
Indeed, $1-\chi$ is supported where the two phases agree.  The remainder
on the right vanishes on $B_2$.

Fix a nonnegative radial function $\chi_{\rm tr}\in C_c^\infty(B_2)$ which is one on
$B_1$, and put
\[
 (\mathsf M_{\rm tr})_{jk}:=\operatorname{Re}\int_{B_2}
        \chi_{\rm tr}Z_k\overline{Z_j}\,dy.
\]
By radial symmetry, $\mathsf M_{\rm tr}$ is a positive scalar multiple of
the identity.
Define
\begin{equation}\label{translationdefectfirstimprovement}
\begin{split}
 b_j(s,t)&:=\sum_{k=1}^2(\mathsf M_{\rm tr}^{-1})_{jk}
 \operatorname{Re}\int_{B_{R_\epsilon}}
 \widehat H_\epsilon^{(1)}\overline{Z_k}\,dy,\\
 H_{\epsilon,\mathrm{sol}}^{(1)}
 &:=\widehat H_\epsilon^{(1)}
 -\chi_{\rm tr}\sum_{j=1}^2b_jZ_j.
\end{split}
\end{equation}

\begin{lemma}[Orthogonality after core regularisation]
\label{regularisedsourceorthogonalitylemma}
The solvable source in \eqref{translationdefectfirstimprovement} satisfies
\begin{equation}\label{regularisedsourceorthogonality}
 \operatorname{Re}\int_{B_{R_\epsilon}}
 H_{\epsilon,\mathrm{sol}}^{(1)}\overline{Z_j}\,dy=0,
 \qquad j=1,2.
\end{equation}
For the corrected curve, the coefficients satisfy, at every fixed order,
\begin{equation}\label{translationdefectboundfirstimprovement}
 |D_{\rm in}^qD_s^aD_t^b(b_1+ib_2)|
 \leq C_m\epsilon\Lambda_\epsilon^{C_m},
 \qquad a+b+q\leq m,\quad q\leq2.
\end{equation}
Thus the subtracted term contributes only
$O(\epsilon^2\Lambda_\epsilon^{C_m})$ to the physical residual.
\end{lemma}

\begin{proof}
We verify both the smooth core bounds and the size of the new projection.
Put $v=\gamma_t^\perp$, $h=1-\epsilon\vec\kappa\cdot y$, and
\[
 B_\epsilon f:=h^{-1}\mathscr D_s(h^{-1}\mathscr D_s f).
\]
In a fixed material coordinate, write
$D_t^x=\mathscr D_t+b_\epsilon\mathscr D_s
-\epsilon^{-1}v\cdot\nabla_y$.  The coefficient $b_\epsilon$ is smooth
and uniformly bounded with its derivatives on $B_3$; it contains the
tangential velocity and the tangential component of frame motion.
Writing $\omega=\widehat\Phi_\epsilon=\beta_\epsilon+\epsilon p_\epsilon$,
direct use of the tubular Laplacian and the vortex equation gives
\begin{equation}\label{exact-core-source-expansion}
\begin{split}
 H_\epsilon^{\rm c}
 &=-i\Lambda_\epsilon v\cdot\nabla_yW
   -\vec\kappa\cdot\nabla_yW
   +\mathcal L_W(ip_\epsilon W)+\epsilon\mathcal R_\epsilon^{\rm c},\\
 \mathcal R_\epsilon^{\rm c}
 &=-\frac{\vec\kappa\cdot y}{h}\vec\kappa\cdot\nabla_yW
   +W\left[-h^{-2}(\mathscr D_s\omega-\mathfrak a)^2
                    -|\nabla_y p_\epsilon|^2\right]\\
 &\quad+\Lambda_\epsilon W\left[
   \mathfrak c+\mathfrak a b_\epsilon-\partial_t\beta_\epsilon
   -b_\epsilon\partial_s\beta_\epsilon+v\cdot\nabla_y p_\epsilon
   -\epsilon(\mathscr D_t+b_\epsilon\mathscr D_s)p_\epsilon\right]\\
 &\quad+iW\left[B_\epsilon(\theta+\omega)
                 -h^{-1}\vec\kappa\cdot\nabla_y p_\epsilon\right].
\end{split}
\end{equation}
Although $\theta$ is not defined at zero, the expression occurring here is
\[
 B_\epsilon(\theta+\omega)
 =h^{-2}(\partial_s^2\beta_\epsilon-\partial_s\mathfrak a
                         +\epsilon\mathscr D_s^2p_\epsilon)
 +\epsilon h^{-3}(D_s\vec\kappa\cdot y)
       (\partial_s\beta_\epsilon-\mathfrak a
                         +\epsilon\mathscr D_sp_\epsilon),
\]
which is smooth.  This proves
$\|\mathcal R_\epsilon^{\rm c}\|_{C^m(B_3)}
+\|H_\epsilon^{\rm c}\|_{C^m(B_3)}\leq
C_m\Lambda_\epsilon^{C_m}$, also after the stated mixed derivatives.

For the projection calculation, the harmonicity of $\theta+\Phi_\epsilon$
and \eqref{laplaceperiodicframe} give on the punctured disc
\begin{equation}\label{core-phase-harmonic-identity}
 \Delta_yP_\epsilon
 =h^{-1}\vec\kappa\cdot\nabla_y\theta
  +\epsilon h^{-1}\vec\kappa\cdot\nabla_yP_\epsilon
  -\epsilon B_\epsilon(\theta+\Phi_\epsilon).
\end{equation}
Also, \eqref{firstinnersourcegeometric} gives exactly
\[
 H_\epsilon^{(1)}+\mathcal L_W(iq_\epsilon W)
 =-i\Lambda_\epsilon v\cdot\nabla_yW
   -\vec\kappa\cdot\nabla_yW+\mathcal L_W(ip_\epsilon W)
   +iW(\vec\kappa\cdot\nabla_y\theta-\Delta_yP_\epsilon).
\]
Consequently
\begin{equation}\label{core-source-projection-comparison}
 \widehat H_\epsilon^{(1)}
 =H_\epsilon^{(1)}+\mathcal L_W(iq_\epsilon W)+d_\epsilon,
 \qquad
 \|D_{\rm in}^q\partial_t^a\partial_s^b d_\epsilon\|_{L^1(B_3)}
 \leq C_m\epsilon\Lambda_\epsilon^{C_m},
\end{equation}
where
$d_\epsilon=\chi\{
\epsilon\mathcal R_\epsilon^{\rm c}
-iW(\vec\kappa\cdot\nabla_y\theta-\Delta_yP_\epsilon)\}$.
Here only an $L^1$ estimate is needed.  To see it directly, use
\eqref{core-phase-harmonic-identity}, $1-h^{-1}=O(\epsilon r)$,
$\nabla_yP_\epsilon=O(\Lambda_\epsilon+|\log r|)$ and the phase jet.
The factor $W=O(r)$ makes every resulting logarithmic term integrable.

The function $iq_\epsilon W$ is compactly supported and is
$O(r^2(\Lambda_\epsilon+|\log r|))$ at zero, with gradient
$O(r(\Lambda_\epsilon+|\log r|))$.  Green's formula on a punctured disc
and $\mathcal L_WZ_j=0$ therefore imply
\[
 \operatorname{Re}\int_{B_{R_\epsilon}}
       \mathcal L_W(iq_\epsilon W)\overline{Z_j}\,dy=0.
\]
The inner boundary term tends to zero like
$O(\rho^2(\Lambda_\epsilon+|\log\rho|))$.
The original source has translation projection
$O(\epsilon\Lambda_\epsilon)$ by the computation preceding
\eqref{coreconstantasymptotics}.  Equation
\eqref{core-source-projection-comparison} proves
\eqref{translationdefectboundfirstimprovement}.  The same argument applies
after ordinary parameter differentiation in the fixed tube; expansion
of the covariant derivatives gives the displayed form.  Finally,
\eqref{translationdefectfirstimprovement} gives exact orthogonality by
matrix multiplication.  The subtraction includes the translation
projection of the compact term $d_\epsilon$.
\end{proof}

To construct the inner correction, we require estimates for the linearised
problem $\mathcal L_W[\phi]=-h$ on a ball of large radius.  The following
lemma gives a convenient general estimate, using the Fourier-mode theory
developed in \cite{delpinojunemanmussolinear2025}.

\begin{lemma}[Planar inverse up to a matching radius]
\label{planarinversefirstimprovementlemma}
Put $\ell(r)=1+\log(2+r)$.  Fix $0<c<1$.  Let $R\geq4$ and let
$h\in C^{m+3}(B_R;\mathbb C)$ satisfy
\begin{equation}\label{hlocalfirstimprovement}
 \sum_{j=0}^{m+3}\|D_y^jh\|_{L^\infty(B_2)}\leq C_h.
\end{equation}
For $1\leq r<R$, write uniquely
\begin{equation*}
 h=iW(g_1+ig_2),
 \qquad g_1,g_2\ \hbox{real-valued},
\end{equation*}
and assume that, for $0\leq j\leq m+3$,
\begin{equation}\label{hweightedfirstimprovement}
 |D_y^jg_1(y)|
 \leq C_h\frac{\ell(r)^2}{(1+r)^{2+j}},
 \qquad
 |D_y^jg_2(y)|
 \leq C_h\frac{\ell(r)^2}{(1+r)^{1+j}}.
\end{equation}
If
\begin{equation}\label{hinverseorthogonality}
 \operatorname{Re}\int_{B_R}h\overline{Z_j}\,dy=0,
 \qquad j=1,2,
\end{equation}
then there is a solution of
\begin{equation}\label{finiteballinnerproblem}
    \mathcal L_W\phi=-h\quad\hbox{in }B_R
\end{equation}
which extends as a Cartesian $C^{m+4}$ function at $r=0$ (and is smooth
there if $h$ is smooth), and has the following matching property: when its
Fourier modes are continued homogeneously to $r>R$, no exponentially
growing branch and no algebraically growing non-zero mode is present.  The
logarithmic zero-mode phase branch forced by the source is allowed.  The
solution may moreover be centred and phase-normalised by
\begin{equation}\label{innernormalisationfirstimprovement}
 \phi(0)=0,
 \qquad
 \operatorname{Re}\int_{B_2}\phi\overline{iW}\,dy=0.
\end{equation}
It satisfies
\begin{equation}\label{planarinverseestimatefirstimprovement}
 \sum_{j=0}^{m+2}
 \sup_{0<r<cR}
 \frac{(1+r)^j|D_y^j\phi(y)|}
      {\ell(r)^3}
 \leq C_{m,c} C_h\ell(R)^2.
\end{equation}
The constant is independent of $R$.  The inverse is defined by extending
the source by zero outside $B_R$.  The differentiated estimate applies
on $B_{cR}$.
\end{lemma}

\begin{proof}
Put
\[
 \phi=iW\psi,
 \qquad
 \mathfrak g:=-(iW)^{-1}h=\mathfrak g_1+i\mathfrak g_2.
\]
Then \eqref{finiteballinnerproblem} is equivalent on $r>0$ to
\begin{equation}\label{psiequationfiniteinverse}
 \Delta\psi+2\frac{w'}w\partial_r\psi
 +\frac{2i}{r^2}\partial_\theta\psi
 -2iw^2\operatorname{Im}\psi
 =\mathfrak g_1+i\mathfrak g_2.
\end{equation}
We now make explicit the Fourier decomposition used in
\cite{delpinojunemanmussolinear2025}.  For $r>0$, write
\[
\begin{aligned}
 \psi&=\psi^0+\sum_{k\geq1}
       \big(\psi_k^{\mathrm c}+\psi_k^{\mathrm s}\big),
 &\psi^0&=\psi_1^0+i\psi_2^0,\\
 \psi_k^{\mathrm c}
 &=a_{k,1}^{\mathrm c}(r)\cos k\theta
   +ia_{k,2}^{\mathrm c}(r)\sin k\theta,
 &\psi_k^{\mathrm s}
 &=a_{k,1}^{\mathrm s}(r)\sin k\theta
   +ia_{k,2}^{\mathrm s}(r)\cos k\theta.
\end{aligned}
\]
Decompose $\mathfrak g$ in the same form, with radial coefficients
$\mathfrak g_j^0$ and $\mathfrak g_{k,j}^{\mathrm c},
\mathfrak g_{k,j}^{\mathrm s}$.  The zero mode satisfies
\[
\begin{aligned}
 (\psi_1^0)''+\left(2\frac{w'}w+\frac1r\right)(\psi_1^0)'
 &=\mathfrak g_1^0,\\
 (\psi_2^0)''+\left(2\frac{w'}w+\frac1r\right)(\psi_2^0)'
      -2w^2\psi_2^0&=\mathfrak g_2^0.
\end{aligned}
\]
The first equation is integrated twice.  For the second, let
$z_{1,0},z_{2,0}$ be the homogeneous solutions in Lemma~2.2 of
\cite{delpinojunemanmussolinear2025}, with $z_{1,0}$ regular at the
origin and $z_{2,0}$ exponentially decaying at infinity.  Variation of
parameters gives
\begin{equation}\label{zeromodefiniteinverse}
\begin{split}
 \psi_{1}^{0}(r)
 &:=\int_0^r\frac1{w(s)^2s}
       \int_0^s w(\rho)^2\rho\,\mathfrak g_1^0(\rho)
       \,d\rho\,ds,\\
 \psi_{2}^{0}(r)
 &:=z_{1,0}(r)\int_r^R w(s)^2s\,\mathfrak g_2^0(s)z_{2,0}(s)\,ds\\
 &\hspace{2em}
   +z_{2,0}(r)\int_0^r w(s)^2s\,\mathfrak g_2^0(s)z_{1,0}(s)\,ds.
\end{split}
\end{equation}
The integration limits make $\psi_2^0$ regular at the origin and remove
its exponentially growing continuation beyond $R$.  A nonzero total
integral in the first line produces the admissible logarithmic phase
branch.

For $k\geq1$, define
\[
\begin{aligned}
 \mathbf a_k^{\mathrm c}
  &=\big(a_{k,1}^{\mathrm c},a_{k,2}^{\mathrm c}\big),
 &\widetilde{\mathbf g}_k^{\mathrm c}
  &=\big(\mathfrak g_{k,1}^{\mathrm c},
         \mathfrak g_{k,2}^{\mathrm c}\big),\\
 \mathbf a_k^{\mathrm s}
  &=\big(a_{k,1}^{\mathrm s},-a_{k,2}^{\mathrm s}\big),
 &\widetilde{\mathbf g}_k^{\mathrm s}
  &=\big(\mathfrak g_{k,1}^{\mathrm s},
        -\mathfrak g_{k,2}^{\mathrm s}\big).
\end{aligned}
\]
Both parity sectors then satisfy the same real system
\[
 \mathbf a''+\left(2\frac{w'}w+\frac1r\right)\mathbf a'
 -\frac1{r^2}
 \begin{pmatrix}
  k^2&2k\\
  2k&k^2+2w^2r^2
 \end{pmatrix}\mathbf a
 =\widetilde{\mathbf g}.
\]
Let $\mathbf z_{1,k},\ldots,\mathbf z_{4,k}$ be the homogeneous
solutions classified in Proposition~3.4 of
\cite{delpinojunemanmussolinear2025} for $k=1$ and Proposition~4.4
there for $k\geq2$.  To obtain a representation formula for $\mathbf a$,
write $\mathbf a=\sum_{j=1}^4c_j\mathbf z_{j,k}$ and impose
\[
 \sum_{j=1}^4c_j'\mathbf z_{j,k}=0,
 \qquad
 \sum_{j=1}^4c_j'\mathbf z_{j,k}'=\widetilde{\mathbf g}.
\]
The Wronskian computation in Propositions~3.8 and 4.8 of that paper,
after normalising the homogeneous basis, gives
\[
\begin{aligned}
 c_1'&=w^2r\,\widetilde{\mathbf g}\cdot\mathbf z_{2,k},
 &c_2'&=-w^2r\,\widetilde{\mathbf g}\cdot\mathbf z_{1,k},\\
 c_3'&=w^2r\,\widetilde{\mathbf g}\cdot\mathbf z_{4,k},
 &c_4'&=-w^2r\,\widetilde{\mathbf g}\cdot\mathbf z_{3,k}.
\end{aligned}
\]
Choosing the antiderivatives at $0$ on the branches singular there and
at $R$ on the branches growing at infinity gives, in either parity
sector,
\begin{equation}\label{modevariationconstantsfirstimprovement}
\begin{split}
 \mathbf a_k(r)={}&
 \mathbf z_{1,k}(r)\int_0^r w(\rho)^2\rho\,
      \widetilde{\mathbf g}_k(\rho)\cdot\mathbf z_{2,k}(\rho)\,d\rho\\
 &+\mathbf z_{2,k}(r)\int_r^R w(\rho)^2\rho\,
      \widetilde{\mathbf g}_k(\rho)\cdot\mathbf z_{1,k}(\rho)\,d\rho\\
 &-\mathbf z_{3,k}(r)\int_r^R w(\rho)^2\rho\,
      \widetilde{\mathbf g}_k(\rho)\cdot\mathbf z_{4,k}(\rho)\,d\rho\\
 &-\mathbf z_{4,k}(r)\int_0^r w(\rho)^2\rho\,
      \widetilde{\mathbf g}_k(\rho)\cdot\mathbf z_{3,k}(\rho)\,d\rho.
\end{split}
\end{equation}
For $k\geq2$, the asymptotics of the four homogeneous solutions show
directly that this choice is regular at the origin and has no growing
continuation beyond $R$.  For $k=1$, Proposition~3.8 initially gives the second
coefficient in the form
\[
 -\int_0^r w(\rho)^2\rho\,
       \widetilde{\mathbf g}_1(\rho)
       \cdot\mathbf z_{1,1}(\rho)\,d\rho.
\]
For the two parity sectors, the conditions
\eqref{hinverseorthogonality} are exactly
\begin{equation}\label{translationmomentfiniteinverse}
 \int_0^R w(\rho)^2\rho\,
       \widetilde{\mathbf g}_1(\rho)
       \cdot\mathbf z_{1,1}(\rho)\,d\rho=0.
\end{equation}
Thus the preceding coefficient equals the integral from $r$ to $R$ in
\eqref{modevariationconstantsfirstimprovement}.  This is the precise point
at which the two translation orthogonality conditions remove the
non-admissible $k=1$ branch.

Formulae \eqref{zeromodefiniteinverse} and
\eqref{modevariationconstantsfirstimprovement} solve every Fourier-mode
equation exactly.  If the source is set equal to zero for $r>R$, then
\eqref{modevariationconstantsfirstimprovement} contains only
$\mathbf z_{1,k}$ and $\mathbf z_{4,k}$ there.  Their asymptotics in
Propositions~3.4 and 4.4 show that they are the decaying algebraic and
decaying exponential branches.  In the zero mode, the imaginary component
decays exponentially, while the real component may have the admissible
logarithmic phase growth described after (1.14) of that paper.  This proves
the asserted matching property.

We next justify the estimate uniformly in $R$.  Let $h_R$ denote $h$ on
$B_R$ and zero outside $B_R$.  In the notation of (1.15)--(1.16) of
\cite{delpinojunemanmussolinear2025}, assumptions
\eqref{hlocalfirstimprovement}--\eqref{hweightedfirstimprovement} give
\begin{equation*}
 \|h_R\|_{**}\leq C C_h\ell(R)^2.
\end{equation*}
The estimates obtained directly from the zero-mode formula in
Proposition~2.1, from (3.18)--(3.19) for $k=1$, and from
(4.17)--(4.18) for $k\geq2$ therefore apply with this right-hand side.
Proposition~5.1 of the same paper sums the modes $k\geq3$ with a constant
independent of the number of modes.  Since the bounded mode solutions in
that proposition are unique, they agree with
\eqref{modevariationconstantsfirstimprovement}.  We consequently obtain
\begin{equation}\label{zerothmode-summed-finiteinverse}
 \sup_{r<2}|\phi|+
 \sup_{2<r<R}\frac{|\phi|}{\ell(r)^3}
 \leq C C_h\ell(R)^2.
\end{equation}

We give an interior argument for the differentiated estimate.  This
also explains the role of the fixed gap between $cR$ and $R$.
Write $U=\operatorname{Re}\psi$, $V=\operatorname{Im}\psi$ and
$A=\Delta+2(w'/w)\partial_r$.  On an annulus of radius $\rho$ contained
in $B_R\setminus B_2$, their equations are
\begin{equation}\label{planar-annular-system}
 AU=\mathfrak g_1+2r^{-2}\partial_\theta V,
 \qquad
 (2w^2-A)V=-\mathfrak g_2+2r^{-2}\partial_\theta U.
\end{equation}
Here and below annuli in the variable $x=y/\rho$ have fixed inner and
outer radii and fixed positive separation.  If $K\Subset K'$ are two
such annuli, put
\[
 X_j(K):=\|U(\rho\,\cdot)\|_{C^{j+1,1/2}(K)}
             +\|V(\rho\,\cdot)\|_{C^{j,1/2}(K)},
 \qquad j\geq1.
\]
For any fixed number of nested annuli the following estimate holds,
with constants depending on their separations:
\begin{equation}\label{annular-bootstrap-estimate}
 X_j(K)\leq C_j\mathcal M_j+C_j\rho^{-2}X_j(K'),
\end{equation}
where $\mathcal M_j$ can be taken to be the sum, on a slightly larger
annulus, of the two $C^0$ norms and
\[
 \|\rho^2\mathfrak g_1(\rho\,\cdot)\|_{C^{j-1,1/2}}
       +\|\mathfrak g_2(\rho\,\cdot)\|_{C^{j,1/2}}.
\]

For completeness, the scalar massive estimate used to prove
\eqref{annular-bootstrap-estimate} is
\begin{equation}\label{massive-annular-estimate}
 \|v\|_{C^{j,1/2}(K)}
 \leq C_j\|F\|_{C^{j,1/2}(K')}
       +C_{j,N}\rho^{-N}\|v\|_{C^0(K')},
\end{equation}
for $(2w(\rho|x|)^2-\rho^{-2}A_\rho)v=F$, where
$A_\rho=\rho^2A$ in the $x$ variable and $\rho$ is sufficiently
large.  To prove it, extend $F$ and the coefficients across the fixed
annulus, retaining a positive lower bound for $2w^2$, and solve the
scalar equation on the whole plane for a bounded particular solution.
The maximum principle, its differentiated equations and difference
quotients give a $C^{j,1/2}$ bound by that of $F$, uniformly for
$\rho^{-1}\leq1$.  For the difference from this particular solution,
the positive zeroth-order coefficient gives the barrier
$C e^{-c_0\rho\operatorname{dist}(x,\partial K')}$.  Interior estimates
on balls of radius $\rho^{-1}$ bound each derivative by a power of
$\rho$ times this exponential.  On $K$ that quantity is bounded by
$C_{j,N}\rho^{-N}$ for every fixed $N$.  This proves
\eqref{massive-annular-estimate}.
Apply it to the second equation of \eqref{planar-annular-system}:
the contribution of $U$ is at most
$C_j\rho^{-2}\|U\|_{C^{j+1,1/2}}$ in scaled variables.  The ordinary
interior Schauder estimate for the first equation bounds
$\|U\|_{C^{j+1,1/2}}$ by its $C^0$ norm, the displayed norm of
$\rho^2\mathfrak g_1$, and $\|V\|_{C^{j,1/2}}$ on a larger annulus.
Using one intermediate annulus gives
\eqref{annular-bootstrap-estimate}.

On a still larger annulus, unit-ball estimates for the original smooth
Cartesian system $\mathcal L_W\phi=-h$ give the rough bound
$X_j\leq C_j\rho^{j+2}\mathcal M_j$.  Iterate
\eqref{annular-bootstrap-estimate} along $n$ fixed nested annuli, with
$2n>j+2$.  The last term is bounded by
$C_{j,n}\rho^{j+2-2n}\mathcal M_j$, and hence
\begin{equation}\label{uniform-annular-derivative-estimate}
 X_j(K)\leq C_j\mathcal M_j.
\end{equation}
The finite number of separations can depend on $j$ and $c$, but not on
$\rho$ or $R$.  Taking $j=m+2$, the source derivatives through $m+3$
control the required H\"older norms.  Equation
\eqref{zerothmode-summed-finiteinverse} and the source assumptions give
$\mathcal M_j\leq C_m C_h\ell(R)^2\ell(\rho)^3$.
Rescaling back, and differentiating $iW\psi$, proves
\eqref{planarinverseestimatefirstimprovement} on $2<r<cR$.
The bounded range of $\rho$ and the core are covered by ordinary
Cartesian interior estimates; all these neighbourhoods remain strictly
inside $B_R$.  This proves the estimate, including the possible massive
boundary contribution, whose effect at distance comparable to $R$ is
exponentially small.

It remains only to impose \eqref{innernormalisationfirstimprovement}.  Put
$Z_0=iW$.  Since $w'(0)>0$ and
\[
 Z_1(0)=w'(0),\qquad Z_2(0)=iw'(0),\qquad Z_0(0)=0,
\]
the real-linear map
$(a_1,a_2)\mapsto a_1Z_1(0)+a_2Z_2(0)$ is an isomorphism from
$\mathbb R^2$ to $\mathbb C$.  Starting from the solution
$\widetilde\phi$ just constructed, there are therefore unique real
$a_1,a_2$ for which
$\widetilde\phi+a_1Z_1+a_2Z_2$ vanishes at the origin.  There is then a
unique real $a_0$ for which
\[
 \phi=\widetilde\phi+a_1Z_1+a_2Z_2+a_0Z_0
\]
satisfies the second condition in
\eqref{innernormalisationfirstimprovement}; adding $a_0Z_0$ does not
disturb the first condition.  Point evaluation on $B_2$ and the fixed
$L^2(B_2)$ pairing bound these three coefficients by the left-hand side
of \eqref{zerothmode-summed-finiteinverse}.  Since $Z_0,Z_1,Z_2$ obey the
same admissible bounds, the estimate and matching property are unchanged.
Thus centring is a bounded finite-rank modification of the mode inverse.
The resulting operator depends only on $R$ and commutes with ordinary
parameter differentiation at fixed normal coordinates.  For angular
derivatives the equivariance generator on Cartesian perturbations is
$\mathcal A=\partial_\theta-i$, for which
$\mathcal A(iW\psi)=iW\partial_\theta\psi$ and
$[\mathcal L_W,\mathcal A]=0$.  The translation-orthogonal source space
and both normalisations are rotation invariant.  Thus angular
differentiation of the relative equation is also legitimate.
Horizontal derivatives are obtained by expanding
$\mathscr D_s=\partial_s-\mathfrak a\partial_\theta$ and
$\mathscr D_t=\partial_t-\mathfrak c\partial_\theta$ after ordinary
parameter differentiation, retaining the derivatives of their
coefficients.  This gives the same estimates at each fixed total order,
including the first two curve variations.
\end{proof}

We apply Lemma \ref{planarinversefirstimprovementlemma} at each $(s,t)$
with $R=R_\epsilon$ and
$h=H_{\epsilon,\mathrm{sol}}^{(1)}$.  The expansion
\eqref{phistargradientexpansion}, the curve equation
\eqref{firstgammaeq}, and the vortex asymptotics give, for $r\geq3$,
the component formula
\begin{equation}\label{first-source-exact-components}
\begin{split}
 H_{\epsilon,\mathrm{sol}}^{(1)}
 &=iW(G_{1,\epsilon}+iG_{2,\epsilon}),\\
 G_{1,\epsilon}
 &=\frac{w'}w e_r\cdot V_\epsilon,\\
 G_{2,\epsilon}
 &=\frac1r e_\theta\cdot V_\epsilon
   +\frac{w'}w\vec\kappa\cdot e_r.
\end{split}
\end{equation}
The compact translation subtraction and the core regularisation are
supported in $B_3$ and hence do not alter these formulae.  Since
$w'/w=O(r^{-3})$, cancellation of
$\log(\ell(t)/\epsilon)$ by the coefficient in the curve velocity yields,
for $a+b+|\alpha|+q\leq m$ and $q\leq2$,
\begin{equation}\label{innersourceweightedbound}
\begin{split}
 |D_{\rm in}^q\mathscr D_t^a\mathscr D_s^b
       D_y^\alpha G_{1,\epsilon}|
 &\leq C_m\frac{\ell(r)^{C_m}}
                   {(1+r)^{3+|\alpha|}},\\
 |D_{\rm in}^q\mathscr D_t^a\mathscr D_s^b
       D_y^\alpha G_{2,\epsilon}|
 &\leq C_m\frac{\ell(r)^{C_m}}
                   {(1+r)^{1+|\alpha|}}.
\end{split}
\end{equation}
The angular mean of the leading part of $G_{1,\epsilon}$ is zero.  More
precisely, if a superscript zero denotes angular average, then for
$a+b+j+q\leq m$ and $q\leq2$,
\begin{equation}\label{first-source-zero-mode-bound}
 |D_{\rm in}^q\mathscr D_t^a\mathscr D_s^b
       \partial_r^jG_{1,\epsilon}^{0}(r)|
 \leq C_m\epsilon\Lambda_\epsilon^{C_m}
       (1+r)^{-2-j}.
\end{equation}
For $r\geq3$, every displayed term in
\eqref{phistargradientexpansion} has zero scalar product with $e_r$ after
angular integration.  Its $O(\epsilon r|\log(\epsilon r)|)$ remainder,
multiplied by $w'/w$, has the size in
\eqref{first-source-zero-mode-bound}.
The core requires a separate calculation.  If
$H_\epsilon^{\rm c}=iW(G_1^{\rm c}+iG_2^{\rm c})$, the exact phase
component of \eqref{exact-core-source-expansion} is
\begin{equation}\label{core-source-phase-component}
\begin{split}
 G_1^{\rm c}
 ={}&-\Lambda_\epsilon\frac{w'}w v\cdot e_r
      +\Delta_y p_\epsilon+2\frac{w'}w\partial_r p_\epsilon
      -h^{-1}\vec\kappa\cdot\nabla_y\theta\\
 &\quad-\epsilon h^{-1}\vec\kappa\cdot\nabla_y p_\epsilon
      +\epsilon B_\epsilon(\theta+\widehat\Phi_\epsilon).
\end{split}
\end{equation}
The leading jet of $P_\epsilon$ in
\eqref{phasejetfirstimprovement} is a first angular harmonic.  Since
$\zeta$ is radial and vanishes on $B_1$, this gives
$\|(p_\epsilon)^0\|_{C^m([0,3])}
\leq C_m\epsilon\Lambda_\epsilon^{C_m}$.
The first term in \eqref{core-source-phase-component} has mean zero, and
\[
 h^{-1}\vec\kappa\cdot\nabla_y\theta
   =-\frac1{\epsilon r^2}\partial_\theta\log h
\]
has mean zero exactly.  The last two terms carry an explicit $\epsilon$
and have bounded angular averages and radial derivatives on $[0,3]$,
as is also clear from the smooth expression for $B_\epsilon$ above.
Thus $(G_1^{\rm c})^0$ has the required small bound.  Multiplication by
the radial $\chi$ preserves it.  Finally, the phase component of
$Z_j/(iW)$ is a first angular harmonic; the radial translation
subtraction therefore has zero phase mean.  This proves
\eqref{first-source-zero-mode-bound} through the core.  All calculations
can be differentiated twice in the curve parameter and to each fixed
horizontal and time order in the common tame family.

The extension in \eqref{regularisedfirstinnersource} supplies the uniform
core bounds \eqref{hlocalfirstimprovement}, and Lemma
\ref{regularisedsourceorthogonalitylemma} supplies the two compatibility
conditions.  We therefore obtain a smooth periodic family
$\phi_{1,\epsilon}$ satisfying
\begin{equation}\label{firstinnercorrectionequation}
 \mathcal L_W\phi_{1,\epsilon}
 =-H_{\epsilon,\mathrm{sol}}^{(1)},
 \qquad \phi_{1,\epsilon}(s,0,t)=0,
 \qquad |y|<R_\epsilon,
\end{equation}
and, for every fixed $m$,
\begin{equation}\label{firstinnercorrectionestimate}
 \sum_{a+b+|\alpha|\leq m}
 \sup_{s,t,\,|y|<7R_\epsilon/8}
 \frac{(1+r)^{|\alpha|}
 |\mathscr D_t^a\mathscr D_s^bD_y^\alpha
             \phi_{1,\epsilon}|}
 {(1+\log(2+r))^3}
 \leq C_m\Lambda_\epsilon^{C_m}.
\end{equation}
For later use we record, and prove, the sharper componentwise estimate
which is needed on the fixed matching annulus.  Write
\begin{equation*}
 \phi_{1,\epsilon}=iW(\psi_{1,\epsilon}+i\psi_{2,\epsilon}),
 \qquad \psi_{1,\epsilon},\psi_{2,\epsilon}\ \hbox{real-valued}.
\end{equation*}
Then, for $2\leq r\leq7R_\epsilon/8$,
\begin{equation}\label{firstinnercorrectiontail}
\begin{split}
 |D_{\rm in}^q\mathscr D_t^a\mathscr D_s^bD_y^\alpha
       \psi_{2,\epsilon}|
 &\leq C_m\Lambda_\epsilon^{C_m}(1+r)^{-1-|\alpha|},\\
 |D_{\rm in}^q\mathscr D_t^a\mathscr D_s^bD_y^\alpha
       \psi_{1,\epsilon}|
 &\leq C_m\Lambda_\epsilon^{C_m}
   \left((1+r)^{-1-|\alpha|}
    +\epsilon\ell(r)^2(1+r)^{-|\alpha|}\right),
\end{split}
\end{equation}
whenever $a+b+|\alpha|+q\leq m$ and $q\leq2$.

We supply the decay argument before applying the interior derivative
estimate.  Extend the solvable source by zero for $r>R_\epsilon$ and
continue the mode solution as in Lemma
\ref{planarinversefirstimprovementlemma}.  Set
$U=\psi_{1,\epsilon}$ and $V=\psi_{2,\epsilon}$.
Angular differentiation preserves the two translation conditions by
rotation equivariance.  The global zeroth-order inverse estimate,
applied also to finitely many angular derivatives, gives
\[
 |\partial_\theta^j U|+|\partial_\theta^j V|
 \leq M_\epsilon\ell(r)^C,\qquad r\geq2,
 \qquad M_\epsilon\leq C_m\Lambda_\epsilon^{C_m}.
\]
This use of angular differentiation is unaffected by the radial jump
of the extended source at $R_\epsilon$.
The second equation of \eqref{planar-annular-system} and the positive
potential $2w^2$ improve this to
\begin{equation}\label{global-amplitude-tail-bound}
 |\partial_\theta^j V(r,\theta)|
 \leq C_mM_\epsilon\ell(r)^C/r,
 \qquad r\geq2,
\end{equation}
for every needed $j$.  Indeed, its right-hand side is bounded by
$C M_\epsilon\ell(r)^C/r$, including for $r>R_\epsilon$ where the
source vanishes.  For sufficiently large fixed $r_0$,
$(2w^2-A)(r^{-1}\ell(r)^C)\geq c r^{-1}\ell(r)^C$.
Apply the scalar maximum principle to each differentiated equation
on $r_0<r<L$, adding an outer boundary barrier and then sending
$L\to\infty$.  Polynomially bounded boundary values have exponentially
vanishing influence on compact subannuli.  The bounded range
$2\leq r\leq r_0$ is absorbed in the constant.

The first equation of \eqref{planar-annular-system} now has a right-hand
side of size $M_\epsilon\ell(r)^C/r^3$ for the non-zero angular part.
Interior $W^{2,p}$ estimates on rescaled annuli, with $p>2$, give
$|\partial_r(U-U^0)|\leq C M_\epsilon\ell(r)^C/r$.
They require only bounded forcing and hence apply also across
$r=R_\epsilon$.  Consequently
\[
 \Delta(U-U^0)=F,\qquad F^0=0,
 \qquad |F|\leq C M_\epsilon\ell(r)^C/r^3,
 \qquad r\geq r_0,
\]
because $w'/w=O(r^{-3})$.
For either sine or cosine coefficient $u_k$ and the corresponding $F_k$,
the absence of a growing harmonic branch gives the formula
\begin{equation}\label{phase-tail-poisson-formula}
 u_k(r)=c_k r^{-k}-\frac1{2k}\left[
 r^{-k}\int_{r_0}^r \rho^{k+1}F_k(\rho)\,d\rho
 +r^k\int_r^\infty \rho^{1-k}F_k(\rho)\,d\rho\right].
\end{equation}
The $k=1$ integrals are bounded by
$C M_\epsilon r^{-1}\ell(r)^{C+1}$.  For $k\geq2$ they are bounded
by $C M_\epsilon r^{-1}\ell(r)^{C+1}/(k(k-1))$.
The coefficients $c_k$ are determined by the bounded data at $r_0$,
and their contribution is summable for $r\geq2r_0$.
Thus
$|U-U^0|+|V|\leq C M_\epsilon r^{-1}\ell(r)^{C+1}$.
Translation orthogonality has already entered the global inverse; without
it a growing $k=1$ phase branch would invalidate this argument.

Apply \eqref{uniform-annular-derivative-estimate} to the pair
$(U-U^0,V)$.  On an annulus of radius $\rho\leq7R_\epsilon/8$, all its
input norms now have size
$C M_\epsilon\rho^{-1}\ell(\rho)^{C_m}$, by
\eqref{innersourceweightedbound}.  The annuli used in that estimate
are chosen strictly inside $B_{R_\epsilon}$.  This proves the
$r^{-1-|\alpha|}$ estimates for all spatial derivatives of the
non-zero phase part and of the amplitude part.

For the phase mean, the exact equation is
\begin{equation}\label{small-phase-zero-mode-primitive}
 (U^0)'(r)=-\frac1{w(r)^2r}
            \int_0^r w(\rho)^2\rho\,G_{1,\epsilon}^0(\rho)\,d\rho.
\end{equation}
The regular origin condition fixes this primitive.  The phase
normalisation in \eqref{innernormalisationfirstimprovement} fixes its
additive constant, which is
$O(\epsilon\Lambda_\epsilon^{C_m})$ by the core part of
\eqref{first-source-zero-mode-bound}.  The same bound and
\eqref{small-phase-zero-mode-primitive} imply
\[
 |U^0(r)|\leq C_m\epsilon\Lambda_\epsilon^{C_m}\ell(r)^2,
 \qquad |\partial_r^jU^0(r)|
 \leq C_m\epsilon\Lambda_\epsilon^{C_m}\ell(r)^2(1+r)^{-j},
 \quad j\geq1.
\]
No orthogonality condition against $iW$ is required for the source.
Ordinary parameter derivatives commute with the fixed-radius inverse;
expanding the horizontal derivatives as in its proof establishes all
the mixed versions.  On $r\leq7R_\epsilon/8$ the harmless additional
powers of $\ell(r)$ in the non-zero modes are absorbed in
$\Lambda_\epsilon^{C_m}$.  This proves
\eqref{firstinnercorrectiontail}.
In particular,
$\epsilon\phi_{1,\epsilon}=o(1)$ uniformly on every fixed inner disc.

We record also the location of the large core logarithm, since the
modulus estimates later use this information.  Set
\begin{equation}\label{core-logarithm-splitting}
 \sigma_\epsilon
 :=\frac{i\Lambda_\epsilon}{2}(\zeta(r)-1)
                    (J\vec\kappa\cdot y)W,
 \qquad
 \bar h_\epsilon:=H_{\epsilon,\mathrm{sol}}^{(1)}
                          -\mathcal L_W\sigma_\epsilon.
\end{equation}
The function $\sigma_\epsilon$ is smooth and compactly supported in
$B_2$, vanishes at zero and has zero pairing with $iW$.
For a constant vector $a$, the identity
$\mathcal L_W(i(a\cdot y)W)=2ia\cdot\nabla_yW$, together with
$\Lambda_\epsilon\gamma_t^\perp
=\Lambda_\epsilon J\vec\kappa+O(1)$, cancels every leading
$\Lambda_\epsilon$ term in the core of $\bar h_\epsilon$.
Equations \eqref{phasejetfirstimprovement} and
\eqref{exact-core-source-expansion} give
\[
 \|\bar h_\epsilon\|_{C^m(B_3)}
 \leq C_m+C_m\epsilon\Lambda_\epsilon^{C_m}.
\]
Outside $B_3$ its components obey \eqref{innersourceweightedbound},
so its global $**$ norm and the required angular versions are bounded
by the same quantity, independently of $R_\epsilon$.
Green's formula preserves its two translation conditions.  If $T_R$
denotes the normalised inverse of the lemma, then exactly
\[
 \phi_{1,\epsilon}=-\sigma_\epsilon+T_{R_\epsilon}\bar h_\epsilon.
\]
Thus the term of size $\Lambda_\epsilon$ on the fixed core is a pure
phase multiple of $W$.  The remaining core profile is bounded by
$C_m+C_m\epsilon\Lambda_\epsilon^{C_m}$ in every fixed Cartesian norm.
The decay proof above, now with this $**$ bound, also gives on
$3\leq r\leq7R_\epsilon/8$
\begin{equation}\label{first-tail-separated-logarithms}
\begin{split}
 \epsilon\big|D_y^\alpha(\psi_{1,\epsilon}-\psi_{1,\epsilon}^0)\big|
 +\epsilon|D_y^\alpha\psi_{2,\epsilon}|
 &\leq C_m\epsilon\ell(r)^{C_m}(1+r)^{-1-|\alpha|}\\
 &\quad+C_m\epsilon^2\Lambda_\epsilon^{C_m}(1+r)^{-|\alpha|}.
\end{split}
\end{equation}
The small zero-phase term obeys the second bound as well.  This
separation is valid after the indicated finite material derivatives.

We now compute the residual using the exact conjugated operator already
defined above.
Direct expansion gives
\begin{equation}\label{exact-linear-inner-residual-identity}
\begin{split}
 &\mathscr C_{\epsilon,\gamma,\omega}[W+\phi]
 -\mathscr C_{\epsilon,\gamma,\omega}[W]
 -\mathcal L_W\phi\\
 ={}&i\epsilon^2\Lambda_\epsilon D_t^x\phi
 +(\epsilon^2\Delta_{\gamma,\epsilon}-\Delta_y)\phi
 +2i\epsilon^2\nabla_{\gamma,\epsilon}\omega\cdot
             \nabla_{\gamma,\epsilon}\phi\\
 &+\big(i\epsilon^2\Delta_x\omega
        -\epsilon^2|\nabla_x\omega|^2
        -\epsilon^2\Lambda_\epsilon D_t^x\omega\big)\phi
 +\mathcal Q_W(\phi),
\end{split}
\end{equation}
where the nonlinear remainder is
\begin{equation}\label{exact-inner-nonlinearity}
 \mathcal Q_W(\phi)
 =-2\operatorname{Re}(\overline W\phi)\phi
   -|\phi|^2W-|\phi|^2\phi.
\end{equation}
The second term in \eqref{exact-linear-inner-residual-identity} is, by
\eqref{laplaceperiodicframe},
\begin{equation*}
\begin{split}
 (\epsilon^2\Delta_{\gamma,\epsilon}-\Delta_y)\phi
 ={}&-\frac{\epsilon\vec\kappa}{h}\cdot\nabla_y\phi
 +\frac{\epsilon^2}{h^2}\mathscr D_s^2\phi
 +\frac{\epsilon^3(D_s\vec\kappa)\cdot y}{h^3}
       \mathscr D_s\phi,\\
 &\hspace{8em} h=1-\epsilon\vec\kappa\cdot y.
\end{split}
\end{equation*}
Thus \eqref{exact-linear-inner-residual-identity} separately records the
normal geometric, horizontal, time, phase and nonlinear terms.

Take $\omega=\widehat\Phi_\epsilon$ and
$\phi=\epsilon\phi_{1,\epsilon}$ in the core where
$\eta_{\rm in}=1$.  Identity \eqref{exact-core-source-identity} gives
\begin{equation*}
 \mathscr C_{\epsilon,\gamma,\widehat\Phi_\epsilon}[W]
 =\epsilon\widehat H_\epsilon^{(1)}
   +\epsilon^2\mathscr R_{0,\epsilon},
 \qquad
 \mathscr R_{0,\epsilon}
 :=(1-\chi)\epsilon^{-2}\mathcal E_\epsilon.
\end{equation*}
This is smooth and identically zero on $B_2$.
On $r\geq2$, the exact expression
\eqref{exact-first-error-remainder-identity} and the differentiated phase
expansion give, for $a+b+|\alpha|+q\leq m$, $q\leq2$,
\begin{equation}\label{base-residual-radial-bound}
 |D_{\rm in}^q\mathscr D_t^a\mathscr D_s^bD_y^\alpha
      \mathscr R_{0,\epsilon}|
 \leq C_m\Lambda_\epsilon^{C_m}
       \ell(r)^{C_m}(1+r)^{-|\alpha|}.
\end{equation}
Using \eqref{firstinnercorrectiontail}, the five remaining groups in
\eqref{exact-linear-inner-residual-identity} obey, respectively,
\begin{equation}\label{first-inner-residual-ledger}
\begin{array}{c|c}
\text{term}&\text{bound after }D_{\rm in}^q\mathscr D_t^a
 \mathscr D_s^bD_y^\alpha\\ \hline
-\epsilon\vec\kappa h^{-1}\!\cdot\nabla_y\phi
&C_m\epsilon^2\Lambda_\epsilon^{C_m}
  \ell(r)^{C_m}(1+r)^{-2-|\alpha|}\\
\epsilon^2h^{-2}\mathscr D_s^2\phi
+\epsilon^3h^{-3}(D_s\vec\kappa\cdot y)\mathscr D_s\phi
&C_m\epsilon^3\Lambda_\epsilon^{C_m}
  \ell(r)^{C_m}(1+r)^{-1-|\alpha|}\\
i\epsilon^2\Lambda_\epsilon D_t^x\phi
&C_m\epsilon^2\Lambda_\epsilon^{C_m}
  \ell(r)^{C_m}(1+r)^{-1-|\alpha|}\\
\text{the two phase-coefficient terms}
&C_m\epsilon^2\Lambda_\epsilon^{C_m}
  \ell(r)^{C_m}(1+r)^{-1-|\alpha|}\\
\mathcal Q_W(\phi)
&C_m\epsilon^2\Lambda_\epsilon^{C_m}
  \ell(r)^{C_m}(1+r)^{-2-|\alpha|}.
\end{array}
\end{equation}
Here ``the two phase-coefficient terms'' means the third line and the
parenthesised fourth line of
\eqref{exact-linear-inner-residual-identity}; their bounds use that the
polar phase is contained in $W$.  The regular phase gradient obeys
$|\nabla_x\widehat\Phi_\epsilon|\leq C\Lambda_\epsilon^{C}$ in
the tube, with the corresponding radial derivative bounds, and has no
polar $1/(\epsilon r)$ term.  Factors $\epsilon\ell(r)^2$ from the
phase mean are absorbed using $\epsilon r\leq\delta$ and the permitted
logarithmic powers.

Since $\mathcal L_W\phi_{1,\epsilon}
=-H_{\epsilon,\mathrm{sol}}^{(1)}$, the order-$\epsilon$ terms cancel.
The remaining compact coefficient is
$\epsilon\chi_{\rm tr}\sum_jb_jZ_j
=O(\epsilon^2\Lambda_\epsilon^{C_m})$.
To pass from the additive calculation to \eqref{innergluedprofile},
subtract $W+\phi$ where $\eta_{\rm in}=1$.  The difference is exactly
\[
 D_\epsilon=(1-\eta_0)W
       \left(e^{\phi/W}-1-\phi/W\right),
\]
defined as zero on $B_1$.  Taylor's formula on $r\geq1$ and
\eqref{firstinnercorrectiontail} bound its mixed normal derivatives by
$C_m\epsilon^2\Lambda_\epsilon^{C_m}
\ell(r)^{C_m}(1+r)^{-2-|\alpha|}$.  The exact residual identity
applied once more to this difference preserves the required
$O(\epsilon^2\Lambda_\epsilon^{C_m})$ bound; derivatives of $\eta_0$
occur only in the fixed annulus $1<r<2$.
Therefore, define
\[
 \mathscr R_{\epsilon}^{\rm in}
 :=S_\epsilon(U_\epsilon^{\rm in})
 -\epsilon e^{i\widehat\Phi_\epsilon}
       \chi_{\rm tr}\sum_jb_jZ_j.
\]
\begin{equation}\label{innerimprovederrorbound}
\begin{split}
 &\max_{q\leq2}\ \max_{a+|\alpha|\leq m-q}
 \epsilon^{|\alpha|+q}(\epsilon\Lambda_\epsilon)^a
 \big|D_{\rm in}^q\partial_t^aD_x^\alpha
       \mathscr R_{\epsilon}^{\rm in}\big|\\
 &\hspace{8em}\leq C_m\epsilon^2\Lambda_\epsilon^{C_m},
 \qquad |z|\leq\frac13\delta.
\end{split}
\end{equation}
On the support of $D\eta_{\rm in}$ one has $r\simeq\epsilon^{-1}$ and
\eqref{firstinnercorrectiontail} gives, for
$a+|\alpha|+q\leq m$ and $q\leq2$,
\begin{equation}\label{first-cutoff-tail-bounds}
 |D_{\rm in}^q\partial_t^aD_x^\alpha
      (\epsilon\phi_{1,\epsilon})|
 \leq C_m\epsilon^2\Lambda_\epsilon^{C_m},
 \qquad \tfrac12\delta\leq|z|\leq\tfrac34\delta.
\end{equation}
The annulus is contained in $r\leq3R_\epsilon/4$, strictly within
the domain of the interior derivative bounds.  At fixed positive
physical distance, each factor $\epsilon^{-1}$ from differentiating
$y=z/\epsilon$ is cancelled by an extra $r^{-1}$ in the tail.
This applies also to fixed-$x$ time derivatives and curve variations;
the remaining coordinate and cutoff derivatives are uniformly bounded.
The linear cutoff commutator is
\begin{equation*}
\begin{split}
 &i\epsilon^2\Lambda_\epsilon(D_t^x\eta_{\rm in})
       \epsilon\phi_1
 +2\epsilon^2\nabla\eta_{\rm in}\cdot
       \nabla(\epsilon\phi_1)
 +\epsilon^2(\Delta\eta_{\rm in})\epsilon\phi_1\\
 &\qquad
 +2i\epsilon^2(\nabla\widehat\Phi_\epsilon\cdot
       \nabla\eta_{\rm in})\epsilon\phi_1.
\end{split}
\end{equation*}
Consequently it is $O(\epsilon^4\Lambda_\epsilon^{C_m})$, including
the indicated physical derivatives, and hence in the scaled family
norm on this annulus.  The exponential-minus-additive difference there
is $O(\epsilon^4\Lambda_\epsilon^{C_m})$ with the same derivative
bounds and contributes no larger residual.  This also records its support and all its
parameter-differentiated versions.

Finally, the construction leading to \eqref{regularphaseextension}, the
smoothness of \eqref{regularisedfirstinnersource}, and the regular
Frobenius choices in Lemma \ref{planarinversefirstimprovementlemma} show
that $U_\epsilon^{\rm in}$ is smooth at $z=0$.

\subsection{The outer wave correction}

Take a smooth tubular cutoff $\zeta_\delta$ satisfying
\begin{equation*}
 \zeta_\delta=0\quad\hbox{for }|z|\leq\tfrac14\delta,
 \qquad
 \zeta_\delta=1\quad\hbox{for }|z|\geq\tfrac13\delta.
\end{equation*}
On the support of $\zeta_\delta$, $U_\epsilon^{\rm in}$ is bounded away
from zero.  Define smooth real functions
\begin{equation}\label{outererrorsfirstimprovement}
 E_{1,\epsilon}+iE_{2,\epsilon}
 :=\zeta_\delta
   \frac{S_\epsilon(U_\epsilon^{\rm in})}
        {iU_\epsilon^{\rm in}}.
\end{equation}
The omitted part of the normalised residual is supported in
$|z|<\delta/3$ and is controlled by
\eqref{innerimprovederrorbound}.  On the outer region write
\begin{equation*}
 U_\epsilon^{\rm in}=\rho_\epsilon e^{i\Theta_\epsilon},
 \qquad
 \frac{\nabla U_\epsilon^{\rm in}}{U_\epsilon^{\rm in}}
 =\nabla\log\rho_\epsilon+i\nabla\Theta_\epsilon.
\end{equation*}
Only derivatives of $\Theta_\epsilon$ occur, so this definition is
independent of the local branch of its phase.

To derive the equations for the outer correction, we use the following
identity.  If $U$ does not vanish and $\psi=\psi_1+i\psi_2$, then direct
differentiation gives
\begin{equation}\label{multiplicativeerroridentity}
\begin{split}
 \frac{S_\epsilon(Ue^{i\psi})}{iUe^{i\psi}}
 ={}&\frac{S_\epsilon(U)}{iU}
 +i\epsilon^2\Lambda_\epsilon\partial_t\psi
 +\epsilon^2\Delta_x\psi
 +2\epsilon^2\frac{\nabla_xU}{U}\cdot\nabla_x\psi
 -2i|U|^2\psi_2\\
 &+i\epsilon^2(\nabla_x\psi)^2
 -i|U|^2\big(1-e^{-2\psi_2}-2\psi_2\big).
\end{split}
\end{equation}
The identity extends across the zeros after multiplication by $iUe^{i\psi}$.
The first line after $S_\epsilon(U)/(iU)$ is the outer linearised operator;
the second line is quadratic in $\psi$.

Dropping for a moment terms containing $\nabla\rho_\epsilon$ and replacing
$\rho_\epsilon$ by one, the real and imaginary parts of the linear problem
associated with \eqref{multiplicativeerroridentity} are
\begin{equation}\label{outerlinearsystemfirstimprovement}
\begin{split}
 -\epsilon^2\Lambda_\epsilon\partial_t\psi_2
 +\epsilon^2\Delta\psi_1
 -2\epsilon^2\nabla\Theta_\epsilon\cdot\nabla\psi_2
 +E_{1,\epsilon}&=0,\\
 \epsilon^2\Lambda_\epsilon\partial_t\psi_1
 +\epsilon^2\Delta\psi_2
 +2\epsilon^2\nabla\Theta_\epsilon\cdot\nabla\psi_1
 -2\psi_2+E_{2,\epsilon}&=0.
\end{split}
\end{equation}
In the second equation of \eqref{outerlinearsystemfirstimprovement}, the
spatial-derivative terms
\[
 \epsilon^2\Delta\psi_2
 \quad\hbox{and}\quad
 2\epsilon^2\nabla\Theta_\epsilon\cdot\nabla\psi_1
\]
are lower order relative to the mass term $-2\psi_2$.  We therefore set
\begin{equation}\label{outerpsi2definition}
 \psi_{2,\epsilon}^{\rm out}
 :=\frac12\left(E_{2,\epsilon}
 +\epsilon^2\Lambda_\epsilon
       \partial_t\psi_{1,\epsilon}^{\rm out}\right).
\end{equation}
Introduce the fast time
\begin{equation*}
    \tau=\frac{\sqrt2}{\epsilon\Lambda_\epsilon}t.
\end{equation*}
Substituting \eqref{outerpsi2definition} into the first equation of
\eqref{outerlinearsystemfirstimprovement} shows that its leading part
vanishes if
\begin{equation}\label{outerwaveequationfirstimprovement}
\left\{
\begin{aligned}
 -\partial_{\tau\tau}\psi_{1,\epsilon}^{\rm out}
 +\Delta_x\psi_{1,\epsilon}^{\rm out}
 &=F_\epsilon^{\rm out}(x,\tau),\\
 \psi_{1,\epsilon}^{\rm out}(x,0)
 &=\partial_\tau\psi_{1,\epsilon}^{\rm out}(x,0)=0,
\end{aligned}
\right.
\end{equation}
where all coefficients on the right are evaluated at
$t=\epsilon\Lambda_\epsilon\tau/\sqrt2$ and
\begin{equation}\label{outerwavesourcefirstimprovement}
 F_\epsilon^{\rm out}
 :=\frac{1}{\sqrt2\epsilon}\partial_\tau E_{2,\epsilon}
 +\nabla\Theta_\epsilon\cdot\nabla E_{2,\epsilon}
 -\epsilon^{-2}E_{1,\epsilon}.
\end{equation}
We solve this three-dimensional wave equation in the physical variable $x$.

We make the cancellation explicit.  With
$\psi_1=\psi_{1,\epsilon}^{\rm out}$ and
$\psi_2=\psi_{2,\epsilon}^{\rm out}$, the part of the real linearised
error retained in \eqref{outerlinearsystemfirstimprovement} equals
\begin{equation}\label{outerrealcancellation}
 -\sqrt2\epsilon^3
 \nabla\Theta_\epsilon\cdot\nabla\partial_\tau\psi_1,
\end{equation}
while the corresponding imaginary part is
\begin{equation}\label{outerimaginaryremainder}
 \epsilon^2\Delta\psi_2
 +2\epsilon^2\nabla\Theta_\epsilon\cdot\nabla\psi_1.
\end{equation}
The terms suppressed in deriving \eqref{outerlinearsystemfirstimprovement}
are
\begin{equation}\label{outeramplituderemainder}
 2\epsilon^2\nabla\log\rho_\epsilon\cdot\nabla\psi
 -2i(\rho_\epsilon^2-1)\psi_2,
\end{equation}
together with the quadratic second line of
\eqref{multiplicativeerroridentity}.  On the support of
$\zeta_\delta$, one has
\begin{equation}\label{outeramplitudeestimates}
 |\rho_\epsilon-1|+|\nabla\rho_\epsilon|
 +|\Theta_\epsilon-\varphi|
 +|\nabla(\Theta_\epsilon-\varphi)|
 \leq C\epsilon^2\Lambda_\epsilon^{C};
\end{equation}
hence all terms in \eqref{outeramplituderemainder} are lower order.

The real part of the normalised error in the matching region is important
because it is divided by $\epsilon^2$ in
\eqref{outerwavesourcefirstimprovement}.  We therefore record the complete
calculation.  Put
\begin{equation*}
 \mu_\epsilon:=\rho_\epsilon-1,
 \qquad
 \vartheta_\epsilon:=\Theta_\epsilon-\varphi.
\end{equation*}
On the fixed annulus
\begin{equation*}
 \mathcal A_{\delta,t}
 :=\{x:\tfrac14\delta\leq
       \operatorname{dist}(x,\Gamma_t)\leq\delta\}
\end{equation*}
one has $r\simeq\epsilon^{-1}$ and $r\geq2$, so
\eqref{innergluedprofile} gives
\begin{equation*}
 \mu_\epsilon
 =w(r)e^{-\epsilon\eta_{\rm in}\psi_{2,\epsilon}}-1,
 \qquad
 \vartheta_\epsilon
 =\epsilon\eta_{\rm in}\psi_{1,\epsilon}.
\end{equation*}
The differentiated vortex asymptotics and
\eqref{firstinnercorrectiontail}, used only on the support
$|z|\leq3\delta/4$ of $\eta_{\rm in}$, therefore imply, for
$a+|\alpha|+q\leq m$, $q\leq2$,
\begin{equation}\label{matching-amplitude-phase-estimate}
 |D_{\rm in}^q\partial_t^aD_x^\alpha\mu_\epsilon|
 +|D_{\rm in}^q\partial_t^aD_x^\alpha\vartheta_\epsilon|
 \leq C_m\epsilon^{2-q}\Lambda_\epsilon^{C_m}.
\end{equation}
Indeed, whenever a fixed-$x$ derivative falls on $y=z/\epsilon$, its
factor $\epsilon^{-1}$ is compensated by the additional factor
$r^{-1}\simeq\epsilon$ in the differentiated tails.  Derivatives falling
on the tubular coordinates or the cutoffs are harmless; those containing
derivatives of $\eta_{\rm in}$ are supported where
$\delta/2\leq\operatorname{dist}(x,\Gamma_t)\leq3\delta/4$.

Since $\Delta_x\varphi=0$ away from the curve, the continuity component of
\eqref{amplitudephaseidentity} is the exact expression
\begin{equation*}
\begin{split}
 E_{1,\epsilon}
 =\zeta_\delta\epsilon^2\bigg[
 &\Lambda_\epsilon
       \frac{\partial_t\mu_\epsilon}{1+\mu_\epsilon}
 +2\frac{\nabla\mu_\epsilon}{1+\mu_\epsilon}\cdot
       (\nabla\varphi+\nabla\vartheta_\epsilon)
 +\Delta\vartheta_\epsilon\bigg].
\end{split}
\end{equation*}
Here $\nabla\varphi$ and its derivatives are bounded on
$\mathcal A_{\delta,t}$.  Applying
\eqref{matching-amplitude-phase-estimate} with $m$ replaced by $m+2$ and
using the Leibniz rule gives
\begin{equation*}
 \|D_{\rm in}^q\partial_t^aD_x^\alpha E_{1,\epsilon}
       \|_{L^\infty(\mathcal A_{\delta,t})}
 \leq C_m\epsilon^{4-q}\Lambda_\epsilon^{C_m}
\end{equation*}
for $q\leq2$ and $a+|\alpha|\leq m-q$.  Since
$(\epsilon\Lambda_\epsilon)^a\leq1$, it follows that
\begin{equation}\label{matching-E1-differentiated-estimate}
 \max_{q\leq2}\max_{a+|\alpha|\leq m-q}
 \epsilon^q(\epsilon\Lambda_\epsilon)^a
 \|D_{\rm in}^q\partial_t^aD_x^\alpha
       E_{1,\epsilon}\|_{L^\infty(\mathcal A_{\delta,t})}
 \leq C_m\epsilon^4\Lambda_\epsilon^{C_m}.
\end{equation}
Thus $-\epsilon^{-2}E_{1,\epsilon}$ in the wave source is
$O(\epsilon^2\Lambda_\epsilon^{C_m})$, also after the spatial, fast-time
and first two curve derivatives used below.

To obtain a suitable weighted estimate for the solution of
\eqref{outerwaveequationfirstimprovement}, we use the far-field
cancellation furnished by the conservation law at the end of Section
\ref{normalgraphsection}.

\begin{lemma}[Far-field cancellation]
\label{farfieldfirstimprovementlemma}
Uniformly for $0\leq t\leq T$,
\begin{equation*}
 \varphi(x,t)
 =-\frac12\frac{\mathcal A_{\rm vec}(\gamma(t))\cdot x}{|x|^3}
 +O(|x|^{-3}),
 \qquad |x|\to\infty.
\end{equation*}
Consequently,
\begin{equation}\label{phasefarfieldderivativesfirstimprovement}
 |\nabla_x\varphi(x,t)|\leq C\langle x\rangle^{-3},
 \qquad
 |\partial_t\varphi(x,t)|\leq C\langle x\rangle^{-3}.
\end{equation}
Outside the fixed tube, where $U_\epsilon^{\rm in}=e^{i\varphi}$,
\begin{equation}\label{outerEfarfieldfirstimprovement}
 E_{1,\epsilon}=0,
 \qquad
 E_{2,\epsilon}
 =\epsilon^2\big(\Lambda_\epsilon\varphi_t
                 +|\nabla\varphi|^2\big)
 =O(\epsilon^2\Lambda_\epsilon\langle x\rangle^{-3}).
\end{equation}
The same estimate, with one additional factor $\Lambda_\epsilon$ at
worst, holds after one $t$ derivative.
\end{lemma}

\begin{proof}
Expand the kernel in \eqref{globalphase} for $|x|$ larger than twice the
diameter of the chosen spanning surface.  Since its induced boundary
orientation is opposite to that of $\gamma$,
\[
 \int_{\Sigma_t}\vec n_{\Sigma_t}\,d\Sigma
 =-\frac12\oint_\gamma\gamma\times d\gamma
 =-\mathcal A_{\rm vec}(\gamma(t)).
\]
The first term in the expansion is therefore
$-\frac12\mathcal A_{\rm vec}\cdot x/|x|^3$ and the next term is
$O(|x|^{-3})$.  Differentiation in $x$ gives the first estimate in
\eqref{phasefarfieldderivativesfirstimprovement}.  On differentiating in
$t$, the apparently larger term of order $|x|^{-2}$ is proportional to
$\partial_t\mathcal A_{\rm vec}$.  It vanishes by Lemma
\ref{vectorareaconservationlemma}; the remaining terms are
$O(|x|^{-3})$.  Formula \eqref{outerEfarfieldfirstimprovement} is the
normalised form of \eqref{outerfirsterror}.
\end{proof}

We use the following wave estimate.  The logarithmic loss in
its statement is more than sufficient for the present construction.

\begin{lemma}[A three-dimensional outer wave estimate]
\label{threeDwavefirstimprovementlemma}
Suppose that $F\in C^m(\mathbb R^3\times[0,\tau_*])$ and, for
$a+|\alpha|\leq m$,
\begin{equation}\label{wavesourceboundfirstimprovement}
 |\partial_\tau^aD_x^\alpha F(x,\tau)|
 \leq M_m\left(\langle x\rangle^{-3-|\alpha|}
 +\sum_{\nu=1}^{N_m}\mathbf 1_{K_{\nu,\tau}}(x)\right),
\end{equation}
where the sets $K_{\nu,\tau}$ have uniformly bounded diameter and remain
in a fixed bounded region.  The zero-data solution of
$-\psi_{\tau\tau}+\Delta\psi=F$ satisfies
\begin{equation*}
 \sum_{a+|\alpha|\leq m-1}
 \|\langle x\rangle\partial_\tau^aD_x^\alpha
       \psi(\cdot,\tau)\|_{L^\infty}
 \leq C_mM_m\log(2+\tau_*).
\end{equation*}
\end{lemma}

\begin{proof}
Kirchhoff's formula reads
\begin{equation*}
 \psi(x,\tau)
 =-\frac1{4\pi}\int_0^\tau(\tau-q)
   \int_{\mathbb S^2}F(x+(\tau-q)\omega,q)
   \,d\omega\,dq.
\end{equation*}
The compactly supported terms may be absorbed into
$C M_m\langle x\rangle^{-3}$ because their supports remain in a fixed
bounded set.  Put $r=|x|$ and $s=\tau-q$.  For $r>0$,
\begin{equation*}
 s\int_{\mathbb S^2}\langle x+s\omega\rangle^{-3}\,d\omega
 =\frac{2\pi}{r}\int_{|r-s|}^{r+s}
       \frac{\rho}{\langle\rho\rangle^3}\,d\rho.
\end{equation*}
If $2\leq r\leq2\tau_*+2$, Fubini's theorem gives
\begin{align*}
 \int_0^\tau s\int_{\mathbb S^2}
       \langle x+s\omega\rangle^{-3}\,d\omega\,ds
 &\leq\frac{C}{r}\left(
       \int_0^r\frac{\rho^2}{\langle\rho\rangle^3}\,d\rho
       +r\int_r^\infty\frac{\rho}{\langle\rho\rangle^3}\,d\rho
       \right)\\
 &\leq\frac{C\log(2+\tau_*)}{r}.
\end{align*}
For $r>2\tau_*+2$, the same integral is at most
$C\tau_*^2r^{-3}\leq Cr^{-1}$, while the region $r\leq2$ is bounded
directly.  This proves the asserted weighted estimate for spatial
derivatives.  Time differentiation of the Duhamel formula produces the
same integral with time derivatives of $F$, together with homogeneous
spherical means of their initial traces; the latter satisfy the same
$\langle x\rangle^{-1}$ bound.  This proves the result.
\end{proof}

The direct estimates from \eqref{outererrorsfirstimprovement} and Lemma
\ref{farfieldfirstimprovementlemma} give for each fixed $m$
\begin{equation*}
 |\partial_\tau^aD_x^\alpha F_\epsilon^{\rm out}|
 \leq C_m\epsilon^2\Lambda_\epsilon^{C_m}
 \left(\langle x\rangle^{-3-|\alpha|}
 +\mathbf 1_{K_t}(x)\right),
 \qquad a+|\alpha|\leq m,
\end{equation*}
where $K_t$ is a fixed-width tube about the compact curve.  Here the
matching estimate, including its support and differentiated family form,
is \eqref{matching-E1-differentiated-estimate}; the
$O(\epsilon^2\Lambda_\epsilon^C)$ amplitude defect enters
$E_{2,\epsilon}$ instead.  Since
\[
 0\leq\tau\leq
 \tau_\epsilon:=\frac{\sqrt2T}{\epsilon\Lambda_\epsilon},
 \qquad
 \log(2+\tau_\epsilon)\leq C\Lambda_\epsilon,
\]
Lemma \ref{threeDwavefirstimprovementlemma} yields a solution of
\eqref{outerwaveequationfirstimprovement} satisfying
\begin{equation}\label{outerpsi1estimatefirstimprovement}
 \sum_{a+|\alpha|\leq m-1}
 \|\langle x\rangle\partial_\tau^aD_x^\alpha
       \psi_{1,\epsilon}^{\rm out}\|_{L^\infty}
 \leq C_m\epsilon^2\Lambda_\epsilon^{C_m}.
\end{equation}
Equation \eqref{outerpsi2definition} then gives
\begin{equation}\label{outerpsi2estimatefirstimprovement}
 \sum_{a+|\alpha|\leq m-2}
 \|\langle x\rangle\partial_\tau^aD_x^\alpha
       \psi_{2,\epsilon}^{\rm out}\|_{L^\infty}
 \leq C_m\epsilon^2\Lambda_\epsilon^{C_m}.
\end{equation}

\subsection{Conclusion of the first improvement}

We collect the preceding construction in a proposition.  To record the
dependence of the residual on the initial curve, for a global function
$G=G(x,t;\Gamma_{\rm in})$ define the scaled family seminorm
\begin{equation}\label{first-global-family-norm}
 \|G\|_{\mathfrak R_m^{(2)}}
 :=\max_{0\leq q\leq2}\
   \max_{a+|\alpha|\leq m-q}
   \sup_{\Gamma_{\rm in},x,t}
   \epsilon^{|\alpha|+q}
   (\epsilon\Lambda_\epsilon)^a
   |D_{\rm in}^q\partial_t^aD_x^\alpha G|.
\end{equation}
The factors in this definition compensate for the losses caused by
differentiating the vortex core, the moving tubular coordinates and the
fast-time scaling.

\begin{proposition}[First global improvement]
\label{firstglobalimprovementproposition}
Let $\gamma$ be a smooth embedded solution of the corrected curve equation
\eqref{geometriccurveequation} on $[0,T]$, with the uniform geometric
bounds in Theorem \ref{uniformlocalcurvetheorem}.  Then, for all sufficiently
small $\epsilon$, there is a smooth approximation
$u_\epsilon^{(1)}$ of the form \eqref{firstglobalgluingansatz} such that
\begin{equation}\label{globalfirstimprovementestimate}
 \|S_\epsilon(u_\epsilon^{(1)})\|_{L^\infty(\mathbb R^3\times[0,T])}
 +\epsilon
 \|\nabla_xS_\epsilon(u_\epsilon^{(1)})\|_{L^\infty(\mathbb R^3\times[0,T])}
 \leq C\epsilon^2\Lambda_\epsilon^C.
\end{equation}
The construction may be performed simultaneously for
$\Gamma_{\rm in}\in\mathscr O_{\rm in}$.  For every fixed $m$ the resulting
family satisfies
\begin{equation}\label{global-first-improvement-family-estimate}
 \big\|S_\epsilon(u_\epsilon^{(1)})\big\|_{\mathfrak R_m^{(2)}}
 \leq C_m\epsilon^2\Lambda_\epsilon^{C_m}.
\end{equation}
In the inner variables, for every fixed $R>0$,
\begin{equation}\label{profilefirstimprovementestimate}
\begin{split}
 u_\epsilon^{(1)}\big(\gamma(s,t)
 +\epsilon y_\alpha\vec e_\alpha(s,t),t\big)
 =e^{i\beta_\epsilon(s,t)}
 \left[W(y)+\epsilon\Phi_{1,\epsilon}(s,y,t)
 +O_R(\epsilon^2\Lambda_\epsilon^C)\right],
\end{split}
\end{equation}
where $\Phi_{1,\epsilon}$ is smooth in $(s,y,t)$ and
\begin{equation*}
 \sup_{0\leq t\leq T}\sup_{s\in\mathbb R/\ell(t)\mathbb Z}
 \sum_{a+b+|\alpha|\leq2}
 |\mathscr D_t^a\mathscr D_s^bD_y^\alpha
       \Phi_{1,\epsilon}(s,y,t)|
 \leq C_R\Lambda_\epsilon^C,
 \qquad |y|\leq R.
\end{equation*}
In addition, the approximation is exactly centred:
\begin{equation}\label{first-improvement-exact-centring}
 u_\epsilon^{(1)}(\gamma(s,t),t)=0.
\end{equation}
This identity holds for the whole initial-curve family, and hence also
after its first two full parameter derivatives when the evaluation point
is varied with the curve.
Moreover,
\begin{equation}\label{firstimprovementouterbounds}
 \max_{q\leq2}\max_{a+|\alpha|\leq m-q}
 \epsilon^q\big(
 \|D_{\rm in}^q\partial_\tau^aD_x^\alpha
       \psi_{1,\epsilon}^{\rm out}\|_{L^\infty}
 +\|D_{\rm in}^q\partial_\tau^aD_x^\alpha
       \psi_{2,\epsilon}^{\rm out}\|_{L^\infty}\big)
 \leq C_m\epsilon^2\Lambda_\epsilon^{C_m}
\end{equation}
for every fixed $m\geq4$.  Moreover, outside a fixed tubular neighbourhood,
\begin{equation}\label{outerremainderfirstimprovement}
 |S_\epsilon(u_\epsilon^{(1)})(x,t)|
 \leq C\Lambda_\epsilon^C\left(
       \epsilon^4\langle x\rangle^{-4}
       +\epsilon^5\langle x\rangle^{-1}\right).
\end{equation}
\end{proposition}

\begin{proof}
The curve equation gives the two compatibility conditions, up to the
$O(\epsilon\Lambda_\epsilon^{C_m})$ projection defect recorded in
Lemma \ref{regularisedsourceorthogonalitylemma}, including the compact
contribution from the exact core residual.  Lemma
\ref{planarinversefirstimprovementlemma} then constructs
$\phi_{1,\epsilon}$, and \eqref{firstinnercorrectionequation} cancels the
solvable order-$\epsilon$ coefficient in the core; the projected defect is
already of order $\epsilon^2\Lambda_\epsilon^{C_m}$ in physical variables.  The
identity \eqref{exact-linear-inner-residual-identity}, the radial
ledger \eqref{first-inner-residual-ledger}, and the displayed cutoff
commutator prove \eqref{innerimprovederrorbound}.

In the outer region, insert \eqref{outerpsi2definition} and
\eqref{outerwaveequationfirstimprovement} into the identity
\eqref{multiplicativeerroridentity}.  The terms involving
$E_{1,\epsilon}$ and $E_{2,\epsilon}$ cancel, and the remaining linear
terms are precisely \eqref{outerrealcancellation},
\eqref{outerimaginaryremainder}, and
\eqref{outeramplituderemainder}.  Estimates
\eqref{outeramplitudeestimates},
\eqref{outerpsi1estimatefirstimprovement}, and
\eqref{outerpsi2estimatefirstimprovement} bound them by
$C\epsilon^4\Lambda_\epsilon^C$ away from the core.  The quadratic terms
in \eqref{multiplicativeerroridentity} obey the same bound.  Outside a
fixed tubular neighbourhood one has $\rho_\epsilon=1$,
$\Theta_\epsilon=\varphi$, and
 \begin{equation*}
  \psi_{2,\epsilon}^{\rm out}
  =\frac12E_{2,\epsilon}
  +\frac{\epsilon}{\sqrt2}
       \partial_\tau\psi_{1,\epsilon}^{\rm out}.
 \end{equation*}
Lemma \ref{farfieldfirstimprovementlemma} and the weighted estimate
\eqref{outerpsi1estimatefirstimprovement} give
\begin{equation*}
 |\Delta E_{2,\epsilon}|
 \leq C\epsilon^2\Lambda_\epsilon^C\langle x\rangle^{-5},
 \qquad
 |D_x^\alpha\partial_\tau^a\psi_{1,\epsilon}^{\rm out}|
 \leq C\epsilon^2\Lambda_\epsilon^C\langle x\rangle^{-1}
 \quad(a+|\alpha|\leq3).
\end{equation*}
Consequently, the real remainder \eqref{outerrealcancellation} is
$O(\epsilon^5\Lambda_\epsilon^C\langle x\rangle^{-4})$.  The two terms
in \eqref{outerimaginaryremainder} are respectively
\begin{equation*}
 O\left(\epsilon^4\Lambda_\epsilon^C\langle x\rangle^{-5}
       +\epsilon^5\Lambda_\epsilon^C\langle x\rangle^{-1}\right)
 \quad\hbox{and}\quad
 O\left(\epsilon^4\Lambda_\epsilon^C\langle x\rangle^{-4}\right).
\end{equation*}
The terms in \eqref{outeramplituderemainder} vanish in this region, while
the quadratic terms are
$O(\epsilon^4\Lambda_\epsilon^C\langle x\rangle^{-6}
 +\epsilon^6\Lambda_\epsilon^C\langle x\rangle^{-2})$.
These estimates prove \eqref{outerremainderfirstimprovement}.

Inside $|z|<\delta/3$, multiplication by the outer factor changes the
inner residual by at most $C\epsilon^2\Lambda_\epsilon^C$, which is the
size allowed in \eqref{globalfirstimprovementestimate}.  The fixed annuli
where the cutoffs vary are controlled by
\eqref{first-cutoff-tail-bounds} and
\eqref{matching-E1-differentiated-estimate}.  This
proves \eqref{globalfirstimprovementestimate}.

All maps used above are defined on the common chart ball
$\mathscr O_{\rm in}$.  Differentiation with respect to
$\Gamma_{\rm in}$ commutes with the Fourier-mode inverse and with the
zero-data wave solution.  Derivatives of tubular coordinates and cutoffs
are exactly the losses included in
\eqref{first-global-family-norm}.  Equations
\eqref{innersourceweightedbound}, \eqref{firstinnercorrectiontail},
\eqref{first-inner-residual-ledger} and
\eqref{matching-E1-differentiated-estimate} therefore give
\eqref{global-first-improvement-family-estimate}, while differentiated
Kirchhoff estimates give \eqref{firstimprovementouterbounds}.  This proves
the first and second initial-curve variation estimates uniformly in
$\epsilon$.

Finally, in every fixed inner disc the regularised phase and the additive
part of \eqref{innergluedprofile} give
\eqref{profilefirstimprovementestimate}.  Its coefficient is smooth by
construction, and the stated estimates follow from
\eqref{firstinnercorrectionestimate} and
\eqref{firstimprovementouterbounds}.  At $y=0$ both $W$ and
$\phi_{1,\epsilon}$ vanish by
\eqref{firstinnercorrectionequation}; the regular phase and the outer
factor in \eqref{firstglobalgluingansatz} are finite exponentials and
therefore nonzero.  This proves
\eqref{first-improvement-exact-centring}, including its differentiated
family version.  This completes the first global improvement.
\end{proof}

\section{Approximation to arbitrary order}
\label{arbitraryapproxsection}

We now iterate the construction of the preceding section.  The main result
is that, for any prescribed $M$ and $N$, one can construct a smooth family
of approximations whose residual and first two variations with respect to
the initial curve are $O(\epsilon^M)$ in $H^N(\mathbb R^3)$.  The
approximations remain centred on corrected curves and are obtained from
the preceding construction by finitely many small corrections.

The construction is triangular.  At each stage, earlier outer corrections
are kept fixed as functions of physical space and time.  A curve correction,
obtained using the same reference linearisation at every stage, cancels the
leading translation contribution.  A finite batch of planar inversions then
reduces the remaining inner error in a fixed physical tube, allowing
controlled polynomial growth in the normal variable.  A new outer
correction, with phase obtained from a three-dimensional wave equation and
amplitude determined algebraically, reduces the resulting outer error.
The remaining errors are transferred to the next stage.

Throughout this section put
\[
 \Lambda_\epsilon=|\log\epsilon|,\qquad
 c_\epsilon=\epsilon^2\Lambda_\epsilon,\qquad
 \tau=\frac{\sqrt2\,t}{\epsilon\Lambda_\epsilon},\qquad
 \tau_\epsilon=\frac{\sqrt2 T}{\epsilon\Lambda_\epsilon}.
\]
The corresponding material time derivative is
\begin{equation}\label{fasttimehorizontalderivative}
 \mathscr D_\tau
 =\frac{\epsilon\Lambda_\epsilon}{\sqrt2}\mathscr D_t.
\end{equation}
All orders of differentiation and all numbers of corrections below are
finite and are chosen before $\epsilon$ is made small.  Constants and
powers of $\Lambda_\epsilon$ may depend on these orders.

\subsection{Statement of the construction}

Recall the initial-curve neighbourhood $\mathscr O_{\rm in}$ from Section
\ref{firstglobalimprovementsection}.  We write $D_{\rm in}$ for
 differentiation with respect to the initial curve.

\begin{theorem}[Approximation to arbitrary order]
\label{arbitraryapproximationtheorem}
Let $M\geq4$ and $N\geq0$ be integers.  There exist finite integers $K$ and
$s_*$, a neighbourhood $\mathscr O_{\rm in}$ bounded in $H^{s_*}$, and
$\epsilon_0>0$, such that, for $0<\epsilon<\epsilon_0$, one can construct
smooth families of constant-speed curves and global approximations
\begin{equation}\label{curvefamilyarbitrary}
 \gamma_\epsilon^{(k)}[\Gamma_{\rm in}],\qquad
 u_\epsilon^{(k)}[\Gamma_{\rm in}],\qquad 1\leq k\leq K,
\end{equation}
on $[0,T]$.  The first pair is the pair constructed in Sections
\ref{normalgraphsection} and \ref{firstglobalimprovementsection}.  All the
curves have initial geometric curve $\Gamma_{\rm in}$.  The terminal
approximation satisfies
\begin{equation}\label{arbitraryresidualtheorem}
 \sup_{0\leq t\leq T}\sum_{q=0}^2
 \big\|D_{\rm in}^qS_\epsilon(u_\epsilon^{(K)})(t)
                  \big\|_{H^{N-q}(\mathbb R^3)}
 \leq C_{M,N}\epsilon^M.
\end{equation}
The norm of a differential is its operator norm on unit directions in
$H^{s_*}$.

The zero set of $u_\epsilon^{(K)}$ in the fixed tube is precisely
$\gamma_\epsilon^{(K)}(t)$.  In particular,
\begin{equation}\label{arbitrary-exact-centring-identity}
 u_\epsilon^{(K)}[\Gamma_{\rm in}]
       (\gamma_\epsilon^{(K)}(s,t),t)=0,
\end{equation}
including the first two derivatives of this identity with respect to
$\Gamma_{\rm in}$.  On the fixed material circle
$\rho\in\mathbb T=\mathbb R/\mathbb Z$, put
$s_k(\rho,t)=\ell_\epsilon^{(k)}(t)\rho$ and remove the value of the phase
at the centre to define $u_{\epsilon,\mathrm{cen}}^{(k)}(\rho,y,t)$.
Then
\begin{equation}\label{arbitrary-centred-normal-derivative}
 \sup_{\rho,t}
 \big|D_yu_{\epsilon,\mathrm{cen}}^{(K)}(\rho,0,t)-DW(0)\big|
 \leq C\epsilon\Lambda_\epsilon^C,
\end{equation}
and, for every fixed $R$,
\begin{equation}\label{arbitraryinnerprofile}
 \sup_{\rho,t}
 \big\|u_{\epsilon,\mathrm{cen}}^{(K)}(\rho,\cdot,t)-W
                  \big\|_{C^N(B_R)}
 \leq C_{R,N}\epsilon\Lambda_\epsilon^{C_{R,N}}.
\end{equation}
The same bound holds for the difference of two successive centred
approximations:
\begin{equation}\label{successive-approximation-bound}
 \sup_{\rho,t}
 \big\|u_{\epsilon,\mathrm{cen}}^{(k)}
       -u_{\epsilon,\mathrm{cen}}^{(k-1)}\big\|_{C^N(B_R)}
 \leq C_{R,N,k}\epsilon\Lambda_\epsilon^{C_{R,N,k}}.
\end{equation}
Finally,
\begin{equation}\label{sharperarbitrarycurvebound}
 \sup_t\big\|
 \gamma_\epsilon^{(K)}(\ell_\epsilon^{(K)}(t)\,\cdot,t)
 -\gamma_\epsilon^{(1)}(\ell_\epsilon^{(1)}(t)\,\cdot,t)
                    \big\|_{H^N(\mathbb T)}
 \leq C_N\epsilon\Lambda_\epsilon^{C_N}.
\end{equation}
All constants are uniform on $\mathscr O_{\rm in}$.
\end{theorem}

\subsection{The global field and the residual used in the iteration}

We first give the recursive form of the objects in Theorem
\ref{arbitraryapproximationtheorem}.  Fix $d>0$ with $10d<\delta$ and
choose radial cutoffs
\[
 \chi_d=1\quad\hbox{on }|z|\leq2d,\qquad
 \chi_d=0\quad\hbox{on }|z|\geq3d,
\]
\[
 \chi_{\rm inv}=1\quad\hbox{on }|z|\leq4d,\qquad
 \chi_{\rm inv}=0\quad\hbox{on }|z|\geq5d.
\]
Their rescaled arguments are $z=\epsilon y$.  Thus
\begin{equation}\label{matchingcutoffderivatives}
 |D_y^a\chi_d(\epsilon y)|+|D_y^a\chi_{\rm inv}(\epsilon y)|
 \leq C_{a,d}\epsilon^a,
\end{equation}
and on a derivative annulus these quantities are also bounded by
$C_{a,d}(1+r)^{-a}$.

We use the constant-speed chart $\mathfrak C_t$ of Section
\ref{normalgraphsection}, expressed on $\mathbb T$.  For a chart path
$z(t)$ write
\[
 \Gamma^{\rm mat}(\rho,t)=\mathfrak C_t(z(t))(\rho)
      =\gamma(\ell(t)\rho,t),\qquad
 V^{\rm mat}=\partial_t\Gamma^{\rm mat}.
\]
The tubular map and its normal coordinate are
\[
 \mathcal X_\Gamma(\rho,y,t)
 =\Gamma^{\rm mat}(\rho,t)+\epsilon y_\alpha e_\alpha(\rho,t),
 \qquad r=|y|.
\]
Time differentiation in the tubular formula always uses this material
velocity, including its tangential component.

Let $B_\epsilon[\Gamma]$ be the inner factor of Section
\ref{firstglobalimprovementsection}, with its first correction computed
using the basic corrected velocity $\mathscr G_\epsilon(t,z)$.
In the exact core source \eqref{regularisedfirstinnersource}, all time
jets used for this seed are those induced by the basic corrected
evolution from that curve (including the frame and phase jets).
This defines the factor for any curve in the chart, whether or not its
actual velocity equals that basic velocity.  The discrepancy of the two
velocities remains in its residual.  Write $\phi_{1,\epsilon}[\Gamma]$
for this first, centred correction.  For an additional material profile
$\Pi(\rho,y,\tau)$, define
\begin{equation}\label{matchedglobalinnercorrection}
 \Pi^{\rm g}(x,t)
 =\chi_d(z_\Gamma(x,t))\,
   \Pi(\rho_\Gamma(x,t),z_\Gamma(x,t)/\epsilon,\tau).
\end{equation}
Every profile is centred: $\Pi(\rho,0,\tau)=0$.  With the fixed core
cutoff $\eta_0$ from Section \ref{firstglobalimprovementsection}, put
\begin{equation}\label{exact-global-inner-factor}
 \mathfrak A_W(\Xi)
 =\eta_0(W+\Xi)+(1-\eta_0)W e^{\Xi/W},\qquad
 I_\epsilon[\Gamma,\Pi]
 =e^{i\widehat\Phi_\epsilon[\Gamma]}
   \mathfrak A_W(\epsilon\phi_{1,\epsilon}[\Gamma]+\Pi^{\rm g})
\end{equation}
in $|z|<3d$.  Outside that region set
$I_\epsilon[\Gamma,\Pi]=B_\epsilon[\Gamma]$.  The formulae agree where
$\Pi^{\rm g}=0$, so they define a smooth global field.
For real global functions $A,B$, our approximations have the form
\begin{equation}\label{stagewise-approximation-form}
 U=I_\epsilon[\Gamma,\Pi]e^{iA-B}.
\end{equation}
For $k\geq2$, $A$ and $B$ will be finite sums of phase and amplitude
corrections $a_k,b_k$.  Both are retained as physical space--time
functions.  Higher inner profiles are retained as material functions;
when the curve is changed, they are inserted using the new tubular map.
The first correction $\phi_{1,\epsilon}[\Gamma]$ and the solid-angle
phase are recomputed from the new curve, using the periodic normal frame
and real relative phase supplied by Lemma
\ref{longitudinal-framing-lemma}.  The first-correction estimates are
used here only for $|z|<5d<\delta/2$, inside their domain
$r<7R_\epsilon/8$.

For residual calculations choose a radial $\eta_*$ equal to one on
$B_4$ and zero outside $B_5$, and put, in the tube,
\[
 \Xi=\epsilon\phi_{1,\epsilon}[\Gamma]+\chi_d\Pi,
 \quad
 \nu_\Xi=\eta_*+(1-\eta_*)e^{\Xi/W},
 \quad
 \mathcal R_\Gamma(\Pi;A,B)
 =\frac{S_\epsilon(U)}{e^{i\widehat\Phi_\epsilon+iA-B}\nu_\Xi}.
\]
The denominator is bounded above and below for the profiles constructed
below.  This residual is smooth at the centre.  For $r\geq5$,
\begin{equation}\label{fixed-tube-normalised-residual}
 \mathcal R_\Gamma=iW(E_1+iE_2),\qquad
 E_1+iE_2=\frac{S_\epsilon(U)}{iU}.
\end{equation}
On bounded inner discs we estimate $\mathcal R_\Gamma$ as a complex
field; outside them we estimate its two components relative to $iW$.

An elliptic substep is defined from this actual residual.  If $\mathcal K$
is a planar right inverse chosen as in Lemma
\ref{polynomial-planar-inverse-lemma} below, set
\begin{equation}\label{complete-inner-coefficient-recursion}
 \phi=-\mathcal K\{\chi_{\rm inv}\mathcal R_\Gamma(\Pi;A,B)\},
 \qquad \Pi_{\rm new}=\Pi+\phi.
\end{equation}
The residual after the substep is evaluated using
\eqref{exact-global-inner-factor}; thus this definition includes all
terms of the additive--multiplicative map.
More explicitly, for a smooth residual map $\mathcal F$ and an increment
$\phi$,
\begin{equation}\label{exact-inner-composite-increment}
 \begin{split}
 \mathcal F(\Pi+\phi)={}&\mathcal F(\Pi)+\mathcal L_W\phi
 +(D\mathcal F(\Pi)-\mathcal L_W)\phi\\
 &+\int_0^1(1-s)D^2\mathcal F(\Pi+s\phi)[\phi,\phi]ds.
 \end{split}
\end{equation}
This identity will also be used after parameter differentiation.

\subsection{Polynomial inverses and differentiated radial bounds}

For a source $h=iW(h_1+ih_2)$, write
\[
 h_j=O_m(\epsilon^a\langle r\rangle^p)
\]
when every indicated material, $\tau$ and normal derivative satisfies
\[
 |\partial_\rho^b\partial_\tau^cD_y^\alpha h_j|
 \leq C\epsilon^a\Lambda_\epsilon^C
         \langle r\rangle^{p-|\alpha|},\qquad
 b+c+|\alpha|\leq m,
\]
for $r\geq1$.  Logarithmic factors in $r$ are included in the displayed
power of $\Lambda_\epsilon$, since $r\leq5d/\epsilon$.  On $B_5$ the
corresponding assertion is the $C^m$ bound for the smooth field $h$.
Different radial exponents are used for different terms.  For inner functions, family derivatives are taken at fixed
$(\rho,y,\tau)$; for outer functions they are taken at fixed
$(x,\tau)$.  The conversion between these conventions uses the
tubular chain rule.  For families, put
\begin{equation}\label{weighted-initial-curve-derivative}
 \mathcal D_{\rm in}^q=\epsilon^qD_{\rm in}^q,
 \qquad 0\leq q\leq2.
\end{equation}
A superscript $(2)$ includes these variations, with their directions
measured in the chosen initial-curve Sobolev space.
For outer functions the notation is
\begin{equation}\label{differentiatedouternorms}
 \|F\|_{Y_m}^{(2)}
   =\max_{q\leq2}\|\mathcal D_{\rm in}^qF\|_{Y_{m-q}},
 \qquad
 \|F\|_{\mathcal H_m}^{(2)}
   =\max_{q\leq2}\|\mathcal D_{\rm in}^qF\|_{\mathcal H_{m-q}}.
\end{equation}
We also use $\mathcal Y_m^{(2)}$ and $\mathcal H_m^{(2)}$ as
subscripts for these family norms.  Only nonnegative derivative orders
are used.  The same convention applies to $X_m$ and to the joint
$W^{m,\infty}_{x,\tau}$ norm.
The norms $Y_m,\mathcal H_m,X_m$ are defined in Section
\ref{threeDwaveestimatesection}.

\begin{lemma}[Planar inversion with polynomial growth]
\label{polynomial-planar-inverse-lemma}
Let $p\geq1$ be an integer.  Suppose that $h$ is smooth, is cut off at
$r=5d/\epsilon$, and its components have bounds
\[
 h_1=O_{4m+20}(A\langle r\rangle^{p-2}),\qquad
 h_2=O_{4m+20}(A\langle r\rangle^q),\qquad q\leq p.
\]
There is a linear map $h\mapsto\phi$, independent of the curve
parameters, such that
\[
 \mathcal L_W\phi=h,\qquad \phi=iW(\psi_1+i\psi_2),
\]
with
\begin{equation}\label{polynomial-planar-inverse-bounds}
 \psi_1=O_m(A\langle r\rangle^p),\qquad
 \psi_2=O_m(A\langle r\rangle^{\max(q,p-2)}).
\end{equation}
It can be chosen with $\phi(0)=0$ and
$\operatorname{Re}\int_{B_2}\phi\overline{iW}=0$.
The constants are independent of $\epsilon$ apart from the indicated
logarithms.  The same conclusions hold with two weighted family
variations.

For $h_1=O_{4m+20}(A\langle r\rangle^{-2})$ and
$h_2=O_{4m+20}(A)$, vanishing of the two translation moments gives instead
$\psi_1=O_m(A)$, $\psi_2=O_m(A)$, with logarithmic factors allowed.
A compact translation-moment source of size $a$ has an inverse bounded
by $O_m(a\langle r\rangle)$ in phase and $O_m(a)$ in amplitude.
A radial imaginary source $iW(ih_2)$ with
$h_2=O_{4m+20}(A\langle r\rangle^{-2})$ has a radial purely imaginary
inverse with the same radial decay.
\end{lemma}

\begin{proof}
Use the mode formulae of Propositions 3.8 and 4.8 of
\cite{delpinojunemanmussolinear2025}.  The zero phase mode is obtained by
two integrations and the zero amplitude mode by the massive Green
formula.  In a nonzero mode $k$, the algebraic homogeneous fields have
phase orders $r^{-k},r^k$ and amplitude orders $r^{-k-2},r^{k-2}$;
there are also an exponentially increasing and an exponentially
decreasing field.  The latter two are paired with integrals towards
infinity and towards zero, respectively, so no increasing exponential
is introduced.

For $2\leq k\leq p$, choose the primitive multiplying the algebraically
increasing field with a fixed finite reference point.  That field is
regular at the origin.  For $k>p$, integrate its coefficient towards the
outer endpoint.  The integral bounds
\[
 \int_1^r s^{a-1}(1+\log s)^bds
 =O\big(r^{\max(a,0)}(1+\log r)^{b+1}\big)
\]
and the corresponding convergent tail bound give the phase order $p$.
The amplitude order is the larger of its source order $q$ and the
phase-induced order $p-2$.  At $k=p$ a logarithm is allowed.  In mode
one the regular formula integrates the coefficient of the increasing
algebraic field from zero; its possible $r$ growth is admitted because
$p\geq1$.  These choices allow the prescribed polynomial growth.

For completeness, the large angular modes can be summed uniformly.
After changing the sign of the second component, the radial system is
cooperative.  On $r\geq R_0$, where $w^2\geq3/4$, use a phase barrier
$C\langle r\rangle^p$ and an amplitude barrier
\[
 C\big(\langle r\rangle^q+4k\langle r\rangle^{p-2}\big).
\]
For $R_0$ sufficiently large and
$k>2(p+|q|+m)+20$, the diagonal $k^2/r^2$ term absorbs the angular
coupling in the first equation.  In the second equation the massive
term applied to $4k\langle r\rangle^{p-2}$ absorbs the coupling from
the phase barrier; the remaining massive term controls the amplitude
source.  On $r\leq R_0$, replace the amplitude barrier smoothly by a
constant multiple of the phase barrier.  The $k^2/r^2$ diagonal terms
absorb the derivatives of this interpolation.  The comparison principle
on a finite annulus, with zero outer boundary values, gives
\[
 |\psi_{1,k}|\leq C A_k\langle r\rangle^p,\qquad
 |\psi_{2,k}|\leq C(1+k)A_k
                 \langle r\rangle^{\max(q,p-2)}.
\]
The regular physical mode vanishes at the inner endpoint as it tends
to zero.  Passing to the limit at both endpoints gives the same
polynomially bounded solution selected by the mode formula when $k>p$.
There is no homogeneous polynomial ambiguity in this range.

Differentiating the radial equations gives the derivative estimates.
At derivative order $j$, the potential commutators have bounds
\[
 C_jk^2r^{-2-a}
\]
and multiply lower radial derivatives.  Applying the
same barriers with their radial degrees lowered by $j$, and using
interior elliptic estimates on $r\leq2R_0$, gives constants bounded by
$C_m(1+k)^{2m+4}$ through order $m$.  Fourier coefficients of the
$4m+20$ differentiated source decay faster than this power, also after
the remaining material, time and angular derivatives.  Their series
are therefore absolutely summable with the claimed radial weights.
The finitely many excluded modes are handled by their explicit
formulae.  This proves the full differentiated estimate.
The two zero-order bounds with translation orthogonality are Theorem
1.1 of \cite{delpinojunemanmussolinear2025}.  Its mode-one formula gives
the stated bound for a compact nonorthogonal source.  The radial
amplitude assertion follows directly from the exponentially decaying
Green kernel and $w=1+O(r^{-2})$.  Adding translations and $iW$ imposes
the two normalisations without changing the admissible bounds.  The
inverse, the source cutoffs and the normalisations are fixed in $y$;
therefore they commute with material, time and family differentiation.
\end{proof}

\subsection{The elliptic batches in a fixed tube}

We record the coefficient estimates used in the batches.  Write the
inner logarithmic correction as $\alpha+i\beta$ outside $B_5$.
For the finite histories below, one has
\begin{equation}\label{fixed-tube-background-profile-bounds}
 \alpha=O_m(\epsilon),\qquad
 \beta=O_m(\epsilon\langle r\rangle^{-1}).
\end{equation}
with the first bound for $\alpha$ including its constant phase.
Here and below the physical outer waves evaluated in the tube are
included in the coefficients, not in the functions to be inverted.
The regular part of their phase gradient is bounded by a logarithmic
power.  These estimates follow from the first-correction tails, the
bounds for each newly constructed profile below, and the outer-wave
bounds.  They are thus preserved by the finite construction.

For $r\geq5$ the increment $\psi=\phi/(iW)$ is governed by the exact
multiplicative identity
\begin{equation}\label{fixed-tube-multiplicative-identity}
 \begin{split}
 \frac{S_\epsilon(Ue^{i\psi})}{iUe^{i\psi}}
 ={}&\frac{S_\epsilon(U)}{iU}
 +ic_\epsilon D_t^x\psi+\epsilon^2\Delta_x\psi
 +2\epsilon^2\frac{\nabla U}{U}\cdot\nabla\psi
 -2i|U|^2\psi_2\\
 &+i\epsilon^2(\nabla\psi)^2
 +i|U|^2(e^{-2\psi_2}-1+2\psi_2).
 \end{split}
\end{equation}
The flat normal part is $\widetilde{\mathcal L}_W$.
The remaining linear terms are generated by
\begin{equation}\label{fixed-tube-operator-list}
 \begin{gathered}
 \epsilon C\,D_y,\qquad
 \epsilon^2 C\,\mathscr D_s^2,\qquad
 \epsilon^2 C\,\mathscr D_s,\qquad
 \epsilon^2 C\,yD_y,\qquad
 i\sqrt2\epsilon\partial_\tau,\qquad
 iC_0\psi_2,
 \end{gathered}
\end{equation}
where every differentiated $C$ has at most logarithmic growth and
$C_0=O_m(\epsilon\langle r\rangle^{-1})$.
The first-order terms containing the regular phase gradient exchange
phase and amplitude.  The geometric scalar first-order terms do not.
There is no zeroth-order term acting on the phase component.
Formula \eqref{laplaceperiodicframe} supplies the longitudinal terms;
the material chain rule supplies the velocity, rotation and $yD_y$
terms.  On $B_5$, the difference between the derivative of the actual
map \eqref{exact-global-inner-factor} and its linear part is a
differential operator with coefficients $O_m(\epsilon)$, in addition
to the time and longitudinal terms already displayed.

\begin{lemma}[Elliptic batch]
\label{fixed-tube-elliptic-batch-lemma}
Assume the coefficient bounds stated above at the required finite
orders.  The following two procedures are possible.

\textup{(a)} Suppose the data have bounded physical-time derivatives,
up to logarithmic factors, and the initial residual satisfies
\[
 R_1=O_m(\epsilon^2\langle r\rangle^{-2}),\qquad
 R_2=O_m(\epsilon^2),
\]
with translation moments $O(\epsilon^3\Lambda_\epsilon^C)$.
There are centred corrections $\phi_n$, $2\leq n\leq J$, such that
\begin{equation}\label{slow-elliptic-profile-bounds}
 \psi_{n,1},\psi_{n,2}
 =O_m(\epsilon^n\langle r\rangle^{n-2}),
 \qquad \phi_n=iW(\psi_{n,1}+i\psi_{n,2}),
\end{equation}
and the resulting residual obeys
\begin{equation}\label{slow-elliptic-residual-bounds}
 R_1=O_m(\epsilon^{J+1}\langle r\rangle^{J-3}),\qquad
 R_2=O_m(\epsilon^{J+1}\langle r\rangle^{J-1}).
\end{equation}

\textup{(b)} Let $A\leq\epsilon$ be a stage bound.  Suppose the initial
residual has the form
\begin{equation}\label{radiative-inner-input}
 \begin{aligned}
 R_1&=O_m(\epsilon^2 A\langle r\rangle^{-1}),\\
 R_2&=\sqrt2\epsilon c(\rho,\tau)(1-w^2)+O_m(\epsilon^2 A),
 \qquad c=O_m(A).
 \end{aligned}
\end{equation}
First invert the displayed radial term.  Thereafter there are centred
corrections $\phi_n$, $2\leq n\leq J$, with
\begin{equation}\label{iteratedinnerprofilebounds}
 \psi_{n,1}=O_m(\epsilon^n A\langle r\rangle^{n-1}),\qquad
 \psi_{n,2}=O_m(\epsilon^n A\langle r\rangle^{n-2}),
\end{equation}
and final residual
\begin{equation}\label{finite-inner-Taylor-remainder-assumption}
 R_1=O_m(\epsilon^{J+1}A\langle r\rangle^{J-2}),\qquad
 R_2=O_m(\epsilon^{J+1}A\langle r\rangle^{J-1}).
\end{equation}
The radial correction is $O_m(\epsilon A\langle r\rangle^{-2})$ in
amplitude and has zero phase.  All statements hold after two weighted
family variations.  An output order $m$ uses at most
\begin{equation}\label{elliptic-derivative-reserve}
 \mathfrak r_J(m):=8^{J+3}(m+20)
\end{equation}
material, normal and fast-time derivatives of the input fields and
coefficients.
\end{lemma}

\begin{proof}
In (a), split off the compact translation moments in the first source.
The orthogonal source has an inverse of size $\epsilon^2$ up to
logarithms, by Lemma \ref{polynomial-planar-inverse-lemma}.  The
remaining moment source has size $\epsilon^3$ and an inverse of size
$\epsilon^3\langle r\rangle$.  On $r\leq5d/\epsilon$ this satisfies
the same differentiated bound as $\epsilon^2$.  This constructs
$\phi_2$.
In (b), the radial amplitude inverse constructs the first correction.
Its time derivative, the geometric terms, and the terms containing an
old correction are respectively of sizes
$\epsilon^2A\langle r\rangle^{-2}$,
$\epsilon^2A\langle r\rangle^{-3}$ and
$O_m(\epsilon^2A)$ on bounded discs.  They fit the source for $n=2$.

For each subsequent substep use \eqref{complete-inner-coefficient-recursion}
and recompute the complete residual.  In (b), the input at level $n$ is
\[
 (R_1,R_2)
 =O_m\big(\epsilon^nA(
       \langle r\rangle^{n-3},\langle r\rangle^{n-2})\big).
\]
Lemma \ref{polynomial-planar-inverse-lemma}, with $p=n-1$, gives
\eqref{iteratedinnerprofilebounds}.  The following list gives the
largest output orders of the non-flat linear terms; logarithms are
suppressed:
\[
\begin{array}{c|cc}
 &\hbox{phase source}&\hbox{amplitude source}\\ \hline
 \epsilon D_y&\epsilon^{n+1}A\langle r\rangle^{n-2}
                   &\epsilon^{n+1}A\langle r\rangle^{n-2}\\
 i\epsilon\partial_\tau&\epsilon^{n+1}A\langle r\rangle^{n-2}
                   &\epsilon^{n+1}A\langle r\rangle^{n-1}\\
 \epsilon^2\mathscr D_s^2,\ \epsilon^2\mathscr D_s,
 \ \epsilon^2yD_y
       &\epsilon^{n+2}A\langle r\rangle^{n-1}
                   &\epsilon^{n+2}A\langle r\rangle^{n-1}\\
 iC_0\psi_2&0&\epsilon^{n+1}A\langle r\rangle^{n-3}
\end{array}
\]
Every entry belongs to level $n+1$, since
$\epsilon\langle r\rangle\leq C_d$.
The real quadratic derivative term is bounded by
$\epsilon^{2n}A^2\langle r\rangle^{2n-5}$; the imaginary derivative
and amplitude squares are bounded by
$\epsilon^{2n}A^2\langle r\rangle^{2n-4}$.
They also belong to level $n+1$.  Higher powers are smaller because the
complete logarithmic amplitude is uniformly small.  On $B_5$ the
analytic increment identity \eqref{exact-inner-composite-increment}
gives the same conclusion from the factor $\epsilon$ in its coefficient
difference and the quadratic size of its remainder.  This proves the
induction in (b), including the nonlinear terms.

In (a), the level-$n$ source has orders
$\epsilon^n(\langle r\rangle^{n-4},
\langle r\rangle^{n-2})$.
The inverse has order \eqref{slow-elliptic-profile-bounds}.
Here the material time operator has size $\epsilon^2\Lambda_\epsilon$;
its transport part is still $\epsilon\Lambda_\epsilon D_y$.
The phase output orders are
$\epsilon^{n+1}\langle r\rangle^{n-3}$ and
$\epsilon^{n+2}\langle r\rangle^{n-2}$.
The amplitude output is no larger than
$\epsilon^{n+1}\langle r\rangle^{n-1}$.
The nonlinear derivative product has phase order
$\epsilon^{2n}\langle r\rangle^{2n-6}$, and the amplitude square has
order $\epsilon^{2n}\langle r\rangle^{2n-4}$.
These give the next level.  Physical-time differentiation preserves this
count, so the slow-time assertion is retained.

We verify that the fixed cutoffs do not change either induction.
Where $\chi_d=1$, the source is cancelled by the planar equation.
Where $1-\chi_d$ is nonzero within the inversion tube,
$\epsilon\langle r\rangle\geq c_d$.  A surviving level-$n$ term
therefore belongs to level $n+1$, with a constant depending on $d$.
The differentiated cutoff terms have two normal derivatives in total,
or one material time derivative.  Their orders follow from
\eqref{matchingcutoffderivatives}.  For example, in (b), the phase
commutators have size
$\epsilon^n A\langle r\rangle^{n-3}$ on such an annulus and hence
are bounded by
$C_d\epsilon^{n+1}A\langle r\rangle^{n-2}$.
The products containing two cutoff derivatives have the same quadratic
orders just listed.  Outside the support of $\chi_d$ the residual is
unchanged, but the same annular comparison applies.

Each use of the exact operator consumes at most two spatial or one
fast-time derivative.  The reserve $4m+20$ in the planar lemma allows both the
angular summation and interior estimates.  Differentiating twice in
the parameter distributes derivatives among these same factors; the
weighted Leibniz rule preserves their $\epsilon$ powers.  At each substep an input index $8(m+10)$ dominates this
loss and the two derivatives of the operator.  Iterating that rule,
and allowing for two family variations, is bounded by
$\mathfrak r_J(m)$.  This also includes the core normalisations and
coordinate compositions.  This proves the differentiated conclusions.
\end{proof}

\subsection{Curve increments and preparation of the second approximation}

Only linear curve equations are needed for the additional increments.
Let $z_*(t)$ be a prescribed smooth path in the chart, with uniformly
bounded physical-time derivatives at the finite orders in use, and set
\[
 L_*g=\partial_tg-D_z\mathscr G_\epsilon(t,z_*(t))g.
\]
The spatial estimates \eqref{NashMoser-low-spatial-estimate} and
\eqref{NashMoser-high-spatial-estimate} give a solution of
\begin{equation}\label{fixed-base-curve-increment}
 L_*g=f,\qquad g(0)=0,
\end{equation}
with $\sup_t\|g(t)\|_{H^m}\leq C_m\sup_t\|f(t)\|_{H^{m+1}}$, after
including the fixed high norms of $z_*$ in the constant.
Repeatedly differentiating the equation, and using
\eqref{iterated-linearised-composition-bound}, with an odd spatial output
rounded up once to the next even grade, gives
\begin{equation}\label{fixed-base-fast-time-curve-bound}
 \sup_t\|g(t)\|_{H^m}
 +\sum_{j=0}^m\sup_t\|\partial_\tau^j\partial_tg(t)\|_{H^{m-j}}
 \leq C_m\sum_{j+a\leq2m+20}
        \sup_t\|\partial_\tau^j f(t)\|_{H^a}.
\end{equation}
Here each derivative of a coefficient of $L_*$ is accompanied by
$\epsilon\Lambda_\epsilon/\sqrt2$; the fixed physical interval is
$[0,T]$.  For the first two family variations one may replace the input
index by $2m+280$: \eqref{polarised-inverse-reserves} gives
$\mu_2=200$, and $2m_*=80$ covers the fixed low data and direction
indices.  The background and its family derivatives are measured
at these same finite enlarged levels.  Thus no smallness of high unweighted
physical-time derivatives of a wave trace is required.

We first choose a slowly varying reference curve for the higher
construction.  Let $z_1$ be the chart of $\gamma_\epsilon^{(1)}$ and
let $H_1$ be the residual of its inner factor before its outer wave is
inserted, with the normalisation \eqref{fixed-tube-normalised-residual}.
The first-source calculation and the componentwise tails of
$\phi_{1,\epsilon}$ give
\begin{equation}\label{slow-seed-source-bounds}
 (H_{1,1},H_{1,2})
 =O_m\big(\epsilon^2(
       \langle r\rangle^{-2},1)\big).
\end{equation}
To verify the first component, write the leading factor outside the
core as $w(r)e^{i\varphi}$, with
$\varphi=\theta+\Phi_\epsilon$.  Rotation of the normal frame does not
change $r$, so at fixed Cartesian $x$,
\[
 \partial_t r=-\epsilon^{-1}V^\perp\cdot e_r,
 \qquad \nabla_x\log w=\epsilon^{-1}(w'/w)e_r.
\]
Since $\Delta_x\varphi=0$ there, its real normalised residual is exactly
\[
 E_{1,0}=\epsilon\frac{w'}w e_r\cdot
                    (-\Lambda_\epsilon V^\perp+2\nabla_z\Phi_\epsilon).
\]
This is the real component of the source cancelled by the first planar
correction.  Its compact translation subtraction is retained in the
core.  For the remaining real terms, the first-correction tails give, up to
logarithms,
\[
\begin{array}{c|c}
\text{term}&\text{radial bound}\\ \hline
\epsilon D_y(\epsilon\psi_1,\epsilon\psi_2)
   &\epsilon^2 r^{-2}+\epsilon^3r^{-1}\\
\begin{gathered}
\epsilon^2\mathscr D_s^2(\epsilon\psi_1,\epsilon\psi_2),\\
\epsilon^2\partial_t^{\rm mat}(\epsilon\psi_1,\epsilon\psi_2)
\end{gathered}
   &\epsilon^3r^{-1}+\epsilon^4\\
\text{quadratic derivative terms}
   &\epsilon^2r^{-4}+\epsilon^3r^{-3}+\epsilon^4r^{-2}.
\end{array}
\]
Each is at most $C\epsilon^2r^{-2}\Lambda_\epsilon^C$ on the fixed
tube, since $\epsilon r\leq C_d$.  The imaginary component also
contains $\epsilon^2|\nabla\Phi_\epsilon|^2$ and the regular part of
$c_\epsilon\varphi_t$, bounded by $\epsilon^2\Lambda_\epsilon^C$.
The analytic core terms are smooth and $O_m(\epsilon^2)$.
This proves \eqref{slow-seed-source-bounds}.  Differentiating these
identities gives the same normal weights and the stated material,
physical-time and two family bounds; derivatives of the exact tubular
coefficients satisfy the bounds used in
\eqref{fixed-tube-operator-list}.
Let $Z_j=\partial_{y_j}W$ and put
\begin{equation}\label{completeprojectionmap}
 P_0=\left(\operatorname{Re}\int
       \chi_{\rm inv}H_1\overline{Z_j}\,dy\right)_{j=1,2},
 \qquad
 m_\epsilon=2\pi\int_0^\infty
                 \chi_{\rm inv}(\epsilon r)w(r)w'(r)\,dr.
\end{equation}
Then $m_\epsilon=\pi+O(\epsilon^2)$ and $P_0=O_m(\epsilon^2)$.
The latter integral converges uniformly: the phase contribution is
bounded by $C\epsilon^2\langle r\rangle^{-2}dr$ and the amplitude
contribution by the same integrable weight.  If $J(a_1,a_2)=(-a_2,a_1)$,
the leading change of this projection under an additional normal velocity
$A$ is $-\epsilon\Lambda_\epsilon m_\epsilon JA$.
Let $\mathscr L_1$ be the constant-speed chart lift of a normal velocity
along the first curve.  Solve
\begin{equation}\label{seed-curve-correction}
 L_1g_0=\mathscr L_1
       \left[(\epsilon\Lambda_\epsilon m_\epsilon)^{-1}J^{-1}P_0\right],
 \qquad g_0(0)=0,
\end{equation}
where $L_1=\partial_t-D_z\mathscr G_\epsilon(t,z_1)$, and set
$z_*=z_1+g_0$.  This curve is closed and of constant speed by the chart
recovery.  Its distance from the first curve is
$O_m(\epsilon\Lambda_\epsilon^C)$, with physical-time derivatives.

Recompute the first inner factor at this curve.  Its two source moments
are $O_m(\epsilon^3)$.
Indeed, along paths satisfying the basic velocity, the derivative of
\eqref{slow-seed-source-bounds} is still $O_m(\epsilon^2)$, with the
finite curve-derivative reserve from Section \ref{normalgraphsection}.
The change in actual velocity is the right-hand side of
\eqref{seed-curve-correction}, up to a quadratic remainder.  Its leading
projection cancels $P_0$.  The remaining velocity terms act on an
$O_m(\epsilon)$ correction or on a higher tubular coefficient, and their
moments are $O_m(\epsilon^3)$ by the same radial integral.  The change
of chart and frame and the quadratic remainder have this size as well.
This argument is unchanged after the stated derivatives.

Apply part (a) of Lemma \ref{fixed-tube-elliptic-batch-lemma} at $z_*$.
After choosing its depth, denote the resulting inner factor by $U_2^-$.
For the next subsection choose a fixed radial core cutoff
\[
 \chi_{\rm c}=1\quad(r\leq6),\qquad
 \chi_{\rm c}=0\quad(r\geq7).
\]
The two smooth real functions
\begin{equation}\label{whole-tube-outer-source-components}
 H_1+iH_2=(1-\chi_{\rm c})\frac{S_\epsilon(U_2^-)}{iU_2^-}
\end{equation}
are defined as zero near the centre.  The source contains the remaining
error throughout the tube outside the core cutoff.

\subsection{The phase wave and the algebraic amplitude}

For any pre-wave approximation $U$ define $H_1,H_2$ by
\eqref{whole-tube-outer-source-components}, with $U_2^-$ replaced by $U$,
and write
\[
 q=|U|^2,\qquad p=\operatorname{Re}(\nabla U/U),\qquad
 v=\operatorname{Im}(\nabla U/U).
\]
Only the values away from the core enter the source.  Set
\begin{equation}\label{higherouterwavesource}
 F=\frac1{\sqrt2\epsilon}\partial_\tau H_2
             +v\cdot\nabla_xH_2-\epsilon^{-2}H_1.
\end{equation}
The phase increment $a$ and the amplitude increment $b$ are defined by
\begin{equation}\label{higherouteralgebraiccorrection}
 \begin{gathered}
 -a_{\tau\tau}+\Delta_xa=F,\qquad a(0)=a_\tau(0)=0,\\
 b=\tfrac12(H_2+\sqrt2\epsilon a_\tau),\qquad
 U^+=Ue^{i(a+ib)}.
 \end{gathered}
\end{equation}
Thus $a$ is real and $-b$ is the increment of logarithmic amplitude.

Applying Lemma \ref{exactfreewavecancellationlemma} below with the
source \eqref{whole-tube-outer-source-components} gives
\begin{equation}\label{higherouterexactremainders}
 S_\epsilon(U^+)
 =e^{i(a+ib)}\chi_{\rm c}S_\epsilon(U)+iU^+(T_1+iT_2),
\end{equation}
where $T_1,T_2$ are defined in
\eqref{exactfreewaverealremainder}--\eqref{exactfreewaveimagremainder}.
The combined right-hand side extends smoothly through the zero of $U$.

We use the wave norm $X_m$ of \eqref{threeDwaveXnorm} and the source
norms $Y_m,\mathcal H_m$ of Section \ref{threeDwaveestimatesection}.
With the fixed dipole functions $f_j$ defined in Subsection
\ref{dipole-wave-estimates-subsection}, write
\begin{equation}\label{nextwavesourcesplitting}
 F=F_{\rm loc}+\sum_{j=1}^3c_j(\tau)f_j+F_{\rm g}.
\end{equation}
For families we use \eqref{differentiatedouternorms} and the dipole
coefficient norm defined in that subsection.  The estimates below
identify and bound these three source contributions; their wave
solutions satisfy \eqref{fixed-tube-wave-estimate}.

\begin{lemma}[Transfer of an elliptic batch to the outer source]
\label{innertooutertransferlemma}
Let the required output derivative order be $m$ and take $J\geq m+8$.
For a slow batch as in part (a) of Lemma
\ref{fixed-tube-elliptic-batch-lemma}, the source in
\eqref{higherouterwavesource} satisfies
\begin{equation}\label{slow-seed-wave-source-bound}
 \|F_{\rm loc}\|_{Y_m}^{(2)}+
       \sum_j\|c_j\|_{C_\tau^m}^{(2)}
 \leq C\epsilon^2\Lambda_\epsilon^C,
 \qquad
 \|H_2\|_{\mathcal H_m}^{(2)}\leq C\epsilon^2\Lambda_\epsilon^C.
\end{equation}
There is no general global term in this initial source.
For a radiative batch of size $A$, its contribution supported in the
fixed tube obeys
\begin{equation}\label{fixed-tube-source-transfer-bound}
 \|F_{\rm tube}\|_{Y_m}^{(2)}
 \leq C_{m,J,d}\epsilon A\Lambda_\epsilon^C,
 \qquad
 \|H_{2,\rm tube}\|_{\mathcal H_m}^{(2)}
 \leq C_{m,J,d}\epsilon^2 A\Lambda_\epsilon^C.
\end{equation}
\end{lemma}

\begin{proof}
We give the derivative count in the physical variables.  The largest
normal-derivative contribution to $\epsilon^{-2}H_1$ in a radiative
batch is
\[
 \epsilon^{-2-a}\epsilon^{J+1}A
          \langle r\rangle^{J-2-a}
 =\epsilon A(\epsilon\langle r\rangle)^{J-2-a}.
\]
The material-time part of
$\epsilon^{-1}\partial_\tau H_2$ is bounded by
\[
 \epsilon^{-1-a}\epsilon^{J+1}A
          \langle r\rangle^{J-1-a}
 =\epsilon A(\epsilon\langle r\rangle)^{J-1-a}.
\]
The material-to-Cartesian time conversion adds a normal derivative with
a logarithmically bounded coefficient and a longitudinal derivative
with a factor $\epsilon\Lambda_\epsilon$.  These are smaller than the
displayed terms.  The singular part of $v$ is
$\epsilon^{-1}r^{-1}e_\theta$; its product with $\nabla H_2$ has bound
$C\epsilon^2A(\epsilon\langle r\rangle)^{J-3-a}$.
The regular part of $v$ gives no larger bound.

On the core-cutoff annulus $r\simeq1$, derivatives of its physical
argument cost powers of $\epsilon^{-1}$, but the numerator retains
$\epsilon^{J+1}$.  On the physical cutoff annuli all such derivatives
are bounded by constants depending on $d$.  Longitudinal derivatives
consume only their reserved material orders.  Two weighted curve
variations have the same count: in a tubular coordinate composition
$\epsilon D_{\rm in}$ is a bounded combination of $D_y$, a material
variation, and $\epsilon\partial_\rho$.  The additional frame terms
have coefficients bounded by $C\epsilon\langle r\rangle$.
Thus $J\geq m+8$ pays for every indicated derivative.
The tube has uniformly bounded physical volume, which proves both
norm bounds in \eqref{fixed-tube-source-transfer-bound}.

For the slow batch the real radial order is one lower and material time
differentiation gains $\epsilon\Lambda_\epsilon$.  The same calculation
therefore gives $\epsilon^2$ instead of $\epsilon A$.  On the fixed annuli between the inversion tube and the original
$\delta$ tube, the first-correction tails give $E_1=O_m(\epsilon^4)$
and $E_2=O_m(\epsilon^2)$ in physical derivatives.  Their slow-time
source has the same $O_m(\epsilon^2)$ bound and has bounded support.
Beyond the original tube the field has unit leading modulus and solid-angle phase.  The
far-field expansion of the phase is a dipole plus a remainder with
$|x|^{-3}$ decay.  Its time derivatives are bounded on the slow curve.
Hence the only $|x|^{-2}$ contribution to $F$ is the finite dipole term;
all other terms belong to $Y_m$.  The imaginary residual itself has
$|x|^{-2}$ decay and is in every required $H^m$.  This proves
\eqref{slow-seed-wave-source-bound}.
\end{proof}

Applying Lemma \ref{innertooutertransferlemma} and the wave estimate
\eqref{fixed-tube-wave-estimate} to $U_2^-$ defines $a_2,b_2$ and
$u_\epsilon^{(2)}=U_2^-e^{ia_2-b_2}$, with
\begin{equation}\label{prepared-second-approximation-bounds}
 \|a_2\|_{X_m}^{(2)}\leq C_m\epsilon^{3/2}\Lambda_\epsilon^{C_m},
 \qquad
 \|b_2\|_{\mathcal H_m}^{(2)}
       +\|b_2\|_{W^{m,\infty}_{x,\tau}}^{(2)}
 \leq C_m\epsilon^2\Lambda_\epsilon^{C_m}.
\end{equation}
The depth of the slow batch is chosen for the highest derivative level
needed subsequently.  This preparation uses the inner factor of
Section \ref{firstglobalimprovementsection}; the approximation
$u_\epsilon^{(1)}$ remains the reference in the theorem, while the
outer sums for stages $k\geq2$ begin with $a_2,b_2$.

\subsection{A complete improvement cycle}

Fix the slow path $z_*$ constructed above and write
\[
 r_*:=\partial_tz_*-\mathscr G_\epsilon(t,z_*),\qquad
 L_*:=\partial_t-D_z\mathscr G_\epsilon(t,z_*).
\]
For a real physical phase $A$, let $\mathscr V(z,A)$ be the reduced
constant-speed lift of the normal velocity
\[
 \frac{2}{\Lambda_\epsilon}
       \Pi_{t_\Gamma}^{\perp}\nabla_x A(\Gamma^{\rm mat}(\rho,t),t),
       \qquad \Gamma^{\rm mat}=\mathfrak C_t(z).
\]
The mean tangential velocity in this lift is zero.  Thus
$\mathscr V(z,A)$ is linear in $A$ and is a smooth function of the curve,
its frame and the trace of $\nabla A$.  The lift is the one specified by
\eqref{constant-speed-normal-lift}.

Suppose a post-wave field $u_\epsilon^{(k)}$ has been constructed, with
chart $z_k$, material inner profile $\Pi_k$ and cumulative phases
\[
 A^{(k)}=\sum_{j=2}^k a_j,\qquad B^{(k)}=\sum_{j=2}^k b_j.
\]
Define its modulation defect by
\begin{equation}\label{fixed-base-modulation-defect}
 D_k:=\partial_tz_k-\mathscr G_\epsilon(t,z_k)
                       -r_*-\mathscr V(z_k,A^{(k)}).
\end{equation}
Solve the prescribed linear equation
\begin{equation}\label{higher-translation-Newton-velocity}
 L_*g_k=-D_k,\qquad g_k(0)=0,\qquad z_{k+1}=z_k+g_k.
\end{equation}
Insert $\Pi_k$ in the new material tube, recompute
$\phi_{1,\epsilon}[\Gamma_{k+1}]$ and the solid-angle phase, and keep
$A^{(k)},B^{(k)}$ fixed as physical functions.  Apply the radiative
elliptic batch to this field.  The resulting field is denoted by
$U_{k+1}^-$.  Equations \eqref{higherouterwavesource} and
\eqref{higherouteralgebraiccorrection} then define $a_{k+1},b_{k+1}$ and
\[
 u_\epsilon^{(k+1)}=U_{k+1}^-e^{ia_{k+1}-b_{k+1}}.
\]
Every source and every curve increment is therefore prescribed before
its linear equation is solved.

For the induction, the incoming core remainder is retained with its
scaled estimate
\[
 \|\chi_{\rm c}\mathcal R_k^-\|_{C^r_{\rho,y,\tau}}^{(2)}
 \leq C_r\epsilon^{J_k+1}A_{k-1}\Lambda_\epsilon^{C_r};
\]
at the second stage $A_{k-1}$ is replaced by one.  The imaginary source
of the preceding wave has the full radial estimate
\eqref{slow-elliptic-residual-bounds} or
\eqref{finite-inner-Taylor-remainder-assumption} on the tube and
$\|H_{2,k}\|_{\mathcal H_r}^{(2)}\leq
C\epsilon^2A_{k-1}\Lambda_\epsilon^C$ globally.  These are part of the
induction data, together with the exact post-wave identity.  

\begin{lemma}[One improvement cycle]
\label{fixed-tube-cycle-lemma}
Let the cumulative low norms satisfy
\[
 \|z_k-z_*\|+\|A^{(k)}\|
 \leq C\epsilon^{5/4}\Lambda_\epsilon^C,
 \qquad \|B^{(k)}\|\leq C\epsilon^2\Lambda_\epsilon^C,
\]
where the phase norm controls its physical derivatives and its value in
$L^\infty$, and let the more precise amplitude and fixed-tube bounds
from the preceding construction hold.  Suppose the newest phase has
size $A_k\leq\epsilon^{5/4}$ at the required higher derivative level and its
amplitude has size $C\epsilon A_k\Lambda_\epsilon^C$.  At $k=2$ the
amplitude bound is instead $C\epsilon^2\Lambda_\epsilon^C$, as in
\eqref{prepared-second-approximation-bounds}.
Suppose also that $A_{k-1}\leq\epsilon^{-1/4}A_k$ for $k\geq3$
and that the incoming core and source estimates just stated hold.
All elliptic depths are at least eight.
Then the above cycle, with elliptic depth chosen for its output index,
has the following estimates:
\begin{equation}\label{nextwavesourcebounds}
 \begin{aligned}
 \|g_k\|_{C_tH^m}
   +\sum_{j=0}^m\|\partial_\tau^j\partial_tg_k\|_{C_tH^{m-j}}
     &\leq C_m A_k\Lambda_\epsilon^{C_m},\\
 \|F_{k+1}^{\rm loc}\|_{Y_m}^{(2)}
   +\sum_{\ell=1}^3\|c_{k+1,\ell}\|_{C_\tau^m}^{(2)}
     &\leq C_m\epsilon A_k\Lambda_\epsilon^{C_m},\\
 \|F_{k+1}^{\rm g}\|_{\mathcal H_m}^{(2)}
   +\|H_{2,k+1}\|_{\mathcal H_m}^{(2)}
     &\leq C_m\epsilon^2 A_k\Lambda_\epsilon^{C_m},\\
 \|a_{k+1}\|_{X_{m-4}}^{(2)}
     &\leq C_m\epsilon^{1/2} A_k\Lambda_\epsilon^{C_m},\\
 \|b_{k+1}\|_{\mathcal H_{m-6}}^{(2)}
   +\|b_{k+1}\|_{W^{m-6,\infty}_{x,\tau}}^{(2)}
     &\leq C_m\epsilon^{3/2} A_k\Lambda_\epsilon^{C_m}.
 \end{aligned}
\end{equation}
The source decomposition is
\begin{equation}\label{fixed-tube-cycle-source-decomposition}
 F_{k+1}=F_{k+1}^{\rm loc}
       +\sum_{\ell=1}^3 c_{k+1,\ell}(\tau)f_\ell
       +F_{k+1}^{\rm g},
\end{equation}
with the fixed dipole functions of Subsection
\ref{dipole-wave-estimates-subsection}.  All bounds have the two weighted
family versions.  Only finitely many higher input derivatives are used.
\end{lemma}

\begin{proof}
We separate the calculation into modulation, the tubular residual,
and the far field.
For the product terms we use the differentiated estimates
\begin{equation}\label{outerclassproductestimates}
 \begin{split}
 \|fg\|_{Y_m}^{(2)}
  &\leq C_m\|f\|_{Y_m}^{(2)}\|g\|_{W^{m,\infty}_{x,\tau}}^{(2)},\\
 \|fg\|_{\mathcal H_m}^{(2)}
  &\leq C_m\big(
       \|f\|_{\mathcal H_m}^{(2)}\|g\|_{W^{m,\infty}_{x,\tau}}^{(2)}
      +\|g\|_{\mathcal H_m}^{(2)}\|f\|_{W^{m,\infty}_{x,\tau}}^{(2)}\big).
 \end{split}
\end{equation}
They follow by Leibniz' rule and the Sobolev product estimate; the two
weighted parameter derivatives obey the same product rule.

\smallskip
\noindent\emph{Modulation.}
Write $D_k^-$ for the modulation defect immediately before $a_k$ was
inserted.  By linearity in the phase,
\[
 D_k=D_k^- -\mathscr V(z_k,a_k).
\]
At $k=2$, $D_2^-=0$.  The estimate from the preceding cycle is
$D_k^-=O(\epsilon A_{k-1}\Lambda_\epsilon^C)$, which is smaller than
the permitted next stage size.  The trace and the constant-speed lift,
followed by \eqref{fixed-base-fast-time-curve-bound}, give the first
line of \eqref{nextwavesourcebounds}.  Before the new wave is
inserted, the new defect has the exact expression
\begin{equation}\label{exact-post-Newton-remainder}
 \begin{split}
 D_{k+1}^-={}&
 -\int_0^1\{D_z\mathscr G_\epsilon(t,z_k+s g_k)
                    -D_z\mathscr G_\epsilon(t,z_*)\}g_k\,ds\\
 &-\int_0^1D_z\mathscr V(z_k+s g_k,A^{(k)})g_k\,ds.
 \end{split}
\end{equation}
The first integrand contains $z_k-z_*$ or $g_k$ in addition to $g_k$;
the second contains a derivative of the accumulated phase.  The tame
bounds of Section \ref{normalgraphsection} and the physical trace
bounds therefore yield
\begin{equation}\label{post-Newton-translation-remainder-estimate}
 D_{k+1}^-=O_m(\epsilon A_k\Lambda_\epsilon^C).
\end{equation}
The same formula applies after the two weighted parameter variations.

\smallskip
\noindent\emph{The tubular residual.}
The exact post-wave identities
\eqref{higherouterexactremainders}
are used before any component is estimated.  On $r\geq5$, put
\[
 p_0=\epsilon^{-1}\frac{w'}w e_r,
 \qquad v_0^{\rm sing}=\epsilon^{-1}r^{-1}e_\theta.
\]
The terms involving the trace of the newest phase gradient are
\[
 2\epsilon\frac{w'}w e_r\cdot\nabla a_k(\Gamma_k,t),
 \qquad
 \frac{2\epsilon}{r}e_\theta\cdot\nabla a_k(\Gamma_k,t).
\]
A normal velocity increment contributes the same pair with
$2\nabla a_k$ replaced by $-\Lambda_\epsilon\delta V^\perp$.
Equation \eqref{higher-translation-Newton-velocity} and
\eqref{post-Newton-translation-remainder-estimate} cancel these leading
terms, with an error $O_m(\epsilon^2 A_k)$ in the velocity source.
Their radial weights are $r^{-3}$ and $r^{-1}$, respectively.

For the remaining terms, Taylor's formula in the physical normal
variable gives
\[
 D_{x,\tau}^{\alpha}a_k(\mathcal X_{\Gamma_k})
 -D_{x,\tau}^{\alpha}a_k(\Gamma_k)
 =\epsilon y_\ell\int_0^1
   \partial_{x_\ell}D_{x,\tau}^{\alpha}a_k
       (\Gamma_k+s\epsilon y)\,ds.
\]
For example, the remaining $p_0\cdot\nabla a_k$ term has size
$\epsilon^2 A_k r^{-2}$, and the remaining
$v_0^{\rm sing}\cdot\nabla a_k$ term has size
$\epsilon^2 A_k$.  The singular transport of $a_{k,\tau}$ has size
$\epsilon^2 A_k/r$ in the real component.  The leading amplitude term is
\[
 \sqrt2\epsilon(1-w^2)a_{k,\tau}(\Gamma_k,t).
\]
Subtracting this radial term leaves size
$\epsilon^2 A_k/r$; it is estimated using the same Taylor formula.

Here is the remaining component count.  Each row includes the indicated
physical derivatives before passage to material coordinates; normal
weights are written only for $r\geq5$ and logarithmic factors are
suppressed.
\[
\begin{array}{c|cc}
\text{contribution}&\text{real}&\text{imaginary}\\ \hline
 \begin{gathered}(p-p_0)\nabla a_k,\\
 (v-v_0^{\rm sing})\nabla a_k\end{gathered}
     &\epsilon^2 A_k r^{-2}&\epsilon^2 A_k\\
 (|U|^2-w^2)\epsilon a_{k,\tau}
     &0&\epsilon^2 A_k r^{-1}\\
 \begin{gathered}\epsilon^2\Delta(\epsilon a_{k,\tau}),\\
 \epsilon^2p\cdot\nabla(\epsilon a_{k,\tau})\end{gathered}
     &0&\epsilon^2 A_k\\
 \text{quadratic terms}
     &\epsilon^2 A_k r^{-1}&\epsilon^2 A_k\\
 \text{change of curve}
     &\epsilon^2 A_k r^{-1}&\epsilon^2 A_k
\end{array}
\]
The first two rows use
$p-p_0=O_m(r^{-2})$,
$v-v_0^{\rm sing}=O_m(\Lambda_\epsilon^C)$ and
$|U|^2-w^2=O_m(\epsilon/r)$.
In the third row the $p_0$ part has the stronger bound
$\epsilon^2 A_k/r^3$; the other parts have an additional $\epsilon$.
For the product row, one high derivative is placed on one factor and the
other factors use their established low norms.  The nonlinear amplitude
remainder starts quadratically in $b_k$.

We explain separately the algebraic source part $H_{2,k}/2$ of $b_k$.
It is zero near the preceding core.  The preceding elliptic batch gives
thereafter, with all needed derivatives,
\[
 |D_y^jH_{2,k}|\leq
 C A_{k-1}\epsilon^{J_k+1}
       \langle r\rangle^{J_k-1-j}\Lambda_\epsilon^C.
\]
At $k=2$ the same bound has $A_{k-1}=1$ and is the slow estimate.
Consequently multiplication by $1-w^2$, or application of
$\epsilon^2\Delta_x$, has physical size at most
$C\epsilon^4 A_{k-1}\Lambda_\epsilon^C$ on the fixed tube and a much
higher order on every bounded inner disc.  The terms involving $p_0$
have still better weights.  These contributions belong to the real and
imaginary bounds in the table, after the derivative reserve is used.
We retain the radial weight when differentiating the core cutoff.
The shift of the core in a later radiative cycle is $o(\epsilon)$, so the
same conclusions hold in the new tube.  Fixed cutoff annuli are included
by Taylor's formula and the cutoff bounds.

For the last row of the table, write the two tubular residuals in the
same material variables.  Here is a componentwise verification of this
comparison.  The sum of all higher inner profiles satisfies
\[
 D_y^j(\alpha_{\rm hi},\beta_{\rm hi})
       =O_m(\epsilon^2\langle r\rangle^{-j}).
\]
For a slow profile this is
$\epsilon^n\langle r\rangle^{n-2-j}
 \leq C_d\epsilon^2\langle r\rangle^{-j}$;
for a radiative profile it follows from
$\sum_i A_i\leq C\epsilon^{5/4}\Lambda_\epsilon^C$ and
$\langle r\rangle\leq C_d/\epsilon$.
Let $f_z$ denote the actual normal velocity minus the basic corrected
normal velocity, written in the chosen frame.  The amplitude--phase
identity gives the material-rebuild difference in the form
\[
\begin{aligned}
 \delta R_1={}&\epsilon\frac{w'}w e_r\cdot
       \big(-\Lambda_\epsilon\delta f_z
                         +2\delta[\nabla A^{(k)}]\big)
       +O_m(\epsilon^2\Delta_g\langle r\rangle^{-1}),\\
 \delta R_2={}&\frac{\epsilon}{r} e_\theta\cdot
       \big(-\Lambda_\epsilon\delta f_z
                         +2\delta[\nabla A^{(k)}]\big)
       +O_m(\epsilon^2\Delta_g),
\end{aligned}
\]
where $\Delta_g$ is the required finite norm of $g,\partial_tg$.
To check the remainders, use
$E_1=c_\epsilon\partial_t\log|U|
 +2\epsilon^2\nabla\log|U|\cdot\nabla\arg U
 +\epsilon^2\Delta\arg U$.
The solid-angle Laplacian is zero, and the first correction cancels
the basic order-$\epsilon$ source.  Differences acting on higher
profiles are bounded by $\epsilon\Delta_g D_y(\alpha_{\rm hi},
\beta_{\rm hi})$ or $\epsilon^2\Delta_g$ times their longitudinal
or rotational derivatives.  They have real size at most
$C\epsilon^3\Delta_g/r+C\epsilon^4\Delta_g$.
Differences of old outer coefficients satisfy
$\delta D_{x,\tau}^{\alpha}A^{(k)}=O(A^{(k)}_{\rm low}\Delta_g)$
and $\delta D_{x,\tau}^{\alpha}B^{(k)}=O(\epsilon^2\Delta_g)$;
when the time derivative of the latter occurs in $E_1$, it has the
additional factor $\epsilon$.  Since
$A^{(k)}_{\rm low}\leq C\epsilon^{5/4}\Lambda_\epsilon^C$,
these contributions have size at most $C\epsilon^3\Delta_g$.
Both bounds fit $\epsilon^2\Delta_g/r$ on the fixed tube.
For $E_2$, the massive term contributes $O(\epsilon^2\Delta_g)$,
and the derivative and phase terms are no larger.  The basic residual
itself has differentiated orders
$\epsilon^2(r^{-2},1)$, so its variation obeys the same estimates.
This proves the displayed rebuild comparison.

Here the change in evaluating an old physical wave is estimated by
integrating its physical gradient against the curve and frame variations.
The fixed forcing $r_*$ changes its normal identification by $O(g_k)$ and
contributes $\epsilon\Lambda_\epsilon O(r_*g_k)$, which has the same size.

On a bounded inner disc these calculations are made for the smooth
complex residual.  The leading two terms combine as
$2i\epsilon\nabla a_k(\Gamma_k)\cdot\nabla W$, and cancel against the
velocity term before division by $W$.  The radial amplitude contribution
is $iW$ times the displayed imaginary component.  The other terms have
smooth $C^m$ bounds of order $\epsilon^2 A_k\Lambda_\epsilon^C$.
The retained core remainder is compact in $y$.  Since $J_k\geq8$
and $A_{k-1}\leq\epsilon^{-1/4}A_k$, it satisfies the same input
bound, at every reserved scaled derivative index; the slow remainder
at $k=2$ does so as well.  Thus the input is exactly
\eqref{radiative-inner-input}.  Part (b) of the elliptic-batch lemma applies.
Lemma \ref{innertooutertransferlemma} proves the tubular part of
the source estimates and the $\epsilon^2 A_k$ bound for $H_2$.

\smallskip
\noindent\emph{The far field.}
On the fixed annuli between $4d$ and the outer boundary of the
original tube, the first-profile tails give $q-1,p=O_m(\epsilon^2)$
in physical derivatives, while $v$ and its derivatives are bounded
by $C_{m,d}\Lambda_\epsilon^C$.  The same residual formula therefore
gives a compactly supported source of size $C\epsilon A_k
\Lambda_\epsilon^C$.  Changes of the curve are estimated by the
physical Taylor formula there; the phase-time term gives
$\partial_{\tau\tau}\delta\varphi=O_m(\epsilon A_k)$, and
the amplitude and spatial terms are smaller.

Beyond the original tube, the base has unit modulus.  There
\begin{equation}\label{smalloutercoefficientsarbitrary}
 \begin{aligned}
 q-1,\ p&=O_m(\epsilon^2\Lambda_\epsilon^C),\\
 v&=v_\Gamma+\widetilde v,\qquad
 \widetilde v=O_m(\epsilon^{5/4}\Lambda_\epsilon^C).
 \end{aligned}
\end{equation}
Here $v_\Gamma=\nabla\varphi_\Gamma$ has $|x|^{-3}$ decay with its
spatial derivatives.  More precisely, on the exterior and with weighted
family variations,
\begin{equation}\label{differentiatedfarfieldphaseclass}
 \begin{aligned}
 |D_x^\alpha v_\Gamma|&\leq C_\alpha\Lambda_\epsilon^C
                        \langle x\rangle^{-3-|\alpha|},\\
 |\partial_\tau^jD_x^\alpha v_\Gamma|
   &\leq C_{j,\alpha}\epsilon\Lambda_\epsilon^C
                        \langle x\rangle^{-3-|\alpha|},\quad j\geq1.
 \end{aligned}
\end{equation}
These follow by differentiating the solid-angle kernel and using the
bounds for $\partial_\tau^j\partial_tg_i$.  Inserting the complete pair
of residuals
\eqref{higherouterexactremainders} into
\eqref{higherouterwavesource} produces the following types:
\[
\begin{array}{c|c|c}
\text{type}&\text{norm}&\text{bound}\\ \hline
 \epsilon v_\Gamma\nabla a_{k,\tau},\\
 \epsilon(\partial_\tau v_\Gamma)\nabla a_k
       &Y_m& C\epsilon A_k\Lambda_\epsilon^C\\
 \epsilon^2\Delta a_{k,\tau\tau}
       &\mathcal H_m&C\epsilon^2 A_k\Lambda_\epsilon^C\\
 \epsilon^{-1}(q-1)b_{k,\tau},\\
 p\nabla a_k,\ \epsilon\widetilde v\nabla a_{k,\tau}
       &\mathcal H_m&C\epsilon^2 A_k\Lambda_\epsilon^C\\
 \text{terms with at least two new correction factors}
       &\mathcal H_m&C\epsilon A_k^2\Lambda_\epsilon^C
\end{array}
\]
The two terms containing the unperturbed phase gradient and one fast
time derivative of the newest phase are, precisely,
\begin{equation}\label{leadinglocalisedwavesource}
 2\sqrt2\epsilon v_\Gamma\cdot\nabla a_{k,\tau}
 +\sqrt2\epsilon(\partial_\tau v_\Gamma)\cdot\nabla a_k.
\end{equation}
Terms with derivatives on the other displayed coefficients have the same
bounds.  Products with $v_\Gamma$ are placed in $Y_m$, and all purely
radiative products use one high Sobolev norm and bounded low derivatives.
At $k=2$ the pre-wave field has unit modulus outside the fixed tube,
so its $q-1$ and $p$ contributions vanish there.  The remaining
algebraic $H_{2,2}$ terms contain spatial derivatives of that source
and have $|x|^{-3}$ or better decay.  Its slow time derivatives give
$\epsilon\Delta\partial_\tau H_{2,2}=O(\epsilon^4
\Lambda_\epsilon^C)$ in the required norm.  The terms supported in
fixed tubular collars have the same or smaller order and belong to
$Y_m$.  Thus the exceptional first amplitude is covered by the stated
local bound.  In subsequent cycles
the old residual is subtracted as a whole, so only its newest increment
is present.

Changing the curve introduces the phase difference
$\delta\varphi=\varphi_{\Gamma_{k+1}}-\varphi_{\Gamma_k}$.
Its exact contribution, relative to the same old outer factor, is
\[
 \delta E_1=2\epsilon^2p\cdot\nabla\delta\varphi,\qquad
 \delta E_2=\sqrt2\epsilon\partial_\tau\delta\varphi
       +2\epsilon^2v\cdot\nabla\delta\varphi
       +\epsilon^2|\nabla\delta\varphi|^2.
\]
Here $\Delta\delta\varphi=0$ outside the tube.  The solid-angle
expansion is a dipole whose coefficient is a constant multiple of the
vector area, followed by a $|x|^{-3}$ remainder.  By the curve increment
estimate,
\[
 \partial_\tau^j\delta\varphi=O_m(\epsilon A_k|x|^{-2}
                              \Lambda_\epsilon^C),\qquad j\geq1,
\]
at the finite orders used here.  The $|x|^{-2}$ part of
$\partial_{\tau\tau}\delta\varphi$ is assigned to the three fixed
dipole functions; its coefficients are $O_m(\epsilon A_k\Lambda_\epsilon^C)$.
The remainder is in $Y_m$.  Every term containing a spatial derivative
of $\delta\varphi$ has $|x|^{-3}$ decay, and the remaining products
belong to the global Sobolev class with the stated smaller bound.
No conservation of the area of the adjusted curve is used.

Combining the three regions proves the source bounds.  The same
calculation before applying the source operator gives
$\|H_{2,k+1}\|_{\mathcal H_m}^{(2)}\leq
C\epsilon^2A_k\Lambda_\epsilon^C$: in the tube this is the
transfer lemma, and outside it every term of the imaginary residual
is bounded either by $\epsilon^2v_\Gamma\cdot\nabla a_k$,
by the $\epsilon^2A_k$ dipole change, or by a smaller product.  The wave
estimate \eqref{fixed-tube-wave-estimate} and the algebraic formula for
$b_{k+1}$ give the last two lines of
\eqref{nextwavesourcebounds}.  All calculations are finite
Leibniz expansions.  Weighted family differentiation acts on the actual
cutoffs, normalisations, curve lift and retained histories; the same
radial derivative count and two additional input variations apply.
This proves the lemma.
\end{proof}

\subsection{Completion of the construction}

We now choose all the finite orders.  Let
\begin{equation}\label{explicit-number-of-stages}
 K=4(M+N+20)+2,\qquad m_K=N+20.
\end{equation}
For $k=K,K-1,\ldots,3$, define successively
\begin{equation}\label{explicit-stage-orders}
 J_k=m_k+20,\qquad
 m_{k-1}=8^{J_k+6}(m_k+40).
\end{equation}
Set $J_2=m_2+20$ and choose
\[
 s_*=8^{J_2+10}(m_2+40).
\]
Increase $s_*$ by the fixed low indices in Section
\ref{normalgraphsection} and the finite initial-data and family
reserves in Corollary \ref{forcedNashMoserfamilycorollary}, if needed.
These choices are finite.
The small initial neighbourhood is taken in the constant-speed chart,
with bounded $H^{s_*}$ norm and with its low norm small enough for
Section \ref{normalgraphsection}.

The input index for a cycle with output phase norm $X_{m_k}$ can be
bounded by $8^{J_k+5}(m_k+40)$.  Indeed, the source-to-wave estimate
uses four additional derivatives, the two family variations and
physical coordinate compositions use a further fixed reserve, each
elliptic inverse uses $4m+20$ at level $m$, and a residual computation
uses at most two further derivatives.  The recursion
\eqref{elliptic-derivative-reserve} bounds this finite sequence of
losses.  Including its two family variations, the fixed-base curve estimate
uses at most twice the relevant spatial order plus $280$.  This is
bounded by $8(m+40)$ at output order $m$, and is covered by the
displayed cycle reserve.  The products require only one factor at
the highest level.  Thus \eqref{explicit-stage-orders} leaves more
than the entire required input index at the preceding stage.  The
choice of $s_*$ supplies the slow seed and all derivatives of the fixed
reference operator at every index used later.

Put
\begin{equation}\label{quantitativeinductionassertions}
 A_k=\epsilon^{(k+3)/4},\qquad 2\leq k\leq K.
\end{equation}
We prove, at stage $k$, the phase and amplitude bounds
\begin{equation}\label{outerinductivebounds}
 \begin{aligned}
 \|a_k\|_{X_{m_k}}^{(2)}&\leq A_k,\\
 \|b_k\|_{\mathcal H_{m_k-2}}^{(2)}
       +\|b_k\|_{W^{m_k-2,\infty}_{x,\tau}}^{(2)}
       &\leq C_k\epsilon A_k\Lambda_\epsilon^{C_k}
                       &&(k\geq3),\\
 \|b_2\|_{\mathcal H_{m_2-2}}^{(2)}
       +\|b_2\|_{W^{m_2-2,\infty}_{x,\tau}}^{(2)}
       &\leq C\epsilon^2\Lambda_\epsilon^C.
 \end{aligned}
\end{equation}
The stage also retains the source bounds and the scaled core remainder
from the elliptic batch that produced it, as specified before Lemma
\ref{fixed-tube-cycle-lemma}.

The seed estimate \eqref{prepared-second-approximation-bounds} gives
\eqref{outerinductivebounds} for $k=2$, since
$\epsilon^{3/2}\Lambda_\epsilon^C\leq\epsilon^{5/4}$ after decreasing
$\epsilon_0$.  Suppose the assertions hold through $k$.  The one-cycle
lemma gives
\[
 \|a_{k+1}\|_{X_{m_{k+1}}}^{(2)}
 \leq C_k\epsilon^{1/2}A_k\Lambda_\epsilon^{C_k}
 \leq\epsilon^{1/4}A_k=A_{k+1}.
\]
Its amplitude estimate gives the corresponding next line of
\eqref{outerinductivebounds}, and the elliptic batch supplies the next
radial and core bounds.  At each step, decreasing $\epsilon_0$ absorbs
the finitely many constants and logarithmic factors.  There are only
$K-2$ such decreases.

The low coefficient bounds used in the lemmas close at the same time.
On the fixed tube, each slow profile has size
\[
 \epsilon^n\langle r\rangle^{n-2}
 \leq C_{n,d}\epsilon^2,
\]
and each radiative phase profile has size
\[
 \epsilon^n A_k\langle r\rangle^{n-1}
 \leq C_{n,d}\epsilon A_k.
\]
The corresponding amplitude is bounded by
$C_{n,d}\epsilon^2A_k$, apart from the explicitly retained radial term
$C\epsilon A_kr^{-2}$.  The differentiated estimates give exactly
\eqref{fixed-tube-background-profile-bounds}.  The new global amplitudes
have total size $O(\epsilon^2\Lambda_\epsilon^C)$, while the global
phase derivatives and their values have total size
$O(\epsilon^{5/4}\Lambda_\epsilon^C)$.
The radiative curve increments have the same latter bound by
\eqref{nextwavesourcebounds}.  All these bounds hold at the low
indices and their required two weighted variations; the higher indices
are supplied by \eqref{explicit-stage-orders}.
The first geometric curve increment is
$O(\epsilon\Lambda_\epsilon^C)$ and is part of the slow reference.
Consequently every curve remains in the common tubular and
constant-speed neighbourhood.  The normalisations and denominators in
the construction are bounded uniformly away from zero wherever a
quotient is used.

We verify the terminal residual in the physical norm.  In a fixed tube,
the exact core-plus-outer identity gives, with two weighted variations,
\begin{equation}\label{terminal-physical-residual-bound}
 \begin{split}
 &\sup_t\|\mathcal D_{\rm in}^qS_\epsilon(u_\epsilon^{(K)})(t)
                                      \|_{H^N}\\
 &\qquad\leq C_{M,N}\Lambda_\epsilon^{C_{M,N}}
       \left(\epsilon^{-N-2}A_K
             +\epsilon^{J_K-N-2}A_{K-1}\right),\qquad q\leq2.
 \end{split}
\end{equation}
Here the smooth physical form of the identity is used at the zero.
In particular, $U(p+iv)=\nabla U$, so the terms with an apparent
singular quotient contain derivatives of the smooth field $U$.
Every physical derivative of its vortex or a core cutoff costs at
most $\epsilon^{-1}$.  The newest wave has the physical derivative
bounds in \eqref{outerinductivebounds}; its algebraic source part has
the bounds retained from the preceding batch.  The term involving
$\chi_{\rm c}S_\epsilon(U_K^-)$ has the second bound in
\eqref{terminal-physical-residual-bound}.  All other terms in
\eqref{higherouterexactremainders} contain a newest phase or amplitude
factor.  These facts give the first bound; the additional two negative
powers are a reserve for the normalisation and mixed coordinate
terms.  The tube has bounded volume, so the same pointwise bounds
imply its $L^2$ estimate.

Outside the fixed tube, terms involving $v_\Gamma$ have an integrable
spatial weight, and the remaining terms contain a spatial derivative of
a wave or an algebraic amplitude in $H^N$.  The dipole residual has
$|x|^{-2}$ decay and belongs to $L^2(\mathbb R^3)$.  The Sobolev product
estimates and \eqref{outerinductivebounds} therefore give the same
bound there.  This proves \eqref{terminal-physical-residual-bound}
globally.

Removing the two parameter weights costs at most $\epsilon^{-2}$.
Since $J_K\geq N+20$ and
\[
 A_K=\epsilon^{M+N+85/4},
\]
the right-hand side is bounded by $C_{M,N}\epsilon^M$ after removing
those weights and absorbing the finite logarithmic power.  We have
proved the stronger $H^N$ bound for each differential; the stated
$H^{N-q}$ estimate follows by the Sobolev inclusion, also when the
latter index is negative.

We conclude with the geometry and the core profile.  All curve
increments have zero initial value, so every stage has the same
initial geometric curve.  The chart makes each curve closed and
constant-speed.  Summing the slow and radiative increments gives
\eqref{sharperarbitrarycurvebound}.
Every added planar field is centred at $y=0$.  There the gluing map is
additive, and all outer factors are nonzero exponentials.  Therefore
\eqref{arbitrary-exact-centring-identity} holds identically as a
parameter-dependent equality.  Its two family derivatives include
the derivatives of the point at which the field is evaluated.

For definiteness, the core phase used in the theorem is the value of
the regular prefactor,
\[
 \beta_k(\rho,t)=\widehat\Phi_\epsilon[\Gamma_k](\rho,0,t)
                  +A^{(k)}(\Gamma_k^{\rm mat}(\rho,t),t).
\]

The centred expression is
\[
 u_{\epsilon,\mathrm{cen}}^{(k)}(\rho,y,t)
 =e^{-i\beta_k(\rho,t)}
   u_\epsilon^{(k)}(\mathcal X_{\Gamma_k}(\rho,y,t),t).
\]
For $k=1$ use the corresponding prefactor of Section
\ref{firstglobalimprovementsection}.  The first correction is
$O(\epsilon\Lambda_\epsilon^C)$ on each fixed inner ball.  The slow
higher profiles, the radiative profiles and the outer amplitudes are
smaller there by the established estimates.  Taylor's formula for the
regular phase prefactor gives
\eqref{arbitrary-centred-normal-derivative} and
\eqref{arbitraryinnerprofile}.  Applying these estimates at two
successive stages proves \eqref{successive-approximation-bound}.

Choose $r_0>0$ so that
$\|DW(y)-DW(0)\|\leq w'(0)/4$ on $B_{r_0}$.  The centred $C^1$
estimate and the fundamental theorem of calculus give
\[
 |u_{\epsilon,\mathrm{cen}}^{(k)}(\rho,y,t)|
 \geq \tfrac12 w'(0)|y|,\qquad |y|\leq r_0.
\]
On a fixed annulus around that disc, uniform closeness to $W$ excludes
additional zeros.  In the rest of the tube the field is a nonvanishing
vortex factor multiplied by exponentials with uniformly small
logarithmic amplitudes; the fixed transition annulus is covered by the
same smallness.  There are no further zeros.  This completes the proof
of Theorem \ref{arbitraryapproximationtheorem}.

\section{Three-dimensional wave estimates for the outer iteration}
\label{threeDwaveestimatesection}

This section proves the coupled phase--amplitude cancellation identity
and the wave estimates used in Section \ref{arbitraryapproxsection}.
The wave sources are decomposed into a part with $\langle x\rangle^{-3}$
decay, a finite-dimensional dipole contribution, and a general global
part controlled in Sobolev norms.  The first two contributions have
energy growth of order $(1+\tau_\epsilon)^{1/2}$, up to logarithmic
factors; the general energy
estimate has growth of order $1+\tau_\epsilon$.

\subsection{Coupled phase--amplitude cancellation}

For a smooth field $U$, write on $\{U\ne0\}$
\[
 E_1+iE_2=\frac{S_\epsilon(U)}{iU},\qquad q=|U|^2,\qquad
 p=\operatorname{Re}(\nabla U/U),\qquad
 v=\operatorname{Im}(\nabla U/U).
\]

\begin{lemma}[Residual after a coupled outer correction]
\label{exactfreewavecancellationlemma}
\label{exact-residual-ledger-lemma}
Let $H_1,H_2$ be prescribed smooth real functions, and suppose that the
real functions $a,b$ satisfy
\begin{equation}\label{freewavecancellationequations}
 -a_{\tau\tau}+\Delta_xa
 =\frac1{\sqrt2\epsilon}\partial_\tau H_2
   +v\cdot\nabla_xH_2-\epsilon^{-2}H_1,
 \qquad b=\tfrac12(H_2+\sqrt2\epsilon a_\tau).
\end{equation}
Put $U^+=Ue^{i(a+ib)}$.  Then
\begin{equation}\label{partial-outer-residual-identity}
 \begin{split}
 S_\epsilon(U^+)
 ={}&e^{i(a+ib)}\{S_\epsilon(U)-iU(H_1+iH_2)\}\\
    &+iU^+(T_1+iT_2),
 \end{split}
\end{equation}
where
\begin{equation}\label{exactfreewaverealremainder}
 T_1=2\epsilon^2p\cdot\nabla a
     -\sqrt2\epsilon^3v\cdot\nabla a_\tau
     -2\epsilon^2\nabla a\cdot\nabla b
\end{equation}
and
\begin{equation}\label{exactfreewaveimagremainder}
 \begin{split}
 T_2={}&-2(q-1)b+\epsilon^2\Delta b
       +2\epsilon^2p\cdot\nabla b+2\epsilon^2v\cdot\nabla a\\
       &+\epsilon^2(|\nabla a|^2-|\nabla b|^2)
       -q(1-e^{-2b}-2b).
 \end{split}
\end{equation}
In particular, when $H_j=E_j$, the normalised residual of $U^+$ is
$T_1+iT_2$.  When $H_j=(1-\chi)E_j$ for a smooth cutoff $\chi$, the
first term in \eqref{partial-outer-residual-identity} is
$e^{i(a+ib)}\chi S_\epsilon(U)$.
If $U$ has zeros, assume that $H_1,H_2$ vanish in a neighbourhood of
them and that $a,b$ are smooth.  The wave source is defined to be zero
there in the term $v\cdot\nabla H_2$, and the combined right-hand side
of \eqref{partial-outer-residual-identity} extends smoothly across the
zero set.
\end{lemma}

\begin{proof}
Apply \eqref{multiplicativeerroridentity} with the increment $a+ib$.
Since $\epsilon^2\Lambda_\epsilon\partial_t
=\sqrt2\epsilon\partial_\tau$, the wave and algebraic equations give
\[
 \begin{split}
 &E_1-\sqrt2\epsilon b_\tau+\epsilon^2\Delta a
             -2\epsilon^2v\cdot\nabla b\\
 &\quad=E_1-\frac\epsilon{\sqrt2}\partial_\tau H_2
       +\epsilon^2(-a_{\tau\tau}+\Delta a)
       -\epsilon^2v\cdot\nabla H_2
       -\sqrt2\epsilon^3v\cdot\nabla a_\tau\\
 &\quad=E_1-H_1-\sqrt2\epsilon^3v\cdot\nabla a_\tau.
 \end{split}
\]
Adding $2\epsilon^2p\cdot\nabla a$ and the real nonlinear term
$-2\epsilon^2\nabla a\cdot\nabla b$ gives $E_1-H_1+T_1$.
For the imaginary component,
\[
 E_2+\sqrt2\epsilon a_\tau-2qb
 =E_2-H_2-2(q-1)b.
\]
The remaining spatial and nonlinear terms give $E_2-H_2+T_2$.
Multiplication by $iU^+$ proves
\eqref{partial-outer-residual-identity} and its two special cases.

Near a zero of $U$, one has $H_2=0$ and
$b=\epsilon a_\tau/\sqrt2$.  The terms in $iU^+(T_1+iT_2)$ containing
$p$ or $v$ then combine as
\[
 2\epsilon^2e^{i(a+ib)}\nabla U\cdot(i\nabla a-\nabla b),
\]
using $U(p+iv)=\nabla U$.  Every remaining term is a product of smooth
fields.  This proves the asserted extension through the zero set.
\end{proof}

\subsection{Pointwise estimates from Kirchhoff's formula}

For a smooth function $F$ on
$\mathbb R^3\times[0,\tau_\epsilon]$, define
\begin{equation}\label{threeDwaveYnorm}
 \|F\|_{Y_m}
 :=\sup_{0\leq\tau\leq\tau_\epsilon}
 \sup_{x\in\mathbb R^3}
 \langle x\rangle^3
 \sum_{j+|\alpha|\leq m}
 |\partial_\tau^jD_x^\alpha F(x,\tau)|.
\end{equation}
For a source controlled only in global Sobolev norms define also
\begin{equation*}
 \|F\|_{\mathcal H_m}
 :=\sup_{0\leq\tau\leq\tau_\epsilon}
   \sum_{b=0}^m
   \|\partial_\tau^bF(\tau)\|_{H^{m-b}}.
\end{equation*}
Thus the norm $\mathcal H_m^{(2)}$ in
\eqref{differentiatedouternorms} is the corresponding norm after up to two
weighted curve variations.  The $Y_m$ class consists of sources with
$\langle x\rangle^{-3}$ decay, whether or not they are compactly supported.
Sources supported in a moving tube of uniformly bounded diameter belong to
this class uniformly because all curves remain in a common bounded set.

\begin{lemma}[Weighted spherical means]
\label{weightedsphericalmeanslemma}
For $T\geq2$,
\begin{equation}\label{weightedsphericalmeanestimate}
 \sup_{x\in\mathbb R^3}
 \langle x\rangle
 \int_0^T r\int_{\mathbb S^2}
       \langle x+r\omega\rangle^{-3}\,d\omega\,dr
 \leq C\log(2+T).
\end{equation}
\end{lemma}

\begin{proof}
The spherical-mean calculation in the proof of Lemma
\ref{threeDwavefirstimprovementlemma}, with $\tau_*$ replaced by $T$,
gives \eqref{weightedsphericalmeanestimate}.
\end{proof}

\begin{lemma}[Pointwise three-dimensional wave estimate]
\label{pointwisethreeDwavelemma}
Let $\psi$ solve
\begin{equation}\label{zerodatawavegeneral}
 -\psi_{\tau\tau}+\Delta\psi=F,
 \qquad \psi(\cdot,0)=\psi_\tau(\cdot,0)=0.
\end{equation}
For $m\geq3$,
\begin{equation}\label{pointwisethreeDwaveestimate}
 \sup_{\tau\leq\tau_\epsilon}
 \sum_{j+|\alpha|\leq m-2}
 \|\langle x\rangle
       \partial_\tau^jD_x^\alpha\psi(\tau)\|_\infty
 \leq C_m\Lambda_\epsilon^{C_m}\|F\|_{Y_m}.
\end{equation}
\end{lemma}

\begin{proof}
Kirchhoff's formula is
\begin{equation}\label{Kirchhoffarbitraryorder}
 \psi(x,\tau)=-\frac1{4\pi}
 \int_0^\tau r\int_{\mathbb S^2}
 F(x+r\omega,\tau-r)\,d\omega\,dr.
\end{equation}
Spatial derivatives and the first time derivative are estimated directly
by Lemma \ref{weightedsphericalmeanslemma}.  Repeated time derivatives are
handled by differentiating \eqref{Kirchhoffarbitraryorder}; alternatively,
use $\psi_{\tau\tau}=\Delta\psi-F$ after two derivatives.  The homogeneous
initial terms which then occur are spherical means of traces of derivatives
of $F$ and satisfy the same estimate.  Finally,
$\log(2+\tau_\epsilon)\leq C\Lambda_\epsilon$.
\end{proof}

\subsection{Energy estimates}

For $m\geq4$, define
\begin{equation}\label{threeDwaveXnorm}
\begin{split}
 \|\psi\|_{X_m}:={}&
 \sup_{0\leq\tau\leq\tau_\epsilon}
 \sum_{j=0}^{m}
 \left(
  \|\partial_\tau^{j+1}\psi\|_{H^{m-j}}
  +\|\nabla\partial_\tau^j\psi\|_{H^{m-j}}
 \right)\\
 &+\sup_{0\leq\tau\leq\tau_\epsilon}
 \sum_{j+|\alpha|\leq m-2}
 \|\partial_\tau^jD_x^\alpha\psi\|_\infty.
\end{split}
\end{equation}
For a family depending smoothly on the initial curve, set
\begin{equation*}
 \|\psi\|_{X_m}^{(2)}
 :=\max_{0\leq q\leq2}
   \|\mathcal D_{\rm in}^q\psi\|_{X_{m-q}}.
\end{equation*}
Only terms with nonnegative Sobolev indices are included.  Notice that no
$L^2$ norm of $\psi$ itself is required; a forced phase may have a
$|x|^{-1}$ tail while its derivatives have finite energy.

\begin{lemma}[Square-root energy growth]
\label{squarerootwaveenergylemma}
Under the assumptions of Lemma \ref{pointwisethreeDwavelemma},
\begin{equation*}
 \sup_{\tau\leq\tau_\epsilon}
 \sum_{j=0}^{m-3}
 \left(
  \|\partial_\tau^{j+1}\psi(\tau)\|_{H^{m-3-j}}
  +\|\nabla\partial_\tau^j\psi(\tau)\|_{H^{m-3-j}}
 \right)
 \leq C_m\epsilon^{-1/2}
 \Lambda_\epsilon^{C_m}\|F\|_{Y_m}.
\end{equation*}
\end{lemma}

\begin{proof}
Put
$\psi_{j,\alpha}=\partial_\tau^jD_x^\alpha\psi$ and
$F_{j,\alpha}=\partial_\tau^jD_x^\alpha F$.  The commuted energy identity is
\begin{equation*}
 \frac{d}{d\tau}\int
 \big(|\partial_\tau\psi_{j,\alpha}|^2
      +|\nabla\psi_{j,\alpha}|^2\big)
 =-2\int F_{j,\alpha}\partial_\tau\psi_{j,\alpha}.
\end{equation*}
The definition of $Y_m$ and the weighted pointwise estimate
\eqref{pointwisethreeDwaveestimate} give
\[
 \left|\int F_{j,\alpha}\partial_\tau\psi_{j,\alpha}\right|
 \leq C_m\Lambda_\epsilon^{C_m}\|F\|_{Y_m}^2
       \int_{\mathbb R^3}\langle x\rangle^{-4}\,dx
 \leq C_m\Lambda_\epsilon^{C_m}\|F\|_{Y_m}^2.
\]
Integrating and taking a square root gives
\[
 C_m\tau_\epsilon^{1/2}\Lambda_\epsilon^{C_m}\|F\|_{Y_m}
 \leq C_m\epsilon^{-1/2}\Lambda_\epsilon^{C_m}\|F\|_{Y_m}.
\]
Time-commuted equations contain only lower initial traces of $F$; their
energy is bounded directly by the $Y_m$ norm.
\end{proof}

\begin{lemma}[Unweighted energy]
\label{unweightedwaveenergylemma}
If \eqref{zerodatawavegeneral} holds and $\|F\|_{\mathcal H_m}<\infty$, then
\begin{equation}\label{unweightedwaveenergyestimate}
 \|\psi\|_{X_{m-3}}
 \leq C_m\epsilon^{-1}\Lambda_\epsilon^{C_m}
 \|F\|_{\mathcal H_m}.
\end{equation}
\end{lemma}

\begin{proof}
For $0\leq j\leq m-3$, commute the equation with
$\partial_\tau^jD_x^\alpha$, where $|\alpha|\leq m-3-j$.  The standard
energy identity gives
\[
 \|\partial_\tau^{j+1}\psi(\tau)\|_{H^{m-3-j}}
 +\|\nabla\partial_\tau^j\psi(\tau)\|_{H^{m-3-j}}
 \leq C_m(1+\tau)\|F\|_{\mathcal H_m}.
\]
For $j\geq1$, the initial traces on the left are finite sums of spatial
derivatives of $\partial_\tau^\ell F(0)$ with
$\ell\leq j-2$; this follows recursively from
$\psi_{\tau\tau}=\Delta\psi-F$.  They are bounded by
$C_m\|F\|_{\mathcal H_m}$.  Since
$\tau_\epsilon\leq C\epsilon^{-1}\Lambda_\epsilon^{-1}$, this proves the
energy part of \eqref{unweightedwaveenergyestimate}.  The pointwise part
follows from
\begin{equation*}
 \|h\|_\infty
 \leq C\|\nabla h\|_2^{1/2}\|\nabla^2h\|_2^{1/2}
\end{equation*}
and its differentiated versions.  The three spare derivatives in the
statement control every term with
$j+|\alpha|\leq m-5$, including the terms obtained by using the equation to
replace two time derivatives.  This is precisely the pointwise component of
$X_{m-3}$.
\end{proof}

\begin{proposition}[Mixed outer wave estimate]
\label{mixedthreeDwaveestimate}
Let
\begin{equation*}
 -\psi_{\tau\tau}+\Delta\psi=F_{\rm loc}+F_{\rm g},
 \qquad \psi(0)=\psi_\tau(0)=0,
\end{equation*}
Assume that $F_{\rm loc}$ belongs to the weighted source class and that
both sources depend smoothly on
$\Gamma_{\rm in}\in\mathscr O_{\rm in}$.  Then
\begin{equation}\label{mixedthreeDwavebound}
 \|\psi\|_{X_{m-3}}^{(2)}
 \leq C_m\Lambda_\epsilon^{C_m}
 \left(
  \epsilon^{-1/2}\|F_{\rm loc}\|_{\mathcal Y_m^{(2)}}
  +\epsilon^{-1}\|F_{\rm g}\|_{\mathcal H_m^{(2)}}
 \right).
\end{equation}
If the sources have no curve dependence, the same statement holds with
$X_{m-3}$, $Y_m$ and $\mathcal H_m$ in place of their differentiated
counterparts.
\end{proposition}

\begin{proof}
Split $\psi=\psi_{\rm loc}+\psi_{\rm g}$ by linearity.  Apply Lemmas
\ref{pointwisethreeDwavelemma} and
\ref{squarerootwaveenergylemma} to
$\mathcal D_{\rm in}^q\psi_{\rm loc}$ and Lemma
\ref{unweightedwaveenergylemma} to
$\mathcal D_{\rm in}^q\psi_{\rm g}$, for
$0\leq q\leq2$, with $m$ replaced by $m-q$.  Curve derivatives commute
with the constant-coefficient wave operator and preserve the zero initial
conditions.  Taking the maximum over $q$ gives
\eqref{mixedthreeDwavebound}.
\end{proof}

\subsection{Wave estimates for dipole sources}
\label{dipole-wave-estimates-subsection}

Choose a smooth radial function $\eta$ of compact support and integral
one, put
\[
 V=\eta*(4\pi|x|)^{-1},\qquad f_j=\partial_jV,\quad 1\leq j\leq3.
\]
Outside a fixed ball, $f_j$ is a constant multiple of $x_j/|x|^3$.
These functions are independent of time and of the initial-curve
parameter.  For a scalar coefficient family $c=c(\tau;\Gamma_{\rm in})$
we use
\[
 \|c\|_{C_\tau^m}^{(2)}
 :=\max_{0\leq q\leq\min\{2,m\}}
   \max_{0\leq b\leq m-q}\sup_{0\leq\tau\leq\tau_\epsilon}
       |\partial_\tau^b\mathcal D_{\rm in}^q c(\tau)|.
\]
The parameter differentials are measured in the family norms of
Section \ref{arbitraryapproxsection}.

\begin{lemma}[Wave estimate with a dipole source]
\label{fixed-tube-wave-estimate-lemma}
Let $a$ be the zero-data solution of
\[
 -a_{\tau\tau}+\Delta_xa
 =F_{\rm loc}+\sum_{j=1}^3c_j(\tau)f_j+F_{\rm g}.
\]
For $m\geq10$,
\begin{equation}\label{fixed-tube-wave-estimate}
 \|a\|_{X_{m-4}}^{(2)}
 \leq C_m\Lambda_\epsilon^{C_m}
 \left\{\epsilon^{-1/2}
       \left(\|F_{\rm loc}\|_{Y_m}^{(2)}
               +\sum_j\|c_j\|_{C_\tau^m}^{(2)}\right)
       +\epsilon^{-1}\|F_{\rm g}\|_{\mathcal H_m}^{(2)}\right\}.
\end{equation}
\end{lemma}

\begin{proof}
The local and general global terms are covered by
\eqref{mixedthreeDwavebound}.  For a dipole,
$\widehat f_j(\xi)=i\xi_j|\xi|^{-2}\widehat\eta(\xi)$ up to the
Fourier normalisation.  Consequently
\[
 \sup_{r>0}r^2(1+r^2)^m
       \int_{S^2}|\widehat f_j(r\omega)|^2d\omega<\infty.
\]
For $F=c(\tau)f_j$, the Fourier formula for the zero-data solution
expresses the energy in terms of
$\widehat f_j(\xi)\int_0^\tau e^{\pm i|\xi|s}c(s)ds$.
Polar integration and one-dimensional Plancherel give
\[
 \|a_\tau(\tau)\|_{H^m}^2+\|\nabla a(\tau)\|_{H^m}^2
 \leq C_m\int_0^\tau|c(s)|^2ds
 \leq C_m\tau\|c\|_\infty^2.
\]
Time differentiation adds forced solutions with differentiated $c$ and
homogeneous solutions whose data are derivatives of $c$ at zero times
spatial derivatives of $f_j$.  Their energies have the same bound,
with $1+\tau$ in place of $\tau$.  The Fourier splitting into
$|\xi|<1$ and $|\xi|\geq1$ gives
$\|a\|_\infty\leq C(\|\nabla a\|_2+\|D_x^2a\|_2)$.
Sobolev embedding gives the other pointwise terms of $X_{m-4}$.
Finally $\tau_\epsilon^{1/2}\leq C\epsilon^{-1/2}$.
Parameter differentiation commutes with this fixed wave operator and
with the fixed dipole functions.
\end{proof}

\section{Linear theory in the global rescaled variables}
\label{curvedlinearsection}

We now prove an $H^1$ inverse estimate for the rescaled Gross--Pitaevskii
equation linearised at the accurate approximation of
Section~\ref{arbitraryapproxsection}.

\subsection{The rescaled equation and the linear estimate}
\label{linear-sec-problem}

Put
\begin{equation}\label{globalrescaledvariable}
 {\Lambda_\epsilon}=|\log\epsilon|,\qquad {c_\epsilon}=\epsilon^2{\Lambda_\epsilon},
 \qquad X=x/\epsilon,\qquad
 \tau=\frac{\sqrt2t}{\epsilon{\Lambda_\epsilon}}.
\end{equation}
Inner products of complex fields are real:
\[
 \langle f,h\rangle=\operatorname{Re}\int_{{\mathbb R}^3}f\overline h\,dX.
\]
The spaces $H^m$ below are complex Sobolev spaces regarded as real
Hilbert spaces.  Constants are uniform in $0<\epsilon<\epsilon_0$
and $0\leq t\leq T$.  Intermediate profile constants and logarithmic
powers may depend on the finite construction orders.  The constants
in the linear estimate below can be chosen independently of those
orders, with a threshold $\epsilon_0$ that may depend on them;
Lemma~\ref{linear-lem-uniformgrowth} proves this uniformity.
Dependence on the reference curve, its low geometric neighbourhood,
the fixed cutoffs and $T$ is allowed.

Fix $\Gamma_{\rm in}\in\mathscr O_{\rm in}$ and let
$U_\epsilon=u_\epsilon^{(K)}[\Gamma_{\rm in}]$ be the terminal
approximation of Theorem~\ref{arbitraryapproximationtheorem}, with
accuracy and regularity orders to be chosen sufficiently large later.
Set
\begin{equation}\label{linear-eq-background}
 \begin{split}
  V(X,t)=V_\epsilon(X,t)&:=U_\epsilon(\epsilon X,t),\\
  \mathscr R=R_\epsilon&:=i{c_\epsilon} V_t+\Delta_X V+(1-|V|^2)V.
 \end{split}
\end{equation}
Equivalently, the rescaled Gross--Pitaevskii operator and its residual
satisfy
\begin{equation}\label{rescaledresidualidentity}
 \begin{split}
  \mathscr G_\epsilon(v)&:=ic_\epsilon\partial_tv+\Delta_Xv
                                      +(1-|v|^2)v,\\
  R_\epsilon(X,t)&=\mathscr G_\epsilon(V_\epsilon)(X,t)
                    =S_\epsilon(U_\epsilon)(\epsilon X,t).
 \end{split}
\end{equation}
Define
\begin{equation*}
 H^X_t h=-\Delta_Xh+(|V|^2-1)h
                   +2\operatorname{Re}(\overline Vh)V,
 \qquad \mathcal L_{\epsilon,t}^{X}:=-H^X_t.
\end{equation*}
We consider the Gross--Pitaevskii equation linearised at $V_\epsilon$
in the rescaled variables $(X,t)$, with forcing $F$ and zero initial data:
\begin{equation}\label{curvedlinearproblem}
 i{c_\epsilon}\Phi_t-H^X_t\Phi=F,\qquad \Phi(0)=0.
\end{equation}
The objective is a bound in $C([0,T];H^1)$ with a fixed algebraic
loss in $\epsilon$.  No orthogonality condition is imposed on $F$
or on $\Phi$.

We use the finite-stage formulae of
Sections~\ref{normalgraphsection}--\ref{threeDwaveestimatesection} to obtain
the differentiated bounds for $U_\epsilon$ needed below; these are
established in Subsection~\ref{linear-sec-backgroundbounds}.

The main result is the following linear estimate.
\begin{theorem}[Rescaled linear estimate about the constructed filament]\label{curvedfilamentlineartheorem}
Choose the approximation $U_\epsilon$ with sufficiently high accuracy
and regularity orders.  There are positive constants $C_T$ and
$C_{\rm lin}$, independent of these orders, such that, for sufficiently
small $\epsilon>0$ and every
\[
 F\in L^\infty((0,T);H^1({\mathbb R}^3_X)),
\]
equation \eqref{curvedlinearproblem} has a unique solution
$\Phi\in C([0,T];H^1({\mathbb R}^3_X))$, satisfying
\begin{equation}\label{curvedfilamentlinearestimate}
 \|\Phi\|_{L^\infty_tH^1_X}
 \leq C_T\epsilon^{-4-C_{\rm lin}T}\|F\|_{L^\infty_tH^1_X}.
\end{equation}
\end{theorem}

The proof constructs a positive quadratic functional in coordinates
moving with the filament.  The quadratic form associated with the
transformed linearised equation need not be nonnegative, so it is
augmented by terms controlling the translation components and the
squared $L^2$ norm of the perturbation.  The resulting functional
controls the $H^1$ norm.  Controlling its time derivative requires
delicate cancellations in the differentiated quadratic form and the
translation equations.  These yield an evolution inequality with
growth rate of order $|\log\epsilon|$, from which Gronwall's inequality
gives the stated estimate.

\subsection{The planar form and regional inequalities}
\label{linear-sec-planar}

\subsubsection{Angular blocks}
Let $W(y)=w(r)e^{i\theta}$ be the degree-one vortex.  We use
\begin{equation}\label{linear-eq-wproperties}
 \begin{gathered}
 w''+r^{-1}w'-r^{-2}w+(1-w^2)w=0,\\
 0<w<1,\quad w'>0,\quad w(r)\asymp r\ (r\leq1),\\
 w=1-\tfrac12r^{-2}+O(r^{-4}),\qquad w'=r^{-3}+O(r^{-5}).
 \end{gathered}
\end{equation}
The profile properties follow from
\cite{chenelliottqi1994,herveherve1994}.  Write
\begin{equation}\label{linear-eq-planarform}
 B_W[z]=\int_{{\mathbb R}^2}\left\{
  |\nabla_yz|^2-(1-w^2)|z|^2+2\operatorname{Re}(\overline Wz)^2\right\}dy
\end{equation}
and
\begin{equation*}
 \|z\|_{\mathcal H_W}^2
 =\int_{{\mathbb R}^2}\left\{|\nabla_yz|^2+\frac{|z|^2}{1+r^2}
                              +\operatorname{Re}(\overline Wz)^2\right\}dy.
\end{equation*}
By \cite{delpinofelmerkowalczyk2004}, $B_W$ is nonnegative on
$\mathcal H_W$, and its kernel is
$\operatorname{span}_{{\mathbb R}}\{Z_1,Z_2\}$, where $Z_j=\partial_{y_j}W$.
The field $iW$ does not belong to $\mathcal H_W$.

For $r>0$, write $\psi=z/(iW)=\psi_1+i\psi_2$ and expand its real
and imaginary parts in angular Fourier series.  Let $f(r)$ and $g(r)$
be the angular averages of $\psi_1$ and $\psi_2$, respectively.  Then
\[
 z=iW(f+ig)+z^1+z^{\geq2},
\]
where $z^1$ is obtained by retaining the $\cos\theta$ and $\sin\theta$
terms in $\psi$ and multiplying by $iW$, while $z^{\geq2}$ is obtained
by retaining all harmonics $\cos(k\theta)$, $\sin(k\theta)$ with
$k\geq2$ and multiplying by $iW$.

For tests away from zero,
\begin{equation}\label{linear-eq-planarpsi}
 B_W[iW\psi]
 =\int w^2\left\{|\nabla\psi|^2
       -\frac4{r^2}\psi_2\partial_\theta\psi_1+2w^2\psi_2^2\right\}dy.
\end{equation}
Integration by parts using \eqref{linear-eq-wproperties} proves the formula.
For a zero-mode field it extends through the origin and gives
\begin{equation}\label{linear-eq-Ezero}
 E_0(f,g):=2\pi\int_0^\infty
      \left(w^2f_r^2+w^2g_r^2+2w^4g^2\right)r\,dr
       =B_W[iW(f+ig)].
\end{equation}
For mode one, the separate terms in \eqref{linear-eq-planarpsi} need not
be integrable at zero.  The original form \eqref{linear-eq-planarform}
is used there.

Choose a nonnegative smooth radial $\chi_*$ with fixed compact
support and positive near zero.  Set
\[
 M_{ij}=\operatorname{Re}\int\chi_*Z_j\overline{Z_i},\qquad
 \beta_i=M^{-1}_{ij}\operatorname{Re}\int\chi_*z^1\overline{Z_j},\qquad
 z^{1,\perp}=z^1-\beta_jZ_j.
\]
The matrix $M$ is positive definite.
\begin{lemma}[Fixed-core coercivity]\label{linear-lem-planarcoercivity}
There are fixed $k,C>0$ such that
\begin{equation}\label{linear-eq-planarcoercivity}
 B_W[z]\geq E_0(f,g)+k\|z^{1,\perp}\|_{\mathcal H_W}^2
                             +k\|z^{\geq2}\|_{\mathcal H_W}^2.
\end{equation}
Also
\begin{equation}\label{linear-eq-amplitudecontrol}
 \|Wg\|_{\mathcal H_W}^2+\|Wg\|_2^2\leq C E_0(0,g).
\end{equation}
On the complement of the zero phase, equivalently,
$B_W[z]\geq k\|z\|_{\mathcal H_W}^2-C|\beta|^2$.
\end{lemma}
\begin{proof}
The operator preserves the real Fourier blocks.  On the exterior
of a sufficiently large fixed disc, use
\[
 \frac4{r^2}|\psi_2\partial_\theta\psi_1|
 \leq \frac1{2r^2}|\partial_\theta\psi_1|^2
                                   +\frac8{r^2}|\psi_2|^2.
\]
The last term is absorbed by the amplitude potential.  The nonzero
phase is controlled with weight $r^{-2}$ by angular Poincar\'e.
A fixed radial partition, with its transition sufficiently wide,
therefore gives
\begin{equation}\label{linear-eq-compactcoercivity}
 \|z\|_{\mathcal H_W}^2
       \leq C B_W[z]+C_{R_*}\|z\|_{L^2(B_{R_*})}^2
\end{equation}
on the complement of the zero phase.  The partition errors in its
exterior are absorbed by the angular term; those on the fixed
interior remain in the second term.

If coercivity failed under the two fixed-core measurements, a
sequence of unit norm and vanishing form would converge weakly in
$\mathcal H_W$ and strongly on the fixed disc.  Form
Cauchy--Schwarz and the kernel characterisation in
\cite{delpinofelmerkowalczyk2004} imply that its limit is a translation.  The two
measurements exclude that kernel, so the limit is zero, contradicting
\eqref{linear-eq-compactcoercivity}.  Subtracting $\beta_jZ_j$ does not
change the form.  This proves \eqref{linear-eq-planarcoercivity}.

For the amplitude, a radial cutoff and integration of $(r^2g^2)'$
give
\begin{equation}\label{linear-eq-coreHardy}
 \int_0^2g^2r\,dr\leq C\int_0^3w^2g_r^2r\,dr
                              +C\int_1^3g^2r\,dr.
\end{equation}
Indeed the cutoff version is $\int rg^2\leq\int r^3g_r^2$.
On $[1,\infty)$, the term $w^4g^2$ controls the required mass.
The product rule applied to $Wg$ and \eqref{linear-eq-wproperties} prove
\eqref{linear-eq-amplitudecontrol}.
\end{proof}

\subsubsection{Free boundary values on a large normal disc}
For $R_0$ sufficiently large define
\begin{equation}\label{linear-eq-diskform}
 p_{R_0}[z]=B_W^{B_{R_0}}[z]
       -\frac{w'(R_0)}{w(R_0)}\int_{\partial B_{R_0}}|z|^2.
\end{equation}
The superscript means integration over the disc.  No boundary
condition is imposed.  The zero block of \eqref{linear-eq-diskform} is
exactly \eqref{linear-eq-Ezero} integrated from $0$ to $R_0$.
\begin{lemma}[Uniform disc comparison]\label{linear-lem-disk}
On fields with no zero phase,
\begin{equation}\label{linear-eq-diskcoercivity}
 p_{R_0}[z]\geq k\|z\|_{\mathcal H_W(B_{R_0})}^2-C|\beta|^2,
\end{equation}
with constants independent of sufficiently large $R_0$.
\end{lemma}
\begin{proof}
The boundary subtraction in \eqref{linear-eq-diskform} is exactly the term
arising when the modulus is integrated by parts.  On $r\geq R_*$,
the resulting phase--amplitude density is positive by the argument
of Lemma~\ref{linear-lem-planarcoercivity}.  Consequently
\eqref{linear-eq-compactcoercivity} holds on $B_{R_0}$ uniformly, with
$p_{R_0}$ in place of $B_W$.  If coercivity failed for a sequence
$R_0\to\infty$ under the two zero measurements, local compactness
and positivity of the exterior density give a whole-plane limit
with nonpositive form.  Nonnegativity and the measurements make
that limit zero, contradicting the uniform compact-remainder
estimate.  Removing the measurements subtracts $\beta_jZ_j$.
Their truncated norms and the forms $p_{R_0}[Z_i,Z_j]$ are uniformly
bounded; continuity and Young's inequality give
\eqref{linear-eq-diskcoercivity}.  The estimates extend by density to the
disc form domain.
\end{proof}

\subsubsection{The radial phase and exterior trace}
The following inequality will be used with the coefficient $c/r$,
not with a uniform mass weight.
\begin{lemma}[Radial form bound]\label{linear-lem-radial}
For real radial $f$ and $R_0>3$,
\begin{align*}
 c\int_0^{R_0}\frac{w^2}{\langle r\rangle}f^2r\,dr
 &\leq \eta\int_0^{R_0}w^2f_r^2r\,dr
       +C_\eta c^2\int_0^{R_0}w^2f^2r\,dr\\
 &\quad+cR_0w(R_0)^2f(R_0)^2.
\end{align*}
If $f(R_0)=0$, the boundary term is absent.
\end{lemma}
\begin{proof}
Since $(rw^2)'\geq w^2$,
\[
 \int_0^{R_0}w^2f^2\,dr
 \leq R_0w(R_0)^2f(R_0)^2
                     -2\int_0^{R_0}rw^2ff_r\,dr.
\]
Cauchy--Schwarz, Young's inequality and
$r/\langle r\rangle\leq1$ prove the claim.  At zero the boundary
term vanishes for smooth Cartesian zero-mode fields; the weighted
Hardy inequality gives the extension to the form domain.
\end{proof}

For a fixed physical tube of radius $d$, write $\Omega_{\epsilon,t}$
for its exterior and $\Sigma_{\epsilon,t}$ for its boundary in $X$.
For smooth compactly supported $F$ on the closure of this exterior,
\begin{equation}\label{linear-eq-exteriortrace}
 \|F\|_{L^6(\Omega_{\epsilon,t})}^2
 \leq C_d\|\nabla F\|_2^2,\qquad
 \|F\|_{L^2(\Sigma_{\epsilon,t})}^2
 \leq C_d\epsilon^{-1}\|\nabla F\|_2^2.
\end{equation}
To prove the physical versions, first control $F$ on a large fixed
outer annulus by integration along rays from infinity.  Poincar\'e
and the trace theorem on the bounded connected region between that
annulus and the tube give the interior trace.  Extend into the
bounded hole using a bounded extension on this fixed smooth domain;
its $H^1$ norm is bounded by the same gradient norm.  The ordinary
homogeneous Sobolev inequality on ${\mathbb R}^3$ gives the first estimate.
The geometry lies in a uniformly smooth compact family.  Rescaling
$x=\epsilon X$ leaves the Sobolev inequality unchanged and gives
the factor $\epsilon^{-1}$ in the trace inequality.

The same argument at a fixed \emph{inner} rescaled radius provides
an extension of a phase function into each core disc that preserves
constants.  Using averages and a fixed radial reflection on an
annulus gives a gradient bound uniform in the length of the
rescaled filament.  In particular, if $f$ is a zero-mode phase,
write
\[
 E_{\rm ph,in}(f)=2\pi\int\!\int_0^{R_0}
                      w^2|f_r|^2r\,dr\,d\sigma,\qquad
 E_{\rm ph,out}(f)=\int_{\Omega_{\epsilon,t}}|\nabla F|^2 .
\]
Here $F$ is the exterior phase of the same field, whose boundary
angular average is $f(R_0)$; the outer notation includes this
specified continuation.  The radius $R_0$ is the interface radius
under consideration, and uniformly comparable metric weights
may be inserted.  Then
\begin{equation}\label{linear-eq-phasecoremass}
 \int_{r<r_*}w^2f^2\,dy\,d\sigma
 \leq C{\Lambda_\epsilon}\left(E_{\rm ph,in}(f)+E_{\rm ph,out}(f)\right).
\end{equation}
Indeed, for $r_*<r<R_0=d/\epsilon$, integrate $f_r$ and use
$\int_{r_*}^{R_0}dr/r\leq C{\Lambda_\epsilon}$; the value at $R_0$ is controlled
by \eqref{linear-eq-exteriortrace} and angular Jensen.  Estimate
\eqref{linear-eq-coreHardy} completes the bounded-core part.  These
extensions do not require the zero phase to vanish at the interface.

\begin{samepage}
\subsection{Canonical coordinates, phase and exact forms}
\label{linear-sec-coordinates}

We now use a unitary change of coordinates and phase to keep the
filament fixed.  We identify the transformed operator and its quadratic
form, including the terms generated by the motion of the coordinates.

\subsubsection{The global unitary map}
Let $\gamma_\epsilon$ be the arclength representative of the terminal
curve $\gamma_\epsilon^{(K)}[\Gamma_{\rm in}]$, with physical length
$\ell(t)$.  Use the material curve
$\Gamma(\varrho,t)=\gamma_\epsilon(\ell(t)\varrho,t)$,
$\varrho\in{\mathbb R}/{\mathbb Z}$, and its oriented periodic normal frame $e_1,e_2$.
\end{samepage}
Put
\[
 \mathcal T_t(\varrho,z)=\Gamma(\varrho,t)+z_ae_a(\varrho,t),
 \qquad
 \widetilde\gamma_\epsilon(\sigma,t)
       =\epsilon^{-1}\gamma_\epsilon(\epsilon\sigma,t).
\]
Here $0\leq\sigma<\ell(t)/\epsilon$ is rescaled arclength.  The
rescaled tubular map is
\begin{equation}\label{rescaledtubularcoordinates}
 \begin{split}
 \mathcal F_t(\varrho,y)
    &=\epsilon^{-1}\mathcal T_t(\varrho,\epsilon y)\\
    &=\widetilde\gamma_\epsilon(\ell(t)\varrho/\epsilon,t)
                                      +y_ae_a(\varrho,t).
 \end{split}
\end{equation}
Choose a fixed physical $d>0$ so that all the tubes of radius $8d$
are embedded.  A radial cutoff $\chi_{\rm fl}$ equals one on the
$6d$ tube and vanishes outside the $7d$ tube.  Define
\[
 Y^{\rm ph}(\mathcal T_t(\varrho,z),t)
       =\chi_{\rm fl}(z)\partial_t\mathcal T_t(\varrho,z),
 \qquad Y(X,t)=\epsilon^{-1}Y^{\rm ph}(\epsilon X,t).
\]
Let $\Phi_t^{\rm amb}$ be the flow of $Y^{\rm ph}$, initially the
identity, and put
\begin{equation*}
 \Xi_t(Z)=\epsilon^{-1}\Phi_t^{\rm amb}(\epsilon Z),\qquad
 \Xi_t(\mathcal F_0(\varrho,y))=\mathcal F_t(\varrho,y)
                  \quad(|\epsilon y|<6d).
\end{equation*}
The second identity follows from uniqueness of the flow.  Both
material and normal coordinates are preserved in this region.

Let $P_t$ be the leading regularised vortex field, equal to
$e^{i\widehat\Phi_\epsilon}W$ in the tube and continued by the
leading solid-angle ansatz.  Off its zero curve define
\begin{equation}\label{linear-eq-gauge}
 n_t(X)=\frac{P_t(X)}{|P_t(X)|}e^{iA_\epsilon(\epsilon X,t)},\qquad
 g_t(X)=n_t(X)\overline{n_0(\Xi_t^{-1}X)}.
\end{equation}
In the material tube put
\[
 \Theta_{\rm reg}(\varrho,y,t)=\widehat\Phi_\epsilon(\varrho,y,t)
             +A_\epsilon(\mathcal T_t(\varrho,\epsilon y),t).
\]
The angular factors cancel, giving
\[
 g_t(\mathcal F_t(\varrho,y))
       =e^{i(\Theta_{\rm reg}(\varrho,y,t)
                      -\Theta_{\rm reg}(\varrho,y,0))}.
\]
This provides the smooth extension through the filament.  The
reference winding is not removed.

With $J_t=\det D\Xi_t$, define
\begin{equation*}
 (\mathscr S_th)(Z)=J_t(Z)^{1/2}
         \overline{g_t(\Xi_tZ)}h(\Xi_tZ),\qquad
 q=\mathscr S_t\Phi,\quad G=\mathscr S_tF.
\end{equation*}
The map is complex-linear and unitary on $L^2$.  The low-order
bounds proved in Subsection~\ref{linear-sec-backgroundbounds} give
\begin{equation}\label{linear-eq-H1map}
 C^{-1}\|h\|_{H^1_X}\leq\|\mathscr S_th\|_{H^1_Z}
                              \leq C\|h\|_{H^1_X}.
\end{equation}
In fact, $D\Xi,(D\Xi)^{-1}$ are bounded,
$D^2\Xi,\nabla_Z\log J=O(\epsilon)$, and
$\nabla_Xg=O(\epsilon{\Lambda_\epsilon}^P)$.  The flow equation proves the first
two assertions.  The local formula for $g$ proves its core bound;
on the complement both phases in \eqref{linear-eq-gauge} stay away from
their respective cores.  A change of variables and the product rule
then prove \eqref{linear-eq-H1map}, including the inverse estimate.

Put
\[
 K_t=(\partial_t\mathscr S_t)\mathscr S_t^{-1},\qquad
 \vartheta=-i\overline g(\partial_t+Y\cdot\nabla_X)g.
\]
The operator $K_t$ is skew-adjoint.  The exact moving equation is
\begin{equation}\label{linear-eq-movingequation}
 i{c_\epsilon} q_t-H_tq=G,\qquad
 H_t=\mathscr S_tH^X_t\mathscr S_t^{-1}+i{c_\epsilon} K_t.
\end{equation}
Its symmetric real form is denoted by $B_t$.  If
$h=\mathscr S_t^{-1}q$, its Eulerian expression is
\begin{equation}\label{linear-eq-eulerianform}
 \begin{split}
 \mathfrak b_t[h]:=B_t[q]=\int\big\{&|\nabla h|^2+(|V|^2-1)|h|^2
                       +2\operatorname{Re}(\overline Vh)^2\\
              &-{c_\epsilon} Y\cdot j(h)+{c_\epsilon}\vartheta|h|^2\big\},
 \quad j(h)=\operatorname{Im}(\overline h\nabla h).
 \end{split}
\end{equation}
Indeed
\[
 \mathscr S_t^{-1}K_t\mathscr S_t
       =Y\cdot\nabla_X+\tfrac12\operatorname{div} Y-i\vartheta.
\]
The half-divergence term has zero real quadratic pairing after
multiplication by $i$.  This proves \eqref{linear-eq-eulerianform} and fixes
both the transport and gauge signs.

\subsubsection{Material metric and strain}
Write $\sigma_0=\ell(0)\varrho/\epsilon$,
$a_\ell(t)=\ell(t)/\ell(0)$, and
\[
 H=1-\epsilon\kappa\cdot y,\qquad j_0=a_\ell H.
\]
Here $\kappa=(\kappa_1,\kappa_2)$ is expressed in the current
periodic frame.  If $\mathfrak a$ is its physical normal connection,
\begin{equation}\label{linear-eq-metric}
 \begin{aligned}
 \partial_{\sigma_0}\mathcal F_t
   &=a_\ell Ht_\Gamma+\epsilon a_\ell\mathfrak a\,y^\perp,
       &\partial_{y_a}\mathcal F_t&=e_a,\\
 dX&=j_0\,dy\,d\sigma_0,
       &y^\perp&=-y_2e_1+y_1e_2.
 \end{aligned}
\end{equation}
The $\sigma_0$ derivative is a fixed-coordinate derivative.  Current
rescaled arclength is $\sigma=a_\ell\sigma_0$ at a fixed time.
For scalars,
\[
 \nabla_X=e_a\partial_{y_a}
   +t_\Gamma(a_\ell H)^{-1}
            (\partial_{\sigma_0}-\epsilon a_\ell\mathfrak a\partial_\theta).
\]
The form on the tube may be described in the fixed material field
\begin{equation}\label{linear-eq-materialfield}
 z=\sqrt{j_0}\,e^{-i\Theta_{\rm reg}}h\circ\mathcal F_t.
\end{equation}
We denote this local transform by $\mathscr S_{\rm mat}$ and put
$\widetilde G=\mathscr S_{\rm mat}F$ for the transformed forcing.
For $q$ fixed, this is a time-independent transform of its restriction
to the initial tube.  In particular, coefficient differentiation in
$z$ is the same fixed-field differentiation as in $q$.

Here is the full operator in this representation.  Indices $A,B$ range
over $\sigma_0,y_1,y_2$, and $Y^A$ are the contravariant components of
$Y\circ\mathcal F_t$.  Put
\[
 M=(a_\ell H)^{-2},\qquad
 \omega=\epsilon a_\ell\mathfrak a,\qquad
 b^a=\omega(y^\perp)^a,\qquad \theta_A=\partial_A\Theta_{\rm reg}.
\]
The inverse metric and the real drift are
\begin{equation}\label{linear-eq-materialcoefficients}
 \begin{gathered}
 G^{\sigma_0\sigma_0}=M,\qquad G^{\sigma_0a}=-Mb^a,\qquad
 G^{ab}=\delta_{ab}+Mb^ab^b,\\
 d^A=2G^{AB}\theta_B-c_\epsilon Y^A .
 \end{gathered}
\end{equation}
In the $dy\,d\sigma_0$ measure,
\begin{equation}\label{linear-eq-fullmaterialoperator}
 \mathbb H_tz=-\partial_A(G^{AB}\partial_Bz)
 -i\bigl(d^A\partial_A+\tfrac12\partial_Ad^A\bigr)z
 +V_{\rm sc}z+2\operatorname{Re}(\overline vz)v,
\end{equation}
where $v=e^{-i\Theta_{\rm reg}}V\circ\mathcal F_t$ and
\[
\begin{split}
 V_{\rm sc}={}&\rho^2-1
 +\tfrac12\partial_A(G^{AB}\partial_B\log j_0)
 +\tfrac14G^{AB}\partial_A\log j_0\,\partial_B\log j_0\\
 &+G^{AB}\theta_A\theta_B
 +c_\epsilon\partial_t\Theta_{\rm reg}
 -c_\epsilon Y^A\theta_A .
\end{split}
\]
The time derivative here holds the material coordinates fixed.
Indeed conjugation of $-j_0^{-1}\partial_A(j_0G^{AB}\partial_B)$
by $\sqrt{j_0}$ gives the two displayed density potentials.
Conjugation by the phase gives the first part of $d$ and
$G^{AB}\theta_A\theta_B$.  The moving map contributes
$ic_\epsilon(Y^A\partial_A+\frac12\partial_AY^A)$ and
$c_\epsilon(\partial_t\Theta_{\rm reg}-Y^A\theta_A)$.
The rank-one term contains $v$, rather than $\sqrt{j_0}v$.
All zeroth-order terms in \eqref{linear-eq-fullmaterialoperator}
commute with multiplication by a real scalar.

Let $S(Y)=(\nabla Y+\nabla Y^T)/2$.  The identity
$\dot g_{AB}=2\partial_A\mathcal F_t\cdot S(Y)\partial_B\mathcal F_t$
and \eqref{linear-eq-metric} give
\begin{equation}\label{linear-eq-strain}
 \begin{split}
 S(Y)_{ab}&=0,\\
 S(Y)_{ta}&=\frac{\epsilon}{2H}
      \left(\mathfrak a_t+\frac{\dot\ell}{\ell}\mathfrak a\right)
                                            (y^\perp)_a,\\
 S(Y)_{tt}&=\frac{\dot\ell}{\ell}+\frac{H_t}{H}.
 \end{split}
\end{equation}
To check the middle identity, differentiate
$g_{0a}=\epsilon a_\ell\mathfrak a(y^\perp)_a$ and use
$S_{ab}=0$.  Differentiation of
$g_{00}=a_\ell^2H^2+\epsilon^2a_\ell^2\mathfrak a^2r^2$
then gives the last identity; the connection contributions cancel.
Thus the normal strain vanishes exactly and the mixed strain is
$O(\epsilon r)$ in the material tube.  Also
\begin{equation*}
 \operatorname{div} Y=\partial_t\log j_0
                     =\dot\ell/\ell+H_t/H.
\end{equation*}

\subsubsection{The complete phase identities}
Off the filament set $\rho=|V|$, $V=\rho e^{i\Theta}$ locally and
$\alpha=\Theta-\Theta_{\rm lead}$, where $n=e^{i\Theta_{\rm lead}}$.
Only phase derivatives occur below.  Write
\begin{equation}\label{linear-eq-defects}
 \begin{split}
 D&=2\nabla\Theta-{c_\epsilon} Y,\qquad
 \mathcal V={D_t^{\rm mat}} V-i\vartheta V,\\
 e&=-{c_\epsilon}{D_t^{\rm mat}}\alpha-\rho^{-2}\operatorname{Re}(\mathscr R\overline V),\\
 s&=\rho^{-2}\operatorname{div}(\rho^2D)
       =2\rho^{-2}\operatorname{Im}(\mathscr R\overline V)
                  -{c_\epsilon}{D_t^{\rm mat}}\log\rho^2-{c_\epsilon}\operatorname{div} Y.
 \end{split}
\end{equation}
In the tube, $V\circ\mathcal F_t=e^{i\Theta_{\rm reg}}v$ and
\begin{equation*}
 \mathcal V\circ\mathcal F_t=e^{i\Theta_{\rm reg}}v_t.
\end{equation*}
For compactly supported real $f,g$ away from the curve,
\begin{equation}\label{linear-eq-fullphaseform}
 \mathfrak b_t[iV(f+ig)]
 =\int\rho^2\{ |\nabla f|^2+|\nabla g|^2
      -2D\cdot\nabla f\,g+2\rho^2g^2-sfg+e(f^2+g^2)\}.
\end{equation}
In particular the pure phase sees only its gradient and $e$.
To prove the formula, expand the kinetic term and integrate the
$\nabla\rho$ terms by parts.  The real residual identity
\[
 \rho^{-2}\operatorname{Re}(\mathscr R\overline V)
   =-{c_\epsilon}\Theta_t+\Delta\rho/\rho-|\nabla\Theta|^2+1-\rho^2
\]
eliminates its scalar part.  Integrate
$\rho^2D\cdot(f\nabla g-g\nabla f)$ by parts to obtain the
$s$ term.  The imaginary residual identity gives the last equality
in \eqref{linear-eq-defects}.

If $\mathfrak H_t$ is the operator of \eqref{linear-eq-eulerianform}, then
\begin{equation}\label{linear-eq-phaseJacobi}
 \mathfrak H_t(iV)=-{c_\epsilon}\mathcal V
        -\tfrac12{c_\epsilon}(\operatorname{div} Y)V-i\mathscr R=V(s/2+ie).
\end{equation}
The first expression is smooth at the filament and is always used
there.  It follows by differentiating the phase invariance of the
nonlinear equation, or by direct substitution.

\subsection{Differentiated estimates for the complete approximation}
\label{linear-sec-backgroundbounds}

We next obtain the differentiated bounds for the approximation needed
to estimate the coefficients of the transformed operator.

\subsubsection{Finite records and derivative reserve}
We distinguish the three contributions to the inner logarithmic
correction: the first correction, the finite slow elliptic batch,
and the radiative batches.  On $r\geq r_*$ their phase and amplitude
coefficients have the following bounds.  A fixed finite power of
${\Lambda_\epsilon}$ may multiply every term except where a logarithm is explicitly
kept as a function of $r$:
\begin{equation}\label{linear-eq-profiletable}
\begin{array}{c|c|c}
 &\text{phase}&\text{logarithmic amplitude}\\ \hline
 \text{first}&
 \epsilon\ell_r^P/r+\epsilon^2{\Lambda_\epsilon}^P&
 \epsilon\ell_r^P/r+\epsilon^2{\Lambda_\epsilon}^P\\
 \text{slow, }n\geq2&
 \epsilon^n\langle r\rangle^{n-2}&
 \epsilon^n\langle r\rangle^{n-2}\\
 \text{radiative, }n\geq2&
 \epsilon^n A_k\langle r\rangle^{n-1}&
 \epsilon^n A_k\langle r\rangle^{n-2}\\
 \text{radial radiative}&0&
 \epsilon A_k\langle r\rangle^{-2}
\end{array}
\end{equation}
Here $\ell_r=1+\log(2+r)$, $A_k=\epsilon^{(k+3)/4}$, and
$\sum_kA_k\leq C\epsilon^{5/4}$.  Each normal derivative lowers the
radial power by one.  Physical-time derivatives of the slow batch,
which is retained in material variables, preserve its displayed order
up to logarithms.  The first row has this property through two actual
material time derivatives, as verified below.  Radiative terms and
higher time derivatives of the recomputed first profile may be
estimated with a loss $(\epsilon\Lambda_\epsilon)^{-1}$ per derivative.
The formulae are used for $r\leq C_d/\epsilon$.
On $r\leq r_*$ the assertion is the corresponding smooth Cartesian
estimate, with all centred corrections vanishing at $y=0$.
We use $L(s)=1+|\log s|$ for $0<s\leq8d$.
Expressions involving $L(\epsilon r)$ are estimated on $r\geq r_*$;
on the fixed core their smooth replacements have size
$C\Lambda_\epsilon$.

The slow and radiative rows, including the derivative rule, follow
from Lemma~\ref{fixed-tube-elliptic-batch-lemma} and the bounds
\eqref{slow-elliptic-profile-bounds} and
\eqref{iteratedinnerprofilebounds}.  The first
row requires a slightly more precise separation than a uniform
${\Lambda_\epsilon}^P/r$ bound.  The leading source in the first normal inverse has
phase and amplitude orders $\ell_r^P/r^3$ and $\ell_r^P/r$.  Indeed
$-{\Lambda_\epsilon} V^\perp+2\nabla_z\Phi_\epsilon$ has its constant
$\log(1/\epsilon)$ part cancelled by the corrected velocity.
The remaining leading coefficients are combinations of $1,\log r$
and angular harmonics of degree one.  The Taylor remainder is
$O(\epsilon r(1+|\log(\epsilon r)|))$.

For these two powers, \eqref{first-tail-separated-logarithms} proves
the first row and its normal derivatives at each instantaneous curve,
using the basic seed velocity prescribed before
\eqref{exact-global-inner-factor}.  Its proof uses the global
zeroth-order inverse, the scalar massive estimate and the Poisson
formula \eqref{phase-tail-poisson-formula}; all differentiated
estimates are taken a fixed distance inside the matching radius.
The small phase mean is controlled through the core by
\eqref{core-source-phase-component} and
\eqref{small-phase-zero-mode-primitive}, which give
$O(\epsilon^2\Lambda_\epsilon^P)$ after restoring the first-correction
factor $\epsilon$.
On the fixed core, \eqref{core-logarithm-splitting} isolates the only
leading $\Lambda_\epsilon$ term as the smooth pure phase
$-\sigma_\epsilon$.  Its linear contribution to the modulus is zero;
its quadratic contribution is $O(\epsilon^2\Lambda_\epsilon^P)$.
The compact remainder from the exact core residual has this same
allowed physical order.  These facts justify the separation of radial
logarithms and global logarithmic remainders in the first row.

To transfer the time estimates, write the centred first profile as
$\mathcal F_\epsilon(t,z)$ in the curve chart, with all seed time jets
induced by the basic corrected evolution.  At every required finite
spatial level, \eqref{nextwavesourcebounds} gives
\[
 \sum_{k=2}^{K-1}\|\partial_tg_k\|_{C_tH^s}
       \leq C\epsilon^{5/4}\Lambda_\epsilon^P,\qquad
 \sum_{k=2}^{K-1}\|\partial_{tt}g_k\|_{C_tH^s}
       \leq C\epsilon^{1/4}\Lambda_\epsilon^P.
\]
The second estimate uses one fast derivative of $\partial_tg_k$.
Thus the actual path $z=z_K$ has uniformly bounded $z_t,z_{tt}$,
since $z_*$ is slow.  At fixed material coordinates the chain rule gives
\[
 \frac{d^2}{dt^2}\mathcal F_\epsilon(t,z(t))
 =\mathcal F_{\epsilon,tt}
  +2D_z\mathcal F_{\epsilon,t}[z_t]
  +D_z^2\mathcal F_\epsilon[z_t,z_t]
  +D_z\mathcal F_\epsilon[z_{tt}].
\]
Here partial $t$ derivatives hold $z$ fixed.  Differentiating the
seed source and its normalised planar inverse preserves the radial
bounds in \eqref{first-tail-separated-logarithms}; on the core use
\eqref{core-logarithm-splitting}.  This proves the first-row rule
through two material time derivatives, also with the mixed normal
and longitudinal derivatives at the reserved levels.  Differentiating
the same source and inverse formulas gives the two weighted family
variations, with the reduced indices of Section~\ref{arbitraryapproxsection}.
The same curve and frame estimates bound $Y^{\rm ph}$ and
$D_t^{\rm mat}Y^{\rm ph}$ in every fixed $W_x^{s,\infty}$ norm used here.

\begin{lemma}[Residual reserve]\label{linear-lem-residualreserve}
For any fixed integers $m,j_*$ and $L$, the finite construction may
be chosen so that
\begin{equation}\label{linear-eq-residualreserve}
 \sup_t\sum_{j\leq j_*}
 \| (\partial_t+Y^{\rm ph}\cdot\nabla_x)^j
                S_\epsilon(U_\epsilon)(t)\|_{H_x^m}
                \leq C\epsilon^L.
\end{equation}
The same assertion holds with fixed-Cartesian time derivatives and
after up to two initial-curve differentials, with directions measured
in the chosen initial-curve Sobolev space.
After rescaling, all the finite residual seminorms used below can
therefore be made $O(\epsilon^{20})$.
\end{lemma}
\begin{proof}
Put $D=\partial_t+Y^{\rm ph}\cdot\nabla_x$.  Derivatives at fixed
$(x,\tau)$, including family derivatives, commute with the wave
operator.  Material derivatives are obtained by the chain rule;
in particular,
\[
 D^2f=f_{tt}+2Y^{\rm ph}\cdot\nabla_xf_t
       +(Y^{\rm ph}\otimes Y^{\rm ph}):D_x^2f
       +(DY^{\rm ph})\cdot\nabla_xf.
\]
Spatial differentiation of this identity uses the coefficient bounds
just proved.  For higher time orders, \eqref{fixed-base-fast-time-curve-bound}
and \eqref{nextwavesourcebounds} control all required fast-time
coefficient derivatives.  Repeated Leibniz expansion then costs at
most one $\epsilon^{-1}$ and a logarithmic power per physical time
derivative.  A physical normal derivative or a derivative of a core
cutoff has the same cost.  Spatial derivatives of fixed physical
collars are bounded; their time derivatives obey this chain rule.
The weighted family derivatives \eqref{weighted-initial-curve-derivative}
have bounded tubular conversion coefficients, since $\epsilon r$
is bounded.

Differentiate the exact terminal residual identities, retaining the
radial estimates of each batch.  The calculation proving
\eqref{terminal-physical-residual-bound} now gives, for $j\leq j_*$
and $q\leq2$,
\[
 \sup_t\|\mathcal D_{\rm in}^q D^jS_\epsilon(U_\epsilon)\|_{H_x^m}
 \leq C\Lambda_\epsilon^P
   \big(\epsilon^{-m-j-2}A_K
       +\epsilon^{J_K-m-j-2}A_{K-1}\big).
\]
The same estimate holds with $D^j$ replaced by $\partial_t^j$.
The two displayed extra powers cover the smooth normalisation and
coordinate terms as in \eqref{terminal-physical-residual-bound}.
In the far field, place a differentiated wave gradient or amplitude
in its global Sobolev norm and the remaining factors in their
pointwise norms; the dipole residual retains its integrable weight.

Choose $N\geq m+j_*+4$ and $M\geq\max\{20,L\}$ in
\eqref{explicit-number-of-stages}--\eqref{explicit-stage-orders}.
The terminal level $m_K=N+20$ includes the product, Sobolev and two
family reserves.  The four-derivative source-to-wave loss and the
planar reserve \eqref{elliptic-derivative-reserve} require at most
$8^{J_k+5}(m_k+40)$ at the preceding stage, below $m_{k-1}$ in
\eqref{explicit-stage-orders}.  Thus these are estimates for the
same finite construction.  Removing the family weights costs at
most $\epsilon^{-2}$.  Since $A_K=\epsilon^{M+N+85/4}$, both terms
give \eqref{linear-eq-residualreserve} after absorbing logarithms.
Rescaling uses
$\|D_X^\alpha f(\epsilon X)\|_2
 =\epsilon^{|\alpha|-3/2}\|D_x^\alpha f\|_2$.

For the residual derivatives used below, one sufficient choice is $m=6$, $j_*=2$,
$M=20$ and $N=12$.  Then $K=210$, $m_K=32$ and $J_K=52$;
the first term above, after two family differentials and $L^2$
rescaling, is $O(\epsilon^{159/4}\Lambda_\epsilon^P)$, and the
second is smaller.  Sobolev embedding supplies four pointwise spatial
derivatives.  All required residual seminorms are therefore
$O(\epsilon^{20})$.  Constants and the smallness threshold may depend
on these fixed construction orders.
\end{proof}

\subsubsection{Time derivatives of the algebraic amplitude}
The cumulative outer phase and amplitude are real functions $A,B$
and satisfy
\begin{equation*}
 -A_{\tau\tau}+\Delta_xA=F_{\rm w},\qquad
 B=\tfrac12H_2+\frac{\epsilon}{\sqrt2}A_\tau.
\end{equation*}
These identities are obtained by summing the stage equations;
$H_2$ denotes the sum of their imaginary source components.
Direct differentiation gives
\begin{align}
 B_t&=\frac{H_{2,\tau}}{\sqrt2\epsilon{\Lambda_\epsilon}}
                      +\frac{\Delta_xA-F_{\rm w}}{{\Lambda_\epsilon}},
                  \label{linear-eq-Bt}\\
 B_{tt}&=\frac{H_{2,\tau\tau}}{\epsilon^2{\Lambda_\epsilon}^2}
           +\frac{\sqrt2}{\epsilon{\Lambda_\epsilon}^2}
                        (\Delta_xA_\tau-F_{{\rm w},\tau}).
                  \label{linear-eq-Btt}
\end{align}
\begin{lemma}[Amplitude derivatives]\label{linear-lem-amplitudetime}
At every fixed finite spatial level required below,
\begin{equation}\label{linear-eq-amplitudetimebounds}
 \|B_t\|_{W_x^{m,\infty}}\leq C\epsilon^{5/4}{\Lambda_\epsilon}^P,
 \qquad
 \|B_{tt}\|_{W_x^{m,\infty}}\leq C\epsilon^{1/4}{\Lambda_\epsilon}^P.
\end{equation}
The same bounds hold for the first two material derivatives in the
tube.  Also $\|B\|_{W_x^{m,\infty}}\leq C\epsilon^2{\Lambda_\epsilon}^P$.
\end{lemma}
\begin{proof}
The bound \eqref{outerinductivebounds} at wave level $m+6$, the
definition \eqref{threeDwaveXnorm} and $H_x^2\hookrightarrow L_x^\infty$ give
$\|A_\tau\|_{W^{m+2,\infty}}+
\|\nabla_xA\|_{W^{m+2,\infty}}\leq C\epsilon^{5/4}{\Lambda_\epsilon}^P$.
This level includes two weighted family variations.
The initial wave source and its first fast derivative are
$O(\epsilon^2{\Lambda_\epsilon}^P)$; subsequent sources are
$O(\epsilon A_k{\Lambda_\epsilon}^P)$.  The later $H_{2,k}$ and their two fast
derivatives are $O(\epsilon^2A_{k-1}{\Lambda_\epsilon}^P)$.
These are $W_x^{m,\infty}$ bounds: use the $Y_{m+6}$ and
$\mathcal H_{m+6}$ source estimates in \eqref{nextwavesourcebounds},
with the slow source estimates \eqref{slow-seed-wave-source-bound}
at the first stage.

For the initial slow $H_{2,2}$, let $a=|\alpha|$, $j\leq2$ and
$a+j\leq m+2$.
The slow elliptic remainder and the tubular chain rule, including
the core cutoff, give on the fixed tube
\[
 |\partial_t^jD_x^\alpha H_{2,2}|
 \leq C\Lambda_\epsilon^P
       \epsilon^{J+1-a-j}\langle r\rangle^{J-1-a-j}
 =C\epsilon^2\Lambda_\epsilon^P
                   (\epsilon\langle r\rangle)^{J-1-a-j}.
\]
The source vanishes on the core, and $J\geq m+5$ suffices also
for the two weighted family variations.  Outside the tube the slow
phase expansion gives the same $W_x^{m,\infty}$ bound.  Hence
$H_{2,2,\tau}=O(\epsilon^3{\Lambda_\epsilon}^P)$ and
$H_{2,2,\tau\tau}=O(\epsilon^4{\Lambda_\epsilon}^P)$.
Insert these estimates in \eqref{linear-eq-Bt}--\eqref{linear-eq-Btt} and sum the
finite geometric sequence $A_k$.  The largest terms are respectively
$\epsilon^{5/4}{\Lambda_\epsilon}^P$ and $\epsilon^{1/4}{\Lambda_\epsilon}^P$.
For material derivatives use the expansion in the proof of
Lemma~\ref{linear-lem-residualreserve}.  The additional terms are
bounded in $W_x^{m,\infty}$ by $C\epsilon^2\Lambda_\epsilon^P$ at
the first derivative and by $C\epsilon^{5/4}\Lambda_\epsilon^P$
at the second.  The same source and wave levels control $B$ in
$W_x^{m+2,\infty}$ and $B_t$ in $W_x^{m+1,\infty}$, as required
by that expansion.
\end{proof}

\subsubsection{Low-order coefficient bounds}
The finite profile records, the phase expansion and
Lemma~\ref{linear-lem-amplitudetime} imply
\begin{equation}\label{linear-eq-lowcoeff}
 \begin{gathered}
 \rho\asymp w\quad\text{in the tube},\qquad
 |v-W|\leq C(\epsilon{\Lambda_\epsilon}^P/\langle r\rangle+\epsilon^2{\Lambda_\epsilon}^P),\\
 |\nabla_X\Theta_{\rm reg}|\leq C\epsilon{\Lambda_\epsilon}^P,
 \qquad |{c_\epsilon}\vartheta|\leq C\epsilon^2{\Lambda_\epsilon}^P,\\
 |{c_\epsilon} Y|\leq C\epsilon{\Lambda_\epsilon},\qquad
 \|Y^{\rm ph}\|_{W^{4,\infty}}+
 \|D_t^{\rm mat}Y^{\rm ph}\|_{W^{2,\infty}}\leq C.
 \end{gathered}
\end{equation}
The first line on a fixed core disc follows from the centred $C^1$
estimate: both $v$ and $W$ vanish simply, so their modulus ratio is
bounded above and below.  Off that disc use the logarithmic ansatz
and \eqref{linear-eq-profiletable}.  Every slow higher profile is at most
$C\epsilon^2{\Lambda_\epsilon}^P$ on the fixed tube, and each radiative phase is at
most $C\epsilon A_k{\Lambda_\epsilon}^P$.  The real outer amplitude is retained in
$\rho$.  The phase estimates follow by differentiating the
regularised solid-angle expansion; its normal derivative is
$O(\epsilon L(\epsilon r))$ outside the fixed core and
$O(\epsilon{\Lambda_\epsilon})$ on it.  The smooth longitudinal phase at the centre
has bounded physical derivatives.  The curve and frame bounds follow
from Section~\ref{normalgraphsection} and the slow and radiative
curve-increment equations \eqref{fixed-base-curve-increment} and
\eqref{higher-translation-Newton-velocity}.  One additional physical time derivative of a
radiative velocity costs $(\epsilon{\Lambda_\epsilon})^{-1}$, but its size is
$O(A_k{\Lambda_\epsilon}^P)$, leaving $O(\epsilon^{1/4}{\Lambda_\epsilon}^P)$.

On the complement of any fixed positive physical tube,
\begin{equation*}
 \rho\geq\tfrac12,\quad |\rho^2-1|\leq C\epsilon^2{\Lambda_\epsilon}^P,
 \quad |\nabla_X\Theta|\leq C\epsilon{\Lambda_\epsilon}^P.
\end{equation*}
In particular $V$ has no zero there.  The same finite estimates and
\eqref{linear-eq-defects} give
\begin{equation}\label{linear-eq-edefect}
 |s|+\langle r\rangle|e|\leq C{c_\epsilon}\quad (r\geq r_*),
 \qquad |\mathcal V|\leq Cw\quad(r\leq r_*).
\end{equation}
For clarity, the slow contribution to ${c_\epsilon}{D_t^{\rm mat}}\alpha$ is bounded by
$C{\Lambda_\epsilon}^P(\epsilon^3/\langle r\rangle+\epsilon^4)$.
A radiative term of index $n$ contributes
$C\epsilon^{n+1}A_k\langle r\rangle^{n-1}{\Lambda_\epsilon}^P$.
Multiplication by $\langle r\rangle\leq C/\epsilon$ bounds the
latter by $C\epsilon A_k{\Lambda_\epsilon}^P=o({c_\epsilon})$.  The leading radial
radiative correction has no phase.  The residual is covered by
Lemma~\ref{linear-lem-residualreserve}.  The density defect uses
$|{D_t^{\rm mat}}\log\rho^2|+|\operatorname{div} Y|\leq C$, since ${D_t^{\rm mat}} w=0$.
This proves \eqref{linear-eq-edefect} without a sign assumption.

\subsection{The unaugmented form on adjoining regions}
\label{linear-sec-regions}

We now combine the planar estimates near the filament with estimates
in the exterior, keeping track of the boundary terms at their interface.
This gives the lower bound for the unaugmented quadratic form needed
in the construction of the positive quadratic functional.

\subsubsection{A phase-preserving core identification}
There is a smooth invertible real-linear matrix $\mathcal C$ with
\begin{equation}\label{linear-eq-coreC}
 \mathcal C(iW)=iv,\qquad
 \mathcal Cz=(v/W)z\quad(r\geq r_*).
\end{equation}
To construct it near zero, write the smooth centred field $iv(y)$
as $F(iW(y))$ on a disc where $y\mapsto iW(y)$ is a diffeomorphism.
Then
\[
 \mathcal C(y)=\int_0^1DF(\lambda iW(y))\,d\lambda
\]
has the first property.  Interpolate radially on a fixed annulus to
multiplication by $v/W$.  Both matrices have the first property and
are close to the identity there; their interpolation is invertible.
The finite differentiated profile estimates give bounded derivatives
of $\mathcal C^{\pm1}$.  We use this matrix only to identify test fields.

For a field $h$ in the physical $4d$ tube put
\begin{equation}\label{linear-eq-eta}
 \eta=\mathcal C^{-1}e^{-i\Theta_{\rm reg}}h\circ\mathcal F_t
       =iW(f+ig)+\eta^1+\eta^{\geq2}.
\end{equation}
Then $f,g$ are real radial functions of $(\sigma_0,r)$ and the pure
zero-phase field is exactly $iVf$.  Let $\eta_+=-Wg+\eta^1+\eta^{\geq2}$.
The corresponding physical field is denoted by $h_+$.  Thus
$h=iVf+h_+$.

\subsubsection{Interface terms}
Set $R_I=4d/\epsilon$, denote the open tube $r<R_I$ in $X$ by
$\mathcal D_t$, and let $\Omega_t$ be its exterior.  For a subdomain
$\mathcal D$ away from the zero, the domain version of
\eqref{linear-eq-fullphaseform} has boundary term
\begin{equation*}
 \mathfrak F_{\partial\mathcal D}[h]
 =\int_{\partial\mathcal D}\left\{
      \frac{\partial_\nu\rho}{\rho}|h|^2
                      +\rho^2(D\cdot\nu)fg\right\}dS,
       \qquad h=iV(f+ig).
\end{equation*}
The first term comes from integrating the modulus; the second
comes from integrating the current.  Define $b_{\rm in}^\circ$
by subtracting this flux from the original form on the tube and
$b_{\rm out}^\circ$ by adding the same flux, with the normal pointing
out of the tube.  Then exactly
\begin{equation}\label{linear-eq-interface}
 \mathfrak b_t[h]=b_{\rm in}^\circ[h]+b_{\rm out}^\circ[h].
\end{equation}
The two regional forms are evaluated on the same field $h$, so their
interface fluxes cancel.  At the zero, define the inner form from its
original smooth Cartesian expression.  All uses of phase variables can
first be made away from the zero and passed to that expression.

The signed zero-phase coefficient is
\begin{equation}\label{linear-eq-Kphase}
 K_0(\sigma_0,r)=\int_0^{2\pi}j_0\rho^2e\,d\theta
 =-\int_0^{2\pi}j_0
       \{{c_\epsilon}\operatorname{Im}(\overline V\mathcal V)+\operatorname{Re}(\mathscr R\overline V)\}\,d\theta.
\end{equation}
Its second expression is smooth at the centre.  By
\eqref{linear-eq-edefect} and $\rho\asymp w$,
\begin{equation*}
 |K_0|\leq C{c_\epsilon} w^2/\langle r\rangle
                              +C\delta_{\rm res} w\,\mathbf1_{r<r_*},
 \qquad \delta_{\rm res}=\sup_t\|\mathscr R(t)\|_\infty.
\end{equation*}
Lemma~\ref{linear-lem-radial}, followed by \eqref{linear-eq-coreHardy}, gives
\begin{align}
 \int[K_0]_-f^2r\,dr\,d\sigma_0
 &\leq\eta E_{\rm ph,in}(f)+C_\eta{c_\epsilon}^2\|iVf\|_2^2\notag\\
 &\quad+C{c_\epsilon} R_I\int f(\sigma_0,R_I)^2\,d\sigma_0.
                    \label{linear-eq-phaseinterfacebound}
\end{align}
Here the residual contribution is absorbed using
$\delta_{\rm res}\ll{c_\epsilon}^2$.  At the interface, \eqref{linear-eq-coreC} is complex
multiplication, so $f$ is the angular average of the exterior phase.
Angular Jensen and \eqref{linear-eq-exteriortrace} imply
\begin{equation}\label{linear-eq-smalltrace}
 {c_\epsilon} R_I\int f(R_I)^2d\sigma_0
 \leq C_d\frac{{c_\epsilon}}{\epsilon}
                         \int_{\Omega_t}|\nabla\psi_1|^2
 =o(1)\int_{\Omega_t}|\nabla\psi_1|^2.
\end{equation}
The positive longitudinal term in the pure phase is
$\int j_0\rho^2(a_\ell H)^{-2}|f_{\sigma_0}|^2$.

\subsubsection{The exterior and its transition region}
On $\Omega_t$, write $h=iV\psi$ and use
$\rho\geq1/2$, $|D|=o(1)$ and $|s|\leq C{c_\epsilon}$.  Young's
inequality in \eqref{linear-eq-fullphaseform} gives
\begin{equation*}
 b_{\rm out}^\circ[h]\geq k\int_{\Omega_t}\rho^2
       (|\nabla\psi|^2+\rho^2\psi_2^2)
       -C\int_{\Omega_t}\rho^2(|e|+s^2)\psi_1^2.
\end{equation*}
Beyond the inner corrections and ambient velocity, $\alpha=0$ and
\[
 e=-\rho^{-2}\operatorname{Re}(\mathscr R\overline V).
\]
In the bounded physical transition region, the first and slow inner phases are
$O(\epsilon^2{\Lambda_\epsilon}^P)$ with slow time derivatives.  A radiative
phase is $O(\epsilon A_k{\Lambda_\epsilon}^P)$, and its material derivative is
$O(A_k{\Lambda_\epsilon}^P)$.  Therefore
\[
 |{c_\epsilon}{D_t^{\rm mat}}\alpha|\leq a_\epsilon,
 \qquad a_\epsilon=C{\Lambda_\epsilon}^P(\epsilon^4+
                                     \epsilon^2\sum_kA_k).
\]
This coefficient has bounded physical support, and
\[
 \|{c_\epsilon}{D_t^{\rm mat}}\alpha\|_{L_X^{3/2}}\leq C\epsilon^{-2}a_\epsilon=o(1).
\]
The exterior Sobolev inequality in \eqref{linear-eq-exteriortrace} absorbs
this part into the gradient.  The residual and $s^2$ cost only
$C{c_\epsilon}^2\|h\|_2^2$.  Thus
\begin{equation}\label{linear-eq-outbound}
 b_{\rm out}^\circ[h]\geq k E_{\rm out}(h)
                               -C{c_\epsilon}^2\|h\|_{L^2(\Omega_t)}^2,
 \quad
 E_{\rm out}(h)=\int_{\Omega_t}\rho^2
                (|\nabla\psi|^2+\rho^2\psi_2^2).
\end{equation}

\subsubsection{The modewise comparison and its norms}
We use the following precise regional quantities.  The normal disc
norms are truncated at $R_I$ and the phase radial integral is also
truncated there.  In the coordinates \eqref{linear-eq-eta}, put
\begin{align*}
 \mathcal N_*(h)&=\int\{E_0(f,g;R_I)
                   +\|\eta^1+\eta^{\geq2}\|_{\mathcal H_W(B_{R_I})}^2\}
                           d\sigma_0+E_{\rm out}(h),
                     \\
 \mathcal K(h)&=\int_{r<R_I}|\partial_{\sigma_0}\eta|^2dy\,d\sigma_0,
                     \\
 \mathcal X(h)^2&=\mathcal N_*(h)+\mathcal K(h)+{c_\epsilon}^2\|h\|_2^2.
\end{align*}
They are regarded as quantities of $q$ by $h=\mathscr S_t^{-1}q$.
Uniform bounded factors $a_\ell,H,\rho/w$ can be inserted without
changing the estimates.  The outer energy also controls
$\int|f(R_I)|^2d\sigma_0$ by \eqref{linear-eq-exteriortrace}.
These definitions contain the translation directions in their
natural planar norm; only the zero phase is exempt from the
$\langle r\rangle^{-2}$ mass weight.

The following regional comparison is useful in both the static and
the differentiated estimates.
\begin{lemma}[Regional continuity and fixed-core coercivity]
\label{linear-lem-regional}
The form is continuous in $\mathcal X$:
\begin{equation}\label{linear-eq-Xcontinuity}
 |B_t[q_1,q_2]|\leq C\mathcal X(q_1)\mathcal X(q_2).
\end{equation}
On fields whose fixed-core translation measurements vanish,
\begin{equation}\label{linear-eq-fixedcorepositive}
 B_t[q]\geq k(\mathcal N_*(q)+\mathcal K(q))
                                             -C{c_\epsilon}^2\|q\|_2^2.
\end{equation}
In a decomposition involving translations, the longitudinal square
in this estimate may instead be retained as the square of the
\emph{complete} field; the transverse cross terms are the ones
estimated in Subsection~\ref{linear-sec-jacobi}.
\end{lemma}
\begin{proof}
On $r<r_*$ use the smooth Cartesian form and the matrix
$\mathcal C$.  Its coefficients differ from those of the planar
form by $O(\epsilon{\Lambda_\epsilon}^P)$.  The local zero-phase mass is bounded
by \eqref{linear-eq-phasecoremass}; hence these errors are small relative
to the natural phase energy.  The other local components are
controlled by their $\mathcal H_W$ norm.

For $r\geq r_*$, $h=iV\psi$ exactly.  The nonzero phase has zero
angular mean.  Thus $\psi_1/r$, $\nabla\psi_1$, $\psi_2$ and
$\nabla\psi_2$ are controlled by its planar norm.  The perturbations
of the radial planar coefficients are $o(1)$; the drift
$D-2\nabla_X\theta$ is $o(1)$ and pairs a phase derivative with the
amplitude.  The errors $e$ and $s$ are controlled by
\[
 \sup r^2|e|+\sup r|s|\leq C\epsilon{\Lambda_\epsilon}=o(1).
\]
The metric errors are $O(d)$, and the periodic angular connection
is controlled using $\epsilon r\leq4d$.  On the zero phase the
singular angular drift is orthogonal to its radial derivative and
annihilates its angular derivative.  Consequently its diagonal is
exactly its positive gradient plus \eqref{linear-eq-Kphase}.  Its mixed
pairing with the complement follows from
\eqref{linear-eq-phaseJacobi} and one integration by parts: the leading
planar pairing vanishes, the derivative terms are $O(d)+o(1)$ times
the product of the two positive energies, and the remaining
coefficients have amplitude size $Cs$ and nonzero-phase size $Ce$.
Young's inequality costs $C{c_\epsilon}^2\|h\|_2^2$ by
\eqref{linear-eq-edefect}.  In the fixed core use the smooth first
expression in \eqref{linear-eq-phaseJacobi} and \eqref{linear-eq-coreHardy}.

Apply Lemma~\ref{linear-lem-disk}, the phase estimate
\eqref{linear-eq-phaseinterfacebound}, and the exterior estimate
\eqref{linear-eq-outbound} in \eqref{linear-eq-interface}.  The interface traces
are controlled by \eqref{linear-eq-smalltrace}.  Choosing $d$ first and
then $\epsilon$ makes the error fractions absorbable.  This proves
\eqref{linear-eq-fixedcorepositive}.  The same comparisons in absolute
value, using polarisation of the radial estimate, prove
\eqref{linear-eq-Xcontinuity}.  In the longitudinal terms use the exact
phase--amplitude expression and keep the full square.  Connection
terms on nonzero modes are bounded by the small angular fraction;
on the pure phase use its exact identity before estimating.  This
justifies the final assertion.
\end{proof}

\subsection{Symplectic translations and the two projections}
\label{linear-sec-translations}

We next define the translation components that must be controlled to
obtain a positive quadratic functional.  We compare the projection
used in that functional with an auxiliary projection suited to the
planar coercivity estimates.

\subsubsection{The translation fields}
Let $\chi_R(y)=\chi(r/R)$, $R=d/\epsilon$, where $\chi=1$ on
$[0,1]$ and vanishes on $[2,\infty)$.  All its support is inside the
material-flow tube.  Set
\begin{equation*}
 N_a=e_a\cdot\nabla_X=\partial_{y_a},\qquad
 T_a^0=\chi_RN_aV.
\end{equation*}
At each material position and radius define
\begin{equation}\label{linear-eq-ampadjust}
 a_a=\frac{\int_0^{2\pi}H\operatorname{Re}(\overline VT_a^0)\,d\theta}
                {\int_0^{2\pi}H\rho^2\,d\theta},
 \qquad T_a=T_a^0-a_aV.
\end{equation}
The two occurrences of $a$ in this definition are indices; the
functions $a_1,a_2$ are real radial functions.  The numerator and
denominator vanish to order $r^2$ at zero.  Their angular averages
are smooth in $r^2$, and the denominator has a positive leading
coefficient.  Hence $a_a$ extends smoothly there.  It vanishes
outside the support of $\chi_R$.  By construction,
\begin{equation}\label{linear-eq-phaseorthogonality}
 \int_0^{2\pi}H\operatorname{Re}(\overline VT_a)\,d\theta=0.
\end{equation}
For $\rho=w$ the formula is explicitly
\begin{equation}\label{linear-eq-ampradialleading}
 a_a=-\frac{\epsilon\chi_R\kappa_a r w'}{2w}.
\end{equation}
This shows its $O(\epsilon)$ core size and $O(\epsilon r^{-2})$
tail.  The differentiated modulus estimates of
Subsection~\ref{linear-sec-backgroundbounds} give the same orders with
$r$-logarithms for the first profile and smaller higher terms.

In the fixed material representation put
\begin{equation}\label{linear-eq-tau}
 \tau_a=\sqrt{j_0}\,e^{-i\Theta_{\rm reg}}T_a\circ\mathcal F_t,
 \qquad U_{\rm mat}=\sqrt{j_0}\,v.
\end{equation}
For a real vector function $b=b(\sigma_0)$, $\mathscr T_tb$ denotes
the global reference field whose material expression is $\tau_ab_a$;
it vanishes off the support of the translations.  Equation
\eqref{linear-eq-phaseorthogonality} becomes
$\langle iU_{\rm mat}f,i\tau_a\rangle_y=0$ for every real radial $f$.

\subsubsection{Gram matrices and logarithmic norm}
Define
\begin{equation}\label{linear-eq-sympprojection}
 p_i(q)=\langle z,i\tau_i\rangle_y,\qquad
 N_{ij}=\langle \tau_j,i\tau_i\rangle_y,\qquad b=N^{-1}p(q).
\end{equation}
These are real-linear functionals.  The matrix $N$ is
antisymmetric and has the form $N=\beta_R\mathbb J$, where
$\mathbb J=\left(\begin{smallmatrix}0&1\\-1&0\end{smallmatrix}\right)$.
In the planar case
\begin{equation*}
 \beta_R=2\pi\int_0^\infty\chi_R^2ww'\,dr
                  =\pi+O(R^{-2}).
\end{equation*}
For the actual fields its leading term is $a_\ell\pi$, with an
$o(1)$ error.  To check this, in the unadjusted determinant expand
$\nabla\rho$ and $\nabla\Theta$.  The leading term is $ww'/r$;
terms with a regular phase gradient have an extra $\epsilon$
and integrable $w'$; differentiated amplitude tails give
$O(\epsilon\ell_r^P/r^3)$ after the area factor is included.
The $H-1$ first angular harmonic has zero leading average.  Terms
with two corrections and the adjustment \eqref{linear-eq-ampadjust} tend
to zero by the same radial integrals.  Thus
\begin{equation}\label{linear-eq-Gramest}
 \|N\|+\|N^{-1}\|\leq C,\qquad
 |N_t|\leq C,\qquad |N_{\sigma_0}|\leq C\epsilon.
\end{equation}
The time estimate uses the common-real-factor identity in
Lemma~\ref{linear-lem-connection}; alternatively it follows by
differentiating the determinant before estimating the tails.
The slow/radiative derivative rules give the last two bounds.

The tests satisfy $\|\tau_a\|_{L_y^2}\leq C\sqrt{\Lambda_\epsilon}$.
Their natural transverse norm is bounded.  The regular phase tail
$\epsilon L(\epsilon r)$ has a mode-zero part.  Its phase energy
uses its gradient, not its quotient by $r$.  In particular,
\begin{equation}\label{linear-eq-liftbounds}
 \|\mathscr T_t b\|_2^2\leq C{\Lambda_\epsilon}\|b\|_2^2,
 \qquad \mathcal N_*(\mathscr T_tb)\leq C\|b\|_2^2.
\end{equation}
For the measurement of a nonzero phase--amplitude field, the only
borderline integral is the amplitude paired with the translation
phase tail.  Indeed, for $Z_1$,
\[
 \langle iW\psi,i\chi_RZ_1\rangle_y
 =\int\chi_R\left(ww'\psi_1\cos\theta
                     -\frac{w^2}{r}\psi_2\sin\theta\right)dy.
\]
The second term has dual norm $C\sqrt{1+\log R}$, since
$\int_1^{2R}dr/r=\log(2R)$.  The first is integrable without a
logarithm.  The corrected projection annihilates the complete zero
phase exactly.  Its other coefficient differences have the bounded
norms just described.  Hence the same single logarithm applies to
\eqref{linear-eq-sympprojection}.

\subsubsection{A fixed-core projection used only in the proof}
Apply the fixed-core Cartesian measurements of
Lemma~\ref{linear-lem-planarcoercivity} to $\eta$ in \eqref{linear-eq-eta}.
Let $M_*(\sigma_0,t)$ be their matrix on the fields $T_a$ under that
identification.  It is a small perturbation of the fixed positive
matrix $M$, and is uniformly invertible.  Define
\begin{equation}\label{linear-eq-auxprojection}
 \beta=M_*^{-1}\mathcal L_*(\eta),\qquad
 w_*=q-\mathscr T_t\beta.
\end{equation}
Thus $w_*$ has zero fixed-core translation measurements and
$\beta(\mathscr T_ta)=a$.  Both projections annihilate $iVf$ for
radial real $f$.  In particular no orthogonality condition on the
full field has been imposed.

Set
\[
 E_*(w_*)=\mathcal N_*(w_*),\qquad K(q)=\mathcal K(q).
\]
The $\sigma_0$ derivative of the auxiliary coefficients uses the
fixed measurement of $\eta_{\sigma_0}$ and the derivative of the
finite matrix.  It follows that
\begin{equation}\label{linear-eq-betaderivative}
 \|\beta_{\sigma_0}\|_2^2
                  \leq C K(q)+C\epsilon^2\|\beta\|_2^2.
\end{equation}
The single logarithmic measurement estimate, applied to $w_*$,
gives
\begin{equation}\label{linear-eq-twoprojections}
 \|\beta-b\|_2^2\leq C{\Lambda_\epsilon} E_*(w_*).
\end{equation}
It follows from \eqref{linear-eq-liftbounds}, the derivative versions of
those bounds and \eqref{linear-eq-betaderivative} that
\begin{equation}\label{linear-eq-auxX}
 \mathcal X(w_*)^2
 \leq C\{E_*(w_*)+{\Lambda_\epsilon} K(q)+{c_\epsilon}\|\beta\|_2^2
                                      +{c_\epsilon}^2\|q\|_2^2\}.
\end{equation}
To see the longitudinal term, differentiate
$w_*=q-\mathscr T_t\beta$ once.  The coefficient of
$\beta_{\sigma_0}$ has $L_y^2$ norm $C\sqrt{\Lambda_\epsilon}$; the derivative
of the lift has norm $C\epsilon\sqrt{\Lambda_\epsilon}$.  Thus its contribution
is bounded by $C{\Lambda_\epsilon}\|\beta_{\sigma_0}\|_2^2+
C\epsilon^2{\Lambda_\epsilon}\|\beta\|_2^2$.  The same calculation in the
phase-preserving identification eliminates an undifferentiated
zero-phase term.  The mass difference costs
${c_\epsilon}^2{\Lambda_\epsilon}\|\beta\|_2^2=o({c_\epsilon})\|\beta\|_2^2$.
These observations prove \eqref{linear-eq-auxX}.  They explain the use of
the auxiliary projection: applying two large-radius $L^2$ bounds
to the symplectic remainder would introduce an unnecessary second
logarithm.

\subsection{The exact coefficient derivative}
\label{linear-sec-derivative}

We now turn to the time variation of the quadratic form.  The
cancellations in its coefficients give the required estimate on the
remainder left by the auxiliary translation projection.

\subsubsection{Eulerian and phase--amplitude identities}
In differentiating $B_t[q]$, the reference field $q$ is fixed.  Thus
$h=\mathscr S_t^{-1}q$ satisfies
\[
 h_t=-Y\cdot\nabla h-\tfrac12(\operatorname{div} Y)h+i\vartheta h.
\]
Differentiation of \eqref{linear-eq-eulerianform}, with integration by
parts, gives the smooth-field identity
\begin{equation}\label{linear-eq-exactderivative}
\begin{split}
 \dot B_t[q]={}&-2\int(\partial_iY_j)
                       \operatorname{Re}(\partial_ih\,\overline{\partial_jh})\\
 &+\int\left(\tfrac12\Delta\operatorname{div} Y+{c_\epsilon}\vartheta_t
                       +2\operatorname{Re}(\overline V\mathcal V)\right)|h|^2\\
 &+\int(2\nabla\vartheta-{c_\epsilon} Y_t)\cdot j(h)
      +4\int\operatorname{Re}(\overline Vh)\operatorname{Re}(\overline{\mathcal V}h).
\end{split}
\end{equation}
For example, the derivative of the kinetic term is
\[
 -2\int(\partial_iY_j)\operatorname{Re}(\partial_ih\overline{\partial_jh})
 +\tfrac12\int\Delta\operatorname{div} Y\,|h|^2
 +2\int\nabla\vartheta\cdot j(h).
\]
The derivative of $\int Y\cdot j(h)$ is
$\int Y_t\cdot j(h)+\int Y\cdot\nabla\vartheta|h|^2$.
Its second term cancels the transport part of
${c_\epsilon}{D_t^{\rm mat}}\vartheta$.  The rank-one potential gives the last term
of \eqref{linear-eq-exactderivative}.  This proves the formula without
any division by the background.

To expose the phase cancellations put, off the zero,
\begin{equation*}
 \begin{split}
 r_1+ir_2&=\mathscr R/V,\qquad v_\Theta=\nabla\Theta,\\
 \lambda&={c_\epsilon}(\vartheta-{D_t^{\rm mat}}\Theta)=-{c_\epsilon}{D_t^{\rm mat}}\alpha,\qquad
 k={D_t^{\rm mat}}\log\rho^2+\operatorname{div} Y,\\
 F_0&={c_\epsilon} Y_t-2\nabla\vartheta+4S(Y)v_\Theta,\\
 E_0^{\rm c}&={D_t^{\rm mat}}(\lambda-r_1)+{c_\epsilon}^{-1}D\cdot\nabla\lambda,\\
 L_0&=\rho^{-2}\operatorname{div}(\rho^2F_0)+4\rho^2\lambda/{c_\epsilon}.
 \end{split}
\end{equation*}
The superscript on $E_0^{\rm c}$ distinguishes this coefficient
from the planar zero-mode energy.  For $h=iV\psi$,
\begin{equation}\label{linear-eq-psitransport}
 {D_t^{\rm mat}}\psi=-\tfrac12k\psi+i\lambda\psi/{c_\epsilon}
\end{equation}
when $q$ is fixed.  The complete derivative is
\begin{equation}\label{linear-eq-sixterms}
\begin{split}
 \dot B_t[\mathscr S_t(iV\psi)]={}&
 -2\int\rho^2S(Y)_{ij}\operatorname{Re}(\partial_i\psi\overline{\partial_j\psi})\\
 &+2\int\rho^2F_0\cdot\nabla\psi_1\,\psi_2
                       +\int\rho^2L_0\psi_1\psi_2\\
 &+\int\rho^2E_0^{\rm c}|\psi|^2
                -\tfrac12\int\rho^2\nabla k\cdot\nabla|\psi|^2\\
 &+2\int\rho^2{D_t^{\rm mat}}(\rho^2)\psi_2^2.
\end{split}
\end{equation}
This rearrangement is the dimension-independent form of the identity in
\cite[Lemma~7.3]{delpinojunemanmusso2025}.  Here it follows from the
canonical moving form, including its density.

For completeness, use the non-integrated current version of
\eqref{linear-eq-fullphaseform}, whose density is
\[
 \rho^2\{ |\nabla\psi|^2+D\cdot\operatorname{Im}(\overline\psi\nabla\psi)
                   +2\rho^2\psi_2^2+(\lambda-r_1)|\psi|^2\}.
\]
In its material derivative, the density contribution to the
kinetic term cancels the real scalar part of
\eqref{linear-eq-psitransport}.  The imaginary part gives a phase-gradient
term that cancels the corresponding current derivative.  Integrating
the remaining current once gives
\[
 -\int\rho^2F_0\cdot\operatorname{Im}(\overline\psi\nabla\psi)
  =2\int\rho^2F_0\cdot\nabla\psi_1\,\psi_2
                  +\int\operatorname{div}(\rho^2F_0)\psi_1\psi_2.
\]
The amplitude potential gives the last term in
\eqref{linear-eq-sixterms} and $4\rho^4\lambda\psi_1\psi_2/{c_\epsilon}$.
The remaining scalar terms give $E_0^{\rm c}$ and the gradient of
$k$.  This proves \eqref{linear-eq-sixterms}.

The residual continuity identity simplifies its coefficients to
\begin{equation}\label{linear-eq-structuredcoefficients}
\begin{split}
 E_0^{\rm c}&=-{D_t^{\rm mat}} r_1-{c_\epsilon}({D_t^{\rm mat}})^2\alpha-D\cdot\nabla({D_t^{\rm mat}}\alpha),\\
  L_0&={c_\epsilon}{D_t^{\rm mat}} k+D\cdot\nabla k-2{D_t^{\rm mat}} r_2\\
     &\quad+2\rho^{-2}\operatorname{div}(\rho^2\nabla({D_t^{\rm mat}}\alpha))
                           -4\rho^2{D_t^{\rm mat}}\alpha.
\end{split}
\end{equation}
For instance
${c_\epsilon}\partial_t\rho^2+2\operatorname{div}(\rho^2v_\Theta)=2\rho^2r_2$.
Substitution in $\rho^{-2}\operatorname{div}(\rho^2F_0)$ cancels the terms
$2r_2k$ arising from the density derivative.  The result is
\eqref{linear-eq-structuredcoefficients}.  In particular the derivative of
the modulus is not an independent potential acting on the pure
phase.  For real $\psi=f$, only the strain, $E_0^{\rm c}f^2$ and
$-\rho^2\nabla k\cdot\nabla(f^2)/2$ remain.

\subsubsection{Regional estimates for the lower-order coefficients}
All estimates in this subsection concern $r\geq r_*$; the bounded
core is estimated in \eqref{linear-eq-exactderivative}.
Choose a fixed physical cutoff $\chi_{\rm loc}$ equal to one on a
bounded region containing the supports of all inner corrections and
$Y^{\rm ph}$, for every time under consideration.  It vanishes outside
a larger fixed bounded region.  Put
$q_\alpha=D_t^{\rm mat}\alpha$ and
$\widetilde L=L_0+2D_t^{\rm mat}r_2$, and define
\begin{equation}\label{linear-eq-coefficientallocation}
\begin{aligned}
 L_{\rm loc}&=\chi_{\rm loc}(\epsilon X)\widetilde L,&
 L_\infty&=(1-\chi_{\rm loc}(\epsilon X))\widetilde L-2D_t^{\rm mat}r_2,\\
 E_{\rm loc}&=-c_\epsilon D_t^{\rm mat}q_\alpha-D\cdot\nabla q_\alpha,&
 E_\infty&=-D_t^{\rm mat}r_1,\\
 K_{\rm loc}&=\chi_{\rm loc}(\epsilon X)\nabla k,&
 K_\infty&=(1-\chi_{\rm loc}(\epsilon X))\nabla k .
\end{aligned}
\end{equation}
Thus $L_0=L_{\rm loc}+L_\infty$,
$E_0^{\rm c}=E_{\rm loc}+E_\infty$, and
$\nabla k=K_{\rm loc}+K_\infty$.  These are cutoffs of the
computed coefficients; no derivative of $\chi_{\rm loc}$ enters
their definition.  All local pieces have bounded physical support.
We prove
\begin{equation}\label{linear-eq-lowercoeffnorms}
 \begin{gathered}
 \|F_0\|_\infty\leq C,\quad
 \|L_{\rm loc}\|_{L_X^3}+\|K_{\rm loc}\|_{L_X^3}
                 +\|E_{\rm loc}\|_{L_X^{3/2}}\leq C,\\
 \|L_\infty\|_\infty+\|K_\infty\|_\infty\leq Cc_\epsilon,
 \qquad \|E_\infty\|_\infty\leq Cc_\epsilon^2,\qquad
 |D_t^{\rm mat}\rho^2|\leq C\rho^2 .
 \end{gathered}
\end{equation}
The allocation \eqref{linear-eq-coefficientallocation} puts every
residual term in the infinity pieces.  We first estimate the other
terms, including the volume of the long tube.
For a coefficient of size $\epsilon^a r^b$ on
$r_*<r<C/\epsilon$, its $L_X^p$ norm is bounded by
\begin{equation}\label{linear-eq-radialledger}
 C\epsilon^{a-1/p}
      \left(\int_{r_*}^{C/\epsilon}r^{pb+1}\,dr\right)^{1/p}.
\end{equation}
The factor $\epsilon^{-1/p}$ is the rescaled length.  Logarithms
are retained under the integral when the exponent is critical.

Write $A_*=\sum_kA_k\leq C\epsilon^{5/4}$.  The following table
allocates every summand of \eqref{linear-eq-structuredcoefficients}.
Norms in the table are restricted to the bounded physical region
containing $\operatorname{supp}\chi_{\rm loc}$, with $r\geq r_*$.
Each entry in the last column is multiplied by $C\Lambda_\epsilon^P$,
except for the geometric constant in the last row.
\[
\begin{array}{c|c|c}
 \text{coefficient}&\text{norm}&\text{bound}\\ \hline
 c_\epsilon D_t^{\rm mat}q_\alpha&L_X^{3/2}&
                  \epsilon^2+A_*/\epsilon\\
 D\cdot\nabla q_\alpha&L_X^{3/2}&\epsilon^{1/3}+A_*\\
 c_\epsilon D_t^{\rm mat}k&L_X^3&\epsilon\\
 D\cdot\nabla k&L_X^3&\epsilon^{2/3}\\
 \rho^{-2}\operatorname{div}(\rho^2\nabla q_\alpha)&L_X^3&
                  \epsilon^{2/3}\\
 \rho^2q_\alpha&L_X^3&\epsilon^{2/3}+A_*/\epsilon\\
 \nabla k&L_X^3&1+\epsilon^{2/3}+\epsilon^{-1/3}A_*
\end{array}
\]
Here are the radial calculations behind the table.  The slow part
of $q_\alpha$ is bounded by
$C(\epsilon\ell_r^P/r+\epsilon^2\Lambda_\epsilon^P)$.
Its product with the leading angular drift in
$D\cdot\nabla q_\alpha$ is $C\epsilon\ell_r^P/r^3$,
whose $L_X^{3/2}$ norm is $C\epsilon^{1/3}$.
The second normal derivatives in the fifth row have this radial
power in $L_X^3$, giving $C\epsilon^{2/3}$.
The factors $\nabla\log\rho^2$ are bounded by the derivative of
$\log w^2$ and the amplitude rows in
\eqref{linear-eq-profiletable}; multiplying by them gives the same
or a smaller radial integral.

The radiative contributions of index $n$ to $q_\alpha$ and
$c_\epsilon D_t^{\rm mat}q_\alpha$ have respective sizes
\[
 \epsilon^{n-1}A_kr^{n-1}\Lambda_\epsilon^P,\qquad
 \epsilon^nA_kr^{n-1}\Lambda_\epsilon^P .
\]
At fixed physical distance these have sizes $A_k\Lambda_\epsilon^P$
and $\epsilon A_k\Lambda_\epsilon^P$.  Formula
\eqref{linear-eq-radialledger} gives $C(A_k/\epsilon)\Lambda_\epsilon^P$
in the sixth and first rows.  Applying one normal derivative and
the angular drift gives $CA_k\Lambda_\epsilon^P$ in the second row.
Two normal derivatives give a smaller contribution to the fifth row.

For the third, fourth and last rows use
$k=D_t^{\rm mat}\log\rho^2+\operatorname{div}Y$.
The gradient of $\operatorname{div}Y$ is $O(\epsilon)$ with bounded
physical support, hence has $L_X^3$ norm $O(1)$.
Its product with the $r^{-1}$ drift is $O(\epsilon/r)$ and has
$L_X^3$ norm $O(\epsilon^{2/3})$.
The first amplitude gives $\epsilon\ell_r^P/r^2$ in $\nabla k$.
The differentiated radial radiative amplitude gives
$A_*\Lambda_\epsilon^P/r^3$, whose $L_X^3$ norm is
$C\epsilon^{-1/3}A_*\Lambda_\epsilon^P$.
The other radiative amplitudes are estimated by the same formula.
The outer terms $\epsilon\nabla_xB_t$ and $c_\epsilon B_{tt}$
have local $L_X^3$ norms at most
$C\epsilon^{5/4}\Lambda_\epsilon^P$ by
\eqref{linear-eq-amplitudetimebounds}.  These bounds prove the table.
In particular the only entry without a positive power of $\epsilon$
is the geometric $O(1)$ contribution to $\nabla k$.

On the support of $1-\chi_{\rm loc}(\epsilon X)$ one has
$Y=0$, $\alpha=0$, $\rho=e^{-B}$ and $k=-2B_t$.  Consequently
\[
 \widetilde L=-2c_\epsilon B_{tt}-2D\cdot\nabla_XB_t,\qquad
 \nabla k=-2\epsilon\nabla_xB_t .
\]
The drift is bounded there, so both expressions have $L^\infty$
norm $C\epsilon^{9/4}\Lambda_\epsilon^P=o(c_\epsilon)$.
Lemma~\ref{linear-lem-residualreserve}, with division by the
uniformly nonzero $V$ on $r\geq r_*$, gives the remaining infinity
bounds in \eqref{linear-eq-lowercoeffnorms}.
Finally \eqref{linear-eq-strain} bounds $F_0$: the $r^{-1}$ vortex
gradient is multiplied by the $O(\epsilon r)$ mixed strain, and
the normal strain vanishes.  The regular phase and velocity terms
are bounded.  The amplitude estimates also give
$|D_t^{\rm mat}\rho^2|\leq C\rho^2$.

H\"older, the constant-preserving extension discussed after
\eqref{linear-eq-exteriortrace}, and the amplitude control give, for example,
\[
 \begin{split}
 \left|\int\rho^2L_{\rm loc}fg\right|
       &\leq C\|L_{\rm loc}\|_3\|f\|_6\|\rho^2g\|_2,\\
 \left|\int\rho^2E_{\rm loc}|\psi|^2\right|
       &\leq C\|E_{\rm loc}\|_{3/2}\|\psi\|_6^2,\\
 \left|\int\rho^2K_{\rm loc}\cdot\nabla|\psi|^2\right|
       &\leq C\|K_{\rm loc}\|_3\|\psi\|_6\|\rho\nabla\psi\|_2.
 \end{split}
\]
The infinity pieces use
${c_\epsilon}\int\rho^2|fg|\leq\eta\int\rho^4g^2+
C_\eta{c_\epsilon}^2\int\rho^2f^2$ and its derivative analogue.
On the core, the profile-dependent coefficients of
\eqref{linear-eq-exactderivative} are $O(\epsilon{\Lambda_\epsilon}^P)$ after the
stationary planar contribution has cancelled.  Equation
\eqref{linear-eq-phasecoremass} costs at most one further logarithm;
$\epsilon{\Lambda_\epsilon}^{P+1}=o(1)$.  The other Cartesian components are
controlled by their local planar norm.  Thus all lower-order terms
on a fixed-core complementary field satisfy a uniform bound by
$C\mathcal X^2$.

For the strain term, \eqref{linear-eq-strain} yields
\begin{equation*}
 |R_{\rm strain}(h)|
 \leq C\int\rho^2|\nabla_\parallel\psi|^2
       +C\int_{\mathcal D_t}\epsilon^2r^2\rho^2
                                             |\nabla_\perp\psi|^2.
\end{equation*}
At zero use the combined Cartesian expression.  For a translation
quotient, the second integral is bounded by
$C\epsilon^2(1+\int_1^{2R}dr/r)|b|^2\leq C{c_\epsilon}|b|^2$.
For the complement it is part of its natural energy.  We have proved
\begin{equation}\label{linear-eq-complementderivative}
 |\dot B_t[w_*]|\leq C\mathcal X(w_*)^2
\end{equation}
for the fixed-core complementary field of
\eqref{linear-eq-auxprojection}.  The bound has no logarithmic loss.

\subsection{Quantitative translation and connection estimates}
\label{linear-sec-jacobi}

We now estimate the translation terms and their interactions with the
remainder, obtaining the additional bounds needed for the quadratic
form and its time derivative on a general perturbation.

\subsubsection{Exact residual identities for the lifts}
Let $\mathcal P_V=i{c_\epsilon}\partial_t-H^X_t$.  Translation invariance
and the product rule give
\begin{align}
 \mathcal P_V(N_aV)&=N_a\mathscr R+i{c_\epsilon}(e_a)_t\cdot\nabla V
               +(\Delta e_a)\cdot\nabla V
              +2(\partial_je_a)\cdot\partial_j\nabla V,
                    \label{linear-eq-normalJacobi}\\
 \mathcal P_V(fV)&=f\mathscr R-2f\rho^2V
         +(i{c_\epsilon} f_t+\Delta f)V+2\nabla f\cdot\nabla V
                      \quad(f\text{ real}).\notag
\end{align}
The time derivative of $e_a$ in \eqref{linear-eq-normalJacobi} is Eulerian.
For the complete cut-off and adjusted lift,
\begin{equation}\label{linear-eq-fullJacobi}
\begin{split}
 \mathcal P_VT_a={}&\chi_RN_a\mathscr R+\chi_R[ic_\epsilon\partial_t+\Delta_X,N_a]V\\
 &+(i{c_\epsilon}\chi_{R,t}+\Delta\chi_R)N_aV
                   +2\nabla\chi_R\cdot\nabla(N_aV)\\
 &-a_a\mathscr R+2a_a\rho^2V
       -(i{c_\epsilon} a_{a,t}+\Delta a_a)V-2\nabla a_a\cdot\nabla V.
\end{split}
\end{equation}
The commutator in \eqref{linear-eq-fullJacobi} contains only the
differential part of the operator.  The real-linear potential is
already differentiated in \eqref{linear-eq-normalJacobi}.
All expressions are smooth; no inverse of $V$ enters this identity.

Let $\mathbb H_t$ denote the moving operator in the material field
\eqref{linear-eq-materialfield}.  With $\tau_a$ as in \eqref{linear-eq-tau}, set
\begin{equation*}
 \mathscr D_a=\mathbb H_t\tau_a-i{c_\epsilon}\partial_t\tau_a
                  =-\mathscr S_{\rm mat}(\mathcal P_VT_a).
\end{equation*}
The transform on the right includes the density and the leading
phase.  Its derivative is a material coefficient derivative;
large derivatives of the leading phase cancel before taking norms.

\subsubsection{The transverse commutator and the dual norm}
At a fixed time use a local Bishop frame for the calculation.  In
current rescaled arclength $\sigma$, $N_a=\partial_{y_a}$ and
\[
 \Delta_X=\Delta_y-\frac{\epsilon\kappa}{H}\cdot\nabla_y
   +H^{-2}\partial_{\sigma\sigma}
   +\frac{\epsilon(\kappa_\sigma\cdot y)}{H^3}\partial_\sigma.
\]
Consequently
\begin{equation}\label{linear-eq-normalcommutator}
\begin{split}
 [\Delta_X,N_a]V={}&
 \frac{\epsilon^2\kappa_a\kappa_b}{H^2}V_{y_b}
        -\frac{2\epsilon\kappa_a}{H^3}V_{\sigma\sigma}\\
 &-\left(\frac{\epsilon\kappa_{a,\sigma}}{H^3}
       +\frac{3\epsilon^2\kappa_a(\kappa_\sigma\cdot y)}{H^4}\right)V_\sigma.
\end{split}
\end{equation}
In the periodic frame, use covariant normal derivatives and retain
the angular connection.  It has the same estimates.  Notice that
there is no purely normal term of order $\epsilon$ in
\eqref{linear-eq-normalcommutator}.  The coefficient of order $\epsilon$
multiplies a slow longitudinal derivative of the \emph{background}.
No restriction is made on longitudinal frequencies of a test field.

Here is a concrete dual convention for the estimates below.  On
$r\geq r_*$ express a source as $e^{i\Theta}(A+iP)$.  Its amplitude
part is measured in $L_y^2$ and the nonzero angular part of its phase
in $\langle r\rangle L_y^2$.  For the averaged phase part, including
the volume and modulus factors from the pairing, put
\[
 P_0(r)=\int_0^r\int_0^{2\pi}j_0\rho P\,d\theta\,s\,ds.
\]
Measure it by
\begin{equation}\label{linear-eq-dualphase}
 \left(\int_0^{R_I}\frac{|P_0(r)|^2}{w(r)^2r}\,dr\right)^{1/2}
                                                  +|P_0(R_I)|.
\end{equation}
The bounded core is paired in the Cartesian fields; its averaged
phase is estimated by the same radial integration or
\eqref{linear-eq-coreHardy}.  The boundary value in
\eqref{linear-eq-dualphase} is paired with $f(R_I)$ and bounded by
\eqref{linear-eq-exteriortrace}.  Integration in $\sigma_0$ and
Cauchy--Schwarz give a dual norm for $\mathcal X$.  We write
$\|S\|_{\mathcal X'}\leq A$ for the precise operator assertion
\begin{equation}\label{linear-eq-dualoperator}
 \left|\int a_j(\sigma_0)\langle S_j,u\rangle_y\,d\sigma_0\right|
                    \leq A\|a\|_2\mathcal X(u).
\end{equation}
Derivatives of $a$ are not included in this notation; longitudinal
commutators are handled separately in
Subsection~\ref{linear-sec-projectionevolution}.

\begin{lemma}[Coefficient and paired Jacobi estimates]
\label{linear-lem-Jacobiestimates}
The lifts satisfy
\begin{equation}\label{linear-eq-Jacobiestimates}
 \begin{split}
 \|\mathscr D\|_{\mathcal X'}+
       \|\partial_t\mathscr D\|_{\mathcal X'}
       +\|i{c_\epsilon}\tau_t\|_{\mathcal X'}
       +\|i{c_\epsilon}\tau_{tt}\|_{\mathcal X'}&\leq C\epsilon,\\
 |\langle \mathscr D_i,\tau_j\rangle_y|&\leq C{c_\epsilon},\\
 |\langle \partial_t\mathscr D_i,\tau_j\rangle_y|
       +|\langle \mathscr D_i,\partial_t\tau_j\rangle_y|&\leq C{c_\epsilon}{\Lambda_\epsilon}.
 \end{split}
\end{equation}
The coefficient estimates are uniform in the material position.
The first line is understood in \eqref{linear-eq-dualoperator}.
\end{lemma}
\begin{proof}
We first record the differentiated amplitude bounds used in the
calculation.  Put $A_*=\sum_kA_k$ and
\[
 b_0=\epsilon^3,\qquad b_1=\epsilon^{9/4},\qquad b_2=\epsilon^{5/4}.
\]
For $q=0,1,2$ define
\begin{equation}\label{linear-eq-adjustmentenvelope}
\begin{split}
 m_q(r)={}&\frac{\epsilon\ell_r^P}{\langle r\rangle^2}
 +\Lambda_\epsilon^P
       \left(\frac{\epsilon^2}{\langle r\rangle}+\epsilon^3+b_q\right)\\
 &+A_*\Lambda_\epsilon^P
       \left(\frac{\epsilon^{1-q}}{\langle r\rangle^3}
             +\frac{\epsilon^{2-q}}{\langle r\rangle}
             +\epsilon^{3-q}\right).
\end{split}
\end{equation}
The quotient in \eqref{linear-eq-ampadjust}, differentiated at fixed
material coordinates, gives
$|\partial_t^qa_i|\leq Cm_q$.
A normal derivative lowers each displayed radial power by one;
on the $b_q$ term it supplies an additional $\epsilon$.
Each $\sigma_0$ derivative supplies an $\epsilon$.
Cutoff derivatives obey these rules on $R<r<2R$.
For the first profile this follows by differentiating the angular
quotient and using the separated radial logarithms in
\eqref{linear-eq-profiletable}.  The leading planar term is
\eqref{linear-eq-ampradialleading}.
The three terms multiplied by $A_*$ come respectively from the
radial radiative amplitude, index $n=2$, and the higher polynomial
profiles on $r\leq2R$.  The last contribution $b_q$ is
$\epsilon\nabla_x\partial_t^qB$, with the material chain rule
included, and follows from Lemma~\ref{linear-lem-amplitudetime}.
This contribution is estimated as a bounded tail.  On the fixed core
the angular quotient is smooth; the only first-order core logarithm is
the pure phase in \eqref{core-logarithm-splitting}, so it does not enter
its leading amplitude estimate.

Time derivatives in \eqref{linear-eq-fullJacobi} are Eulerian.
Before applying \eqref{linear-eq-adjustmentenvelope} we use
\begin{equation}\label{linear-eq-Jacobitimeconversion}
 \chi_{R,t}=-Y\cdot\nabla\chi_R,\qquad
 a_{i,t}=D_t^{\rm mat}a_i-Y\cdot\nabla a_i,\qquad
 (e_i)_t=D_t^{\rm mat}e_i-Y\cdot\nabla e_i.
\end{equation}
Here $D_t^{\rm mat}\chi_R=0$ on the material tube.
The vector $c_\epsilon Y$ and its first material derivative are
$O(\epsilon\Lambda_\epsilon)$; the Eulerian frame coefficient in
\eqref{linear-eq-Jacobitimeconversion} and its first material
derivative are bounded.  Thus the cutoff and adjustment transport
terms are included in the estimates below.

We use amplitude $L_y^2$, weighted nonzero-phase $L_y^2$, and
the primitive norm \eqref{linear-eq-dualphase} separately.
For reference, the following are bounds for all the groups in
\eqref{linear-eq-fullJacobi}, after applying the material transform.
For a group $S_i$, set $J_{ij}=\langle S_i,\tau_j\rangle_y$.
Dots in the table are material coefficient derivatives.
\[
\begin{array}{c|c|c|c}
 \text{group}&\|S\|_{\mathcal X'}+\|\dot S\|_{\mathcal X'}
       &|J_{ij}|&|\dot J_{ij}|\\ \hline
 \chi_RN_i\mathscr R-a_i\mathscr R&o(\epsilon)&o(c_\epsilon)&o(c_\epsilon)\\
 \chi_R[\Delta_X,N_i]V&C\epsilon&Cc_\epsilon&Cc_\epsilon\Lambda_\epsilon\\
 ic_\epsilon\chi_R(e_i)_t\cdot\nabla V&C\epsilon&Cc_\epsilon&Cc_\epsilon\Lambda_\epsilon\\
 \text{spatial and time cutoff terms}&C\epsilon&Cc_\epsilon&Cc_\epsilon\Lambda_\epsilon\\
 2a_i\rho^2V&C\epsilon&Cc_\epsilon&Cc_\epsilon\Lambda_\epsilon\\
 -(\Delta a_i)V-2\nabla a_i\cdot\nabla V&C\epsilon&Cc_\epsilon&Cc_\epsilon\Lambda_\epsilon\\
 -ic_\epsilon a_{i,t}V&C\epsilon&Cc_\epsilon&Cc_\epsilon\Lambda_\epsilon
\end{array}
\]
We now justify the entries, including the angular cancellations
needed for the paired columns.

For the geometric commutator use \eqref{linear-eq-normalcommutator}
and the bounds
\[
 |\nabla_yV|\leq C(\langle r\rangle^{-1}+\epsilon L(\epsilon r)),
 \qquad |V_\sigma|\leq C\epsilon,\qquad
 |V_{\sigma\sigma}|\leq C\epsilon^2.
\]
After removal of the common leading phase the same estimates
hold for one material derivative, with higher profile terms
estimated by \eqref{linear-eq-profiletable}.
The resulting envelope is
$C\epsilon^2/\langle r\rangle+C\epsilon^3L(\epsilon r)$.
The corresponding weighted squares are bounded by
\[
 C\epsilon^4\int_0^{2R}r\,dr+
 C\epsilon^6\int_{r_*}^{2R}r^3L(\epsilon r)^2\,dr
 \leq C\epsilon^2.
\]
The leading phase contribution has first angular harmonic.  Its
weighted angular mean has the extra geometric or profile factor,
and is $O(\epsilon^3(1+L(\epsilon r)))$.
Its primitive is
$O(\epsilon^3r^2(1+L(\epsilon r)))$, with both terms of
\eqref{linear-eq-dualphase} bounded by $C\epsilon$.
The radial integrals in a translation pairing are instead
\begin{equation}\label{linear-eq-criticaltranslationintegral}
 \epsilon^2\int_1^{2R}\frac{dr}{r}
 +\epsilon^3\int_{r_*}^{2R}L(\epsilon r)\,dr
 +\epsilon^4\int_{r_*}^{2R}L(\epsilon r)^2r\,dr
 \leq Cc_\epsilon .
\end{equation}
These estimates include a material derivative: all leading
geometric coefficients have bounded derivatives, and differentiating
a radiative coefficient retains an additional positive power of
$\epsilon$ in these integrals.

For the frame-time term the phase--amplitude expression is
\[
 ic_\epsilon(e_i)_t\cdot\nabla V
 =e^{i\Theta}\{-c_\epsilon\rho(e_i)_t\cdot\nabla\Theta
                    +ic_\epsilon(e_i)_t\cdot\nabla\rho\}.
\]
The $r^{-1}$ tail is amplitude-directed, with $L_y^2$ norm
$Cc_\epsilon\sqrt{\Lambda_\epsilon}=o(\epsilon)$.
The phase part contains $w'$; its leading angular mean vanishes.
The geometric mean correction has primitive bounded by
$Cc_\epsilon\epsilon\ell_r^{P+1}$.
For a material derivative the additional bounded phase tail from
$\nabla_yB_t$ has dual norm at most
$Cc_\epsilon b_1R^2=o(\epsilon)$.
In the paired column the two leading phase tails have zero real
pairing after multiplication by $i$.  The remaining normal
determinant contains $\nabla\rho\,\nabla\Theta$ and has bounded
absolute integral.  The longitudinal component contains the
additional $\epsilon$ from $V_\sigma$.  This gives $Cc_\epsilon$.
After differentiation, the additional outer-amplitude integral
is bounded by $Cc_\epsilon b_1R=o(c_\epsilon)$; the first and
polynomial profiles are bounded by the same convergent radial
integrals.

On the cutoff annulus the spatial commutator has size $CR^{-3}$,
weighted dual norm $CR^{-1}$, and translation pairing $CR^{-2}$.
Its averaged primitive is also bounded by $CR^{-1}$.
The time cutoff from \eqref{linear-eq-Jacobitimeconversion} has
bounded coefficients.  Its leading rotated translation is
amplitude-directed, with annular norm $Cc_\epsilon$.
Its phase component contains an amplitude derivative of $V$.
The corresponding weighted norm and primitive are bounded by
$Cc_\epsilon\epsilon\Lambda_\epsilon^P$.
After one material derivative add
$Cc_\epsilon b_1R^2\Lambda_\epsilon^P=o(\epsilon)$
in the phase dual norm and
$Cc_\epsilon b_1R\Lambda_\epsilon^P=o(c_\epsilon)$
in the translation pairing.

For the massive adjustment, \eqref{linear-eq-adjustmentenvelope}
gives
\[
 \|a_i\|_{L_y^2}+\|\partial_ta_i\|_{L_y^2}\leq C\epsilon.
\]
For example, the first tail uses
\begin{equation}\label{linear-eq-logtailintegral}
 \int_1^\infty\ell_r^{2P}r^{-3}\,dr<\infty,
\end{equation}
and the outer differentiated tail uses
$b_1R\Lambda_\epsilon^P=O(\epsilon^{5/4}\Lambda_\epsilon^P)=o(\epsilon)$.
The paired estimate uses the angular identity, rather than this
$L^2$ estimate.  With
\[
 \tau_j=\sqrt{j_0}\rho e^{i(\theta+\alpha)}(f_j+ip_j),
 \qquad
 f_j=\chi_R\partial_{y_j}\log\rho-a_j,\qquad
 p_j=\chi_R\partial_{y_j}\Theta ,
\]
one has $\int_0^{2\pi}j_0\rho^2f_j\,d\theta=0$.  Hence
\begin{equation}\label{linear-eq-massivepaircancellation}
 \int j_0a_i\rho^4f_j\,dy
 =\int j_0a_i(\rho^2-w^2)\rho^2f_j\,dy .
\end{equation}
The first amplitude factors on the right have orders
$\epsilon\ell_r^P/r^2$ and $\epsilon\ell_r^P/r$.
The leading $f_j$ is $O(r^{-3})$; its corrections have orders
$m_0$.  Their radial integrals are $O(\epsilon^2)$.
Differentiation of \eqref{linear-eq-massivepaircancellation}
retains the extra amplitude factor.  The new bounded outer
contributions have sizes $b_1$ and $B_t=O(\epsilon^{5/4}\Lambda_\epsilon^P)$;
their products with the remaining first amplitude and the
translation amplitude are $O(\epsilon^{9/4}\Lambda_\epsilon^P)$
or smaller.  This verifies the differentiated paired column for
the massive term.

For the derivative adjustment write, in a Bishop chart,
\[
 \nabla a_i\cdot\nabla V
   =(a_i)_rV_r+H^{-2}(a_i)_\sigma V_\sigma .
\]
The corresponding covariant formula holds in the periodic frame,
since $a_i$ is radial.  Thus the leading angular phase gradient
does not enter this product.  The leading amplitude part of
$(a_i)_rV_r$ is radial and its pairing with $T_j$ vanishes by
the same angular identity.  The nonradial part of $\Delta_Xa_i$
contains an additional $\epsilon$, or two longitudinal derivatives.
The remaining phase gradient in $V_r$ has order
$\epsilon L(\epsilon r)+\epsilon\ell_r^P/r^2+
\epsilon^2\Lambda_\epsilon^P/r$.
The borderline paired integral is therefore bounded by
\[
 C\epsilon^2\int_{r_*}^{2R}
          \frac{\ell_r^P L(\epsilon r)}{r^3}\,dr
       +C\epsilon^2\leq Cc_\epsilon .
\]
The first core phase contributes the single factor
$\Lambda_\epsilon$ identified in
\eqref{core-logarithm-splitting}; other global logarithms here
have an additional positive power of $\epsilon$.
The amplitude dual norms use two further radial powers in
\eqref{linear-eq-adjustmentenvelope} and the convergent integral
\eqref{linear-eq-logtailintegral}.  The phase dual norms and
primitives use the displayed phase gradient.
The differentiated estimate follows with $m_1$ and the first
material derivative of that gradient.  Its additional radiative
terms have a factor $A_*$ and satisfy the paired bound because
$A_*\Lambda_\epsilon^P=o(\epsilon)$.

The last adjustment row has two parts by
\eqref{linear-eq-Jacobitimeconversion}.
The primitive of $c_\epsilon D_t^{\rm mat}a_i$ is estimated from
$m_1$, and that of its material derivative from $m_2$.
The largest new bounded term has dual norm
\[
 Cc_\epsilon b_2R^2\Lambda_\epsilon^P
       =C\epsilon^{5/4}\Lambda_\epsilon^{P+1}=o(\epsilon).
\]
The decaying radial radiative term gives at most
$Cc_\epsilon(A_*/\epsilon)\Lambda_\epsilon^{P+1/2}=o(\epsilon)$.
The transport part has leading size
$Cc_\epsilon/r^3$ and primitive norm
$Cc_\epsilon\sqrt{\Lambda_\epsilon}$.
The bounded tails acquire an additional normal factor
$\epsilon$ before multiplication by $c_\epsilon Y$.
In a translation pairing use $|p_j|\leq C(r^{-1}+
\epsilon L(\epsilon r))$ together with the profile remainders.
In particular the bounded second-time amplitude contributes
$Cc_\epsilon b_2R\Lambda_\epsilon^P=o(c_\epsilon)$.
These calculations justify the last row, including its derivative.
The residual row follows from Lemma~\ref{linear-lem-residualreserve};
even the $c_\epsilon$-weighted mass estimate leaves the stated
positive powers of $\epsilon$.

To separate differentiation of the test lift, put
\begin{equation}\label{linear-eq-zetasplit}
 \tau_t=\omega_\ell\tau+\zeta,\qquad
 \omega_\ell=\tfrac12\partial_t\log a_\ell .
\end{equation}
The profile and metric estimates give
\begin{equation}\label{linear-eq-zetabound}
 \mathcal X(\zeta a)\leq
 C\epsilon\Lambda_\epsilon\|a\|_2+C\|a_{\sigma_0}\|_2.
\end{equation}
For clarity, the logarithmic part of the normal regular-phase
gradient is angularly constant.  Its zero-phase energy is
estimated by
$\epsilon^2\int_1^{2R}dr/r$; its nonconstant angular part is
$O(\epsilon)$ and has the same energy bound.
The $H_t$ term is $\epsilon r$ times a translation.
The first core phase has size $C\epsilon\Lambda_\epsilon$,
and the radiative contributions have smaller natural norm.
Also $\|\zeta\|_{L_y^2}\leq C$ and
$\|\partial_{\sigma_0}\zeta\|_{L_y^2}
+\|\partial_{\sigma_0}\tau\|_{L_y^2}
\leq C\epsilon\sqrt{\Lambda_\epsilon}$.
These prove \eqref{linear-eq-zetabound}, with the longitudinal
derivative falling on $a$ included.
The pointwise transverse version of this estimate and the dual
column of the table imply
$|\langle S_i,\zeta_j\rangle_y|\leq C\epsilon^2\Lambda_\epsilon$.
Together with the paired column this bounds
$\langle S_i,\partial_t\tau_j\rangle_y$.
Finally
$\langle\partial_tS_i,\tau_j\rangle_y
=\partial_tJ_{ij}-\langle S_i,\partial_t\tau_j\rangle_y$,
which proves both differentiated paired estimates.

It remains to estimate the two rotated time derivatives of the lift.
Differentiating the displayed phase--amplitude formula for $\tau_j$
shows that the largest $p_{j,tt}$ term rotates into amplitude.
Its bounded radiative part is $O(A_*/\epsilon)$ and gives
$Cc_\epsilon(A_*/\epsilon)R\Lambda_\epsilon^P
=CA_*\Lambda_\epsilon^{P+1}=o(\epsilon)$.
The $f_{j,tt}$ term rotates into phase and is bounded using $m_2$,
including $c_\epsilon b_2R^2$.
The term $\alpha_{tt}p_j$ is also phase-directed after rotation.
For a radiative phase of index $n$ its weighted norm and primitive
are bounded by
\[
 Cc_\epsilon\epsilon^{n-2}A_k\Lambda_\epsilon^P
       \left(\int_1^{2R}r^{2n-2}r\,dr\right)^{1/2}
 \leq CA_k\Lambda_\epsilon^{P+1}=o(\epsilon).
\]
The slow terms use the separated radial logarithms, and the common
length-density term has amplitude norm
$Cc_\epsilon\sqrt{\Lambda_\epsilon}$.
Products of first time derivatives have an additional small
profile factor and satisfy the same bounds.  The first rotated
time derivative is covered by these estimates with one fewer
time derivative.  All fixed-core estimates use smooth Cartesian
fields and \eqref{linear-eq-coreHardy}.
Integration in $\sigma_0$ and Cauchy--Schwarz now give
\eqref{linear-eq-dualoperator}, without an additional length
factor, and complete the proof of \eqref{linear-eq-Jacobiestimates}.
\end{proof}

\subsubsection{The material connection}
\begin{lemma}[Connection and its derivative]\label{linear-lem-connection}
For a real $b=b(\sigma_0)$ fixed in time, define
\[
 C_t[b]={c_\epsilon}\operatorname{Re}\int i\partial_t(\tau_jb_j)
                                     \overline{\tau_kb_k}\,dy\,d\sigma_0.
\]
Then
\begin{equation}\label{linear-eq-connectionbound}
 |C_t[b]|+|\partial_t C_t[b]|\leq C{c_\epsilon}\|b\|_2^2.
\end{equation}
\end{lemma}
\begin{proof}
Outside the fixed core, write
\[
 \tau_b=\sqrt{j_0}\rho e^{i(\theta+\alpha)}(f_b+ip_b),
 \quad f_b=\chi_R b_a\partial_{y_a}\log\rho-b_aa_a,
 \quad p_b=\chi_R b_a\partial_{y_a}\Theta.
\]
Direct differentiation gives the exact connection
\begin{equation}\label{linear-eq-connectionidentity}
 C_t[b]={c_\epsilon}\int j_0\rho^2
       \{p_b(f_b)_t-f_b(p_b)_t-\alpha_t(f_b^2+p_b^2)\}.
\end{equation}
Both derivatives of the common real factor $\sqrt{j_0}\rho$
drop out.

Sufficient envelopes on $r\geq r_*$ are
\begin{equation*}
\begin{aligned}
 |f_b|&\leq C|b|\{r^{-3}+{\Lambda_\epsilon}^P
                  (\epsilon r^{-2}+\epsilon^2r^{-1}+\epsilon^3)\},\\
 |(f_b)_t|&\leq C|b|{\Lambda_\epsilon}^P(\epsilon r^{-2}+\epsilon^2),\\
 |p_b|&\leq C|b|(r^{-1}+\epsilon{\Lambda_\epsilon}^P),
       &|(p_b)_t|&\leq C|b|\epsilon{\Lambda_\epsilon}^P,\\
 |\alpha_t|&\leq C\epsilon{\Lambda_\epsilon}^P.
\end{aligned}
\end{equation*}
The first two follow from the modulus, including
$\partial_yB_t=\epsilon e_a\cdot\nabla_xB_t$; the phase
bounds follow from the regular phase and
\eqref{linear-eq-profiletable}.  The same estimates differentiated once
more can be bounded by
\[
 |(f_b)_{tt}|\leq C|b|{\Lambda_\epsilon}^P
                (\epsilon^{1/4}r^{-2}+\epsilon^{5/4}),
 \qquad |(p_b)_{tt}|+|b|\,|\alpha_{tt}|
                                      \leq C|b|\epsilon^{1/4}{\Lambda_\epsilon}^P.
\]
These bounds follow from \eqref{linear-eq-amplitudetimebounds}.

Every integral in \eqref{linear-eq-connectionidentity}, divided by
${c_\epsilon}|b|^2$, is bounded by a sum of
\[
 C{\Lambda_\epsilon}^P\{\epsilon+\epsilon^2R+
                 \epsilon^2\log R+\epsilon^3R^2\}=o(1).
\]
When differentiating the bracket, the two products
$(p_b)_t(f_b)_t-(f_b)_t(p_b)_t$ cancel.  The remaining terms have
the same radial integrals with $\epsilon^{1/4}{\Lambda_\epsilon}^P$ in place
of their small coefficient, and are $o(1)$.  The derivative of
$j_0\rho^2$ is bounded by itself times a fixed constant.  On the
fixed core the smooth fields give these conclusions directly from
the centred first and second time bounds.  The two connection
estimates follow.  Integration in $\sigma_0$ adds no length factor.
\end{proof}

\subsubsection{Static and differentiated blocks}
We record the block estimates in the form used in the assembly.
In a longitudinal expression $\mathcal K$ always refers to the
complete field; the coefficient statements below concern the
transverse part after this square has been retained.
\begin{lemma}[Translation blocks]\label{linear-lem-blockestimates}
For $q=w_*+\mathscr T_t\beta$ as in
\eqref{linear-eq-auxprojection},
\begin{equation}\label{linear-eq-staticblocks}
 B_t[q]\geq k(E_*(w_*)+K(q))
                       -C{c_\epsilon}\|\beta\|_2^2-C{c_\epsilon}^2\|q\|_2^2.
\end{equation}
For any real vector $a=a(\sigma_0)$ and instantaneous test $u$,
\begin{align}
 |\dot B_t[\mathscr T_ta]|
   &\leq C{\Lambda_\epsilon}\|a_{\sigma_0}\|_2^2
                           +C{c_\epsilon}{\Lambda_\epsilon}\|a\|_2^2,\label{linear-eq-dotdiagonal}\\
 2|\dot B_t[\mathscr T_ta,u]|
   &\leq C\mathcal X(u)^2+C{\Lambda_\epsilon}\|a_{\sigma_0}\|_2^2
                           +C{c_\epsilon}{\Lambda_\epsilon}\|a\|_2^2.\label{linear-eq-dotmixed}
\end{align}
\end{lemma}
\begin{proof}
For the static transverse pairing, use
$\mathbb H_t\tau=\mathscr D+i{c_\epsilon}\tau_t$ and the first line
of \eqref{linear-eq-Jacobiestimates}.  The mixed coefficient is
$O(\epsilon)$ in the natural dual norm.  Young's inequality costs
$C\epsilon^2|\beta|^2\leq C{c_\epsilon}|\beta|^2$.
The diagonal pairing of $\mathscr D$ is $O({c_\epsilon})$ by the separate
paired bound; that of $i{c_\epsilon}\tau_t$ is
\eqref{linear-eq-connectionbound}.  The background longitudinal
coefficient terms removed to retain the common square obey the
same estimates: each background longitudinal derivative supplies
an $\epsilon$, and the only critical integral is
\eqref{linear-eq-criticaltranslationintegral}.  Apply the fixed-core
regional coercivity, leaving the common positive longitudinal
square intact.  This proves \eqref{linear-eq-staticblocks}.

The differentiated identities are especially useful.  For $a$
held fixed in time, set
$\mathscr D_a=\mathbb H_t(\tau_ja_j)-i{c_\epsilon}\partial_t(\tau_ja_j)$.
Then exactly
\begin{align*}
 \dot B_t[\mathscr T_ta]
  &=\langle \partial_t\mathscr D_a,\tau_ja_j\rangle
       -\langle \mathscr D_a,\partial_t(\tau_ja_j)\rangle
                              +\partial_tC_t[a],\\
 \dot B_t[\mathscr T_ta,u]
  &=\langle \partial_t\mathscr D_a+i{c_\epsilon}\partial_{tt}(\tau_ja_j),u\rangle
                              -B_t[\partial_t(\mathscr T_ta),u].
\end{align*}
The terms $\langle i{c_\epsilon}\tau_t,\tau_t\rangle$ vanish in deriving the
first formula.  Its coefficient terms are bounded by the paired
part of Lemma~\ref{linear-lem-Jacobiestimates} and
Lemma~\ref{linear-lem-connection}, giving $C{c_\epsilon}{\Lambda_\epsilon}\|a\|_2^2$.
In the second formula the first pairing has coefficient size
$C\epsilon\|a\|_2\mathcal X(u)$.  Separate
$\partial_t\tau=\omega_\ell\tau+\zeta$ in the last term.
The length term uses the static mixed estimate.  The remainder
uses \eqref{linear-eq-Xcontinuity} and \eqref{linear-eq-zetabound}; its
transverse contribution is at most
$C\epsilon{\Lambda_\epsilon}\|a\|_2\mathcal X(u)$, bounded by
$C\mathcal X(u)^2+C\epsilon^2{\Lambda_\epsilon}^2\|a\|_2^2$.
Since $\epsilon^2{\Lambda_\epsilon}^2={c_\epsilon}{\Lambda_\epsilon}$, this is the asserted scale.

It remains to estimate the longitudinal terms in these identities.
For the operator \eqref{linear-eq-fullmaterialoperator} and a real
$a=a(\sigma_0)$, the product rule gives
\begin{equation}\label{linear-eq-longcommutator}
 [\mathbb H_t,a]\tau
 =-2G^{\sigma_0B}\partial_B\tau\,a'
   -\partial_A G^{A\sigma_0}\tau\,a'
   -G^{\sigma_0\sigma_0}\tau\,a''-i d^{\sigma_0}\tau\,a'.
\end{equation}
The drift is specified in \eqref{linear-eq-materialcoefficients}.
One integration by parts gives \eqref{linear-eq-longflux} below.
Lemma~\ref{linear-lem-longitudinalkernels} estimates that full flux
and its first coefficient derivative.  In particular
\eqref{linear-eq-fluxcoefficientderivative} bounds a mixed
longitudinal term by
\[
 C\|a'\|_2\{\sqrt{\Lambda_\epsilon}K(u)^{1/2}+\mathcal X(u)\}
 \leq C\mathcal X(u)^2+C\Lambda_\epsilon\|a'\|_2^2.
\]
For the diagonal use \eqref{linear-eq-translationfluxmatrices},
including its time-differentiated lift version.  The resulting
terms have bounds
\[
 C\Lambda_\epsilon\|a'\|_2^2
 +C\epsilon\Lambda_\epsilon\|a'\|_2\|a\|_2
 \leq C\Lambda_\epsilon\|a'\|_2^2+Cc_\epsilon\|a\|_2^2.
\]
These estimates include the density, gauge, drift and angular
connection contributions.  Together with the coefficient estimates
already proved they give
\eqref{linear-eq-dotdiagonal}--\eqref{linear-eq-dotmixed}.
\end{proof}

\subsection{The projected evolution and longitudinal flux}
\label{linear-sec-projectionevolution}

We derive the equation governing the translation coefficients and
estimate the terms involving derivatives along the filament.  This
gives the time-derivative bound for the translation term that will be
added to the quadratic form.

The real symplectic coefficients in \eqref{linear-eq-sympprojection} are
functions on the fixed initial circle.  If $\widetilde G$ denotes
the forcing in the material field, the exact weak equation is
\begin{equation}\label{linear-eq-projectedequation}
 \begin{split}
 {c_\epsilon}\int\partial_t(Nb)\cdot a\,d\sigma_0
   ={}&-\langle z,\mathbb H_t(\tau_ia_i)-i{c_\epsilon}(\partial_t\tau_i)a_i\rangle\\
     &-\langle \widetilde G,\tau_ia_i\rangle
 \end{split}
\end{equation}
for every smooth real vector $a(\sigma_0)$.  Indeed
${c_\epsilon}\partial_t\langle z,i\tau_i\rangle
 =-\langle \mathbb H_tz+\widetilde G,\tau_i\rangle
       +{c_\epsilon}\langle z,i\partial_t\tau_i\rangle$.
Symmetry and multiplication by $a_i$ give
\eqref{linear-eq-projectedequation}.

The coefficient part of the right-hand side is
$-\langle z,\mathscr D_ia_i\rangle$.  Its longitudinal part is given by
\eqref{linear-eq-longcommutator}.  After one integration by parts the
latter can equivalently be written as
\begin{equation}\label{linear-eq-longflux}
 \begin{split}
 \langle z,[\mathbb H_t,a_i]\tau_i\rangle
 =\int a_i'\Big\{&\operatorname{Re}\big[G^{\sigma_0B}
       (\tau_i\partial_B\overline z-(\partial_B\tau_i)\overline z)\big]\\
                   &-d^{\sigma_0}\operatorname{Im}(z\overline{\tau_i})\Big\}
                                    dy\,d\sigma_0.
 \end{split}
\end{equation}
The metric coefficients satisfy
$G^{\sigma_0\sigma_0}=(a_\ell H)^{-2}$ and
$G^{\sigma_0a}=O(\epsilon r)$; the first expression in
\eqref{linear-eq-longflux} is the covariant longitudinal derivative.

\begin{lemma}[Longitudinal kernels and measurements]
\label{linear-lem-longitudinalkernels}
Let $\mathcal F_i(q)(\sigma_0)$ be the integral in $y$ of the
bracket in \eqref{linear-eq-longflux}.  Then
\begin{equation}\label{linear-eq-fullfluxbound}
 \|\mathcal F(q)\|_2\leq
 C\sqrt{\Lambda_\epsilon}\,K(q)^{1/2}
 +C\mathcal X(w_*)+C\epsilon\Lambda_\epsilon\|\beta\|_2 .
\end{equation}
Its first coefficient derivative, with $q$ fixed, satisfies
\begin{equation}\label{linear-eq-fluxcoefficientderivative}
 \|\dot{\mathcal F}(q)\|_2\leq
 C\sqrt{\Lambda_\epsilon}\,K(q)^{1/2}+C\mathcal X(q).
\end{equation}
The symplectic coefficients satisfy
\begin{equation}\label{linear-eq-sympderivative}
 \|b'\|_2\leq C\sqrt{\Lambda_\epsilon}\,K(q)^{1/2}
 +C\epsilon\sqrt{\Lambda_\epsilon}\,E_*(w_*)^{1/2}
 +C\epsilon\sqrt{\Lambda_\epsilon}
                  (\|\beta\|_2+\|b\|_2).
\end{equation}
On translation fields the flux has the representation
\begin{equation}\label{linear-eq-translationfluxmatrices}
 \mathcal F_i(\mathscr T_ta)
       =\mathsf G_{ij}a_j'+\mathsf C_{ij}a_j,\qquad
 |\mathsf G|+|\partial_t\mathsf G|\leq C\Lambda_\epsilon,\quad
 |\mathsf C|+|\partial_t\mathsf C|\leq C\epsilon\Lambda_\epsilon .
\end{equation}
When one lift is replaced by its material time derivative, the
undifferentiated bounds $|\mathsf G|\leq C\Lambda_\epsilon$ and
$|\mathsf C|\leq C\epsilon\Lambda_\epsilon$ also hold.
\end{lemma}
\begin{proof}
Write $\partial_\sigma=\partial_{\sigma_0}$ in this proof.
The coefficients in \eqref{linear-eq-materialcoefficients} give
$G^{\sigma B}\partial_B=M(\partial_\sigma-\omega\partial_\theta)$.
All lift cutoffs are independent of $\sigma,\theta,t$ in these
coordinates.  In particular no radial derivative of $\chi_R$
enters the longitudinal flux.

First work outside the fixed core, where
$\eta=iW\psi$ and $z=iU_{\rm mat}\psi$.  Write
$\tau_i=U_{\rm mat}(f_i+ip_i)$ as in the Jacobi proof, and put
\[
 m=j_0\rho^2M,\qquad
 D^\sigma=2G^{\sigma B}\partial_B\Theta-c_\epsilon Y^\sigma .
\]
In the Wronskian in \eqref{linear-eq-longflux}, the derivatives
of the common real factor $\sqrt{j_0}\rho$ cancel.
The derivatives of its phase $\theta+\alpha$ combine with $d^\sigma$
to give $D^\sigma$.  Substitution and angular integration by parts
therefore give
\begin{equation}\label{linear-eq-fluxkernels}
\begin{split}
 \mathcal F_i(q)
 ={}&\int m(p_i\partial_\sigma\psi_1
                         -f_i\partial_\sigma\psi_2)\,dy\\
 &+\int(\mathcal A_i\psi_1+\mathcal B_i\psi_2)\,dy,\\
 \mathcal A_i={}&-m\partial_\sigma p_i
       +\omega(\partial_\theta m)p_i
       +2\omega m\partial_\theta p_i-j_0\rho^2D^\sigma f_i,\\
 \mathcal B_i={}&m\partial_\sigma f_i
       -\omega(\partial_\theta m)f_i
       -2\omega m\partial_\theta f_i-j_0\rho^2D^\sigma p_i .
\end{split}
\end{equation}
Use a fixed radial partition at $r_*$ to interpret this formula
together with its smooth core expression.  The partition acts on
the coefficients, not on the field.  Its angular and longitudinal
derivatives vanish.  On the core expand
$z=\sqrt{j_0}\mathcal C\eta$ in \eqref{linear-eq-longflux},
and move its angular derivative onto the smooth test.
Thus throughout the disc the first integral in
\eqref{linear-eq-fluxkernels} is a real-linear measurement of
$\eta_\sigma$; every other term is a coefficient paired with $\eta$.
This retains the ordinary derivative in the definition of $K$.

We give bounds for these coefficients and their first material
derivatives.  For $r\geq r_*$ put
\[
 \mathcal E_\epsilon(r)=C\Lambda_\epsilon^P
 \left\{\frac{\epsilon^2\ell_r^P}{r^2}
       +\frac{\epsilon^3}{r}+\epsilon^4
       +\frac{A_*}{r}+\epsilon A_*\right\},
 \qquad A_*=\sum_kA_k .
\]
The profile and metric estimates imply
\begin{equation}\label{linear-eq-fluxkernelenvelopes}
\begin{aligned}
 |\mathcal A_i|+|\partial_t\mathcal A_i|
 &\leq C\left\{\frac{\epsilon}{r}
      +\epsilon^2L(\epsilon r)
      +\frac{\epsilon\Lambda_\epsilon}{r^3}\right\}
      +\mathcal E_\epsilon(r),\\
 |\mathcal B_i|+|\partial_t\mathcal B_i|
 &\leq C\left\{\frac{\epsilon}{r^3}
      +\epsilon^2(1+L(\epsilon r))
      +\frac{\epsilon\Lambda_\epsilon}{r}\right\}
      +\mathcal E_\epsilon(r).
\end{aligned}
\end{equation}
Here the derivatives hold $\psi$ fixed.
To check every term, the leading amplitudes and phases are
$f_i^0=\chi_R(w'/w)(e_i\cdot e_r)$ and
$p_i^0=\chi_R(e_i\cdot e_\theta)/r$.
The modulus corrections in $f_i$ obey $m_0$, and their time
derivatives obey $m_1$ in
\eqref{linear-eq-adjustmentenvelope}.
The regular part of $p_i$ is $\epsilon L(\epsilon r)$;
its $\sigma$ derivative is $\epsilon^2L(\epsilon r)$.
The logarithmic part is independent of $\theta$, so its angular
derivative has size $C\epsilon$.
The remaining inner terms follow by one normal derivative of
\eqref{linear-eq-profiletable}.  The outer phase contributes at
most $C\epsilon A_*\Lambda_\epsilon^P$ to $p_i$ and
$CA_*\Lambda_\epsilon^P$ to $\partial_tp_i$; an additional
$\sigma$ derivative supplies $\epsilon$.

Also $m_\theta=O(\epsilon r)$ plus differentiated amplitude
corrections, and $|D^\sigma|+|\partial_tD^\sigma|
\leq C\epsilon\Lambda_\epsilon$.
Thus the $2\omega m\partial_\theta p_i^0$ term has size
$C\epsilon/r$, and $\omega m_\theta p_i^0$ has size
$C\epsilon^2$.  The phase kernel's drift contribution from $f_i^0$
has size $C\epsilon\Lambda_\epsilon/r^3$.
The amplitude kernel has respectively
$C\epsilon/r^3$ from $\omega\partial_\theta f_i^0$
and $C\epsilon\Lambda_\epsilon/r$ from $D^\sigma p_i^0$.
The regular part of the latter is bounded by the same envelope,
since $\epsilon rL(\epsilon r)$ is bounded.
All amplitude remainders, including $b_1$, fit
$\mathcal E_\epsilon$.  A time derivative of the outer phase
in $\partial_\sigma p_i$ gives $\epsilon A_*$, the largest
bounded remainder in the phase kernel.
This verifies \eqref{linear-eq-fluxkernelenvelopes}.
On the fixed core the lower coefficient is
$O(\epsilon\Lambda_\epsilon)+O(\epsilon^2\Lambda_\epsilon^P)
+O(A_*\Lambda_\epsilon^P)$.
This follows from the smooth $\mathcal C$ and its $\sigma$
derivative, the regular phase, and $c_\epsilon Y^\sigma$.

Here are the dual estimates for the phase kernel.  Its nonzero
angular part is measured in $rL_y^2$.  For its angular mean form
the radial primitive of $\int\mathcal A_i\,d\theta$.
For the three leading envelopes in the first line of
\eqref{linear-eq-fluxkernelenvelopes} these primitives are bounded by
\[
 C\epsilon r,\qquad
 C\epsilon^2r^2(1+L(\epsilon r)),\qquad C\epsilon\Lambda_\epsilon .
\]
Their norms in \eqref{linear-eq-dualphase} and boundary values
are bounded by
$C\epsilon R+C\epsilon^2R^2+
C\epsilon\Lambda_\epsilon^{3/2}\leq C$.
The same bound holds for the nonzero angular parts.
The remainder has dual norm at most
\[
 C\Lambda_\epsilon^P
       \{\epsilon^2\Lambda_\epsilon^{P+3/2}
                           +\epsilon^2+A_*/\epsilon\}=o(1).
\]
The amplitude kernel has $L_y^2$ norm bounded by
$C\epsilon\Lambda_\epsilon^{3/2}
+C\Lambda_\epsilon^P(\epsilon^2+A_*\sqrt{\Lambda_\epsilon})$,
and hence by a constant.  On the core use
\eqref{linear-eq-phasecoremass}; its logarithmic cost is absorbed
by the positive powers just displayed.
At $R_I$ use \eqref{linear-eq-exteriortrace}.
Consequently the lower integral in \eqref{linear-eq-fluxkernels}
has norm at most $C\mathcal N_*(q)^{1/2}$.
The first integral is bounded by
$C\sqrt{\Lambda_\epsilon}K(q)^{1/2}$,
since its test has $L_y^2$ norm $C\sqrt{\Lambda_\epsilon}$.

For the refined translation estimate insert
$\psi_1=p_j$, $\psi_2=-f_j$ in the lower integral.
The only leading logarithmic integral is
$\epsilon\int_1^{2R}dr/r$.
The other phase-kernel terms use
$\epsilon^2\int L(\epsilon r)\,dr=O(\epsilon)$ and
$\epsilon^3\int L(\epsilon r)^2r\,dr=O(\epsilon)$.
The drift is paired with a translation amplitude, so its leading
integral contains $r^{-3}r^{-1}r\,dr$ and is
$O(\epsilon\Lambda_\epsilon)$.
The remainder contributes at most
$C\epsilon^2\Lambda_\epsilon^P+
CA_*\Lambda_\epsilon^{P+1}=o(\epsilon\Lambda_\epsilon)$.
The core has the same bound.
Splitting only this lower integral according to
$q=w_*+\mathscr T_t\beta$, while leaving $\eta_\sigma$ as the
derivative of the complete field, proves
\eqref{linear-eq-fullfluxbound}.

For a coefficient derivative with $z$ fixed, outside the core
\[
 \left.\partial_t\psi\right|_z=-(\mu+i\alpha_t)\psi,\qquad
 \mu=\partial_t\log(\sqrt{j_0}\rho).
\]
Here $\mu$ is bounded.  Differentiating $\psi_\sigma$ adds
$-(\mu_\sigma+i\alpha_{t\sigma})\psi$.
The geometric part of $\mu_\sigma$ is $O(\epsilon^2r)$;
multiplied by $p_i$ it gives the $\epsilon^2L(\epsilon r)$
envelope in \eqref{linear-eq-fluxkernelenvelopes}.
The amplitude and phase remainders obey the same
$\mathcal E_\epsilon$ bound.  On the core use the differentiated
smooth matrix $\sqrt{j_0}\mathcal C$.
The preceding dual calculation therefore proves
\eqref{linear-eq-fluxcoefficientderivative}.

The principal matrix on translations is
$\mathsf G_{ij}=\int M\operatorname{Re}(\tau_i\overline{\tau_j})\,dy$.
It and its time derivative have size $C\Lambda_\epsilon$.
The coefficient $\mathsf C_{ij}$ is $\mathcal F_i(\mathscr T_te_j)$.
The lower-integral estimate just proved and
$\|\partial_\sigma\eta(\mathscr T_te_j)\|_{L_y^2}
\leq C\epsilon\sqrt{\Lambda_\epsilon}$
give $|\mathsf C|\leq C\epsilon\Lambda_\epsilon$.
Its derivative satisfies the same bound: differentiated regular
phases use the integrals above, and the largest new outer
contribution is $CA_*\Lambda_\epsilon^P=o(\epsilon)$.
For a time-differentiated lift use
\eqref{linear-eq-zetasplit} and the transverse and longitudinal
bounds following \eqref{linear-eq-zetabound}.
This proves \eqref{linear-eq-translationfluxmatrices}, including
the stated mixed coefficient bounds.

Finally differentiate the projection in the $\eta$ variables.
There is a smooth real-linear test
\[
 k_i=(\sqrt{j_0}\mathcal C)^*(i\tau_i),
 \qquad p_i(q)=\langle\eta,k_i\rangle_y ,
\]
where the star denotes the adjoint for the real pairing.
The notation $p_i(q)$ denotes the measurement here, whereas
$p_i$ in \eqref{linear-eq-fluxkernels} is the lift's phase coefficient.
One has $\|k_i\|_{L_y^2}\leq C\sqrt{\Lambda_\epsilon}$.
The identity $p_i(iVf)=0$ for every radial real $f$ implies
\[
 \langle iWf,k_i\rangle_y
 =\langle iWf,\partial_\sigma k_i\rangle_y=0.
\]
Thus the lower differentiated test annihilates the complete
zero phase before it is estimated.
Off the core its phase and amplitude kernels are
$\partial_\sigma(j_0\rho^2f_i)$ and
$\partial_\sigma(j_0\rho^2p_i)$.
Their respective envelopes are
\[
\begin{split}
 C\frac{\epsilon^2\ell_r^P}{r^2}
       +C\Lambda_\epsilon^P(\epsilon^3/r+\epsilon^4),\\
 C\epsilon^2(1+L(\epsilon r))
       +C\Lambda_\epsilon^P
            (\epsilon^2\ell_r^P/r^2+\epsilon^3/r+\epsilon^4
                                      +\epsilon^2A_*).
\end{split}
\]
The first has weighted nonzero-phase $L_y^2$ norm and the
second has amplitude $L_y^2$ norm at most
$C\epsilon\sqrt{\Lambda_\epsilon}$.
The core follows from the smooth differentiated identification
and the zero-phase cancellation.  It follows that
\[
 \|\langle\eta,\partial_\sigma k_i\rangle_y\|_{L^2_\sigma}
 \leq C\epsilon\sqrt{\Lambda_\epsilon}
       \{E_*(w_*)^{1/2}+\|\beta\|_2\}.
\]
The exterior energy controls any boundary trace in this estimate.
Now differentiate $b=N^{-1}p(q)$.
The term with $\eta_\sigma$ gives
$C\sqrt{\Lambda_\epsilon}K(q)^{1/2}$,
the lower test gives the displayed bound, and
\eqref{linear-eq-Gramest} gives $C\epsilon\|b\|_2$.
This proves \eqref{linear-eq-sympderivative}.
\end{proof}

\begin{lemma}[Symplectic coefficient energy]\label{linear-lem-projectedenergy}
Along the homogeneous equation,
\begin{equation}\label{linear-eq-projectedenergy}
 {c_\epsilon}\frac d{dt}\|b\|_2^2
 \leq C\{\mathcal X(w_*)^2+{\Lambda_\epsilon} K(q)
        +{c_\epsilon}{\Lambda_\epsilon}\|\beta\|_2^2+{c_\epsilon}{\Lambda_\epsilon}\|b\|_2^2\}.
\end{equation}
All norms are on the fixed material circle or the complete
reference space as previously defined.  In particular the estimate
requires only one longitudinal derivative of the perturbation.
\end{lemma}
\begin{proof}
Take $a=N^{-T}b$ in \eqref{linear-eq-projectedequation}.  The left side is
\[
 \frac{c_\epsilon}{2}\frac d{dt}\|b\|_2^2
       +c_\epsilon\int b\cdot N^{-1}N_tb .
\]
Since $N=\beta_R\mathbb J$ and
$N^{-1}N_t=(\beta_{R,t}/\beta_R)\operatorname{Id}$,
the second term is bounded by $Cc_\epsilon\|b\|_2^2$.
For the coefficient defect use $q=w_*+\mathscr T_t\beta$ and
Lemma~\ref{linear-lem-Jacobiestimates}:
\[
 |\langle z,\mathscr D_ia_i\rangle|
 \leq C\epsilon\|b\|_2\mathcal X(w_*)
                  +Cc_\epsilon\|b\|_2\|\beta\|_2
 \leq C\mathcal X(w_*)^2
                  +Cc_\epsilon(\|b\|_2^2+\|\beta\|_2^2).
\]

For the complete longitudinal term use
Lemma~\ref{linear-lem-longitudinalkernels}.
The derivative of $N^{-T}$ contributes $O(\epsilon)b$, so
\eqref{linear-eq-sympderivative} gives
\[
 \|a'\|_2\leq C\sqrt{\Lambda_\epsilon}K(q)^{1/2}
 +C\epsilon\sqrt{\Lambda_\epsilon}
       \{E_*(w_*)^{1/2}+\|\beta\|_2+\|b\|_2\}.
\]
Multiply this by \eqref{linear-eq-fullfluxbound} and apply
Young's inequality.  The result is
\[
\begin{split}
 \left|\int a_i'\mathcal F_i(q)\,d\sigma_0\right|
 \leq{}&C\Lambda_\epsilon K(q)+C\mathcal X(w_*)^2
       +C\epsilon^2\Lambda_\epsilon E_*(w_*)\\
 &+C\epsilon^2\Lambda_\epsilon^2
                      (\|\beta\|_2^2+\|b\|_2^2).
\end{split}
\]
The third term is bounded by $C\mathcal X(w_*)^2$, and
$\epsilon^2\Lambda_\epsilon^2=c_\epsilon\Lambda_\epsilon$.
Together with the coefficient estimate this proves
\eqref{linear-eq-projectedenergy}.
The kernel proof includes the radial phase primitive and the
interface trace, so no undifferentiated zero-phase mass is used.
Only $\eta_{\sigma_0}$, the derivative in $K(q)$, occurs in the
field-derivative term.
\end{proof}

\begin{samepage}
\subsection{Positive energy and completion of the proof}
\label{linear-sec-assembly}

We now combine the preceding bounds into a positive quadratic
functional and apply Gronwall's inequality to obtain the estimate in
Theorem~\ref{curvedfilamentlineartheorem}.  We also establish existence
and uniqueness and show that the constants in this estimate can be
chosen independently of the approximation orders.

For positive constants $C_{\rm tr}$ and $\gamma$ to be chosen below,
define the quadratic functional
\begin{equation}\label{linear-eq-Qintro}
 Q_t[q]=B_t[q]+C_{\rm tr}{c_\epsilon}\|b_t(q)\|_{L^2_{\sigma_0}}^2
                                      +\gamma{c_\epsilon}^2\|q\|_2^2.
\end{equation}
\end{samepage}

\subsubsection{Static comparisons for the same symplectic projection}
By \eqref{linear-eq-staticblocks},
\[
B_t[q]\geq k(E_*(w_*)+K(q))-C{c_\epsilon}\|\beta\|_2^2
                                      -C{c_\epsilon}^2\|q\|_2^2.
\]
Equation \eqref{linear-eq-twoprojections} gives
\[
 {c_\epsilon}\|\beta\|_2^2\leq2{c_\epsilon}\|b\|_2^2
                                  +C{c_\epsilon}{\Lambda_\epsilon} E_*(w_*).
\]
The second coefficient tends to zero.  Thus, for sufficiently
large fixed $C_{\rm tr}$ and $\gamma$, the quadratic functional
\eqref{linear-eq-Qintro} satisfies
\begin{equation}\label{linear-eq-Qstrong}
 Q_t[q]\geq k_0\{E_*(w_*)+K(q)+{c_\epsilon}\|b\|_2^2
                                      +{c_\epsilon}^2\|q\|_2^2\}.
\end{equation}
In particular, the translation correction is controlled by $Q_t$.

The symplectic remainder $q^\perp=q-\mathscr T_tb$ is fully
represented: $q^\perp=w_*+\mathscr T_t(\beta-b)$.
Equations \eqref{linear-eq-liftbounds} and
\eqref{linear-eq-twoprojections} imply
\[
 \mathcal N_*(q^\perp)\leq C{\Lambda_\epsilon} E_*(w_*).
\]
This logarithmic loss concerns the symplectic remainder $q^\perp$.
For the derivative estimate we use the stronger comparison
\eqref{linear-eq-Qstrong}.

The corresponding $H^1$ comparisons are
\begin{equation}\label{linear-eq-Qintrocomparisons}
 {c_\epsilon}^2\|q\|_{H^1}^2\leq C Q_t[q],\qquad
 Q_t[q]\leq C\|q\|_{H^1}^2.
\end{equation}
Indeed, \eqref{linear-eq-eulerianform}, \eqref{linear-eq-lowcoeff} and
Young's inequality give
\[
 \mathfrak b_t[h]\geq\tfrac12\|\nabla h\|_2^2-C\|h\|_2^2.
\]
Using \eqref{linear-eq-H1map}, multiply the resulting estimate by
${c_\epsilon}^2$ and use \eqref{linear-eq-Qstrong}.  This proves the lower bound in
\eqref{linear-eq-Qintrocomparisons}.  The upper bound follows from the
bounded coefficients of \eqref{linear-eq-eulerianform} and
\[
 {c_\epsilon}\|b\|_2^2\leq C{c_\epsilon}{\Lambda_\epsilon}\|q\|_2^2\leq C\|q\|_2^2.
\]
Thus both norm comparisons and the control of the added form have
constants independent of $\epsilon$.

\subsubsection{The combined coefficient bound}
At a fixed time the decomposition
$q=w_*+\mathscr T_t\beta$ can be used in the bilinear form
$\dot B_t$ without differentiating $w_*$ or $\beta$.
Equations \eqref{linear-eq-complementderivative},
\eqref{linear-eq-dotdiagonal} and \eqref{linear-eq-dotmixed} give
\begin{equation*}
 |\dot B_t[q]|\leq C\{\mathcal X(w_*)^2
       +{\Lambda_\epsilon}\|\beta'\|_2^2+{c_\epsilon}{\Lambda_\epsilon}\|\beta\|_2^2\}.
\end{equation*}
Use \eqref{linear-eq-betaderivative} and \eqref{linear-eq-auxX}.  The right-hand
side is bounded by
\[
 C\{E_*(w_*)+{\Lambda_\epsilon} K(q)+{c_\epsilon}{\Lambda_\epsilon}\|\beta\|_2^2
                                         +{c_\epsilon}^2\|q\|_2^2\}.
\]
Next,
\[
 {c_\epsilon}{\Lambda_\epsilon}\|\beta\|_2^2
 \leq C{c_\epsilon}{\Lambda_\epsilon}\|b\|_2^2+C{c_\epsilon}{\Lambda_\epsilon}^2 E_*(w_*),
 \qquad {c_\epsilon}{\Lambda_\epsilon}^2=\epsilon^2{\Lambda_\epsilon}^3=o(1).
\]
It follows from \eqref{linear-eq-Qstrong} that
\begin{equation}\label{linear-eq-dotBfinal}
                         |\dot B_t[q]|\leq C{\Lambda_\epsilon} Q_t[q].
\end{equation}
Likewise \eqref{linear-eq-projectedenergy}, \eqref{linear-eq-auxX} and the same
comparison prove
\begin{equation}\label{linear-eq-translationratefinal}
 C_{\rm tr}{c_\epsilon}\frac d{dt}\|b\|_2^2\bigg|_{\rm homogeneous}
                                                   \leq C{\Lambda_\epsilon} Q_t[q].
\end{equation}
Both bounds are expressed in terms of $Q_t$ and contain a single
factor $\Lambda_\epsilon$.

\subsubsection{The mass term, forcing, and existence}
Let $A_t[q]=C_{\rm tr}{c_\epsilon}\|b_t(q)\|_2^2$ and
$M[q]=\|q\|_2^2$.  From \eqref{linear-eq-Qstrong}, $A_t\leq C Q_t$.
Polarising $B_t=Q_t-A_t-\gamma{c_\epsilon}^2M$, and applying
Cauchy--Schwarz separately to the nonnegative forms $Q_t,A_t,M$,
gives
\begin{equation*}
 |B_t[q,iu]|\leq C Q_t[q]^{1/2}\|u\|_{H^1}.
\end{equation*}
We used the upper comparison for $Q_t[iu]$ and
$A_t[iu]\leq C Q_t[iu]$.  For $u=q$, the lower $H^1$ comparison
gives
\begin{equation}\label{linear-eq-smallmassrate}
                    {c_\epsilon}|B_t[q,iq]|\leq C Q_t[q].
\end{equation}

Along \eqref{linear-eq-movingequation}, self-adjointness gives
\[
 \frac d{dt}B_t[q(t)]=\dot B_t[q]-\frac2{{c_\epsilon}}B_t[q,iG],
 \qquad
 {c_\epsilon}\frac d{dt}\|q\|_2^2=2B_t[q,iq]+2\langle G,iq\rangle.
\]
The exact derivative of the complete quadratic functional is consequently
\begin{equation*}
 \frac d{dt}Q_t[q(t)]
 =\dot B_t[q]+\frac d{dt}A_t[q(t)]\bigg|_{\rm homogeneous}
            +2\gamma{c_\epsilon} B_t[q,iq]-\frac2{{c_\epsilon}}Q_t[q,iG].
\end{equation*}
The forcing contribution of the projection has combined with the
other forcing terms to give the last pairing.  Equations
\eqref{linear-eq-dotBfinal}, \eqref{linear-eq-translationratefinal},
\eqref{linear-eq-smallmassrate}, and Cauchy--Schwarz for $Q_t$,
followed by \eqref{linear-eq-H1map} for $G$, give
\begin{equation}\label{linear-eq-Qintroderivative}
 \frac d{dt}Q_t[q(t)]\leq C_{\rm lin}{\Lambda_\epsilon} Q_t[q(t)]
                  +\frac C{{c_\epsilon}}Q_t[q(t)]^{1/2}\|F(t)\|_{H^1}.
\end{equation}

For smooth data, regularise $Q_t^{1/2}$ by $(Q_t+\delta)^{1/2}$,
apply an integrating factor, and let $\delta$ decrease to zero.
Since $q(0)=0$,
\[
 Q_t[q(t)]^{1/2}\leq\frac C{{c_\epsilon}}
       \int_0^t e^{C_{\rm lin}{\Lambda_\epsilon}(t-s)}\|G(s)\|_{H^1}\,ds.
\]
The two $H^1$ comparisons give
\begin{equation}\label{curvedfilamentpreciseestimate}
 \|\Phi(t)\|_{H^1}
 \leq\frac C{{c_\epsilon}^2}\int_0^t
       e^{C_{\rm lin}{\Lambda_\epsilon}(t-s)}\|F(s)\|_{H^1}\,ds.
\end{equation}
Finally
${c_\epsilon}^{-2}=\epsilon^{-4}{\Lambda_\epsilon}^{-2}\leq\epsilon^{-4}$ and
$e^{C_{\rm lin}{\Lambda_\epsilon} T}=\epsilon^{-C_{\rm lin}T}$, proving
\eqref{curvedfilamentlinearestimate}.

For fixed $\epsilon$, we construct the solution using Duhamel's formula.
The real-linear zeroth-order part of $-H^X_t$ is
\[
 (1-2|V|^2)h-V^2\overline h,
\]
a bounded time-continuous operator on $H^1$ because
$V\in C_tW_X^{1,\infty}$.  The Duhamel equation based on
$e^{it\Delta_X/{c_\epsilon}}$ is a contraction on a sufficiently short
fixed-$\epsilon$ interval.  Iteration covers $[0,T]$ and also gives
uniqueness.  Smooth forcings give sufficiently regular solutions
for the identities above.

To justify core calculations, a smooth compactly supported field
can be cut off near the zero curve with a logarithmic radial
cutoff.  Its Dirichlet cost is $O(1/|\log\delta|)$ times the fixed
rescaled curve length, which tends to zero for fixed $\epsilon$.
These fields are dense in $H^1$.  All complete form identities have
smooth Cartesian coefficients at the core, so their limits are the
original forms.  For the projected expressions, use
\eqref{linear-eq-longflux} to retain only first derivatives.  Their
quadratic bounds extend continuously on $H^1$.

Approximate $F\in L^1_tH^1_X$ by smooth forcings.  Estimate
\eqref{curvedfilamentpreciseestimate} applies to differences of the corresponding
solutions.  They converge in $C_tH^1_X$ to the fixed-$\epsilon$
Duhamel solution.  Passing to the limit gives the estimates for the quadratic
functional and \eqref{curvedfilamentpreciseestimate} for every such $F$,
with a unique solution $\Phi\in C([0,T];H^1)$.  Taking
$F\in L^\infty_tH^1_X$ also gives \eqref{curvedfilamentlinearestimate}
and completes the proof of
Theorem~\ref{curvedfilamentlineartheorem}.

\subsubsection{Uniformity in the approximation order}

\begin{lemma}[Uniform linear growth constant]
\label{linear-lem-uniformgrowth}
Fix the reference binormal flow, its low geometric neighbourhood and
the cutoffs.  For every $M\geq20$ and $N\geq12$, the construction in
Theorem~\ref{arbitraryapproximationtheorem}, with its finite derivative
reserve, satisfies Theorem~\ref{curvedfilamentlineartheorem} with
$C_T,C_{\rm lin}$ independent of $M,N$.  The constants
$C_{\rm tr},\gamma,C$ in \eqref{linear-eq-Qintro},
\eqref{linear-eq-Qintrocomparisons}, \eqref{linear-eq-Qintroderivative}
and \eqref{curvedfilamentpreciseestimate} have the same uniformity.
The threshold $\epsilon_0(M,N)>0$ may depend on these orders.
The assertion is uniform on the initial-curve neighbourhoods in that
construction, restricted to the fixed low geometric neighbourhood.
\end{lemma}

\begin{proof}
Fix the finite record used by the $H^1$ energy proof: the coefficient
bounds through four spatial and two material time derivatives, the
Jacobi coefficients through one material time derivative, the first
coefficient derivative of the longitudinal flux, and the residual
bounds through two material time and four spatial derivatives, in
the required local and global norms, supplied by
Lemma~\ref{linear-lem-residualreserve}.
All geometric norms required for this record are at fixed finite
indices.  Apply Theorem~\ref{uniformlocalcurvetheorem} at the enlarged
indices required for this record.  The slowly varying curve $z_*$
is obtained from \eqref{seed-curve-correction} before the depth of the
elliptic batches is chosen.  Its required geometric bounds depend
only on this fixed record.

At every spatial level $s$ in the record, the subsequent increments
satisfy
\[
 \sum_{k=2}^{K-1}\|\partial_tg_k\|_{C_tH^s}
 \leq C_{M,N}\epsilon^{5/4}\Lambda_\epsilon^{P_{M,N}},
 \qquad
 \sum_{k=2}^{K-1}\|\partial_{tt}g_k\|_{C_tH^s}
 \leq C_{M,N}\epsilon^{1/4}\Lambda_\epsilon^{P_{M,N}}.
\]
Their zero initial data give the corresponding spatial bound for
$\sum g_k$.  Thus the actual curves, frames and ambient velocities
stay in a common bounded set at the required indices.  Their
geometric bounds can be fixed independently of $M,N$ by decreasing
$\epsilon_0(M,N)$.  Choose $d$ using this common geometric bound.
This fixes all tube and cutoff constants.

Evaluate the first-profile estimates at the fixed output orders
of the record.  The source, its normalised inverse and the two
chain-rule derivatives in Subsection~\ref{linear-sec-backgroundbounds}
then have constants and radial logarithmic powers bounded independently
of $M,N$.  Write $C_{\rm b}$ for a bound on these quantities and the
geometric quantities above, and $P_{\rm b}$ for a fixed radial
logarithmic exponent.  In particular, the integrals
\[
 \int_1^\infty \ell_r^{2P_{\rm b}}r^{-3}\,dr
\]
that occur in the first-profile Jacobi estimates have fixed finite
values.  The leading core logarithm is the pure phase identified
in \eqref{core-logarithm-splitting}.  Its coefficient is part of
this same fixed record.

The remaining slow profiles, radiative profiles and outer amplitudes
can have constants and logarithmic powers depending on the finite
construction.  Collect them in $C_{M,N},P_{M,N}$, increasing these
quantities finitely many times below.  The geometric sequence in
\eqref{quantitativeinductionassertions} gives
\[
 A_*=\sum_{k=2}^K A_k\leq2\epsilon^{5/4}
\]
for sufficiently small $\epsilon$, independently of $K$.  The
construction-dependent errors in the coefficient estimates below
are bounded, at their stated scales, by
\[
 \delta_{M,N}(\epsilon)
 :=C_{M,N}\epsilon^{1/4}\Lambda_\epsilon^{P_{M,N}}
 \longrightarrow0.
\]
We verify this comparison in the norms used by the energy proof.

For the static comparisons, the higher profiles give
$O(\epsilon^2\Lambda_\epsilon^{P_{M,N}})$ modulus corrections and
$O(\epsilon A_*\Lambda_\epsilon^{P_{M,N}})$ radiative phases on
the fixed tube.  Their contributions to $\rho/w$, the core matrix,
the Gram matrices and their required derivatives are small in
the norms of Subsections~\ref{linear-sec-regions} and
\ref{linear-sec-translations}.  In the phase defect
\eqref{linear-eq-edefect}, multiplication by $\langle r\rangle$
and division by $c_\epsilon$ give a slow error
$O(\epsilon\Lambda_\epsilon^{P_{M,N}})$ and a radiative error
$O((A_*/\epsilon)\Lambda_\epsilon^{P_{M,N}})$.
The density defect has its bounded geometric part plus a small
amplitude correction.  The exterior transition terms have the
small $L_X^{3/2}$ norm calculated in
Subsection~\ref{linear-sec-regions}.
Lemma~\ref{linear-lem-residualreserve} gives all the residual
seminorms as $O(\epsilon^{20})$ with a construction-dependent
constant; their ratios to the residual scales in the energy
proof tend to zero.  Hence the regional coercivity, continuity
and projection constants have bounds depending only on $C_{\rm b}$
once $\delta_{M,N}$ is small.

For the differentiated coefficients, the table proving
\eqref{linear-eq-lowercoeffnorms} has a positive power of $\epsilon$
in every profile-dependent entry.  Its largest radiative term is
\[
 C_{M,N}\Lambda_\epsilon^{P_{M,N}}A_*/\epsilon
 \leq2C_{M,N}\epsilon^{1/4}\Lambda_\epsilon^{P_{M,N}}.
\]
The only entry without a positive power is the geometric
$L_X^3$ bound for $\nabla\operatorname{div}Y$.
The core errors are small after the logarithmic phase estimate.
The far-field terms, divided by $c_\epsilon$ or $c_\epsilon^2$
as in \eqref{linear-eq-lowercoeffnorms}, also tend to zero.
Thus the constant in \eqref{linear-eq-complementderivative}
has a bound independent of $M,N$.

For the Jacobi estimates, the relevant scales in
\eqref{linear-eq-Jacobiestimates} are $\epsilon$ in the dual norm,
$c_\epsilon$ in a translation pairing and
$c_\epsilon\Lambda_\epsilon$ in a differentiated pairing.
For example, the higher slow terms in
\eqref{linear-eq-adjustmentenvelope} satisfy, with $R=d/\epsilon$,
\[
 \left\|\Lambda_\epsilon^{P_{M,N}}
       \left(\frac{\epsilon^2}{\langle r\rangle}+\epsilon^3\right)
             \right\|_{L_y^2(r<2R)}
 \leq C\epsilon^2\Lambda_\epsilon^{P_{M,N}+1/2}
 =o(\epsilon).
\]
The bounded first-time amplitude tail gives
$b_1R\Lambda_\epsilon^{P_{M,N}}
=O(\epsilon^{5/4}\Lambda_\epsilon^{P_{M,N}})$.
The second-time tail in the phase dual norm gives
\[
 \epsilon^{-1}c_\epsilon b_2R^2
                  \Lambda_\epsilon^{P_{M,N}}
 \leq C\epsilon^{1/4}\Lambda_\epsilon^{P_{M,N}+1}.
\]
All remaining higher terms have the additional powers displayed
in the proof of Lemma~\ref{linear-lem-Jacobiestimates}.
For the paired estimates use the angular cancellations there
before taking absolute values.  The leading critical integral is
$\epsilon^2\int_1^{2R}dr/r=O(c_\epsilon)$, with coefficient
bounded by $C_{\rm b}$.  The higher-profile terms in that proof
are small relative to the stated paired scales.  The rotated
second lift derivative has largest radiative ratio
$O((A_*/\epsilon)\Lambda_\epsilon^{P_{M,N}})$ in the dual scale.
The connection and its derivative, divided by $c_\epsilon$,
are bounded by the convergent integrals in
Lemma~\ref{linear-lem-connection}.  This proves uniform bounds
for the constants entering the static and differentiated
translation blocks.

For the longitudinal terms, the remainder following
\eqref{linear-eq-fluxkernelenvelopes} has phase dual norm at most
\[
 C_{M,N}\Lambda_\epsilon^{P_{M,N}}
 \left\{\epsilon^2\Lambda_\epsilon^{P_{M,N}+3/2}
                     +\epsilon^2+A_*/\epsilon\right\}=o(1).
\]
Its translation pairing is
$O(\epsilon^2\Lambda_\epsilon^{P_{M,N}}
     +A_*\Lambda_\epsilon^{P_{M,N}+1})
=o(\epsilon\Lambda_\epsilon)$.
These estimates include a coefficient time derivative.
The leading kernels have constants bounded by $C_{\rm b}$;
their critical translation integral is
$\epsilon\int_1^{2R}dr/r=O(\epsilon\Lambda_\epsilon)$.
The principal measurement has its fixed
$O(\sqrt{\Lambda_\epsilon})$ norm, and the lower differentiated
projection test is estimated after annihilating the zero phase.
Thus the constants in Lemma~\ref{linear-lem-longitudinalkernels},
including the matrices divided by their displayed scales, are
bounded independently of the construction orders.

There are finitely many estimates in this fixed energy record.
Choose $\epsilon_0(M,N)$ so that their construction-dependent
remainders are all smaller than the fixed fractions needed in
the coercivity comparisons.  The constants in
\eqref{linear-eq-staticblocks}--\eqref{linear-eq-dotmixed} and
\eqref{linear-eq-projectedenergy} are then bounded in terms of
$C_{\rm b}$ and the fixed planar and cutoff constants.
Choose $C_{\rm tr},\gamma$ and the lower constant in
\eqref{linear-eq-Qstrong} using these bounds.
The algebra proving \eqref{linear-eq-dotBfinal},
\eqref{linear-eq-translationratefinal} and
\eqref{linear-eq-smallmassrate} gives a fixed $C_{\rm lin}$.
The norm comparisons and the forcing estimate have the same
uniformity.  All dependence on $M,N$ in this argument is confined
to the smallness threshold and the intermediate profile constants.
\end{proof}

\section{Solving the full problem}
\label{fullnonlinearsection}

We now combine the approximation and linear estimates in a contraction
argument to complete the proof of Theorem~\ref{mainthm}.  Throughout this
section we use the global rescaled variable \(X\) introduced in
\eqref{globalrescaledvariable}.
Here $m$ is the spatial Sobolev index of the correction in
$C([0,T];H_X^m)$, $M$ is the residual order and $N$ is the regularity
index in Theorem~\ref{arbitraryapproximationtheorem}.  Fix
\[
 C_0:=C_{\rm lin}+1,
\]
where $C_{\rm lin}$ is supplied by
Lemma~\ref{linear-lem-uniformgrowth}.  This constant is chosen before
$m,M,N$.  The added unit absorbs the fixed logarithmic powers in the
higher Sobolev estimate below.  Multiplicative constants and the
smallness threshold may depend on the chosen finite orders.

\subsection{Higher-order linear estimates}

\begin{corollary}[Higher-order linear estimate]
\label{higherordercurvedlinearcorollary}
Fix an integer \(m\geq1\).  If the accuracy and regularity indices in
Theorem \ref{arbitraryapproximationtheorem} are sufficiently large in
terms of \(m\), then the solution operator
\(\mathcal T_\epsilon:F\mapsto\Phi\) for
\eqref{curvedlinearproblem} satisfies
\begin{equation}\label{higherordercurvedlinearestimate}
 \|\mathcal T_\epsilon[F]\|_{L_t^\infty H^m_X}
 \leq C_m\epsilon^{-4m-mC_0T}
 \|F\|_{L_t^\infty H^m_X}.
\end{equation}
Here $C_m$ may depend on the finite construction orders, whereas
$C_0$ is the constant fixed above.
\end{corollary}

\begin{proof}
It is enough to consider smooth $F$ and $\Phi$; the general result then
follows by density.  Commute the original rescaled problem
\eqref{curvedlinearproblem}, whose principal part is the fixed Laplacian.
For $|\alpha|=r\leq m-1$, put
$\Phi_\alpha=\partial_X^\alpha\Phi$.  Then
\begin{equation}\label{commutedcurvedlinearproblem}
 ic_\epsilon\partial_t\Phi_\alpha
 +\mathcal L_{\epsilon,t}^X\Phi_\alpha
 =\partial_X^\alpha F+\mathcal C_\alpha(\Phi),
 \qquad \Phi_\alpha(0)=0,
\end{equation}
where
\[
 \mathcal C_\alpha(\Phi)
 :=[\mathcal L_{\epsilon,t}^X,\partial_X^\alpha]\Phi.
\]
The Laplacian commutes with $\partial_X^\alpha$, and the commutator
vanishes when $r=0$.  For $r\geq1$, every term in
$\mathcal C_\alpha(\Phi)$ contains a derivative of $V_\epsilon$ and at most
$r-1$ derivatives of $\Phi$ before the final $H^1$ derivative is taken.
The differentiated bounds supplied by Theorem
\ref{arbitraryapproximationtheorem} and the standard Moser product estimate
give
\begin{equation}\label{higherorderlinearcommutatorbound}
 \|\mathcal C_\alpha(\Phi)\|_{H^1_X}
 \leq C_m\Lambda_\epsilon^{C_m}\|\Phi\|_{H^r_X}.
\end{equation}

Apply Theorem \ref{curvedfilamentlineartheorem} to
\eqref{commutedcurvedlinearproblem}, sum over $|\alpha|=r$, and use
\eqref{higherorderlinearcommutatorbound}.  This gives
\[
 \|\Phi\|_{L_t^\infty H^{r+1}_X}
 \leq C_m\epsilon^{-4-C_{\rm lin}T}
 \left(\|F\|_{L_t^\infty H^{r+1}_X}
 +\Lambda_\epsilon^{C_m}
  \|\Phi\|_{L_t^\infty H^r_X}\right).
\]
Starting with \eqref{curvedfilamentlinearestimate} and arguing inductively
in $r$, while retaining the logarithmic factors, gives
\[
 \|\Phi\|_{L_t^\infty H^m_X}
 \leq C_m\Lambda_\epsilon^{C_m}
          \epsilon^{-4m-mC_{\rm lin}T}
 \|F\|_{L_t^\infty H^m_X}.
\]
The logarithmic exponent is finite once $m,M,N$ are fixed.
For $T>0$, decrease $\epsilon_0(m,M,N)$ so that
$\Lambda_\epsilon^{C_m}\leq\epsilon^{-mT}$.  The choice
$C_0=C_{\rm lin}+1$ then gives
\eqref{higherordercurvedlinearestimate}.  The same choice works for
every fixed finite range of Sobolev indices after taking the minimum
of the corresponding thresholds.  For $T=0$ the estimate is immediate
from the zero initial condition.
\end{proof}

\subsection{The contraction in the rescaled variables}

Fix an integer $m\geq7$ and choose $N\geq\max\{m,12\}$ sufficiently
large for the finite derivative estimates required.  With $C_0$
already fixed, choose an integer $M$ satisfying
\begin{equation}\label{nonlinearaccuracychoice}
 M\geq
 \left\lceil8m+2mC_0T+4\right\rceil.
\end{equation}
Construct $U_\epsilon$ using these orders and the derivative reserve
of Section~\ref{arbitraryapproxsection}, and put
$V_\epsilon(X,t)=U_\epsilon(\epsilon X,t)$.
Lemma~\ref{linear-lem-uniformgrowth} applies to this construction,
since \eqref{nonlinearaccuracychoice} implies $M\geq20$.
Corollary~\ref{higherordercurvedlinearcorollary} therefore gives the
linear estimate for this same $V_\epsilon$.
All orders are now fixed; $\epsilon$ is subsequently made small.
For the contraction one may take $m=7$, in which case
\eqref{nonlinearaccuracychoice} reads
$M\geq\lceil60+14C_0T\rceil$.

Let
\begin{equation*}
 \mathcal N_\epsilon(\Phi)
 :=-2\operatorname{Re}(\overline{V_\epsilon}\Phi)\Phi
 -|\Phi|^2V_\epsilon-|\Phi|^2\Phi.
\end{equation*}
For \(m\geq2\), the algebra property of \(H^m(\mathbb R^3)\) and the
rescaled coefficient bounds give
\begin{equation}\label{rescalednonlinearestimates}
\begin{split}
 \|\mathcal N_\epsilon(\Phi)\|_{H^m_X}
 &\leq C_m\Lambda_\epsilon^{C_m}
 \big(\|\Phi\|_{H^m_X}^2+\|\Phi\|_{H^m_X}^3\big),\\
 \|\mathcal N_\epsilon(\Phi_1)
 -\mathcal N_\epsilon(\Phi_2)\|_{H^m_X}
 &\leq C_m\Lambda_\epsilon^{C_m}
 \big(\|\Phi_1\|+\|\Phi_2\|\\
 &\hspace{5em}
 +\|\Phi_1\|^2+\|\Phi_2\|^2\big)
 \|\Phi_1-\Phi_2\|.
\end{split}
\end{equation}
The unlabelled norms in the second line are \(H^m_X\) norms.

By \eqref{rescaledresidualidentity} and
\eqref{arbitraryresidualtheorem}, the residual of the rescaled
approximation satisfies
\begin{equation}\label{residualhighordernorm}
 \|R_\epsilon\|_{L_t^\infty H^m_X}
 \leq C_{M,m}\epsilon^{M-3/2}.
\end{equation}
Indeed, if \(H(X)=h(\epsilon X)\), then
\begin{equation*}
 \|\partial_X^\alpha H\|_{L^2_X}
 =\epsilon^{|\alpha|-3/2}
  \|\partial_x^\alpha h\|_{L^2_x}.
\end{equation*}
The $L_X^2$ term therefore gives the largest loss in the
inhomogeneous $H_X^m$ norm.

\begin{proposition}[Exact solution in the rescaled variables]
\label{exactrescaledsolutionproposition}
Let $m\geq7$ and let $V_\epsilon$ be the approximation selected above
using \eqref{nonlinearaccuracychoice}.  For sufficiently small
$\epsilon$, on the prescribed interval $[0,T]$ there is a smooth solution
\[
 v_\epsilon=V_\epsilon+\Phi_\epsilon
\]
of
\begin{equation}\label{exactrescaledGP}
 \mathscr G_\epsilon(v_\epsilon)=0
 \quad\hbox{in }\mathbb R^3_X\times[0,T],
 \qquad
 v_\epsilon(\cdot,0)=V_\epsilon(\cdot,0).
\end{equation}
For the chosen $m$, the correction satisfies
\begin{equation}\label{rescaledcorrectionsmallness}
 \|\Phi_\epsilon\|_{L_t^\infty H^m_X}
 \leq \rho_\epsilon(T),
 \qquad
 \rho_\epsilon(T):=\epsilon^{4m+mC_0T+1}.
\end{equation}
\end{proposition}

\begin{proof}
Fix \(m\geq7\), and put
\begin{equation*}
 a_m(T):=4m+mC_0T,
 \qquad
 \rho_\epsilon(T):=\epsilon^{a_m(T)+1}.
\end{equation*}
The choice \eqref{nonlinearaccuracychoice} gives
$M\geq2a_m(T)+4$ for this approximation.

The exact rescaled expansion is
\begin{equation}\label{exactrescalednonlinearexpansion}
\begin{split}
 &ic_\epsilon\partial_t(V_\epsilon+\Phi)
 +\Delta_X(V_\epsilon+\Phi)
 +(1-|V_\epsilon+\Phi|^2)(V_\epsilon+\Phi)\\
 &\qquad
 =R_\epsilon
 +ic_\epsilon\partial_t\Phi
 +\mathcal L_{\epsilon,t}^{X}\Phi
 +\mathcal N_\epsilon(\Phi).
\end{split}
\end{equation}
Thus a solution with zero correction at \(t=0\) is a fixed point of
\begin{equation}\label{nonlinearfixedpointmap}
 \mathcal A_\epsilon(\Phi)
 :=-\mathcal T_\epsilon
 [R_\epsilon+\mathcal N_\epsilon(\Phi)].
\end{equation}

Consider
\[
 \mathscr B_\epsilon
 :=\left\{\Phi\in C([0,T];H^m_X):
 \|\Phi\|_{L_t^\infty H^m_X}\leq\rho_\epsilon(T)\right\}.
\]
By Corollary \ref{higherordercurvedlinearcorollary},
\eqref{residualhighordernorm}, and
\eqref{nonlinearaccuracychoice},
\[
 \|\mathcal T_\epsilon[R_\epsilon]\|_{L_t^\infty H^m_X}
 \leq C\epsilon^{M-3/2-a_m(T)}
 \leq C\epsilon\,\rho_\epsilon(T).
\]
For \(\Phi\in\mathscr B_\epsilon\),
\eqref{rescalednonlinearestimates} gives
\[
 \|\mathcal T_\epsilon[\mathcal N_\epsilon(\Phi)]\|_
 {L_t^\infty H^m_X}
 \leq C\Lambda_\epsilon^C\epsilon^{-a_m(T)}
 (\rho_\epsilon(T)^2+\rho_\epsilon(T)^3)
 \leq C\Lambda_\epsilon^C\epsilon\,\rho_\epsilon(T).
\]
Hence \(\mathcal A_\epsilon\) maps
\(\mathscr B_\epsilon\) into itself.  Similarly,
\[
 \|\mathcal A_\epsilon(\Phi_1)
 -\mathcal A_\epsilon(\Phi_2)\|_{L_t^\infty H^m_X}
 \leq C\Lambda_\epsilon^C\epsilon
 \|\Phi_1-\Phi_2\|_{L_t^\infty H^m_X}.
\]
Since \(\Lambda_\epsilon^C\epsilon\to0\), this is a contraction.
Its fixed point \(\Phi_\epsilon\) satisfies
\eqref{rescaledcorrectionsmallness}.  Equations
\eqref{exactrescalednonlinearexpansion} and
\eqref{nonlinearfixedpointmap} show that
\(v_\epsilon=V_\epsilon+\Phi_\epsilon\) solves
\eqref{exactrescaledGP}.

The fixed point is smooth.  Fix $\epsilon>0$ and the finite construction.
Its inner and collar terms are smooth and
supported in a bounded region.  In the exterior the solid-angle phase
has spatial gradient $O(|x|^{-3})$ and time derivative $O(|x|^{-2})$.
The smooth source formulas and the
wave energy identity, differentiated at any finite order with
$\epsilon$ fixed, put all derivatives of the outer amplitude and of
the wave phase gradient and time derivative in $L^2$.  Thus, for every
finite $s$, $V_\epsilon\in C_tW_X^{s,\infty}$ and
$R_\epsilon\in C_tH_X^s$.  The constants here may depend on $\epsilon$.
For $s\geq m$, solve the correction equation locally in $H_X^s$ with
zero initial data.  Its $H_X^m$ uniqueness identifies this solution
with $\Phi_\epsilon$.  The skew-adjoint Laplacian and the polynomial
product estimates give
\[
 \frac{d}{dt}\|\Phi_\epsilon\|_{H_X^s}
 \leq C_{\epsilon,s}\|\Phi_\epsilon\|_{H_X^s}
       +c_\epsilon^{-1}\|R_\epsilon\|_{H_X^s},
\]
where $C_{\epsilon,s}$ uses the already bounded $H_X^m$ norm and the
coefficient bounds just stated.  Gronwall's inequality continues this
solution throughout $[0,T]$.  Uniqueness therefore gives all spatial
regularities for the original fixed point, and its equation then gives
all time derivatives.
\end{proof}

Sobolev embedding applied to \eqref{rescaledcorrectionsmallness} gives
\begin{equation}\label{rescaledcorrectionC1}
 \|\Phi_\epsilon\|_{L^\infty}
 +\|\nabla_X\Phi_\epsilon\|_{L^\infty}
 \leq C\rho_\epsilon(T).
\end{equation}
Uniform estimates at stronger finite levels require the orders to be
selected first.  Given integers $j,k\geq0$ and $A>0$, choose $m\geq7$
so that
\[
 m>k+2j+\tfrac32,\qquad
 a_m(T)+1>A+2j+k,
\]
then take $N\geq\max\{m+j+4,12\}$ and $M$ as in
\eqref{nonlinearaccuracychoice}.  For this one construction,
Lemma~\ref{linear-lem-residualreserve} gives
$\|\partial_t^\ell R_\epsilon\|_{L_t^\infty H_X^m}
\leq C\epsilon^{M-3/2}$ for $\ell\leq j$.
The differentiated profile and wave formulas at these reserved levels
also give
$\|\partial_t^\ell V_\epsilon\|_{L_t^\infty W_X^{m,\infty}}
\leq C\epsilon^{-\ell}\Lambda_\epsilon^C$.
Use the correction equation in
\eqref{exactrescalednonlinearexpansion} recursively: a time derivative
costs two spatial derivatives and $c_\epsilon^{-1}\leq\epsilon^{-2}$.
Leibniz' rule and the algebra property at the remaining indices give
\[
 \|\partial_t^\ell\Phi_\epsilon\|_{L_t^\infty H_X^{m-2\ell}}
 \leq C\epsilon^{-2\ell}\Lambda_\epsilon^C
       \big(\rho_\epsilon(T)+\epsilon^{M-3/2}\big),
 \qquad 0\leq\ell\leq j.
\]
Since $M\geq2a_m(T)+4$, Sobolev embedding and the strict exponent
margin above prove
\[
 \sum_{\ell=0}^j
 \|\partial_t^\ell\Phi_\epsilon\|_{L_t^\infty C_X^k}
 \leq C_{j,k,A}\epsilon^A.
\]
The extra margin $k$ also absorbs the factor $\epsilon^{-k}$ when
returning to physical spatial derivatives.  These are estimates for
prescribed finite indices, distinct from the fixed-$\epsilon$
smoothness of the same solution proved above.

\subsection{Return to the original variables and proof of Theorem
\ref{mainthm}}

\begin{proof}[Proof of Theorem \ref{mainthm}]
For the stated theorem take $m=7$, $N\geq12$ and $M$ as in
\eqref{nonlinearaccuracychoice}; increase these finite indices first
if a stronger estimate specified below is desired.  Let $U_\epsilon$
be the terminal approximation and let $v_\epsilon=V_\epsilon+
\Phi_\epsilon$ be the smooth solution supplied by
Proposition~\ref{exactrescaledsolutionproposition}.  Define
\[
 u_\epsilon(x,t):=v_\epsilon(x/\epsilon,t)
 =U_\epsilon(x,t)+\eta_\epsilon(x,t),\qquad
 \eta_\epsilon(x,t):=\Phi_\epsilon(x/\epsilon,t).
\]
The identity
$S_\epsilon(u_\epsilon)(x,t)
 =\mathscr G_\epsilon(v_\epsilon)(x/\epsilon,t)=0$,
divided by $\epsilon^2$, gives \eqref{gpepsintro}.  Moreover,
\begin{equation}\label{correctionreturnestimates}
 \|\partial_x^\alpha\eta_\epsilon\|_{L^\infty_x}
 =\epsilon^{-|\alpha|}
  \|\partial_X^\alpha\Phi_\epsilon\|_{L^\infty_X},
 \qquad
 \|\eta_\epsilon\|_{H^1_{\epsilon,x}}
 =\epsilon^{3/2}\|\Phi_\epsilon\|_{H^1_X},
\end{equation}
where $\|h\|_{H^1_{\epsilon,x}}^2
=\|h\|_2^2+\epsilon^2\|\nabla_xh\|_2^2$.
In particular, $\|\eta_\epsilon\|_{C_x^1}
\leq C\epsilon^{-1}\rho_\epsilon(T)=o(1)$, uniformly in time.
Write the terminal centre on the fixed material circle as
$\Gamma_\epsilon(\varrho,t)
=\gamma_\epsilon^{(K)}(\ell_\epsilon^{(K)}(t)\varrho,t)$.

We first verify finite energy.  Outside the original fixed tube,
\[
 U_\epsilon=e^{-B+i(\varphi_{\Gamma_\epsilon}+A)},
\]
where $A,B$ are the cumulative outer phase and amplitude.  By
\eqref{outerinductivebounds} and \eqref{threeDwaveXnorm},
$B\in H^1\cap L^\infty$ and $\nabla A\in L^2$.
The solid-angle gradient belongs to $L^2$ on this exterior by
\eqref{differentiatedfarfieldphaseclass}.  Hence $\nabla U_\epsilon$
and $1-|U_\epsilon|^2=1-e^{-2B}$ belong to $L^2$ there.
The interior is bounded and $U_\epsilon$ is smooth at fixed
$\epsilon$, so both norms are finite globally and $U_\epsilon$ is
bounded.  Also $\eta_\epsilon\in H^1\cap L^4$ by rescaling and
Sobolev embedding.  The gradient triangle inequality and
\[
 \|1-|U_\epsilon+\eta_\epsilon|^2\|_2
 \leq \|1-|U_\epsilon|^2\|_2
 +2\|U_\epsilon\|_\infty\|\eta_\epsilon\|_2
 +\|\eta_\epsilon\|_4^2
\]
prove $E_\epsilon(u_\epsilon(t))<\infty$ for every $t\in[0,T]$.

Using the smooth oriented normal frame $\vec e_{\alpha,\epsilon}$
of $\Gamma_\epsilon$ and its centre phase $\beta_K$, consider
the real two-component map
\[
 F_\epsilon(\varrho,y,t)
 :=e^{-i\beta_K(\varrho,t)}u_\epsilon\big(
 \Gamma_\epsilon(\varrho,t)
 +\epsilon y_\alpha\vec e_{\alpha,\epsilon}(\varrho,t),t\big).
\]
Put $b_0=w'(0)>0$, so $DW(0)=b_0I_2$.
Choose a fixed small $r_0$ with
$\|DW-b_0I_2\|_{L^\infty(B_{r_0})}\leq b_0/4$.
Exact centring, \eqref{arbitraryinnerprofile} and
\eqref{rescaledcorrectionC1} give
\[
 |F_\epsilon(\varrho,0,t)|\leq C\rho_\epsilon(T),\qquad
 \|D_yF_\epsilon-b_0I_2\|_{L^\infty(B_{r_0})}\leq b_0/2.
\]
Thus $y\mapsto y-b_0^{-1}F_\epsilon(\varrho,y,t)$ is a contraction
on a ball of radius $C\rho_\epsilon(T)$, and $F_\epsilon$ is
injective on $B_{r_0}$.  Its unique zero
$\zeta_\epsilon(\varrho,t)$ is smooth and periodic by the implicit
function theorem, and $|\zeta_\epsilon|\leq C\rho_\epsilon(T)$.
For $r\geq r_0$ in the tube, \eqref{linear-eq-lowcoeff} gives
$|U_\epsilon|\geq c\,w(r_0)>0$; the exterior has a uniform positive
lower bound as well.  Since $\|\eta_\epsilon\|_\infty=o(1)$,
there are no additional zeros.

The physical zero curve is
$\Gamma_\epsilon+h_\epsilon$, where
$h_\epsilon=\epsilon\zeta_{\epsilon,\alpha}
\vec e_{\alpha,\epsilon}$.
Differentiating $F_\epsilon(\varrho,\zeta_\epsilon,t)=0$
quantifies its displacement.  The centred background vanishes
identically at $y=0$ and has bounded mixed $\varrho$ and $y$
derivatives on $B_{r_0}$ at the reserved finite levels.  Each
$\varrho$ derivative of the composed correction costs at most
$\epsilon^{-1}$.  At the zero this gives
\[
 \partial_\varrho\zeta_\epsilon
 =-(D_yF_\epsilon)^{-1}\partial_\varrho F_\epsilon,
 \qquad
 |\partial_\varrho F_\epsilon|
 \leq C\big(|\zeta_\epsilon|+\epsilon^{-1}\rho_\epsilon(T)\big).
\]
The uniform inverse bound, followed by induction in the number of
$\varrho$ derivatives, therefore gives
\begin{equation}\label{exact-zero-curve-displacement}
 \begin{aligned}
 \sup_t\|h_\epsilon\|_{C^0(\mathbb T)}
 &\leq C\epsilon\rho_\epsilon(T),\\
 \sup_t\|h_\epsilon\|_{C^p(\mathbb T)}
 &\leq C_p\epsilon^{1-p}\rho_\epsilon(T)\quad(p\geq1).
 \end{aligned}
\end{equation}
Here any prescribed $p$ is allowed if $m>p+3/2$ and
$N\geq\max\{m+4,p+4,12\}$ were chosen first.  In particular,
the $C^1$ displacement is $O(\rho_\epsilon(T))$.
The curve is embedded by the tubular chart, and its normal winding
is one because the real normal Jacobian has positive determinant.
Theorem~\ref{uniformlocalcurvetheorem} and
\eqref{sharperarbitrarycurvebound} give the claimed
$C/\Lambda_\epsilon$ Hausdorff bound.

Project the old normal frame onto the normal plane of this zero
curve and orthonormalise it.  The resulting oriented frame differs
by $O(\rho_\epsilon(T))$.  The centre shift divided by $\epsilon$
is also $O(\rho_\epsilon(T))$.  Consequently the change to coordinates
centred on the true zero curve is $o(1)$ in $C_y^1(B_R)$ for every
fixed $R$; its small longitudinal displacement is controlled by the
same differentiated profile bounds.  With the centre phase carried
to this parametrisation, \eqref{arbitraryinnerprofile} and
\eqref{rescaledcorrectionC1} prove the asserted inner asymptotics.

Finally let $\mathcal K$ be any compact set disjoint from the
space--time filament of $\gamma^0$.  Its spatial distance from
$\Gamma_\epsilon(t)$ at time $t$ is bounded below by a positive
$\delta$ for small $\epsilon$.  In the part of $\mathcal K$ inside
the fixed tube, $r\geq\delta/\epsilon$.  The first-tail estimate
\eqref{first-tail-separated-logarithms} and the slow and radiative
bounds \eqref{slow-elliptic-profile-bounds}--
\eqref{iteratedinnerprofilebounds} make the inner logarithmic
corrections, including one physical derivative,
$O_{\mathcal K}((\epsilon^2+\epsilon A_*)\Lambda_\epsilon^C)$,
where $A_*:=\sum_{k=2}^K A_k\leq2\epsilon^{5/4}$.
The fixed physical cutoffs preserve these bounds, and
$w(r)-1$ and its physical gradient are $O_{\mathcal K}(\epsilon^2)$.
The outer estimates give $B,\nabla B=o(1)$ and $\nabla A=o(1)$.
It follows, both in this annular part and outside the tube, that
\[
 |U_\epsilon|\longrightarrow1,\qquad
 \nabla\arg U_\epsilon-\nabla\varphi_{\Gamma_\epsilon}
 \longrightarrow0\quad\hbox{uniformly on }\mathcal K.
\]
Here $\nabla\arg U$ denotes $\operatorname{Im}(\overline U\nabla U)/|U|^2$.
On $\mathcal K$, $|U_\epsilon|$ is bounded below and
$\nabla U_\epsilon$ is bounded, so
\eqref{correctionreturnestimates} passes both limits to $u_\epsilon$.
The kernel in \eqref{solidanglegradient} is smooth on this compact
set.  Since $\Gamma_\epsilon\to\gamma^0$ in $C^1$, its integral
converges uniformly to the Biot--Savart expression in
Theorem~\ref{mainthm}.  This proves all its conclusions.

For the sharper assertion in the introduction, prescribe an integer
$q\geq0$ and select $m,N$ to allow $p=\max\{q,1\}$ in
\eqref{exact-zero-curve-displacement}.  Combining that estimate with
\eqref{sharperarbitrarycurvebound} gives, in the compatible material
parametrisations,
\[
 \sup_t\|\gamma_\epsilon^{\rm zero}(t)
       -\gamma_\epsilon^{(1)}(t)\|_{H^q(\mathbb T)}
 \leq C_q\epsilon\Lambda_\epsilon^{C_q},
\]
where $\gamma_\epsilon^{\rm zero}=\Gamma_\epsilon+h_\epsilon$.
This curve is denoted by $\gamma_\epsilon$ in the introduction.
\end{proof}

\section*{Acknowledgements}

The research of M.~del Pino is supported by the Royal Society Research
Professorship grant RP-R1-180114 and by the ERC/UKRI Horizon Europe grant
ASYMEVOL, EP/Z000394/1. The research of J.~Wei is partially supported by
GRF of RGC of Hong Kong entitled ``On critical and supercritical Fujita
equation''. R.~Juneman wishes to thank the Institute of Mathematical
Sciences at CUHK for the excellent working environment.

The authors acknowledge the use of OpenAI tools  
(GPT 5.6 and GPT 6 Astra) in
the preparation of this manuscript, including assistance with its organisation, exposition, and selected calculations. The authors have reviewed the complete manuscript and take full responsibility for all mathematical statements and proofs.

\end{document}